\documentclass[11pt,reqno,a4paper]{amsart}
\usepackage{mathtools,amssymb,amsfonts} %amsmath
\usepackage[foot]{amsaddr}
\usepackage{mathrsfs,bm}
\usepackage{calligra}
\usepackage{textcomp}
\usepackage{physics}
\usepackage{slashed}
\usepackage[only llbracket,rrbracket,longmapsfrom]{stmaryrd}
\usepackage{tikz}
\usetikzlibrary{cd}
\usetikzlibrary{calc}
\usepackage{tcolorbox}
\usepackage{here}
\usepackage{multirow}
\usepackage[
 bookmarks=false,
 draft=false,
 breaklinks,
 colorlinks]{hyperref} %backref
\hypersetup{
  pdftitle={},
  pdfauthor={Taichi Yano},
  pdfsubject={},
  pdfkeywords={},
  pdfcreator={pdfLaTeX}}

\usepackage{fullpage}
\numberwithin{equation}{section}
\numberwithin{figure}{section}
\allowdisplaybreaks
\usepackage{enumitem}
\setenumerate{label=(\arabic*),nosep}
\setitemize{nosep}
\newlist{clist}{enumerate}{1}
\setlist*[clist]{label=(\roman*), nosep}

\usepackage{amsthm}
\usepackage{aliascnt}
\usepackage[capitalise]{cleveref}
\theoremstyle{definition}
\newtheorem{thm}{Theorem}[subsection]
\newaliascnt{dfn}{thm}
\newtheorem{dfn}[dfn]{Definition}
\aliascntresetthe{dfn}
\newaliascnt{prop}{thm}
\newtheorem{prop}[prop]{Proposition}
\aliascntresetthe{prop}
\newaliascnt{lem}{thm}
\newtheorem{lem}[lem]{Lemma}
\aliascntresetthe{lem}

\newaliascnt{cor}{thm}
\newtheorem{cor}[cor]{Corollary}
\aliascntresetthe{cor}

\newaliascnt{conjecture}{thm}
\newtheorem{conjecture}[conjecture]{Conjecture}
\aliascntresetthe{conjecture}

\newaliascnt{ass}{thm}
\newtheorem{ass}[ass]{Assumption}
\aliascntresetthe{ass}

\newaliascnt{mainthm}{thm}
\newtheorem{mainthm}[mainthm]{Theorem}
\aliascntresetthe{mainthm}

\theoremstyle{remark}

\newaliascnt{rem}{thm}
\newtheorem{rem}[rem]{Remark}
\aliascntresetthe{rem}

\newaliascnt{example}{thm}
\newtheorem{example}[example]{Example}
\aliascntresetthe{example}
\crefname{thm}{Theorem}{Theorems}
\crefname{dfn}{Definition}{Definitions}
\crefname{prop}{Proposition}{Propositions}
\crefname{lem}{Lemma}{Lemmas}
\crefname{cor}{Corollary}{Corollaries}
\crefname{ass}{Assumption}{Assumptions}
\crefname{conjecture}{Conjecture}{Conjectures}
\crefname{rem}{Remark}{Remarks}
\crefname{example}{Example}{Examples}
\crefname{figure}{Figure}{Figures}
\crefname{table}{Table}{Tables}
\crefname{section}{\S\!}{\S\S\!}
\crefname{subsection}{\S\!}{\S\S\!}
\crefname{subsubsection}{\S\!}{\S\S\!}
\crefname{appendix}{Appendix}{Appendices}
\crefname{equation}{}{}
\crefname{clisti}{}{}
\crefname{mainthm}{Theorem}{Theorems}

\newcommand{\ol}{\overline}
\newcommand{\wt}{\widetilde}
\newcommand{\ceq}{\coloneqq}
\newcommand{\inj}{\hookrightarrow}
\newcommand{\srj}{\twoheadrightarrow}
\newcommand{\lto}{\longrightarrow}
\newcommand{\mto}{\mapsto}
\newcommand{\sto}{\xrightarrow{\sim}}
\newcommand{\linj}{\lhook\joinrel\lto}
\newcommand{\lmto}{\longmapsto}
\newcommand{\lsto}{\xrr{\sim}}
\newcommand{\bbC}{\mathbb{C}}
\newcommand{\bbP}{\mathbb{P}}
\newcommand{\bbR}{\mathbb{R}}
\newcommand{\bbT}{\mathbb{T}}
\newcommand{\bbU}{\mathbb{U}}
\newcommand{\bbV}{\mathbb{V}}
\newcommand{\bbZ}{\mathbb{Z}}
\newcommand{\clD}{\mathcal{D}}
\newcommand{\clE}{\mathcal{E}}
\newcommand{\clF}{\mathcal{F}}
\newcommand{\clG}{\mathcal{G}}
\newcommand{\clJ}{\mathcal{J}}
\newcommand{\clL}{\mathcal{L}}
\newcommand{\clM}{\mathcal{M}}
\newcommand{\clW}{\mathcal{W}}
\newcommand{\clZ}{\mathcal{Z}}
\newcommand{\frb}{\mathfrak{b}}
\newcommand{\frg}{\mathfrak{g}}
\newcommand{\frh}{\mathfrak{h}}
\newcommand{\frm}{\mathfrak{m}}
\newcommand{\frn}{\mathfrak{n}}
\newcommand{\frw}{\mathfrak{w}}
\newcommand{\frz}{\mathfrak{z}}
\newcommand{\frB}{\mathfrak{B}}
\newcommand{\frX}{\mathfrak{X}}
\newcommand{\scW}{\mathscr{W}}
\newcommand{\ad}{\operatorname{ad}}
\newcommand{\gr}{\operatorname{gr}}
\newcommand{\rk}{\operatorname{rk}}
\newcommand{\res}{\operatorname{res}}
\newcommand{\Act}{\operatorname{Act}}
\newcommand{\End}{\operatorname{End}}
\newcommand{\Hol}{\operatorname{Hol}} 
\newcommand{\Hom}{\operatorname{Hom}}
\newcommand{\Obs}{\operatorname{Obs}}
\newcommand{\Sym}{\operatorname{Sym}}
\newcommand{\supp}{\operatorname{supp}}
\newcommand{\Open}{\operatorname{Open}}
\newcommand{\IA}{\operatorname{InnerAct}}
\newcommand{\PD}{\operatorname{PolyDiff}}
\newcommand{\dbr}[1]{\llbracket #1 \rrbracket} %using stmaryrd package
\newcommand{\dpr}[1]{(\!( #1 )\!)}
\newcommand{\xrr}[1]{\xrightarrow{\ #1 \ }}
\newcommand{\qis}{\stackrel{\text{qis}}{\simeq}}
\newcommand{\qim}{\stackrel{\text{qis}}{\xrr{\sim}}}

\begin{document}

\title{Wakimoto-type realization via BV quantization}
\author{Taichi Yano$^1$}
\address{$^1$ Graduate School of Mathematics, Nagoya University,
 Furo-cho, Chikusa-ku, Nagoya, Japan, 464-8602}
\email{yano.taichi.h5@s.mail.nagoya-u.ac.jp}
\date{November 26, 2026} %June, 2026 -- 
%\subjclass[2020]{Primary ; Secondary }
%\keywords{Wakimoto realization}

\begin{abstract}
We construct a Wakimoto-type free-field realization of Kac--Moody factorization algebras on the complex plane $\mathbb{C}$ associated with simple Lie algebras. 
This realization induces the ordinary Wakimoto realization of affine vertex algebras via the functor introduced by Costello and Gwilliam.
To this end, we develop notions of quantization and equivalence of quantizations for the free degenerate classical field theories introduced by Butson and Yoo.
We also extend our construction to higher Kac--Moody factorization algebras on $\mathbb{C}^d$.
\end{abstract}

\maketitle
\tableofcontents

%%%%%%%%%%%%%%%%%%%%%%%%%%%%%%%%%%%%%%%%%%%%%%%%%%%%%%%%%%%%%%%%%%%%%%%%%%%%%%%%
%%%%%%%%%%%%%%%%%%%%%%%%%%%%%%%%%%%%%%%%%%%%%%%%%%%%%%%%%%%%%%%%%%%%%%%%%%%%%%%%
\section{Introduction}

The \emph{Wakimoto realization} is a free field realization of the universal affine vertex algebra $V^\kappa(\frg)$ associated to a semisimple Lie algebra $\frg$ at the level $\kappa \in (\Sym^2\frg^\vee)^\frg$.
It embeds $V^\kappa(\frg)$ into a tensor product of the $\beta\gamma$-system $\beta\gamma^{\otimes\dim\frn_+}$ (free fields associated with the nilradical $\frn_+$) and the Heisenberg vertex algebra $\pi^{\kappa - \kappa_c}$ associated with the Cartan subalgebra $\frh$ at the shifted level $\kappa - \kappa_c$, where $\kappa_c = -\frac{1}{2}\kappa_{\frg}$ is the critical level and $\kappa_{\frg}$ is the Killing form.
This realization was originally constructed by Wakimoto \cite{Wak86} for $\frg = \mathfrak{sl}_2$ and was subsequently extended to all semisimple Lie algebras by Feigin--Frenkel \cite{FF90, FF92}.
Free field realizations are extremely useful for understanding complicated vertex algebras, and the Wakimoto realization is a fundamental example of a free field realization.

On the other hand, the notion of \emph{factorization algebra} is introduced by Costello and Gwilliam, is a mathematical formulation of the structure of observables in a field theory.
In \cite{CG, CG2}, they formulated the notion of classical field theory and its perturbative quantization based on the \emph{Batalin--Vilkovisky (BV) formalism} and the algebra of observables in terms of factorization algebras.
In \cite[Chapter 5.3]{CG}, they constructed vertex algebras from prefactorization algebras called $S^1$-equivariant holomorphically translation invariant prefactorization algebras on $\bbC$ (see \cref{FA_to_VA}).
Furthermore, Butson and Yoo \cite{Yoo} introduced the notion of a \emph{degenerate classical field theory}, which generalizes the notion of a classical field theory in the sense of \cite{CG2}.
From a geometric perspective, classical field theories are analogous to symplectic geometry, whereas degenerate classical field theories are analogous to Poisson geometry.
This generalization is motivated by the study of \emph{boundary theories} and \emph{boundary conditions}.

The motivation for this note is to understand the Wakimoto realization of vertex algebras in terms of \emph{holomorphic field theories} on $\bbC$ through the above correspondence.
We believe that it is important to understand the Wakimoto realization in the language of holomorphic field theory.
For example, \cite{CG19,BLL+15} propose a relationship between higher-dimensional field theories and vertex algebras.
Beem, Meneghelli, and Rastelli \cite[Chapter 6]{BMR19} conjectured that vertex algebra constructed from four-dimensional supersymmetric field theories \cite{BLL+15} admits free field realizations.
To construct such free field realizations, it should be necessary to formulate and understand four- or two-dimensional field theories mathematically.
As a well-known example of the two-dimensional perspective, it is useful to understand the Wakimoto realization in the language of holomorphic field theory.
Moreover, since the free field realizations they propose are obtained by embedding a field theory into a free field theory (see \cite[Equation (6.1)]{BMR19}), we would like to understand the Wakimoto realization as an embedding between the factorization algebras of observables of the corresponding field theories.

 %The Lie algebra $\frg$ acts on the big cell $U = N_+ \cdot [B_-] \subset G/B_-$ of the flag variety via the adjoint action, and this action lifts to a classical action of the sheaf of Lie algebras $\clG_1 = \Omega^{0,*}_{\bbC} \otimes \frg$ on 
 %the combined $\beta\gamma$-Heisenberg system $\clW_\kappa = \clM_{\frg} \oplus \clM^{\hbar(\kappa-\kappa_c)}_{Heis}$.
%We construct an equivariant action functional $S^{Wak}_\kappa \in \Act(\clG_1, \clW_\kappa)[\hbar]$---a solution to the equivariant quantum master equation in the BV formalism---whose associated Noether map gives a map of factorization algebras
%\[
%  w^{fact}_{\kappa}: 
%  \bbU^\hbar_{\phi(\kappa)}(\clG_1)[\hbar^{-1}] \lto 
%  \Obs^q_{\clW_\kappa}[\hbar^{-1}],
%\]
%where $\bbU^\hbar_{\phi(\kappa)}(\clG_1)$ is the twisted enveloping factorization algebra of $\clG_1$ (whose associated vertex algebra under the correspondence $\bbV$ is the universal affine vertex algebra $V^\kappa(\frg)\dbr{\hbar}$), and $\Obs^q_{\clW_\kappa}$ is the factorization algebra of quantum observables of the free degenerate classical field theory $\clW_\kappa$.
%Applying the functor $\bbV$ and specializing to $\hbar = 2\pi i$, we recover the Wakimoto realization
%\[
%  w^{ch}_\kappa : V^\kappa(\frg) \lto 
%  \beta\gamma^{\oplus \dim \frn_+} \otimes \pi^{\kappa - \kappa_c}.
%\]

In \cref{s.5}, we develop a higher-dimensional analogue of the Wakimoto realization for the \emph{higher Kac--Moody factorization algebra} $\bbU^\hbar_{\phi(\kappa)}(\clG_d)$ on $\bbC^d$.
The higher Kac--Moody factorization algebra on $\bbC^d$ was introduced in \cite{GW} as a higher-dimensional generalization of the affine Kac--Moody vertex algebra, and our construction gives a natural ``free field realization'' in this setting.
Unfortunately, the constructed map is not injective, so it is not appropriate to call it a ``realization''.
This suggests that it is difficult to obtain a realization by simply extending the corresponding results for vertex algebras.

Gui, Wang, and Williams constructed free field realization using a higher-dimensional chiral algebra in \cite{GWW25}.
This construction is based on higher-dimensional holomorphic field theories, but the relationship between their construction and ours is not yet understood.

\subsection*{Main results}

The goal of this note is to construct a map of factorization algebras that recovers the Wakimoto realization and to extend this construction to arbitrary dimensions. To this end, classical field theories in the sense of \cite{CG, CG2} alone are insufficient to realize the Wakimoto realization at a general level.
Therefore, we develop the notions of ``quantization'' and ``equivariant quantization'' of free degenerate classical field theories proposed in \cite{Yoo}.
As a result, we can quantize the two-dimensional holomorphic field theory introduced in \cite{Yoo} as a boundary theory of three-dimensional Chern--Simons theory, and show that the corresponding vertex algebra is an affine vertex algebra. (see \cref{WZW_VA})
\[
\bbV(\Obs^q_{WZW,\kappa})\simeq V^\kappa_\hbar(\frg).
\]
Furthermore, the vertex algebra naturally constructed from the factorization algebra of quantum observables is a $\bbC\dbr{\hbar}$-linear vertex algebra, and we cannot arbitrarily specialize $\hbar$ merely to avoid infinite sums. Here, we formulate conditions under which infinite sums can be avoided in the language of field theory (see \cref{lem_assumption_satisfies}).

In \cref{ss:4.2}, we introduce a holomorphic free degenerate classical field theory $\clM_\frg$ such that the associated $\bbC\dbr{\hbar}$-linear vertex algebra $\bbV(\Obs^q_{\clM_\frg})$ is isomorphic to the $\bbC\dbr{\hbar}$-linear $\beta\gamma$-system $\beta\gamma_\hbar
^{\otimes\dim\frn_+}$.
Using the embedding of $\frg$ into $\mathcal{D}_{\leq 1}(U)$ introduced in \cref{ss:4.1}, we construct a classical action of $\clG_1=\Omega^{0,*}_{\bbC}\otimes\frg$ on $\clM_\frg$. Moreover, this action admits a lift to a quantum action, and we obtain the following theorem.

\begin{mainthm}[\cref{Wakimoto_critical,Noether_map_Wakimoto_critical_level}]
For a simple Lie algebra $\frg$, there exists an outer-quantum action $\{S^{Wak}_{\kappa_c}[L]\}_{L>0}$ of $\clG_1$ on the quantization of the free degenerate field theory $\clM_{\frg}$ that induces a map of amenably holomorphic factorization algebras
\[
  \wt{w}^{fact}_{\kappa_c}\colon\mathbb{U}^\hbar_{\phi(\kappa_c)}(\clG_1)[\hbar^{-1}]\lto\Obs^q_{\clM_\frg}[\hbar^{-1}].
\]
Moreover, the induced map of vertex algebras
\[
  \wt{w}^{ch}_{\kappa_c}:V^{\kappa_c}(\frg)\lto\beta\gamma^{\otimes\dim\frn_+}
\]
under the functor $\bbV$ and specialization $\hbar=2\pi i$, coincides with the Wakimoto realization at critical level $\kappa_c$.
\end{mainthm}

Furthermore, for a 2-cocycle $\kappa\in\Sym^2(\frg^\vee)^\frg$, the action of $\clG_1$ can be extended to the classical field theory $\clW_\kappa$ defined by extending $\clM_\frg$, and we obtain the following theorem.

\begin{mainthm}[\cref{Wakimoto_general_level_FA,Noether_map_Wakimoto_general_level}]
For a simple Lie algebra $\frg$ and a 2-cocyle $\kappa\in\Sym^2(\frg^\vee)^\frg$, there exists an outer-quantum action $\{S^{Wak}_{\kappa}[L]\}_{L>0}$ of $\clG_1$ on the quantization of the free degenerate field theory $\clW_{\kappa}$ that induces a map of amenably holomorphic factorization algebras
\[
  w^{fact}_{\kappa}\colon\mathbb{U}^\hbar_{\phi(\kappa)}(\clG_1)[\hbar^{-1}]\lto\Obs^q_{\clW_\kappa}[\hbar^{-1}].
\]
Moreover, the induced map of vertex algebras
\[
  w^{ch}_{\kappa}:V^{\kappa}(\frg)\lto\beta\gamma^{\otimes\dim\frn_+}\otimes\pi^{\kappa-\kappa_c}
\]
under the functor $\bbV$ and specialization $\hbar=2\pi i$, coincides with the Wakimoto realization at level $\kappa$.
\end{mainthm}

By generalizing the above argument on $\bbC$, we obtain a higher-dimensional analogue.
The main results of this note can be summarized in the following table.

\begin{table}[htbp]
  \centering
  \label{tab:Wakimoto}
  \begin{tabular}{c||cc|cc}
  \hline \hline
  &\multicolumn{2}{|c|}{Critical level}&\multicolumn{2}{|c}{General level}\\ \hline
  Space-time&Existence&Injectivity&Existence&Injectivity\\ \hline
  $\bbC$
  &\cref{Noether_map_Wakimoto_critical_level}
  &\cref{not_inj_critical_level}
  &\cref{Noether_map_Wakimoto_general_level}
  &\cref{thm_inj_C}\\
  $\bbC^{2n}$
  &\cref{Noether_map_higher_Wakimoto_critical_level}
  &\cref{kernel_of_higher_Wakimoto_critical}&×&×\\
  \multirow{2}{*}{$\bbC^{2n+1}$}
  &\multirow{2}{*}{\cref{Noether_map_higher_Wakimoto_critical_level}}&\multirow{2}{*}{\cref{kernel_of_higher_Wakimoto_critical}}&\multirow{2}{*}{\cref{Noether_map_higher_Wakimoto_general_level}}&\cref{not_inj_ngeq2} ($n\ge2$)\\
  &&&&\cref{conj_n_1_inj} ($n=1$)
\end{tabular}
\end{table}

We show that the Wakimoto realization can be understood as a map between factorization algebras of observables (see \cref{BD_FA_morphism_Wakimoto}).
This justifies calling the Wakimoto realization a “\emph{free field realization}” in the sense of quantum field theory.

%%%%%%%%%%%%%%%%%%%%%%%%%%%%%%%%%%%%%%%%%%%%%%%%%%%%%%%%%%%%%%%%%%%%%%%%%%%%%%%%
\subsection*{Organization}

The paper is organized as follows.

In \cref{s.2}, we develop the general theory of free degenerate classical field theories.
In \cref{ss:2.2}, we introduce the notion of factorization algebra of quantum observables for free degenerate classical field theory.
In \cref{ss:2.3}, following Part III of \cite{CG2}, we introduce an action of an $L_\infty$-algebra on a degenerate classical field theory.
We show that this action induces a map between factorization algebras (see \cref{quantum_Noether_map}).

In \cref{s.3}, we recall the notion of vertex algebras and the construction of vertex algebra via amenably holomorphic prefactorization algebra on $\bbC$.
In \cref{ss:3.3}, we review the construction of the $\beta\gamma$-system and the universal affine vertex algebra due to \cite{CG}, and prepare the basic tools concerning compactly supported Dolbeault cohomology.
In \cref{ss:3.4} and \cref{ss:3.5}, we compute the vertex algebra constructed from the factorization algebra of observables of a free degenerate classical field theory.
In \cref{ss:3.6}, we introduce the notion of a BD center of factorization algebra and the associated vertex algebras.

\cref{s.4} is the main part of this note.
In \cref{ss:4.2} and \cref{ss:4.3}, we reconstruct the Wakimoto realizations using the terminology of the equivariant field theory in \cref{ss:2.3}.
In \cref{ss:4.4}, we show that the injectivity of the Wakimoto realization also holds at the level of factorization algebras. 
In \cref{ss:4.5}, we show that the map of factorization algebras constructed in \cref{ss:4.2} and \cref{ss:4.3} lifts to a map between the factorization algebras of quantum observables associated with free degenerate classical field theories.
By comparing their BD centers, we construct a factorization-algebra-level free field realization of the Feigin–Frenkel center.

In \cref{s.5}, we construct the higher dimensional analogue of the construction in \cref{s.4}.
In \cref{ss:5.2}, as a natural higher-dimensional analogue of \cref{ss:4.2}, we construct an analogue of the Wakimoto-type map on complex affine space of arbitrary dimension at the critical level.
In \cref{ss:5.3}, as a natural higher-dimensional analogue of \cref{ss:4.3}, we construct a Wakimoto-type map on complex affine spaces of odd complex dimension. However, we show that this map is never injective in complex dimensions five and higher.

%%%%%%%%%%%%%%%%%%%%%%%%%%%%%%%%%%%%%%%%%%%%%%%%%%%%%%%%%%%%%%%%%%%%%%%%%%%%%%%%
\subsection*{Notations and terminology}
We summarize here the notation and terminology used throughout the main text.
\begin{itemize}
  \item The symbols $\bbZ$ and $\bbC$ denote the set of all integers and complex numbers, respectively.
  \item We denote quasi-isomorphisms by $\qis$ or $\qim$.
  \item Unless otherwise specified, vector spaces are defined over the field $\bbC$, and $\Hom,\End$, and $\otimes$ are understood in this category.
   \item For a vector space $V$, we denote the linear dual $\Hom(V,\bbC)$ by $V^\vee$.
  \item For a Lie algebra $\frg$, we denote by $U(\frg)$ the universal enveloping algebra of $\frg$.
  \item For a topological space $X$, we denote by $\Open(X)$ the category of open subsets of $X$ and open immersions.
  \item For $n\geq1,\ z\in\bbR^n$, and $R>0$, we denote by $D(z,R)\subset\bbR^n$ the open disk of radius $R$ centered at $z$.
  \item We follow \cite{CG} for the terminology of differentiable vector spaces and differentiable (pro-)cochain complexes.
  \item We denote by $\mathrm{DVS}$ the category of differentiable vector spaces and by $\mathrm{pCh(DVS)}$ the category of differentiable pro-cochain complexes (see Appendix C.4 of \cite{CG}).
  \item We denote by $\Sym$ the completed symmetric product in the category $\mathrm{pCh(DVS)}$.
  \item For a vector bundle $E$ on a manifold $X$, we denote by $\clE$ the sheaf of the smooth sections, $\ol{\clE}$ the sheaf of the distributional sections, $\clE_c$ the cosheaf of the compactly supported smooth sections, and $E^!=\mathrm{Dens}_X\otimes_{C^\infty_X} E^\vee$ the ``Verdier dual'' of $E$, where $\mathrm{Dens}_X$ is the density bundle on $X$ (see \S 3.5.5 of \cite{CG}).
  \item For a vector bundle $E$ on a manifold $X$, its continuous dual is denoted by $\clE^*=\ol{\clE}_c^!$ (see Chapter 3.5.5 of \cite{CG}).
\end{itemize}

%%%%%%%%%%%%%%%%%%%%%%%%%%%%%%%%%%%%%%%%%%%%%%%%%%%%%%%%%%%%%%%%%%%%%%%%%%%%%%%%
%%%%%%%%%%%%%%%%%%%%%%%%%%%%%%%%%%%%%%%%%%%%%%%%%%%%%%%%%%%%%%%%%%%%%%%%%%%%%%%%
\section{The BV formalism}
\label{s.2}

In this section, we prepare one of the main tools that will be used throughout this note.
In \cite{Yoo}, the notion of a degenerate classical field theory was introduced.
This notion is a generalization of the notion of a classical field theory defined in \cite{CG2}.
The classical field theory in \cite{CG2} is defined to be equipped with a ``symplectic structure'', whereas a degenerate classical field theory is its ``Poisson'' analogue.
The main part of this section is to define the notion of quantization of a free degenerate classical field theory.
More precisely, in analogy with Chapter 8 of \cite{CG2}, we construct the BD factorization algebra $\Obs^q_{\clM}$ for a free degenerate classical field theory $(\clM,\Pi,Q_{BV})$ (see \cref{FA_str_on_Obs}).

In \cref{ss:2.3}, we define the notion of action of $L_\infty$-algebra on strictly free degenerate field theory.
The goal of this subsection is to show that, given such an action, it induces a map from the factorization algebra defined in \cref{twisted_enveloping_factorization_algebra} to the factorization algebra of observables defined in \cref{ss:2.2} (see \cref{quantum_Noether_map}).
In \cref{s.4} and \cref{s.5}, we use \cref{quantum_Noether_map} to construct a map of factorization algebras.

%%%%%%%%%%%%%%%%%%%%%%%%%%%%%%%%%%%%%%%%%%%%%%%%%%%%%%%%%%%%%%%%%%%%%%%%%%%%%%%%
\subsection{Factorization algebra}
\label{ss:2.1}

In this subsection, we review the notion of a factorization algebra defined in \cite{CG}.
Since in this note, we only consider the factorization algebra valued in either $\mathrm{pCh(DVS)}$ or $\mathrm{DVS}$, in this subsection, let $\mathrm{C}$ denote either $\mathrm{pCh(DVS)}$ or $\mathrm{DVS}$.
However, the results in this subsection should continue to hold if $\mathrm{C}$ is replaced by a general symmetric monoidal category which admits coproduct.
We denote $1_\mathrm{C}\in\mathrm{C}$ by the unit object of the symmetric monoidal category $\mathrm{C}$.

\subsubsection{Basic definitions for factorization algebra}

In this subsection, we review the basic definitions of (pre)factorization algebras.
Our main reference is \cite[Chapter 3]{CG}.

\begin{dfn}[{\cite[Chapter 3.1.1]{CG}}]
\label{def_PFA}
Let $X$ be a topological space.
A \emph{prefactorization algebra} is a tuple $(\clF,m_{\clF},\eta)$ consisting of
\begin{itemize}
\item 
an assignment $\clF(U)\in\mathrm{C}$, for every open subset $U\subset X$,
\item 
an assignment $(m_{\clF})^{U}_V:\clF(U)\to\clF(V)\in\Hom_{\mathrm{C}}(\clF(U),\clF(V))$, for every open subsets $U,V\subset X$ such that $U\subset V$,
\item 
an assignment $(m_{\clF})^{U,V}_W:\clF(U)\otimes\clF(V)\to\clF(W)$, for every open subsets $U,V,W\subset X$ such that $U\sqcup V\subset W$, called the \emph{factorization product},
\item 
a map $\eta:1_\mathrm{C}\to\clF(\emptyset)$,
\end{itemize}
which satisfies the following axioms:
\begin{clist}
\item 
The pair $(\clF,(m_{\clF})^\bullet_\bullet)$ is a precosheaf.
\item 
For open subsets $U_1,U_2,V,W\subset X$ such that $U_1\sqcup U_2\subset V\subset W$, we have
\[
  (m_{\clF})^{V}_W\circ(m_\mathfrak{F})^{U_1,U_2}_V=(m_{\clF})^{U_1,U_2}_W.
\]
Also, for open subsets $U_1,U_2,V_1,V_2,W\subset X$ such that
\[
  U_1\subset V_1,\quad U_2\subset V_2,\quad V_1\sqcup V_2\subset W,
\]
we have
\[
  (m_{\clF})^{V_1,V_2}_{W}\circ\left((m_{\clF})^{U_1}_{V_1}\otimes(m_{\clF})^{U_2}_{V_2}\right)=(m_{\clF})^{U_1,U_2}_W.
\]
\item 
For open subsets $U,V,W\subset X$ such that $U\sqcup V\subset W$, 
the following diagram commutes:
\[
  \begin{tikzcd}
   \clF(U)\otimes\clF(V) 
   \arrow[rr, "\sim"] \arrow[rd, "(m_{\clF})^{U,V}_W"]
   & & \clF(V)\otimes\clF(U) \\
   & \clF(W) \arrow[ru, leftarrow, "(m_{\clF})^{V,U}_W"] 
  \end{tikzcd}
\]
\item 
For open subsets $U_1,U_2,U_3,V_1,V_2,W\subset X$ such that
\[
  U_1\sqcup U_2\subset V_1,\quad U_2\sqcup U_3\subset V_2,\quad U_1\sqcup V_2\subset W,\quad U_3\sqcup V_1\subset W,
\]
the following diagram commutes:
\[
\begin{tikzcd}
(\clF(U_1)\otimes\clF(U_2))\otimes\clF(U_3) \arrow[rr, "\sim"] \arrow[d, "{(m_{\clF})^{U_1,U_2}_{V_1}\otimes1}"] &                & \clF(U_1)\otimes(\clF(U_2)\otimes\clF(U_3)) \arrow[d, "{1\otimes(m_{\clF})^{U_2,U_3}_{V_2}}"] \\
\clF(V_1)\otimes\clF(U_3) \arrow[rd, "{(m_{\clF})^{V_1,U_3}_W}"]                                                        &                & \clF(U_1)\otimes\clF(V_2) \arrow[ld, "{(m_{\clF})^{U_1,V_2}_W}"]                                     \\
                                                        & \clF(W) &  
\end{tikzcd}
\]
\item 
For an open subset $U\subset X$, the following diagram commutes:
\[
  \begin{tikzcd}
      1_{\mathrm{C}}\otimes\clF(U) \arrow[rr, "\eta\otimes1"] \arrow[rd, "\sim"]
        & & \clF(\emptyset)\otimes\clF(U) \\
      & \clF(U) \arrow[ru, leftarrow, "(m_{\clF})^{\emptyset,U}_U"] 
  \end{tikzcd}
\]
\end{clist}
\end{dfn}

By an abuse of notation, we denote the prefactorization algebra $(\clF,m_{\clF},\eta)$ simply by $\clF$.

\begin{dfn}[{\cite[Definition 3.1.6]{CG}}]
\label{def_PFA_mor}
Let $\clF,\clG:\Open(X)\to\mathrm{C}$ be prefactorization algebras on a topological space $X$.
A map $\varphi:\clF\to\clG$ of precosheaves is called a map of prefactorization algebra if the following conditions hold:
\begin{itemize}
\item 
For open subsets $U,V,W\subset X$ such that $U\sqcup V\subset W$, we have
\[
  \varphi_W\circ(m_{\clF})^{U,V}_W=(m_{\clG})^{U,V}_W\circ(\varphi_U\otimes\varphi_V).
\]
\item 
We have $\varphi_\emptyset\circ\eta_{\clF}=\eta_{\clG}$.
\end{itemize}
\end{dfn}

\begin{dfn}[{\cite[Definition 6.3.6]{CG}}]
\label{def_Weiss_cover}
Let $X$ be a topological space and $U\subset X$ be an open subset.
An open cover $\{U_i\}_{i\in I}$ of $U$ is a \emph{Weiss cover} if for finitely many points $x_1,\dots,x_n\in U$, there exists $i\in I$ such that $x_1,\dots,x_n\in U$.
\end{dfn}

\begin{dfn}[{\cite[Definition 6.1.3]{CG}}]
\label{def_FA}
A prefactorization algebra $\clF:\Open(X)\to\mathrm{C}$ is 
a \emph{factorization algebra} if it satisfies the following conditions:
\begin{clist}
\item 
For an open subset $U\subset X$ and Weiss cover $\{U_i\}_{i\in I}$ of $U$, 
the following diagram is a coequalizer diagram.
\[
  \begin{tikzcd}
   \bigoplus_{i,j\in I}\clF(U_i\cap U_j)
   \arrow[r, shift left = 0.2em]
   \arrow[r, shift right = 0.2em, swap]
  &\bigoplus_{i\in I}\clF(U_i)
   \arrow[r]
  &\clF(U)
  \end{tikzcd}
\]
\item 
For open subsets $U,V\subset X$ such that $U\cap V=\emptyset$, 
the factorization product
\[
  (m_{\clF})^{U,V}_{U\sqcup V}: \clF(U)\otimes\clF(V)\lto\clF(U\sqcup V)
\]
is an isomorphism.
\end{clist}
\end{dfn}

%%%%%%%%%%%%%%%%%%%%%%%%%%%%%%%%%%%%%%%%%%%%%%%%%%%%%%%%%%%%%%%%%%%%%%%%%%%%%%%%
\subsubsection{Equivariant (pre)factorization algebras}

Suppose that a group $G$ acts on a topological space $X$.
Then, for an element $a\in G$ and a (pre)factorization algebra $\clF$ on $X$, 
we have a precosheaf
\begin{align*}
  a\clF: \Open(X)\lto\mathrm{C}, \quad U\lmto\clF(aU).
\end{align*}

\begin{dfn}[{\cite[Definition 3.7.1]{CG}}]
\label{def_equivariant_PFA}
A (pre)factorization algebra $\clF$ on $X$ equipped with a map
\[
  \sigma_a=\{\sigma_{a,U}\}_{U\in\Open(X)}: \clF \lto a\clF
\]
of precosheaves for each $a\in G$ is called 
a \emph{$G$-equivariant (pre)factorization algebra} 
if the following conditions hold:
\begin{itemize}
\item 
For $a,b\in G$ and open subset $U\subset X$, 
we have $\sigma_{ab,U}=\sigma_{a,bU}\circ\sigma_{b,U}$.
\item 
For the unit $e\in G$ and an open subset $U\subset X$, 
we have $\sigma_{e,U}=\mathrm{id}_{\clF(U)}$.
\item 
For $a\in G$ and open subsets $U,V,W\subset X$ such that $U\sqcup V\subset W$,
we have
\[
  \sigma_{a,W}\circ(m_{\clF})^{U,V}_W = 
  (m_{\clF})^{aU,aV}_{aW}\circ(\sigma_{a,U}\otimes\sigma_{a,V}).
\]
\item 
For $a\in G$, we have $\sigma_{a,\emptyset}\circ\eta_{\clF}=\eta_{\clF}$.
\end{itemize}
\end{dfn}

The following notion plays an important role on this note.

\begin{dfn}[{\cite[Chapter A.3.2]{CG}}]
\label{BD_algebra}
A tuple $(A,d,\{-,-\})$ is a \emph{Beilinson--Drinfeld (BD) algebra} 
if it satisfies the following conditions:
\begin{itemize}
\item 
$(A,d)$ is a differential graded commutative algebra over $\bbC\dbr{\hbar}$.
\item 
The $\bbC\dbr{\hbar}$-bilinear map $\{-,-\}:A\times A\to A$ 
defines a Poisson bracket on $A$.
\item 
For elements $a,b\in A$ of homogeneous degree $|a|,|b|\in\bbZ$, 
the following holds.
\[
  d(ab)=d(a)b+(-1)^{|a|}ad(b)+\hbar\{a,b\}.
\]
\end{itemize}
\end{dfn}

One of the main objects of this note is 
a factorization algebra equipped with a BD algebra structure.
This notion is introduced informally in \cite{CG2}, 
but a precise formulation should be as follows.

\begin{dfn}
\label{def_BD_FA}
A tuple $(\clF,m_{\clF},\eta,d=(d_U)_{U\in\Open(X)},\{-,-\}=(\{-,-\}_U)_{U\in\Open(X)})$ is a \emph{Beilinson--Drinfeld (BD) factorization algebra} if it satisfies the following conditions:
\begin{clist}
\item 
For an open subset $(\clF(U),d_U,\{-,-\}_U)$ is a BD algebra.
\item 
A tuple $(\clF,m_{\clF},\eta)$ is a factorization algebra on $X$.
\item 
For open subsets $U,V\subset X$ such that $U\subset V$, we have
\[
  (m_{\clF})^U_V \circ \{-,-\}_U = 
  \{-,-\}_V \circ ((m_{\clF})^U_V \otimes (m_{\clF})^U_V).
\]
\end{clist}
\end{dfn}

Let $(\clF,d^{\clF},\{-,-\}^{\clF}),(\clG,d^{\clG},\{-,-\}^{\clG})$ be 
BD factorization algebras on a topological space $X$, 
and $\varphi:\clF\to\clG$ be a map of factorization algebras.
Then the map $\varphi$ is a \emph{map of BD factorization algebras} 
if for an open subset we have
\[
  \varphi_U \circ \{-,-\}^{\clF}_U = 
  \{-,-\}^{\clG}_U \circ (\varphi_U\otimes\varphi_U).
\]

Recall the notion of \emph{local $L_\infty$-algebra} 
on a manifold \cite[\S3.1]{CG2}.
This is a fundamental algebraic concept used both in 
defining actions of $L_\infty$-algebras and in defining field theories.

\begin{dfn}[{\cite[Definition 3.1.3]{CG2}}]
\label{def: local L_infty algebra}
Let $X$ be a manifold.
A \emph{local $L_\infty$-algebra} on $X$ consists of the following data.
\begin{itemize}
\item 
A graded vector bundle $L$ on $X$, whose space of smooth sections will be denoted $\clL$.
\item 
A collection of polydifferential operators $l_n: \clL^{\otimes n}\to\clL$ for $n\geq1$, which are of cohomological degree $2-n$, and endow $\clL$ with the structure of $L_\infty$-algebra.
\end{itemize}
We call a local $L_\infty$-algebra on $X$ a \emph{local dgLA} if for $n\ge3$, the polydifferential operator $l_n:\clL^{\otimes n}\to\clL$ vanishes.
\end{dfn}

For a manifold $X$, let $C^\infty_X$ be the sheaf of smooth functions on $X$, $\mathrm{Dens}_X$ be the sheaf of smooth densities, and $\clD_X$ be the sheaf of differential operators with smooth coefficients.
Also, for a vector bundle $E$ on $X$, we denote by $J(E)$ the jet bundle of $E$, and $\clJ(\clE)$ the sheaf of smooth sections.
Then, the sheaf  $\clJ(\clE)$ has the natural $\clD_X$-module and $C^\infty_X$-module structure.

Let $\clL$ be a local $L_\infty$-algebra on $X$ and $\clE$ be a sheaf on $X$ equipped with a local action of $\clL$ on $\clE$, 
i.e.  we have polydifferential operators
\[
  \clL^{\otimes n}\otimes\clE \lto \clE \quad (n\ge1)
\]
that are compatible with the $L_\infty$-algebra structure on $\clL$ (see Definition 3.4.2 of \cite{CG2}).
Then, for an open subset $U\subset X$, define
\[
  C^*_{loc}(\clL,\clE)(U) := \prod_{n\ge0}
  \Hom_{\clD_X(U)}\bigl(\clJ(\clL)(U)^{\otimes n},\clJ(\clE)(U)\bigr)_{S_n},
\]
where the subscript $S_n$ denotes the coinvariants with respect to the natural symmetric group actions (see Definition 3.5.9 of \cite{CG2}).

For an $L_\infty$-algebra $(L,\{l_k\})$, we denote by $C^{\mathrm{CE}}_\bullet(L)$ the Chevalley--Eilenberg homological complex, i.e.
\[
  C^{\mathrm{CE}}_\bullet(L) := \Sym(L[1]),
\]
with differential $d$ whose restriction to $\Sym^k(L[1])\to L[1]$ are precisely the higher brackets $l_k$.

The following notion will be important when considering the Noether map.
Physical conserved quantities are described by a map from the following factorization algebra (see \cref{quantum_Noether_map}).

\begin{dfn}[{\cite[Chapter 3.4]{CG}}]
\label{twisted_enveloping_factorization_algebra}
Let $X$ be a manifold and $\clL$ be a local dgLA on $X$.
A degree $k+2$ cocycle $\alpha\in C^*_{loc}(\clL)$ determines a $k$-shifted central extension
\[
  0 \lto \underline{\bbC}[k] \lto \wt{\clL}_c \lto \clL_c \lto 0
\]
of precosheaves of $L_\infty$-algebras, where
\begin{itemize}
  \item $l_1$ and $l_2$ are defined by the $L_\infty$-brackets on $\clL_c$.
  \item Add $\alpha$ to the $L_\infty$-brackets of $\clL_c$
\end{itemize}

We define the factorization algebra on $X$ by
\begin{align}\label{eq:U_alpha}
  \bbU_\alpha^{cl}\clL: \left(X\supset V \lmto C_\bullet^{CE}(\wt{\clL}_c(V))\otimes_{\bbC[c]}\bbC_{c=1}\right),
\end{align}
where $c$ is a generator of $\bbC[k]$.
Here, the factorization product on $\bbU^{cl}_\alpha(\clL)$ is defined by the symmetric product.

We denote the associated twisted enveloping factorization algebra on $X$ by
\begin{equation}
\label{eq:U_alphahbar}
  \bbU^\hbar_\alpha\clL:
  \left(X \supset V \lmto C_\bullet^{CE}
        (\wt{\clL}_c(V)\dbr{\hbar}) \otimes_{\bbC\dbr{\hbar}[c]}
        \bbC\dbr{\hbar}_{c=1}\right).
\end{equation}
\end{dfn}

%%%%%%%%%%%%%%%%%%%%%%%%%%%%%%%%%%%%%%%%%%%%%%%%%%%%%%%%%%%%%%%%%%%%%%%%%%%%%%%%
\subsection{Degenerate field theory}
\label{ss:2.2}

First, we review the notion of degenerate classical field theories introduced
in \cite{Yoo}.
Although they formulate the notion in the language of \emph{formal moduli problems}, we describe it here using only the language of $L_\infty$-algebras.

Let $E_1,\dots,E_n,F$ be vector bundles on a manifold $X$, we denote by $\PD(\clE_1\otimes\cdots\otimes\clE_n,\clF)$ the sheaf of polydifferential operators from $\clE_1\otimes\cdots\otimes\clE_n$ to $\clF$ on $X$.
By Lemma 3.5.2 of \cite{CG2}, there exists a natural bijection
\[
\PD(\clE_1\otimes\cdots\otimes\clE_n,\clF)(U)\simeq\Hom_{\clD_X}(\clJ(E_1)(U)\otimes\cdots\otimes\clJ(E_n)(U),\clJ(F)(U)),
\]
for an open subset $U\subset X$.

\begin{dfn}[{\cite[Definition 2.28]{Yoo}}]
\label{polyvector_field_def}
Let $\clM$ be a local $L_\infty$-algebra on a manifold $X$.
Then a sheaf of \emph{$(-1)$-shifted multilocal $j$-polyvector fields} is
\[
\Gamma_{mloc}(\clM,\Sym^j\bbT_{\clM}) \ceq \prod_{k\in\bbZ_{\geq0}}\PD(\clM[1]^{\otimes k}\otimes\clM^![-1]^{\otimes(j-1)},\clM[1])_{S_k\times S_j}
\]
equipped with Chevalley--Eilenberg differential defined below, where the subscript $S_k\times S_j$ denotes the coinvariants with respect to the natural symmetric group actions.
\end{dfn}
%We denote the factors in the direct product defining $\Gamma_{mloc}(\clM,\Sym^j\bbT_{\clM})$ by $\Gamma_{mloc}^k(\clM,\Sym^j\bbT_{\clM})$.

We denote the space of all $(-1)$-shifted multilocal polyvector fields, and its completion, by
\[
\Gamma_{mloc}(\clM,\Sym\bbT_{\clM}) \ceq \bigoplus_{j\geq0}\Gamma_{mloc}(\clM,\Sym^j\bbT_{\clM}),
\]
\[
\ \Gamma_{mloc}(\clM,\widehat{\Sym}\bbT_{\clM}) \ceq \prod_{j\geq0}\Gamma_{mloc}(\clM,\Sym^j\bbT_{\clM}).
\]
In order to define the notion of a multilocal Poisson bivector and to introduce the Chevalley--Eilenberg differential on all the spaces of multilocal polyvector fields, we will need to discuss the Schouten bracket of multilocal vector fields.

Recall the definition of the usual
Schouten--Nijenhuis bracket.

Let $X$ be a smooth manifold. For $k \geq -1$, define
\[
T^k_{\text{poly}}(X) := \Gamma\bigl(X, \wedge^{k+1} TX\bigr).
\]
We further define a graded vector space $\frg = \bigoplus_{k\in\bbZ} \frg^k$ by
\[
\frg^k :=
\begin{cases}
T^k_{\text{poly}}(X) & (k \geq -1), \\
0 & (k < -1).
\end{cases}
\]

The Schouten--Nijenhuis bracket on $\frg$ is the bilinear map
\[
  [-,-]_{\text{SN}} \colon \frg^k \otimes \frg^l \lto \frg^{k+l}
\]
defined as follows.
For $\xi_s, \eta_t \in \Gamma(M,TM)$ and $h \in C^\infty(M)$, we set
\begin{align*}
  [\xi_0 \wedge \cdots \wedge \xi_k,\,
   \eta_0 \wedge \cdots \wedge \eta_l]_{\text{SN}}
  &=
  \sum_{s=0}^k \sum_{t=0}^l (-1)^{s+t}
  [\xi_s,\eta_t] \wedge
  \xi_0 \wedge \cdots \wedge \widehat{\xi_s} \wedge \cdots \wedge \xi_k \\
&\qquad\qquad\wedge
  \eta_0 \wedge \cdots \wedge \widehat{\eta_t} \wedge \cdots \wedge \eta_l, \\[0.5em]
  [\xi_0 \wedge \cdots \wedge \xi_k,\,h]_{\text{SN}}
  &=
  \sum_{s=0}^k (-1)^s \, \xi_s(h)\,
  \xi_0 \wedge \cdots \wedge \widehat{\xi_s} \wedge \cdots \wedge \xi_k .
\end{align*}
Similarly to the above construction, for differential operators $D_i:\clM\to\clM,\ D_i':\clM\to\clM$ and $E_i:\clM^!\to\clM,\ E_i'\in\clM^!\to\clM$, we define
\begin{align*}
  &[D_0\otimes\cdots\otimes D_n\otimes E_1\otimes\cdots\otimes E_{j-1},D_0'\otimes\cdots\otimes D_m'\otimes E_1'\otimes\cdots\otimes E_{l-1}']\\
  &\ceq \sum_{s=0}^n\sum_{t=0}^m(-1)^{s+t}[D_s,D_t']\otimes D_0\otimes\cdots\widehat{D_s}\otimes\cdots\otimes D_n\otimes D_0'\otimes\cdots\otimes\widehat{D_t'}\otimes\cdots\otimes D_m'\\
  &\hspace{7em} 
   \otimes E_1\otimes\cdots\otimes E_{j-1}\otimes E_1'\otimes\cdots\otimes E_{l-1}',\\
  &[D_0\otimes\cdots\otimes D_n\otimes E_1\otimes\cdots\otimes E_{j-1},E_1'\otimes\cdots\otimes E_{l-1}']\\
  & \ceq D_0\otimes\cdots\otimes D_n\otimes E_1\otimes\cdots\otimes E_{j-1}\otimes E_1'\otimes\cdots\otimes E_{l-1}',
\end{align*}
where $[D_s,D_t'] \ceq D_s\circ D_t'-D_t'\circ D_s$ is the commutator of differential operators.

\begin{prop}[{\cite[Proposition 2.29]{Yoo}}]
\label{Schouten_bracket_def}
 The above map $[-,-]$ of sheaves of multilocal polyvector spaces on a manifold $X$
\[
  [-,-]:\Gamma_{mloc}(\clM,\Sym^j\bbT_{\clM})\otimes\Gamma_{mloc}(\clM,\Sym^l\bbT_{\clM})\lto\Gamma_{mloc}(\clM,\Sym^{j+l-1}\bbT_{\clM})
\]
makes $\Gamma_{mloc}(\clM,\Sym\bbT_{\clM})$ a sheaf of graded Lie algebras.
\end{prop}

Note that, as in the case of usual polyvector fields, 
$\Gamma_{mloc}(\clM,\Sym\bbT_{\clM})$ forms a Lie subalgebra of $\Gamma_{mloc}(\clM,\Sym\bbT_{\clM})$ which acts on $\Gamma_{mloc}(\clM,\Sym^j\bbT_{\clM})$ for each $j\geq0$.
Furthermore, by  definition, 
there is an identification of sheaves $C^\infty_M$-modules
\[
  \Gamma_{mloc}(\clM,\Sym\bbT_{\clM})=C^*_{loc}(\clM,\clM[1]).
\]
The above identification gives a graded Lie algebra structure on $C^*_{loc}(\clM,\clM[1])$.

The following is a reformulation of the description in terms of action functionals given in Chapter 4 of \cite{CG2} in the language of vector fields:
\begin{prop}[{\cite[Definition 2.25]{Yoo}}]
\label{equivalence_of_CLFT}
Let $\clM$ be a local $L_\infty$-algebra on $X$.
Then, the following data are equivalent:
\begin{itemize}
\item A collection of polydifferential operators
\[
  \{l_n:\clM^{\otimes n}\lto\clM[2-n]\}_{n\geq1}
\]
defines a local $L_\infty$-algebra structure on $\clM$.
\item The vector field
\[
  Q_{\clM}\in C^*_{loc}(\clM(X),\clM(X)[2])
\]
whose degree 0 term vanishes and satisfies $[Q_{\clM},Q_{\clM}]=0$, where $[-,-]$ denote by the Lie bracket induced by the Schouten bracket on $\Gamma_{mloc}(\clM,\Sym\bbT_{\clM})$  (\cref{Schouten_bracket_def}).
\end{itemize}
\end{prop}

In particular, under this identification, the equivalence of
\cref{equivalence_of_CLFT} determines a vector field
\[
  Q_{\clM}\in C^*_{loc}\bigl(\clM,\clM[2]\bigr)
\]
which encodes the structure of a local $L_\infty$-algebra on $\clM$.

For each $k \geq 0$, the adjoint action induced by
$Q_{\clM}$ on
$\Gamma_{mloc}\bigl(\clM, \Sym^k \bbT_{\clM}\bigr)$
defines a differential of cohomological degree $1$.
We define the resulting differential on
$\Gamma_{mloc}\bigl(\clM, \Sym^\ast \bbT_{\clM}\bigr)$
to be the Chevalley--Eilenberg differential.

In keeping with the usual interpretation of the Schouten bracket, we use the notation
\[
\clL_Q\xi \ceq [Q,\xi]
\]
for $Q\in\Gamma_{mloc}(\clM,\Sym\bbT_{\clM})$, and  $\xi\in\Gamma_{mloc}(\clM,\Sym^k\bbT_{\clM})$.
\begin{dfn}[{\cite[Definition 2.30]{Yoo}}]
\label{Poisson_structure_def}
Let $\clM$ be a local $L_\infty$-algebra on a manifold $X$.
\emph{A $(-1)$-shifted Poisson structure on $\clM$} is a cohomological degree 1 multilocal bivector field
\[
\Pi\in\Gamma_{mloc}(\clM,\Sym^2\bbT_{\clM})[1]
\]
such that
\[
\clL_{Q_{\clM}}\Pi=0,\ [\Pi,\Pi]=0
\]
\end{dfn}
Let $\clM$ be a local $L_\infty$-algebra on a manifold $X$, and $\Pi\in\Gamma_{mloc}(\clM,\Sym^2\bbT_{\clM})[1]$ be a $(-1)$-shifted Poisson structure on $\clM$.
Then we denote by $\Pi^{(n)}$ the component 
\begin{align}
    \label{equ:componentofPoisson}
     \Pi^{(n)}:\clM[1]^{\otimes n}\otimes\clM^![-1]\lto\clM[1],
\end{align}
and call it the \emph{$n$-component} of the $(-1)$-shifted Poisson structure $\Pi$.

\begin{rem}
\label{Poisson_symplectic_generalization}
This notion is a generalization of a $(-1)$-shifted symplectic structure, and hence a generalization of a classical field theory.
Indeed, it is immediate that a $(-1)$-shifted symplectic structure canonically determines a local $(-1)$-shifted Poisson structure.
In particular, a $(-1)$-shifted symplectic structure is concentrated in the $0$-component.
\end{rem}

\begin{dfn}[{\cite[Definition 2.34]{Yoo}}]
\label{degenerate_classical_field_theory_def}
A \emph{degenerate classical field theory} on a manifold $X$ is a triple $(\clM,\Pi,\{l_k\})$ consisting of
\begin{itemize}
  \item a local $L_\infty$-algebra $(\clM,\{l_k\})$ on $X$ such that $l_1:\clM\to\clM$ is an elliptic operator,
  \item a $(-1)$-shifted Poisson structure $\Pi\in\Gamma_{mloc}(\clM,\Sym^2\clM)[1]$ on $\clM$.
\end{itemize}
A degenerate classical field theory is \emph{free} if for $k\geq2$, the $l_k$-bracket is zero.
A free degenerate classical field theory is \emph{strictly free} if the $(-1)$-shifted Poisson structure $\Pi$ belongs to the subsheaf $\PD(\clM^![-1],\clM[1])_{S_2} \subset\Gamma_{mloc}(\clM,\Sym^2\clM)$.
\end{dfn}
If a degenerate classical field theory $\bigl(\clM,\Pi,\{l_k\}\bigr)$ is free, then we simply denote it by $\bigl(\clM,\Pi,l_1=Q_{\mathrm{BV}}\bigr)$.

With these preparations, we define one of the main objects of this note.
\begin{dfn}
\label{classical_observables_of_degenerate_classical_field_theory}
Let $(\clM,\Pi,\{l_k\})$ be a free degenerate classical field theory on a manifold $X$.
Then the factorization algebra of classical observables of $(\clM,\Pi,\{l_k\})$, denoted $\Obs^{cl}_{\clM}$, assigns to each open subset $U\subset X$ the cochain complex
\begin{align*}
  \Obs^{cl}_{\clM}(U) \ceq (\Sym\left(\clM^*(U)[-1]\right),d_{CE}),
\end{align*}
where $d_{CE}$ is the Chevalley--Eilenberg differential of the $L_\infty$-algebra $(\clM(U),\{l_k\})$.
\end{dfn}

\begin{prop}[{\cite[Theorem 2.33]{Yoo}}]
\label{observables_Poisson}
For a degenerate classical field theory $(\clM,\Pi,\{l_k\})$ on a manifold $X$, there exists a $P_0$-factorization algebra $\wt{\Obs}^{cl}_{\clM}$ on $X$, and the $P_0$-factorization algebra structure $\{-,-\}$ is given by
\[
\{f,g\} \ceq Q_f(g),\ \text{for}\ f,g\in\wt{\Obs}^{cl}_{\clM}(U),
\]
for an open subset $U\subset X$.
\end{prop}

As in \cite[Proposition~5.4.5]{CG2}, one can show that a homotopy $\wt{\Obs}^{cl}_{\clM}\to\Obs^{cl}_{\clM}$ exists and the following proposition holds.

\begin{prop}
\label{quasi_isom}
For a degenerate classical field theory $(\clM,\Pi,\{l_k\})$ on a manifold $X$, then for any open subset $U\subset X$, there exists a quasi-isomorphism
\[
  \wt{\Obs}^{cl}_{\clM}(U) \qim \Obs^{cl}_{\clM}(U).
\]
This quasi-isomorphism makes $\Obs^{cl}_{\clM}$ a $P_0$-factorization algebra.
\end{prop}

%\begin{proof}
% The proof of Proposition~5.4.5 of \cite{CG2} works with minor modifications.
%\end{proof}

In Chapter 7.2 of \cite{CG2}, for a classical field theory $(\clM,\omega,S)$ and a gauge fixing operator $Q^{GF}:\clM\to\clM^!\simeq\clM$, a heat kernel is defined by a symmetric element $K_L\in\ol{\clM}\otimes\ol{\clM}$ that satisfies
\[
\omega(K_L,\varphi)=e^{-L[Q,Q^{GF}]}\varphi
\]
for any field $\varphi\in\clM$.
To extend this definition to degenerate classical field theories, we reformulate it in terms of a shifted Poisson structure.
First, the $(-1)$-shifted pairing is identified with the following map
\begin{align*}
  \clM[1] \lto (\clM[1])^![1], \quad
  \varphi \lmto \omega(\varphi,-).
\end{align*}
Since the pairing $\omega$ is nondegenerate, this map has an inverse 
\begin{align*}
  \Pi_\omega\colon (\clM[1])^! \lto \clM[1][-1],
\end{align*}
where the symbol $!$ denotes the Verdier dual in the sense of \cite{CG}.
Thus, the heat kernel $K_L$ is identified with the following map
%\footnote{note that for a vector bundle $E$, $\ol{\clE}$ is the continuous dual to $\clE^!_c$.}
\begin{align*}
  K_L=e^{-L[Q,Q^{GF}]}\circ\Pi_\omega\colon 
  (\clM[1])^!_c \lto \ol{\clM}[1][-1].
\end{align*}
Thus, the following definition is natural.

\begin{dfn}
\label{heat_kernel_of_degenerate_CLFT}
Let $(\clM,\Pi,Q_{BV})$ be a free degenerate classical field theory, and $Q^{GF}\colon\clM\to\clM^!$ be a gauge fixing operator.
Then we define the \emph{heat kernel} $K_L$ of the free degenerate classical field theory $(\clM,\Pi,Q_{BV})$ as follows
\[
  K_L \ceq \sum_{n\ge0}K_L^{(n)},\quad 
  K_L^{(n)} \ceq e^{-L[Q_{BV},Q^{GF}]}\circ\Pi^{(n)}\colon 
  \clM[1]^{\otimes n}\otimes\clM^!_c[-1] \lto \ol{\clM}[1][-1],
\]
where $\Pi^{(n)}$ is the $n$-th component of the $(-1)$-shifted Poisson structure $\Pi$.
Moreover, for $0<\epsilon<L$, we define the \emph{propagator} $P(\epsilon,L)$ of the free degenerate classical field theory $(\clM,\Pi,Q_{BV})$ as follows
\[
  P(\epsilon,L) \ceq \sum_{n\ge0}P^{(n)}(\epsilon,L)
  \ceq \sum_{n\ge0}\int_{\epsilon}^L(Q^{GF}\otimes1)K_t^{(n)} dt.
\]
\end{dfn}

For every $n\ge0$ and $0<\epsilon<L$, these $K_L^{(n)}$, and $P^{(n)}(\epsilon,L)$ define maps $\clM^!_c[-1]\to\ol{\clM}[1][-1]$ upon evaluating a section of $\clM[1]^{\otimes n}$.
Since the operator $Q_{BV}$ is an elliptic operator, these maps admit Schwartz kernel.
Therefore, the operators $K^{(n)}_L$ and $P^{(n)}(\epsilon,L)$ admit Schwartz kernels ``fieldwise'' in this sense, which we also denote by $K^{(n)}_L(x,y)$ and $P^{(n)}(\epsilon,L)(x,y)$, respectively.

In general, for $n\geq0$, given any map
\[
  \alpha\colon \clM^{\otimes n} \lto \clM\otimes\clM,
\]
we define an operator
\[
  \partial_\alpha\colon \Sym(\clM^*[-1]) \lto \Sym(\clM^*[-1]), 
\]
which is the unique differential operator that is zero on $\Sym^{\leq1}(\clM^*[-1])$ and that is given by the pairing induced by $\alpha$ on $\Sym^2(\clM^*[-1])$:
\begin{align*}
  \Sym^2(\clM^*[-1])\ni F \lmto
  \left((\phi_1,\dots,\phi_{n}) \lmto 
  \alpha(\phi_1,\dots,\phi_n)\vee F\right)\in\Sym^n(\clM^*[-1]),
\end{align*}
where we denote by $\vee$ the contraction operator (see Chapter 7.2.7 of \cite{CG2}).

Then we define the following cohomological degree $1$ operators
\begin{align*}
\Delta_L^{(n)}& \ceq \partial_{K_L^{(n)}}:\Sym^k(\clM^*(U)[-1])\lto\Sym^{k+n-2}(\clM^*(U)[-1]),\\
\Delta_L& \ceq \sum_{n\geq0}\Delta_L^{(n)}:\Sym(\clM^*(U)[-1])\lto\Sym(\clM^*(U)[-1]),
\end{align*}
and the following bilinear map 
\begin{gather*}
  \{-,-\}_L\colon 
  \Sym(\clM^*(U)[-1])\times\Sym(\clM^*(U)[-1])\lto\Sym(\clM^*(U)[-1])\\
  (I,J) \lmto \{I,J\}_L \ceq 
  \Delta_L(IJ)-\Delta_L(I)J-(-1)^{|I|}I\Delta_L(J),
\end{gather*}
where $|I|$ is the cohomological degree of the functional $I$.
We call $\Delta_L$ and $\{-,-\}_L$ the scale $L$-BV Laplacian and scale $L$-BV bracket, respectively.

\begin{prop}
\label{quantum_observables_degenerate}
Let $(\clM,\Pi,Q_{BV})$ be a free degenerate classical field theory on a manifold $X$.
For an open subset $U\subset X$ and a positive real number $L>0$, set
\[
  \Obs^q_{\clM,L}(U) \ceq \Sym(\clM^*(U)[-1])\dbr{\hbar}.
\]
Then the map
\[
  Q_{BV}+\hbar\Delta_L\colon \Obs^q_{\clM,L}(U)\lto\Obs^q_{\clM,L}(U)
\]
is square zero, and hence the $(\Obs^q_{\clM,L},Q_{BV}+\hbar\Delta_L)$ is a differentiable pro-cochain complex.
Moreover, for an open subset $U\subset X$, the triple $(\Obs^q_{\clM,L}(U),Q_{BV}+\hbar\Delta_L,\{-,-\}_L)$ is a BD algebra.
\end{prop}

\begin{proof}
  Regarding $\Delta_L$ as an operator-valued functional, substituting fields and proceeding in Chapter 7.2 of \cite{CG2}, a direct computation shows that the following holds.
  \[
  (\Delta_L)^2=0,\quad Q_{BV}\Delta_L=-\Delta_LQ_{BV}.
  \]
  Thus, we have $(Q_{BV}+\hbar\Delta_L)^2=Q_{BV}^2+\hbar(Q_{BV}\Delta_L+\Delta_LQ_{BV})+\hbar^2(\Delta_L)^2=0$.

  The differentiable pro-cochain complex structure on $\Obs^q_{\clM,L}(U)$ is given by the following filtration
  \begin{align}
    \label{diff_pro_cochain_str}
  F^n\Obs^q_{\clM,L}(U) \ceq \prod_{i\geq0}\hbar^i\Sym^{\geq(n-2i)}(\clM^*(U)[-1]).    
  \end{align}
  For an open subset $U\subset X$, by direct computation, we obtain
  \begin{align*}
      (Q_{BV}+\hbar\Delta_L)(IJ)&=Q_{BV}(I)J+(-1)^{|I|}IQ_{BV}(J)+\hbar\Delta_L(IJ)\\
      &=Q_{BV}(I)J+(-1)^{|I|}IQ_{BV}(J)+\hbar\left(\Delta_L(I)J+(-1)^{|I|}I\Delta_L(J)+\{I,J\}_L\right)\\
      &=(Q_{BV}+\hbar\Delta_L)(I)J+(-1)^{|I|}I(Q_{BV}+\hbar\Delta_L)(J)+\hbar\{I,J\}_L,
  \end{align*}
  for $I,J\in\Obs^q_{\clM,L}(U)$.
  Therefore, $(\Obs^q_{\clM,L}(U),Q_{BV}+\hbar\Delta_L,\{-,-\}_L)$ is a BD algebra.
\end{proof}

\begin{rem}
    \label{rem: difference of classical field theory}
    In Chapter 8.2.1 of \cite{CG2}, it is shown that the space of observables of a classical field theory, which is not necessarily free, carries a BD algebra structure.
    In the present construction, when the $(-1)$-shifted Poisson structure is \emph{nondegenerate}, this BD algebra agrees with the one constructed in \cite{CG2} for free classical field theories.
\end{rem}

The following lemma is an analogue of Lemma 8.3.3 of \cite{CG2}.
\begin{lem}
\label{lem_scale_isom}
For $0<L'<L$, and an open subset $U\subset X$, the linear map
\begin{align*}
  W^L_{L'}(U)\colon \Obs^q_{\clM,L'}(U) \lto \Obs^q_{\clM,L}(U), \quad  
  O \lmto e^{\hbar\partial_{P(L',L)}}O
\end{align*}
is an isomorphism of differentiable pro-cochain complexes.
\end{lem}
\begin{proof}
  Since we have
\begin{align*}
  (Q\otimes1)P(0,L)
  = \int_{\epsilon}^L([Q,Q^{GF}]\otimes1)K_tdt
  =-\int_{\epsilon}^L\frac{d}{dt}K_tdt
  =-K_L+K_\epsilon,
\end{align*}
we have
\[
 [Q,\partial_{P(L',L)}]=[Q,\partial_{P(0,L)}-\partial_{P(0,L')}]
 =\Delta_{L'}-\Delta_L.
\]
Therefore, we have
\begin{align*}
  (Q+\hbar\Delta_L)W^L_{L'}(U)(O)
  &=(Q+\hbar\Delta_L)e^{\hbar\partial_{P(L',L)}}O
   =e^{\hbar\partial_{P(L',L)}}(Q+\hbar\Delta_{L'})O\\
  &=W^L_{L'}(U)((Q+\hbar\Delta_{L'})O).
\end{align*}
The above equation implies that the map $W^L_{L'}(U)$ is a map of cochain complexes.
By the definition of the differential operator $\partial_{P(L',L)}$, the map $W^L_{L'}(U)$ preserves the filtration \eqref{diff_pro_cochain_str}.
Since the map $W^{L'}_L(U)$ is the inverse of $W^L_{L'}(U)$, the map $W^L_{L'}(U)$ is an isomorphism of differentiable pro-cochain complexes.
\end{proof}

Using the identification
\[
  \Obs^q_{\clM,L}=\Sym(\ol{\clM}^!_c[-1])\dbr{\hbar},
\]
for open subsets $U\subset V\subset X$, we define a linear map
\[
  (m^{\Obs}_L)^U_V\colon \Obs^q_{\clM,L}(U) \lto \Obs^q_{\clM,L}(V)
\]
by zero extension of the compactly supported distributional sections.

\begin{lem}
\label{lem_degene_fact_alg_str}
Let $(\clM,\Pi,Q_{BV})$ be a free degenerate classical field theory on a manifold $X$.
Then for $L>0$, disjoint open subsets $U,V\subset X$, an open subset $U\sqcup V\subset W\subset X$, and observables $O\in\Obs^q_{\clM,L}(U),O'\in\Obs^q_{\clM,L}(V)$, there exists a limit
\[
  \lim_{\epsilon\to0}W^L_{\epsilon}\left((m^{\Obs}_L)^U_W(W^\epsilon_{L}(O))*(m^{\Obs}_L)^V_W(W^{\epsilon}_L(O'))\right),
\]
and a bilinear map
\begin{align*}
  (m^{\Obs}_L)^{U,V}_W\colon 
  &\Obs^q_{\clM,L}(U)\times\Obs^q_{\clM,L}(V) \lto \Obs^q_{\clM,L}(W),\\
  (O,O')&\lmto
  \lim_{\epsilon\to0}W^L_{\epsilon}\left((m^{\Obs}_L)^U_W(W^\epsilon_{L}(O))*(m^{\Obs}_L)^V_W(W^{\epsilon}_L(O'))\right),
\end{align*}
defines a prefactorization algebra structure on 
\[
  \Obs^q_{\clM,L}\colon \Open(X) \lto\text{pCh(DVS)}, \quad 
  (U \lmto \Obs^q_{\clM,L}(U)).
\]
\end{lem}

\begin{proof}
  Since for $L>0$, there exists a limit $P(0,L)=\lim_{\epsilon\to0}P(\epsilon,L)$, we have
  \begin{align*}
    &\lim_{\epsilon\to0}W^L_{\epsilon}\left((m^{\Obs}_L)^U_W(W^\epsilon_{L}(O))*(m^{\Obs}_L)^V_W(W^{\epsilon}_L(O'))\right)\\
  =&\lim_{\epsilon\to0}e^{\hbar\partial_{P(\epsilon,L)}}\left((m^{\Obs}_L)^U_W\left(e^{-\hbar\partial_{P(\epsilon,L)}}O\right)*(m^{\Obs}_L)^V_W\left(e^{-\hbar\partial_{P(\epsilon,L)}}O'\right)\right)\\
  =&e^{\hbar\partial_{P(0,L)}}\left((m^{\Obs}_L)^U_W\left(e^{-\hbar\partial_{P(0,L)}}O\right)*(m^{\Obs}_L)^V_W\left(e^{-\hbar\partial_{P(0,L)}}O'\right)\right),
\end{align*}
and hence there exists the limit $\lim_{\epsilon\to0}W^L_{\epsilon}\left((m^{\Obs}_L)^U_W(W^\epsilon_{L}(O))*(m^{\Obs}_L)^V_W(W^{\epsilon}_L(O'))\right)$.

For disjoint open subsets $U,V\subset X$, and $O\in\Obs^q_{\clM,L}(U)$ and $O'\in\Obs^q_{\clM,L}(V)$, the condition $\supp{O}\sqcup\supp{O'}=\emptyset$ implies that the scale $L$-BV bracket $\{O,O'\}_L$ vanishes.
Hence, the map $(m^{\Obs}_L)^{U,V}_W$ is a cochain map.
Since it is obvious that the map $(m^{\Obs}_L)^{U,V}_W$ preserves the filtration \eqref{diff_pro_cochain_str}, the map $(m^{\Obs}_L)^{U,V}_W$ is a map of differentiable pro-cochain complexes.
It is straightforward to show that this satisfies the associativity and commutativity properties necessary to define a prefactorization algebra. 
\end{proof}

\begin{rem}
    \label{rem:factorizationprodisdifferfromvol1}
    The factorization product in \cref{lem_degene_fact_alg_str} is a generalization to free degenerate field theories of the factorization product in \cite[Theorem 8.5.1]{CG2} with vanishing \emph{effective interaction}.
    %{red}vol2のfactorization productの定義と本質的に同じだけど、effective interactionが0のときの拡張になっている
\end{rem}

\begin{lem}
\label{gr_of_q_degenerate}
Let $(\clM,\Pi,Q_{BV})$ be a free degenerate classical field theory on a manifold $X$.
For every open subset $U\subset X$, we define a filtration $\{F^p_\hbar\Obs^q_{\clM,L}(U)\}_{p\geq0}$ on $\Obs^q_{\clM,L}(U)$ by saying the $p$-th filter $F^p_\hbar\Obs^q_{\clM,L}(U)$ is the sub-$\bbC\dbr{\hbar}$-module consisting of those observables that are zero modulo $\hbar^k$.
Then, there exists an isomorphism
\[
  \gr^\hbar\Obs^q_{\clM,L} \lsto 
  \Obs^{cl}_{\clM}\otimes_{\bbC}\bbC\dbr{\hbar}
\]
of prefactorization algebras, where $\gr^\hbar\Obs^q_{\clM,L}$ is an associated graded with respect to the above filtration $\{F^p_{\hbar}\Obs^q_{\clM,\omega,S}\}_{p\ge0}$.
This isomorphism makes $\gr^\hbar\Obs^q_{\clM,L}$ into a factorization algebra.
\end{lem}

\begin{proof}
By definition of the differential $Q+\hbar\Delta_L$ on $\Obs^q_{\clM,L}$, for every open subset $U\subset X$, the associated graded $\gr^\hbar\Obs^q_{\clM,L}(U)$ is a pro-cochain complex with the differential $Q_{BV}$ and the filtration
\[
  G^p\gr^\hbar\Obs^q_{\clM,L}(U) \ceq \Sym^{\geq p}(\clM^*(U)[-1]).
\]
Therefore, we have an isomorphism
\[
  \Phi\colon \gr^\hbar\Obs^q_{\clM,L} \xrightarrow{\ \sim \ }
  \Obs^{cl}_{\clM}\otimes_{\bbC}\bbC\dbr{\hbar}
\]
of precosheaf on $X$ valued in the category $\text{pCh(DVS)}$.
It is enough to show that the isomorphism $\Phi$ is a map of prefactorization algebras.
By definition of the factorization product on $\Obs^q_{\clM,L}$, the factorization product on $\gr^\hbar\Obs^q_{\clM,L}$ is nothing but the following map
\begin{align*}
  (\gr^\hbar m^{\Obs}_L)^{U,V}_W\colon
   \gr^\hbar\Obs^q_{\clM,L}(U)\times\gr^\hbar\Obs^q_{\clM,L}(V)&\lto\gr^\hbar\Obs^q_{\clM,L}(W)\\
  (O,O')&\lmto(m^{\Obs}_L)^{U}_W(O)*(m^{\Obs}_L)^V_W(O'),
\end{align*}
for open subsets $U\sqcup V\subset W\subset X$.
Therefore, the map $\Phi$ is a map of prefactorization algebras.
\end{proof}

The following theorem is the main theorem in this subsection.

\begin{thm}
\label{FA_str_on_Obs}
Let $(\clM,\Pi,Q_{BV})$ be a free degenerate classical field theory on a manifold $X$.
Then the prefactorization algebra $\Obs^q_{\clM,L}$ of quantum observables has the following properties.
\begin{itemize}
  \item The tuple $(\Obs^q_{\clM,L},Q_{BV}+\hbar\Delta_L,\{-,-\}_L)$ is a BD factorization algebra.
  \item There exists an isomorphism
  \[
  \Obs^q_{\clM,L}\otimes_{\bbC\dbr{\hbar}}\bbC\lsto\Obs^{cl}_{\clM}
  \]
  of factorization algebras.
\end{itemize}
\end{thm}
\begin{proof}
The proof of Theorem~8.6.1 of \cite{CG2} works in this situation,
where Proposition~8.6.2 in loc.\ cit.\ is replaced by \cref{gr_of_q_degenerate}.
\end{proof}

By \cref{lem_scale_isom}, the BD factorization algebra structure on $\Obs^q_{\clM,L}$ is independent of the scale $L$.
Therefore, we define the “universal” BD factorization algebra $\Obs^q_{\clM}$ as follows.
\begin{align*}
  \Obs^q_{\clM}(U) := 
  \left(\{O[L]\}_{L>0}|O[L]\in\Obs^q_{\clM,L}(U),\ 
        W^{L'}_L(O[L])=O[L'] \text{ for all $L,L'>0$} \right)
\end{align*}
for an open subset $U\subset X$, 
with differential and BD algebra structure is defined scalewise.
Then for $L>0$, the projection
\begin{align*}
  p_L:\Obs^q_{\clM} \lto \Obs^q_{\clM,L}, \quad
  \{O[L']\}_{L'>0} \lmto O[L]
\end{align*}
is a quasi-isomorphism of BD factorization algebras.

%%%%%%%%%%%%%%%%%%%%%%%%%%%%%%%%%%%%%%%%%%%%%%%%%%%%%%%%%%%%%%%%%%%%%%%%%%%%%%%%
\subsection{Equivariant field theory}
\label{ss:2.3}

In this subsection, we consider the action of a local $L_\infty$-algebra on a strictly free degenerate classical field theory.

We now introduce a complex that plays an important role in the study of actions on field theories.
For a local $L_\infty$-algebra $\clL$ on a manifold $X$, we define the presheaf
\begin{align*}
C^*(\clL)& \ceq \Sym(\clL^*[-1]),
\end{align*}
equipped with the Chevalley--Eilenberg differential associated with the $L_\infty$-algebra $\clL(X)$.
Furthermore, we define the presheaf
\begin{align*}
C^*_{red,loc}(\clL)& \ceq \text{Dens}_X\otimes_{\clD_X}\prod_{n>0}\Hom_{C^\infty_X}(\clJ(\clL)^{\otimes n},C^\infty_X)_{S_n}.
\end{align*}
By Corollary 3.5.2 of \cite{CG2}, the sheaf $\clJ(\clL)$ has the local $L_\infty$-algebra structure induced by $\clL$, and hence the presheaf $C^*_{red,loc}(\clL)$ equipped with the Chevalley--Eilenberg differential associated with the $L_\infty$-algebra structure on $\clJ(\clL)(X)$.

Let $(\clM,\Pi,Q_{BV})$ be a strictly free degenerate classical field theory on a manifold $X$, and $\clL$ be a local $L_\infty$-algebra on $X$.

\begin{dfn}
\label{inner_act_outer_act}
We define graded modules
\begin{align*}
 \IA(\clL,\clM) & \ceq C^*_{red,loc}(\clL\oplus\clM), \\ 
 \Act(\clL,\clM) & \ceq
 \IA(\clL,\clM) / C^*_{red,loc}(\clL).
\end{align*}
%where $C^*_{red,loc}(\clL)$ is 
%the central subcomplex of $\IA(\clL,\clM)$.
\end{dfn}
We have the following natural identification
\[
C^*_{red,loc}(\clL\oplus\clM)=C^*_{red,loc}(\clL)\otimes C_{red,loc}^*(\clM)
\]
of sheaves of graded algebras on $X$ (forgetting differentials).
Then, the $P_0$-factorization algebra structure $\{-,-\}$ on $\Obs^{cl}_{\clM}$ which is defined in \cref{quasi_isom} can be extended to the vector spaces $\IA(\clL,\clM)$ and $\Act(\clL,\clM)$ by $C^*(\mathcal{L})$-linearity.
\begin{dfn}
\label{equivariant_action}
A cohomological degree 0 element $S^{\clL}\in\Act(\clL,\clM)$ is an equivariant action functional if the following equality holds
\begin{align}
  \label{equiv_classical_master_equation}
d_{\clL}S^{\clL}+Q_{BV}S^{\clL}+\frac{1}{2}\{S^{\clL},S^{\clL}\}=0,
\end{align}
where $d_{\clL}$ is the Chevalley--Eilenberg differential of the $L_\infty$-algebra $\clL(X)$, and $\{-,-\}$ is the $P_0$-factorization algebra structure on $\Obs^{cl}_{\clM}$.
We call the \eqref{equiv_classical_master_equation} the equivariant classical equation.
\end{dfn}

By definition, we have the following assertion.

\begin{prop}
\label{differential_Inneract_act_classical}
Given an equivariant action functional $S^{\clL}\in\Act(\clL,\clM)$, 
the cohomological degree 1 linear map
\[
  d_{\mathrm{Inner}} \ceq d_{\clL}+Q_{BV}+\{S^{\clL},-\}\colon
  \IA(\clL,\clM) \lto \IA(\clL,\clM)
\]
defines a differential on $\IA(\clL,\clM)$,
with respect to which 
$C^*_{red,loc}(\clL) \subset \IA(\clL,\clM)$ is a subcomplex,
and $d_{\mathrm{Inner}}$ induces a differential on the quotient $\Act(\clL,\clM)$.
We have a short exact sequence of complexes
\[
  0 \lto C^*_{red,loc}(\clL) \lto \IA(\clL,\clM)
    \lto \Act(\clL,\clM) \lto 0.
\]
\end{prop}
\begin{proof}
By \cref{equivariant_action}, the vector space $\IA(\clL,\clM)$ is equipped with the grading given by the cohomological degree of the Chevalley--Eilenberg complex $C^*_{red,loc}(\clL\oplus\clM)$.
It is obvious that the maps $d_\clL$ and $Q_{BV}$ has cohomological degrees 1.
Since the $P_0$-factorization algebra structure $\{-,-\}$ has cohomological degree 1 and $S^\clL$ has cohomological degree 0, the map $\{S^{\clL},-\}$ has cohomological degree 1. Thus, the map $d_{\text{Inner}}$ has cohomological degree 1.
The equivariant classical master equation implies that the map $d_{\text{Inner}}$ is square zero, and hence the map $d_{\text{Inner}}$ makes $\IA(\clL,\clM)$ a cochain complex.
Since the maps $Q_{BV}$ and $\{S^\clL,-\}$ act trivially on $C^*_{red,loc}(\clL)$, the graded subspace $C^*_{red,loc}(\clL)$ is a subcomplex with respect to the differential $d_{\text{Inner}}$.
\end{proof}

\begin{rem}
  \label{rem_inner_act}
  In Chapter 12 of \cite{CG2}, the notion of an inner action is introduced.
  However, since we consider only free theories here, there is no need to consider inner actions.
\end{rem}

\begin{prop}
\label{classical_Noether_map}
In this situation, if the local $L_\infty$-algebra is local dgLA, 
there exists a map
\[
  \bbU_0^{cl}\clL \lto \Obs^{cl}_{\clM}
\]
of factorization algebras, where $\bbU_0^{cl}$ is the factorization algebra \eqref{eq:U_alpha}.
\end{prop}

\begin{proof}
For an open subset $U\subset X$, we define a linear map
\begin{align*}
  \Phi_U\colon C_*(\clL_c)(U) \lto \Obs^{cl}_{\clM}(U), \quad 
  l\lmto\ev*{l,\exp(S^{\clL})}|_{0\times\clM(X)[1]},
\end{align*}
where $\ev*{-,-}\colon C_*(\clL_c)(U)\otimes C^*(\clL)(U)\to\bbC$ is the duality pairing between $C_*(\clL_c)(U)$ and $C^*(\clL)(U)$.
We show that the map $\Phi_U$ is a map of differentiable pro-cochain complexes.
Since the functional $S^{\clL}$ has cohomological degree 0, the equivariant classical master equation \eqref{equiv_classical_master_equation} implies that the following equality
\begin{align*}
  \Phi_U(d_{\clL}l)(m)
  &= \ev*{d_{\clL}l,\exp(S^{\clL})}(m)
   =-\ev*{l,d_{\clL}(\exp(S^{\clL})})(m)\\
  &=\ev*{l,(Q_{BV}+\{S^{\clL},-\}\exp(S^{\clL}))}(m)\\
  &=(Q_{BV}+\{S^{\clL},-\})\ev*{l,\exp(S^{\clL})}(m).
\end{align*}
Therefore, the map $\Phi_U$ is a cochain map.
It is easy to check that the map $\Phi_U$ preserves the filtrations of $\bbU^{P_0}(\clL)(U)$ and $\Obs^{cl}_{\clM}(U)$.
By construction, the map $\Phi=(\Phi_U)_{U\in\Open(X)}$ defines a map of precosheaves.

Finally, we show that the map $\Phi$ is a map of factorization algebras.
It is enough to show that for disjoint open subsets $V_1,V_2\subset X$ and sections $m_{ji}\in \clL_c(V_i),\ (j=1,\dots,k_{i},\ i=1,2)$, the following hold.
\begin{align*}
  &\Phi_{V_1}(m_{11}\cdots m_{k_11})*\Phi_{V_2}(m_{12}\cdots m_{k_22})\\
  &=\left\{\frac{\partial}{\partial m_{11}}\cdots\frac{\partial}{\partial m_{k_11}}\exp(S^{\clL})\right\}\left\{\frac{\partial}{\partial m_{12}}\cdots\frac{\partial}{\partial m_{k_22}}\exp(S^{\clL})\right\}\\
  &=\frac{\partial}{\partial m_{11}}\cdots\frac{\partial}{\partial m_{k_11}}\frac{\partial}{\partial m_{12}}\cdots\frac{\partial}{\partial m_{k_22}}\exp(S^{\clL})\\
  &=\Phi_{V_1\sqcup V_2}(m_{11}\cdots m_{k_11}*m_{12}\cdots m_{k_22}),
\end{align*}
 where for simplicity, we suppress the notation for the zero extension and simply denote the factorization product by $*$.
 Since for $j_1=1,\dots,k_1$ and $j_2=1,\dots,k_2$, the supports of $m_{j_11}\in\clL_c(V_1)$ and $m_{j_22}\in \clL_c(V_2)$ are disjoint, for a functional $F\in C^*(\clL)(V_1\sqcup V_2)$, we have
\[
  \frac{\partial}{\partial m_{j_11}}\frac{\partial}{\partial m_{j_22}}F=0.
\]
Thus, we have the desired equality.
\end{proof}

\begin{rem}
  \label{rem_Noether_map_def}
  Note that the following two observables on $\clM$ are distinct: 
  \begin{align*}
    \ev*{l,\exp(S^{\clL})} \text{ and } \exp(S^{\clL}(l))
  \end{align*}
  for a section $l$ of $\clL$.
\end{rem}

Next, we consider the quantum action of a local $L_\infty$-algebra on a strictly free degenerate classical field theory $(\clM,\Pi,Q_{BV})$ on a manifold $X$.
Note that the $(-1)$-shifted Poisson structure $\Pi$ is concentrated in the $0$-component:
\[
  \Pi\colon (\clM^![1]) \lto \clM[1][-1].
\]
Thus, the heat kernel and the propagator in \cref{heat_kernel_of_degenerate_CLFT} are Schwartz kernels in the usual sense, and hence all the arguments in Chapter 13 of \cite{CG2} carry over directly.

We have the following natural identification
\[
C^*(\clL\oplus\clM)=C^*(\clL)\otimes C^*(\clM)
\]
of sheaves of graded algebras on $X$ (forgetting differentials).
Then, one can extend to the left hand side, by $C^*(\clL)$-linearity, the operations scale $L$-BV Laplacian $\Delta_L$, the scale $L$-BV bracket, and the propagator $P(\epsilon,L)$.

Let
\[
\IA^q(\clL,\clM)\subset C^*(\clL\oplus\clM)(X)\dbr{\hbar}
\]
denote the space of those functionals that satisfy the following conditions:
\begin{itemize}
  \item Their constant parts (degree 0 parts) are zero, and they are zero modulo $\hbar$ when restricted to be functions just on $\clM$.
  \item They have proper support as functionals on $\clL\oplus\clM$, in the sense of Definition 7.2.2 of \cite{CG2}.
  \item They has smooth first derivative, in the sense of Section 5.4 of \cite{CG2}. Note that this condition involves differentiation by elements of both $\clL$ and $\clM$.
\end{itemize}

\begin{dfn}
\label{equivariant_quantum_master}
An element
\[
I\in\IA^q(\clL,\clM)
\]
solves the \emph{scale-$L$ equivariant quantum master equation} if it satisfies
\[
d_{\clL}I+Q_{BV}I+\{I,I\}_L+\hbar\Delta_LI=0.
\]
\end{dfn}
Similarly to the classical case, we define the following space
\begin{align*}
  \Act^q(\clL,\clM)
  =\IA^q(\clL,\clM)/(C^*(\clL)(X)\dbr{\hbar}\cap\IA^q(\clL,\clM)).
\end{align*}

For $\epsilon,L>0$ and $I\in\IA^q(\clL,\clM)$, we define a $C^*(\clL)$-linear map $W(P(\epsilon,L),I)$ by
\begin{align}
 \label{eq:WPI}
 W(P(\epsilon,L),I)=\hbar\log\left(\exp(\hbar\partial_{P(\epsilon,L)})\exp(I/\hbar)\right).
\end{align}

Although the expression in \eqref{eq:WPI} contain negative power of $\hbar$, 
by Lemma 3.4.1 of \cite{C11}, 
%for $\epsilon,L>0$ and $I\in\IA^q(\clL,\clM)$, 
the functional $W(P(\epsilon,L),I)$ belongs to 
$C^*(\clL\oplus\clM)(X)\dbr{\hbar}$.
Also, since the linear map $W(P(\epsilon,L),-)$ is $C^*(\clL)(X)\dbr{\hbar}$-linear, for %$\epsilon,L>0$ and 
the class $[I]\in\Act^q(\clL,\clM)$, the functional $W(P(\epsilon,L),[I])$ is well-defined and belongs to the quotient space $C^*(\clL\oplus\clM)(X)\dbr{\hbar}/C^*(\clL)(X)\dbr{\hbar}$.

\begin{dfn}
\label{equiv_quantum_action}
A \emph{quantum outer-action} of $\clL$ on a quantization of the strictly free degenerate field theory $(\clM,\Pi,Q_{BV})$ is a collection of functionals
\[
I^{\clL}[L]\in\Act^q(\clL,\clM)
\]
for every scale $L>0$ satisfying the following properties:
\begin{itemize}
  \item It satisfies the renormalization group equation
  \[
  W(P(L',L),I^{\clL}[L'])=I^{\clL}[L],
  \]
  for every $L',L>0$. 
  \item Each $I^{\clL}[L]$ satisfies the \emph{equivariant outer-quantum master equation}:
  \[
  d_{\clL}I^{\clL}[L]+Q_{BV}I^{\clL}[L]+\{I^{\clL}[L],I^{\clL}[L]\}_L+\hbar\Delta_LI^{\clL}[L]=0.
  \]
  \item The locality axiom, as explained in Section 7.2.5 of \cite{CG2}, holds: The support of $I^{\clL}[L]$ converges to the diagonal as the scale $L$ tends to zero, with the same bounds explained in Section 7.2.5 of \cite{CG2}.
\end{itemize}
\end{dfn}

This notion is a degenerate and equiavriant analogue of \emph{effective interaction} defined in \cite{CG2}.
In Chapter 13.2.2 of \cite{CG2}, the notion of action of an $L_\infty$-algebra on a quantum field theory in the sense of \cite{CG2} is defined.
\cref{equiv_quantum_action} is a natural extension to our setting.

The vector space $\Act^q(\clL,\clM)$ has a grading
induced by the cohomological grading of the Chevalley--Eilenberg complex $C^*(\clL)(X) \mathbin{\otimes} C^*(\clM)(X)\dbr{\hbar}$.
Moreover, given a quantum outer-action $\{I^{\clL}[L]\}_{L>0}$, it has a  differential $d_{\clL}+Q_{BV}+\hbar\Delta_L+\{I^{\clL}[L],-\}_L$ for each $L>0$, making $\Act^q(\clL,\clM)$ a complex. 
Indeed, the equivariant quantum master equation implies that the linear map $d_{\clL}+Q_{BV}+\hbar\Delta_L+\{I^{\clL}[L],-\}_L$ is square zero.

\begin{dfn}
  \label{equiv_inner_act}
  A \emph{quantum inner-action} of $\clL$ on a quantization of the strictly free degenerate field theory $(\clM,\Pi,Q_{BV})$ is a collection of functionals
\[
\wt{I}^{\clL}[L]\in\IA^q(\clL,\clM)
\]
that satisfies the same three conditions as for an outer-action in this larger space of functionals.
\end{dfn}

Similarly to $\Act^q(\clL,\clM)$, 
%the vector space $\IA^q(\clL,\clM)$ 
%has a natural grading induced by the cohomological grading on the Chevalley--Eilenberg complex $C^*(\clL)(X)\otimes C^*(\clM)(X)\dbr{\hbar}$.
%Moreover, 
given a quantum inner-action $\{\wt{I}^{\clL}[L]\}_{L>0}$, $\IA^q(\clL,\clM)$ is a complex
with differential $d_{\clL}+Q_{BV}+\hbar\Delta_L+\{\wt{I}^{\clL}[L],-\}_L$ for each $L>0$.
%, making the vector space $\IA^q(\clL,\clM)$ a complex.
%Indeed, the equivariant quantum master equation implies that the linear map $d_{\clL}+Q_{BV}+\hbar\Delta_L+\{\wt{I}^{\clL}[L],-\}_L$ is square zero.

\begin{rem}
\label{scale_change}
For $L',L>0$, the equality
\[
  [Q,\partial_{P(L',L)}]=\Delta_{L'}-\Delta_L
\]
implies that if the functional $I[L']\in\IA^q(\clL,\clM)$ satisfies the equivariant quantum master equation for the scale $L'>0$, then the following equality
\begin{align*}
  (Q+\hbar\Delta_L)\exp(W(P(L',L),I[L'])/\hbar)
  &=(Q+\hbar\Delta_L)e^{\hbar\partial_{P(L',L)}}e^{I[L']/\hbar}\\
  &=e^{\hbar\partial_{P(L',L)}}(Q+\hbar\Delta_{L'})e^{I[L']/\hbar}=0,
\end{align*}
and hence the functional $I[L]=W(P(L',L),I[L'])\in\IA^q(\clL,\clM)$ satisfies the equivariant quantum master equation for the scale $L$.
\end{rem}

The following result is a degenerate analogue of \cite[Lemma 11.1.1]{C11}.

\begin{prop}
\label{lift_classical_to_quantum_obstruction_class}
Let $X$ be a manifold, $\clL$ be a local $L_\infty$-algebra on $X$, and $(\clM,\Pi,Q_{BV})$ be a strictly free degenerate classical field theory on $X$.
Suppose that $\clL$ acts on the field theory $(\clM,\Pi,Q_{BV})$.
Then the obstruction class for lifting this classical action to a quantum outer-action on the quantization of $(\clM,\Pi,Q_{BV})$ is an element of 1st cohomology
\[
  H^1(\Act(\clL,\clM)).
\]
Moreover, the obstruction class for lifting this classical action to a quantum inner-action on the quantization of $(\clM,\Pi,Q_{BV})$ is an element of 1st cohomology
\[
  H^1(\IA(\clL,\clM)).
\]
\end{prop}
\begin{proof}
The proof is similar to that of \cite[Lemma 11.1.1]{C11}.

  Let $S^{\clL}\in\Act(\clL,\clM)$ be an equivariant action functional that defines the classical action of $\clL$ on $(\clM,\Pi,Q_{BV})$.
  Let $\{I^{\clL}[L]\in\Act^q(\clL,\clM)/(\hbar^2)\}_{L>0}$ be an arbitrary lift of $S^{\clL}$, i.e. the following holds
  \[
  S^{\clL}=\lim_{L\to0}I^{\clL}[L]\ \bmod \hbar.
  \]
  Then for $L>0$, define
  \[
  \Theta[L] \ceq (d_{\clL}I^{\clL}[L]+Q_{BV}I^{\clL}[L]+\frac{1}{2}\{I^{\clL}[L],I^{\clL}[L]\}_L+\hbar\Delta_LI^{\clL}[L])|_{\Act^q(\clL,\clM)},
  \]
  where the symbol $(\cdots)|_{\Act^q(\clL,\clM)}$ means taking the class in the quotient space $\Act^q(\clL,\clM)$.
  Let us make $\Act^q(\clL,\clM)/(\hbar^2)$ into a dgLA, with bracket $\{-,-\}_L$, and with differential $d_{\clL}+Q_{BV}+\hbar\Delta_L$. 
  Since for every $L>0$, the element $\Theta[L]\in\Act^q(\clL,\clM)/(\hbar^2)$ has cohomological degree 1, by the following \cref{lem_obstructction_class_useful}, we have
  \[
  (d_{\clL}+Q_{BV}+\hbar\Delta_L)\Theta[L]+\{I^{\clL}[L],\Theta[L]\}_L=0.
  \]
  Taking the expression above modulo $\hbar$ and then taking the limit $L\to0$, we obtain the following equation:
  \[
  (d_{\clL}+Q_{BV})\Theta+\{S^{\clL},\Theta\}=0,
  \]
  where we denote $\Theta=\lim_{L\to0}\Theta[L]\ \bmod \hbar$. This implies that the element $\Theta\in\Act(\clL,\clM)$ is a $d_{\mathrm{Inner}}=d_{\clL}+Q_{BV}+\{S^{\clL},-\}$-closed element, and hence the element $\Theta$ defines a cohomology class $[\Theta]\in H^1(\Act(\clL,\clM))$.
  The latter statement follows exactly the same way.
\end{proof}

The following \cref{lem_obstructction_class_useful} is stated in the proof of \cite[Lemma 11.1.1]{C11}, and can be verified by a direct computation.

\begin{lem}
  \label{lem_obstructction_class_useful}
Let $\frg$ be a dgLA, let $X\in\frg$ be a cohomologically odd element, and let 
\[
  O(X)=d_{\frg}X+\frac{1}{2}[X,X].
\]
Then we have
\[
 d_{\frg}O(X)+[X,O(X)]=0.
\]
\end{lem}

The following \cref{obstruction_class_inner_act_to_aut_act_quantum} is crucial, as it explains the origin of the level parameter in the Wakimoto realization studied in \cref{ss:4.2} and subsequent sections, and is a degenerate analogue of \cite[Proposition 13.3.1]{CG2}.

\begin{lem}
\label{obstruction_class_inner_act_to_aut_act_quantum}
Let $X$ be a manifold, $\clL$ be a local $L_\infty$-algebra on $X$, and $(\clM,\Pi,Q_{BV})$ be a strictly free degenerate classical field theory on $X$.
Suppose that the local $L_\infty$-algebra $\clL$ acts on the quantization of $(\clM,\Pi,Q_{BV})$.
Then, there exists a cohomological degree 1 closed element
\[
\alpha\in C^*_{red,loc}(\clL)\dbr{\hbar},
\]
such that the action of $\clL$ lifts an inner-action if and only if $\alpha$ vanishes.
\end{lem}

\begin{proof}
 Let $\{I^{\clL}[L]\in\Act^q(\clL,\clM)\}_{L>0}$ be a quantum outer-action of $\clL$ on the quantization of $(\clM,\Pi,Q_{BV})$.
 For $L>0$, let $\wt{I}^{\clL}[L]\in\IA^q(\clL,\clM)/(\hbar^2)$ be an arbitrary lift of $I^{\clL}[L]$.
 Note that there exists a short exact sequence
\[
  0 \lto C^*_{red,loc}(\clL)\dbr{\hbar}\xrr{i} 
  \IA^q(\clL,\clM) \xrr{p}
  \Act^q(\clL,\clM) \lto 0
\]
of complexes. Thus, we have an exact sequence
\begin{align}
\label{exact_equiv_proof}
  H^0(\IA^q(\clL,\clM))\stackrel{H^0(p)}{\lto} H^0(\Act^q(\clL,\clM))\stackrel{\delta}{\lto} H^1(C^*_{red,loc}(\clL))\dbr{\hbar}
\end{align}
The quantum equivariant master equation for the scale $L$ implies 
that the functional $I^{\clL}[L]$ defines a cohomology class $[I^{\clL}[L]]\in H^0(\Act^q(\clL,\clM))$.
Then the exactness of \eqref{exact_equiv_proof} implies that the lift $\wt{I}^{\clL}[L]\in\IA^q(\clL,\clM)$ defines a cohomology class $H^0(\IA^q(\clL,\clM))$ and the image under the map $H^0(p)$ coincides with $[I^{\clL}[L]]$ if and only if the image $\alpha[L]=\delta([I^{\clL}[L]])$ vanishes.
By \cref{scale_change}, the condition that $\alpha[L]$ vanishes for some $L>0$ is equivalent to the condition that $\alpha[L]$ vanishes for all $L>0$.
Thus, setting $\alpha=\lim_{L\to0}\alpha[L]$, we obtain the desired obstruction class.
\end{proof}

Explicitly, the image $\alpha[L]=\delta([I^{\clL}[L]])$ is given as follows
\[
  \delta([I^{\clL}[L]])=\left.\left(d_{\clL}\wt{I}^{\clL}[L]+Q_{BV}\wt{I}^{\clL}[L]+\hbar\Delta_L\wt{I}^{\clL}[L]+\frac{1}{2}\{\wt{I}^{\clL}[L],\wt{I}^{\clL}[L]\}_L\right)\right|_{\clL\times\{0\}}.
\]
The following theorem is a key result of this note, 
and is a degenerate analogue of \cite[Theorem 13.5.2]{CG2}

\begin{thm}
\label{quantum_Noether_map}
Let $X$ be a manifold, $\clL$ be a local $L_\infty$-algebra on $X$, and $(\clM,\Pi,Q_{BV})$ be a strictly free degenerate classical field theory on $X$.
Suppose that the local $L_\infty$-algebra $\clL$ acts on the quantization of $(\clM,\Pi,Q_{BV})$.

Furthermore, let $\alpha\in H^1(C^*_{red,loc}(\clL))\dbr{\hbar}$ be the obstruction class for lifting this action to an inner-action in \cref{obstruction_class_inner_act_to_aut_act_quantum}.
Then there exists a $\bbC\dbr{\hbar}$-linear map of factorization algebras
\[
\mu\colon\bbU^\hbar_\alpha\clL[\hbar^{-1}]\lto\Obs^q_{\clM}[\hbar^{-1}].
\]
\end{thm}
\begin{proof}[\textbf{Sketch of proof}]
Here we merely construct the desired map and omit the verification that it satisfies the required properties.
 We refer the reader to Theorem 13.5.1 of \cite{CG2} for the details.

Similarly to the construction of the classical Noether map in \cref{classical_Noether_map}, we construct a linear map from the symmetric product $\Sym(\clL_c(V))\dpr{\hbar}$ to $\Obs^q_{\clM}[\hbar^{-1}]$ using duality pairing.

 %Taking an open subset $V\subset X$, we will construct a linear map
 %\[
 %f:C_*^\alpha(\clL_c(V))\to\Obs^q_{\clM}(V)[\hbar^{-1}],
 %\]
 %where $C_*^\alpha(\clL_c(V))$ is a ``twisted'' Chevalley--Eilenberg chain complex whose differential is $d_{\clL}-\alpha$.

 %Using the natural pairing between the symmetric algebra of $\clL_c(V)[1]$ and the space of functionals on $\clL_c(V)[1]$, the linear map
%\[
%f:\Sym(\clL_c(V))\dpr{\hbar}\to C^*(\clM(X)[1])\dpr{\hbar},
%\]
%in turn, is identified with a functional
 %\[
% O_f\in C^*(\clL_c(V)[1]\oplus\clM(X)[1])\dpr{\hbar}.
% \]

Let $\{I^{\clL}[L]\}_{L>0}$ be an outer quantum action of $\mathcal{L}$ on the quantization of $(\mathcal{M},\Pi,Q_{BV})$.
Then, the collection of functionals
\begin{align}
\label{Noether_map_define}
 \exp\left(I^{\clL}[L]/\hbar\right) \in 
 C^*(\clL_c(V)[1]\oplus\clM(X)[1])\dpr{\hbar}, \quad   
\end{align}
satisfies the renormalization group equation and the locality axiom in \cref{equiv_quantum_action}.
Therefore, we have a linear map
\begin{align}
\label{Noether_map}
  f\colon C_*(\clL_c(V)) \lto \Obs^q_{\clM}(V)[\hbar^{-1}], \quad 
  l \lmto 
  \left.\ev*{l,\exp\left(I^{\clL}[L]/\hbar\right)}\right|_{0\times\clM(X)[1]},
\end{align}
where $\ev*{-,-}$ is the duality pairing between $C^*(\clL_c(V))$ and $C_{*}(\clL_c(V))$.
\end{proof}

%%%%%%%%%%%%%%%%%%%%%%%%%%%%%%%%%%%%%%%%%%%%%%%%%%%%%%%%%%%%%%%%%%%%%%%%%%%%%%%%
%%%%%%%%%%%%%%%%%%%%%%%%%%%%%%%%%%%%%%%%%%%%%%%%%%%%%%%%%%%%%%%%%%%%%%%%%%%%%%%%
\section{Costello--Gwilliam functor}
\label{s.3}

In this section, we review the construction of a vertex algebra from a prefactorization algebra introduced in \cite[Chapter 5]{CG}.
We extend the notion of a holomorphic classical field theory introduced in \cite[Chapter 9.2]{CG2} to free degenerate classical field theories. The factorization algebra of quantum observables of a holomorphic field theory is defined over $\bbC\dbr{\hbar}$.
Therefore, the Costello--Gwilliam functor yields a $\bbC\dbr{\hbar}$-linear vertex algebra, and to obtain a vertex algebra over $\bbC$, we need to specialize $\hbar$ to an appropriate value.
Here, we formulate precise conditions for this specialization (see \cref{lem_assumption_satisfies}).

%%%%%%%%%%%%%%%%%%%%%%%%%%%%%%%%%%%%%%%%%%%%%%%%%%%%%%%%%%%%%%%%%%%%%%%%%%%%%%%%
\subsection{Vertex algebra}
\label{ss:3.1}

We review the basic definition of a vertex algebra. For a commutative $\bbC$-algebra $R$, we introduce the definition of an $R$-linear vertex algebra.
The cases $R=\bbC$, $\bbC\dbr{\hbar}$, and $\bbC\dpr{\hbar}$ play a central role in this note, and present some basic examples. We also discuss the conditions under which $\hbar$ can be specialized to a suitable value.

Let $R$ be a $\bbC$-algebra and $V$ be an $R$-module.
A formal series
\[
a(z)=\sum_{n\in\bbZ}a_{(n)}z^{-n-1}\in(\End_RV)\dbr{z^{\pm1}}
\]
is called a field on $V$ if $a(z)v\in V\dpr{z}$ for each $v\in V$.

\begin{dfn}[{\cite[Definition 1.3.1]{FBZ}}]
\label{VA_def}
For a commutative $\bbC$-algebra $R$, an \emph{$R$-linear vertex algebra} is a tuple $(V,\ket*{0},T,Y(-,z))$ consisting of
\begin{itemize}
  \item an $R$ module $V$,
  \item an element $\ket*{0}\in V$,
  \item an $R$-linear map $T:V\to V$,
  \item an $R$-linear map
  \[
    Y(-,z)\colon V\lto(\End_R V)\dbr{z^{\pm1}},\quad 
    a\lmto Y(a,z)=\sum_{n\in\bbZ}a_{(n)}z^{-n-1},
  \]
\end{itemize}
which satisfies the following axioms:
\begin{clist}
  \item For $a\in V$, $Y(a,z)$ is a field.
  \item (Vacuum axiom) $Y(\ket*{0},z)=\mathrm{id}_V$, and $Y(a,z)\ket*{0}=a+O(z)$ for any $a\in V$.
  \item (Translation axiom) $T\ket*{0}=0$ and $[T,Y(a,z)]=\partial_zY(a,z)$ for  any $a\in V$.
  \item (Locality axiom) For $a,b\in V$, there exists a positive integer $N_{ab}\in\bbZ_{\geq0}$ such that
  \[
    (z-w)^{N_{ab}}[Y(a,z),Y(b,w)]=0
  \]
  in $(\End_R V)\dbr{z^{\pm1},w^{\pm1}}$.
\end{clist}
We simply write $V$ to denote an $R$-vertex algebra.
\end{dfn}

A morphism of $R$-linear vertex algebras can be defined in a standard way.
For simply, when $R=\bbC$, then we call a $\bbC$-linear vertex algebra $V$ a vertex algebra.

The following is a standard method for constructing a vertex algebra from another vertex algebra.

\begin{dfn}[{\cite[Chapter 5.7.2]{FBZ}}]
\label{def_center_of_VA}
Let $V$ be an $R$-vertex algebra.
Then the \emph{center} $\clZ^{VA}(V)$ of $V$ is defined by
\[
  \clZ^{VA}(V)\ceq \{v\in V \mid w_{(n)}v=0 \text{ for all } w\in V\}.
\]
\end{dfn}

\begin{example}[\cite{FBZ}]
  The following are well-known examples of vertex algebras.
  \label{example_VA}
  \def\theenumi{\roman{enumi}}
\def\labelenumi{(\theenumi)}
\begin{enumerate}
  \item For a Lie algebra $\frg$ and a 2-cocycle $\kappa\in(\Sym^2\frg^*)^\frg$, let $V^\kappa(\frg)$ be the universal affine vertex algebra 
  of level $\kappa$ \cite[Chapter 2.4]{FBZ}.
  \item For $\kappa\in\bbC$, let $\pi^\kappa$ be the $\kappa$-Heisenberg vertex algebra \cite[Chapter 2.2]{FBZ}.
  \item Let $\beta\gamma$ be the $\beta\gamma$-system \cite[Chapter 12.2.5]{FBZ}.
  \item For a Lie algebra $\frg$, let $\frz(\hat{\frg})=\clZ^{VA}(V^{-h^\vee}(\frg))$ be the \emph{Feigin--Frenkel center} \cite[Chapter 18.4.3]{FBZ}.
\end{enumerate}
\end{example}

Below, we consider the case that $R=\bbC\dbr{\hbar}$.

The following theorem guaranties that a vertex algebra can be defined by a certain set of generators and relations.
\begin{thm}[{\cite[4.4.1. Reconstruction Theorem]{FBZ}}]
\label{reconstruction_thm_Ch}
Let $R$ be a commutative $\bbC$-algebra and $V$ be an $R$-module equipped with a nonzero vector $\ket*{0}$, an endomorphism $T$, a countable ordered set $\{a^\alpha\}_{\alpha\in S}$ of vectors, and fields
\[
a^\alpha(z)=\sum_{n\in\bbZ}a^\alpha_{(n)}z^{-n-1}\in\End_{R}(V)[[z^{\pm}]]
\]
such that
\begin{clist}
\item 
For all $\alpha\in S$, $a^\alpha(z)\ket*{0}=a^\alpha+O(z)$.
\item 
$T\ket*{0}=0$ and $[T,a^\alpha(z)]=\partial_za^\alpha(z)$ for all $\alpha$.
\item 
All fields $a^\alpha(z)$ are mutually local.
\item 
$V$ is spanned by the vectors of the form
$a^{\alpha_1}_{(j_1)}\dots a^{\alpha_m}_{(j_m)}\ket*{0}$
with each $j_i<0$.
\end{clist}
Then, using the formula
\[
Y(a^{\alpha_1}_{(j_1)}\dots a^{\alpha_m}_{(j_m)}\ket*{0},z)=\frac{1}{(-j_1-1)!\dots(-j_m-1)!}:\partial^{-j_1-1}a^{\alpha_1}(z)\dots\partial^{-j_m-1}a^{\alpha_m}(z):
\]
to define a vertex operation, we obtain a unique $R$-linear vertex algebra $(V,\ket*{0},T,Y(-,z))$ such that $Y(a^\alpha,z)=a^\alpha(z)$.
In this situation, we say that $V$ is generated by $\{a^\alpha\}_{\alpha\in S}$.
\end{thm}

We now consider specialization in $\hbar$. For $a\in\bbC\setminus\{0\}$, we would like to specialize $\hbar$ to $a$. 
However, since the formal power series
$-\frac{1}{a}\sum_{n=0}^{\infty}\left(\frac{\hbar}{a}\right)^n$
is the inverse of $\hbar-a$ in $\bbC\dbr{\hbar}$, we have $V/(\hbar-a)V=0$ for any $\bbC\dbr{\hbar}$-linear vertex algebra $V$. 
Therefore, to perform the specialization, we need to impose a finiteness condition on the powers of $\hbar$.
Then, an infinite power series in $\hbar$ appearing in the OPE would be incompatible with the vertex algebra structure and the finiteness condition. 
It is therefore natural to make the following assumption.

\begin{ass}
\label{assumption_for_Chbar_VA}
In this note, we only consider the $\bbC\dbr{\hbar}$-linear vertex algebra $V$ such that the following conditions holds.
\begin{itemize}
  \item There exists a positive integer $N\geq0$, and vectors $a^\alpha\in V/\hbar V,\ (\alpha=1,\dots,N)$ such that $V$ is generated by $\{a^\alpha\}_{\alpha\in\{1,\dots,N\}}$.
  \item For every integers $n,m=1,\dots,N$, there exists a positive integer $M_{n,m}$ such that the operator product expansion takes the following form.
  \begin{align}
    \label{Ch_assumption_siki}
    a^n(z)a^m(w)\sim\sum_{k=0}^{M_{nm}}\sum_{l}\hbar^k\frac{c^{l,k}(w)}{(z-w)^{l+1}},
\end{align}
  for some field $c^{l,k}(w)\in\End_{\bbC}(V/\hbar V)[[z^\pm]]$ on $V/\hbar V$.
\end{itemize}
\end{ass}

Under the above assumption, we have the following proposition.

\begin{prop}
\label{specialize_hbar}
Let $V_\hbar$ be a $\bbC\dbr{\hbar}$-linear vertex algebra that satisfies \cref{assumption_for_Chbar_VA}.
Then, for any $a\in\bbC$, the following quotient space $V_\hbar|_{\hbar=a}$
has a structure of a vertex algebra.
\[
  V_\hbar|_{\hbar=a} \ceq
  \Bigl(\bigoplus_{n\geq0}V_\hbar/\hbar^nV_\hbar\Bigr)/(\hbar-a).
\]
\end{prop}

\begin{proof}
By \cref{assumption_for_Chbar_VA}, the vector space $V_\hbar|_{\hbar=a}$ is spanned by the monomials
\[
  a^{\alpha_1}_{(-n_1-1)}\dots a^{\alpha_l}_{(-n_l-1)}\ket*{0}  
\]
for $l\geq0$, $\alpha_1,\dots,\alpha_l\in\{1,\dots,N\}$, 
and $n_1,\dots,n_l\geq0$.
It is enough to show that for $n,m=1,\dots,N$, there exists $N_{nm}\in\bbZ_{\geq0}$ such that $(z-w)^{N_{nm}}[Y(a^n_{(-1)}\ket*{0},z),Y(a^m_{(-1)}\ket*{0},w)]=0$.
Equivalently, it is enough to show that there exist finitely many fields $c^l$ such that the following holds.
\[
  Y(a^n_{(-1)}\ket*{0},z)Y(a^m_{(-1)}\ket*{0},w) \sim
  \sum_l\frac{c^l(w)}{(z-w)^{l+1}}.
\]
But by \cref{assumption_for_Chbar_VA}, the right hand side of \eqref{Ch_assumption_siki} is a polynomial in $\hbar$.
Thus, setting $c^l(w)=\sum_{k=0}^{M_{nm}}a^kc^{l,k}(w)$, 
we have the desired OPE, and hence the assertion follows.
\end{proof}

\begin{rem}
\label{rem:cdprhbarlinear}
In \cref{ss:4.2} and \cref{ss:4.3}, we use $\bbC\dpr{\hbar}$-linear vertex algebras of the form $V[\hbar^{-1}]$, where $V$ is a $\bbC\dbr{\hbar}$-linear vertex algebra.
Since only finitely many negative powers of $\hbar$ occur in each element of $V[\hbar^{-1}]$, if $V$ satisfies \cref{assumption_for_Chbar_VA}, we can specialize $\hbar$ to $\hbar=a$ for any $a\in\bbC\backslash\{0\}$.
Furthermore, let $V,W$ be $\bbC\dbr{\hbar}$-linear vertex algebras satisfying \cref{assumption_for_Chbar_VA} and $f:V[\hbar
^{-1}]\to W[\hbar^{-1}]$ be a map of $\bbC\dpr{\hbar}$-linear vertex algebras.
We assume that, for each strong generator $v\in V$ of $V$, the element $f(v)$ belongs to $W[\hbar^{\pm1}]$.
Then for $a\in\bbC\backslash\{0\}$, the specialization 
\begin{align*}
  f|_{\hbar=a}\colon V|_{\hbar=a} \lto W|_{\hbar=a}
\end{align*}
is well-defined.
\end{rem}

%\begin{prop}(\cite{Li02}, Proposition 5.6)
%\label{hbar_adic_to_PVA}
%$V$ be a $\bbC\dbr{\hbar}$-linear vertex superalgebra such that $\ker\hbar\subset\hbar V$, and $V/\hbar V$ is a commutative vertex algebra.
%Then $(V/\hbar V,\cdot,Y_-,T)$ carries the structure of a Poisson vertex algebra where
%\begin{align*}
%  (a+\hbar V)\cdot(b+\hbar V)&=a_{(-1)}b+\hbar V,\\
%  Y_-(a+\hbar V)(b+\hbar V)&=\frac{1}{\hbar}Y(a,z)b+\hbar V\dbr{z^{\pm1}},\\
%  T(a+\hbar V)&=(Ta)+\hbar V.
%\end{align*}
%for $a,b\in V$.
%\end{prop}

%\begin{dfn}\cite{Li02}
%\label{def_deformation_of_VA}
%Let $(A,\cdot,Y_-(-,z),T)$ be a Poisson vertex algebra.
%We say that $V$ is a deformation of $A$ if there exists a $\bbC\dbr{\hbar}$-linear vertex algebra $V_\hbar$ such that $V_\hbar/\hbar V_\hbar$ is commutative, $V_\hbar/\hbar V_\hbar\simeq A$, and $V_\hbar/(\hbar-1)V_\hbar\simeq V$.
%Then we call $V_\hbar$ the $\bbC\dbr{\hbar}$-linearization of $(V,A)$.
%\end{dfn}

%For a vertex algebra $V$, we consider the $\bbC\dbr{\hbar}$-linearization of $(V,\gr^FV)$.

Given a vertex algebra $V$, the space $V\dbr{\hbar}$ naturally carries a $\bbC\dbr{\hbar}$-linear vertex algebra structure, which clearly satisfies \cref{assumption_for_Chbar_VA}. Moreover, $V\dbr{\hbar}|_{\hbar=a}\simeq V$.

We now introduce the $\bbC\dbr{\hbar}$-linear vertex algebras used in this note. These are $\bbC\dbr{\hbar}$-linear analogues of the universal affine vertex algebra and the $\beta\gamma$-system in \cref{example_VA}.

Let $\frg$ be a Lie algebra and $\kappa\in(\Sym^2\frg^\vee)^\frg$ a 2-cocycle of $\frg$.
Then we define a $\bbC\dbr{\hbar}$-linear affine Kac--Moody Lie algebra $\widehat{\frg}_\kappa^\hbar$ as a central extension
\[
  0 \lto \bbC\dbr{\hbar}\mathbf{1} \lto 
         \widehat{\frg}^\hbar_\kappa \lto \frg[t^{\pm1}]\dbr{\hbar} \lto 0.
\]
As a vector space, $\widehat{\frg}^\hbar_\kappa\simeq(\frg[t^{\pm1}]\oplus\bbC\mathbf{1})\dbr{\hbar}$, with the commutation relation $[\mathbf{1},-]=0$ and 
\[
  [x\otimes f(t),y\otimes g(t)] = \frac{1}{2\pi i}
  \left(\hbar[x,y]\otimes f(t)g(t)-\hbar\operatorname{Res}_{t=0}(fdg)\kappa(x,y)\mathbf{1}\right)
\]
for $x,y\in\frg$ and $f,g\in\bbC[t^{\pm1}]\dbr{\hbar}$.

Similarly to the usual universal affine vertex algebra, let $\bbC\dbr{\hbar}$ be a representation of $\frg[t]\dbr{\hbar}\oplus\mathbf{C}\dbr{\hbar}\mathbf{1}$ on which $\frg[t]\dbr{\hbar}$ acts trivially and $\mathbf{1}$ acts as the identity, and we define the $\bbC\dbr{\hbar}$-linear vacuum representation
\[
V^\kappa_\hbar(\frg)\ceq\operatorname{Ind}^{\widehat{\frg}^\kappa_\hbar}_{\frg[t]\dbr{\hbar}\oplus\bbC\dbr{\hbar}\mathbf{1}}\bbC\dbr{\hbar}=U(\widehat{\frg}^\kappa_\hbar)\bigotimes_{U(\frg[t]\dbr{\hbar}\oplus\bbC\dbr{\hbar}\mathbf{1})}\bbC\dbr{\hbar}.
\]
\begin{prop}
\label{hbar_adic_universal_affine_VA}
The vector space $V^\kappa_\hbar(\frg)$ has a $\bbC\dbr{\hbar}$-linear vertex algebra structure given as follows. We call $V^\kappa_\hbar(\frg)$ the $\bbC\dbr{\hbar}$-linear universal affine vertex algebra.
\begin{itemize}
\item 
The vacuum is $\ket*{0}=1\otimes1$.
\item 
The translation operator $T\colon V^\kappa_\hbar(\frg)\to V^\kappa_\hbar(\frg)$ is defined by
\[
  T\ket*{0}=0,\quad 
  [T,x\otimes f(t)]=x\otimes\partial_tf(t) \quad 
  (x\in\frg, \, f\in\bbC[t^{\pm1}]\dbr{\hbar}).
\]
The vertex operator $Y(-,z):V^\kappa_\hbar(\frg)\to(\End_{\bbC\dbr{\hbar}}V^\kappa_\hbar(\frg))\dbr{z^{\pm1}}$ is defined by
  \[
  Y(x\otimes t^{-1}\ket*{0},z)=x(z)=\sum_{n\in\bbZ}(x\otimes t^n)z^{-n-1}\ (x\in\frg).
  \]
\end{itemize}
\end{prop}

%\begin{prop}
%  \label{C[[hbar]]_linearization_of_universal_affine_VA}
%The $\bbC\dbr{\hbar}$-linear vertex algebra $V^\kappa_\hbar(\frg)$ is a $\bbC\dbr{\hbar}$-linearization of $(V^\kappa(\frg),\bbC[\clJ_\infty\frg^*])$.
%\end{prop}
%\begin{proof}
%  For $x,y\in\frg$, we have
%  \begin{align*}
%      x(z)y(w)&\sim\frac{\hbar\kappa(x,y)}{(z-w)^2}+\frac{\hbar ([x,y])(w)}{z-w},\\
%      \frac{1}{\hbar}[(x\otimes t^{-1}\ket*{0})_\lambda (y\otimes t^{-1}\ket*{0})]&=([x,y]\otimes t^{-1})\ket*{0}+\lambda\hbar\kappa(x,y)\equiv([x,y]\otimes t^{-1})\ket*{0}\ (\bmod \hbar).
%    \end{align*}
%  Thus, the quotient $V^\kappa_\hbar(\frg)/\hbar V^\kappa_\hbar(\frg)$ is commutative, and $V^\kappa_\hbar(\frg)/\hbar V^\kappa_\hbar(\frg)\simeq\bbC[\clJ_\infty\frg^*]$ is an isomorphism of Poisson vertex algebras.
%  Further, it is obvious that there exists an isomorphism $V^\kappa(\frg)/(\hbar-1)V^\kappa_\hbar(\frg)\simeq V^\kappa(\frg)$.
%\end{proof}

In particular, for an abelian Lie algebra $\frh$ and a bilinear form $\kappa:\frh\times\frh\to\bbC$, we denote
\[
  \pi^\kappa_\hbar \ceq V^\kappa_\hbar(\frh)
\]
and call it the $\bbC\dbr{\hbar}$-linear Heisenberg vertex algebra.

Consider the following Lie algebra
\[
  A_\hbar \ceq \bigoplus_{n\in\bbZ}\bbC\dbr{\hbar}a_n \oplus
  \bigoplus_{n\in\bbZ}\bbC\dbr{\hbar}a_n^* \oplus \bbC\dbr{\hbar}c'
\]
with Lie bracket is given as follows
\[
  [a_m,a_n]=0=[a_n^*,a_m^*],\quad 
  [c',-]=0,\quad 
  [a_m,a_n^*]=\frac{1}{2\pi i}\hbar\delta_{n+m,0}c' 
  \quad (m,n\in\bbZ).
\]
Denote $\bbC\dbr{\hbar}$ by the representation of the Lie subalgebra 
\[
 A_{+,\hbar} \ceq 
 \bigoplus_{n\in\bbZ_{\ge0}}\bbC\dbr{\hbar}a_n   \oplus 
 \bigoplus_{n\in\bbZ_{  >0}}\bbC\dbr{\hbar}a_n^* \oplus 
 \bbC\dbr{\hbar}c'\subset A
\] 
generated by the element $1$ satisfying 
\[
  a_n  \cdot 1=0 \ (n\ge0),\quad 
  a_m^*\cdot 1=0 \ (m  >0),\quad 
  c'\cdot1=1.
\]
Set
\[
  \beta\gamma_\hbar \ceq 
  \operatorname{Ind}^{A_\hbar}_{A_{+,\hbar}}\bbC\dbr{\hbar} = 
  U(A_\hbar) \bigotimes_{U(A_{+,\hbar})} \bbC\dbr{\hbar}.
\]
\begin{prop}
\label{def_hbar_adic_betagamm}
The vector space $\beta\gamma_\hbar$ has a $\bbC\dbr{\hbar}$-linear vertex algebra structure given as follows.
We call $\beta\gamma_\hbar$ the $\bbC\dbr{\hbar}$-linear $\beta\gamma$-system.
\begin{itemize}
  \item The vacuum is $\ket*{0}=1\otimes1$.
  \item The translation operator $T:\beta\gamma_\hbar\to\beta\gamma_\hbar$ is defined by
  \[
  T\ket*{0}=0,\ [T,a_n]=-na_{n-1},\ [T,a_n^*]=-(n-1)a^*_{n-1}.
  \]
  \item The vertex operator $Y(-,z):\beta\gamma_\hbar\to(\End_{\bbC\dbr{\hbar}}\beta\gamma_\hbar)\dbr{z^{\pm1}}$ is defined by
  \begin{align*}
    Y(a_{-1}\ket*{0},z)&=a(z)=\sum_{n\in\bbZ}a_nz^{-n-1},\\
    Y(a_0^*\ket*{0},z)&=a^*(z)=\sum_{n\in\bbZ}a_n^*z^{-n}.
  \end{align*}
\end{itemize}
\end{prop}
Let $V$ be a vector space and $\{v_i\}_{i=1,\dots,\dim V}$ be a basis of $V$.
For $i=1,\dots,\dim V$, $m\in\bbZ_{>0}$, and $n\in\bbZ_{\geq0}$, let $a_{i,m}$ and $a^*_{i,n}$ be the elements of $U(A_\hbar)^{\otimes\dim V}$ whose $i$-th tensor components are $a_m$ and $a^*_n$, respectively, and whose other components are $1$.
Then, the $\bbC\dbr{\hbar}$ vertex algebra $\beta\gamma^{\otimes\dim V}$ is spanned by the following monomials.
\begin{align}
    \label{equ:monomial betagamma}
    a_{\vec{\alpha},K}a^*_{\vec{\beta},L}\ket*{0},
\end{align}
where $\vec{\alpha}\in\{1,\dots,\dim V\}^N,\ \vec{\beta}\in\{1,\dots,\dim V\}^M,\ K\in\bbZ_{>0}^N$, and $L\in\bbZ_{\geq0}^M$ are multiindices.

%\begin{prop}
%  \label{C[[hbar]]_linearization_of_betagamma}
%The $\bbC\dbr{\hbar}$-linear vertex algebra $\beta\gamma_\hbar$ is a $\bbC\dbr{\hbar}$-linearization of $(\beta\gamma,\bbC[\clJ_\infty T^*\mathbb{A}^1])$.
%\end{prop}
%\begin{proof}
%  The proof of \cref{C[[hbar]]_linearization_of_universal_affine_VA} works with minor modifications.
%\end{proof}

\subsection{The Costello--Gwilliam construction of vertex algebras via factorization algebras}
\label{ss:3.2}
In this subsection, we review the construction of the vertex algebra from the factorization algebra on $\bbC$ that was introduced in Chapter 5.3 of \cite{CG}.

We let $(S^1)^n\ltimes\mathbb C^n\subset U(n)\ltimes \bbC^n$ be the subgroup of isometric affine transformations of $\bbC^n$, and consider its natural action on $\mathbb C$, i.e.
for $(\theta_1,\dots,\theta_n)\in (S^1)^n,\ (a_1,\dots,a_n)\in\bbC^n,$ and $(z_1,\dots,z_n)\in\bbC^n$, the action of $(S^1)^{n}\ltimes\bbC^n$ on $\bbC^n$ is given as follows
\begin{align}
  \label{isometric_action_onCn}
((\theta_1,\dots,\theta_n),(a_1,\dots,a_n))\cdot(z_1,\dots,z_n)=(e^{-i\theta_1}z_1+a_1,\dots,e^{-i\theta_n}z_n+a_n).
\end{align}

Recall that we denote by $D(\alpha,R) \subset \bbC$ the open disk  with center $\alpha$ and radius $R$.

\begin{thm}[{\cite[Theorem 5.3.3]{CG}}]
  \label{FA_to_VA}
  Let $\clF$ be an $S^1\ltimes\bbC$-equivariant prefactorization algebra on $\bbC$ satisfies the following conditions:
  \begin{itemize}
    \item For $0<R\leq\infty$, the $S^1$-action on $\clF(D(0,R))$ is tame (see Definition 5.3.2 of \cite{CG}).
    \item For $0<R<R'$ and $\Delta\in\bbZ$, the structure map
    \[
      \clF(D(0,R))_\Delta \lto \clF(D(0,R'))_\Delta
    \]
    is a quasi-isomorphism, where $\clF(D(0,R))_\Delta$ is the weight $\Delta$-subspace with respect to the $S^1$-action.
    \item For $n\neq0$, the cohomology $H^n(\clF(D(0,R)))$ vanishes.
    \item There exists a positive integer $N\in\bbZ_{\geq0}$ such that for $\Delta\geq N$, the vector space $H^0(\clF(D(0,R))_\Delta)=0$ is zero.
    \item For each $\Delta\in\bbZ$ and $R>0$, the differentiable vector space $H^0(\clF(D(0,R))_\Delta)$ is isomorphic to a countable sequantial colimit of finite dimensional graded vector spaces.
  \end{itemize}
  Then the vector space
  \[
  \bbV(\clF)=\bigoplus_{\Delta\in\bbZ} H^0(\clF(D(0,R))_\Delta)
  \]
  admits the structure of a $\bbZ$-graded vertex algebra.
  We call such an $S^1\ltimes\bbC$-equivariant prefactorization algebra $\clF:\Open(\bbC)\to\text{Ch(DVS)}$ an \emph{amenably holomorphic prefactorization algebra}.
\end{thm}

In this note, we consider the amenably holomorphic factorization algebra associated with the free degenerate classical field theory, and its higher dimensional analogue.
We consider the degenerate version of the holomorphically translation-invariant classical field theory defined in Chapter 9.2 of \cite{CG2}.

\begin{dfn}
\label{def_holomorphic_field_theory}
A free degenerate classical field theory $(\clM,\Pi,Q_{BV})$ on $\bbC^n$ is \emph{holomorphic} if the following conditions hold:
\begin{itemize}
  \item There exist graded vector spaces $V_l,\ (l=0,\dots,n)$ such that $\clM\simeq\bigoplus_{l\geq0}\Omega^{l,*}_{\bbC^n}\otimes V_l[-1]$.
  \item The differential $Q_{BV}$ coincides with the Dolbeault differential $\bar{\partial}$.
  \item The $(-1)$-shifted Poisson structure is given by a sum of holomorphic polydifferential operators.
  \item There exist $(S^1)^n$-actions on the vector spaces $V_l$.
\end{itemize}
\end{dfn}
For a holomorphic free degenerate classical field theory $(\clM,\Pi,\bar{\partial})$ on $\bbC^n$.
We can choose a gauge fixing operator as the formal adjoint of the Dolbeault differential.
\[
  Q^{GF}=\bar{\partial}^*\colon \clM \lto \clM^!.
\]
Throughout what follows, whenever we consider a quantization, we always take this gauge fixing operator.

Then we have the following assertions:
\begin{prop}
  \label{Obs_FA_property}
Let $\bigl(\clM=\bigoplus_{l\geq0}\Omega^{l,*}_{\bbC^n}\otimes V_l[-1],\Pi,\bar{\partial}\bigr)$ be a holomorphic free degenerate classical field theory on $\bbC^n$.
Then the factorization algebras $\Obs^{cl}_{\clM},\ \Obs^q_{\clM}$ of observables are $(S^1)^n\ltimes\bbC^n$-equivariant factorization algebras.
Moreover, for the case of $n=1$, if the graded vector spaces $V_0,V_1$ are concentrated in cohomological degree $0$, the factorization algebras $\Obs^{cl}_{\clM},\ \Obs^q_{\clM}$ are amenably holomorphic factorization algebras.
\end{prop}
\begin{proof}
    The proofs of Proposition 9.2.3 and Lemma 9.2.5 of \cite{CG2} works with minor modifications.
\end{proof}

The following lemma gives a sufficient condition for specializing $\hbar$ in a $\bbC\dbr{\hbar}$-linear vertex algebra constructed from a quantum field theory.

\begin{lem}
  \label{lem_assumption_satisfies}
  Let $\bigl(\clM=\bigoplus_{l\geq0}\Omega^{l,*}_{\bbC}\otimes V_l[-1],\Pi,\bar{\partial}\bigr)$ be a holomorphic free degenerate classical field theory on $\bbC$ such that the graded vector spaces $V_0,V_1$ are concentrated in a cohomological degree $0$.
  Suppose that for $n\geq2$, the $n$-th component of the $(-1)$-shifted Poisson structure $\Pi^{(n)}$ is zero.
  Then the associated $\bbC\dbr{\hbar}$-linear vertex algebra $\bbV(\Obs^q_{\clM})$ satisfies \cref{assumption_for_Chbar_VA}, and hence specializing to $\hbar=2\pi i$, we have a vertex algebra $\bbV(\Obs^q_{\clM})|_{\hbar=2\pi i}$.
\end{lem}
\begin{proof}
  Take a positive real number $R>0$.
  By the Atiyah--Bott lemma \cite[Lemma D]{CG2}, we have a quasi-isomorphism
  \begin{align*}
    \Obs^{cl}_{\clM}(D(0,R))
    &=
    \Sym\left(\ol{\left(\Omega^{0,*}_{\bbC} \otimes 
    V_0[-1]\oplus\Omega^{1,*}_{\bbC} \otimes 
    V_1[-1]\right)}^!_c(D(0,R))[-1]\right)\\
    &\qis\Sym\left(\Omega^{1,*}_{\bbC,c}(D(0,R))\otimes 
    V_0^\vee[1]\oplus\Omega^{0,*}_{\bbC,c}(D(0,R))\otimes V_1^\vee[1]\right)
  \end{align*}
  of cochain complexes.
  Here by Poincar\'{e}'s lemma, we have an isomorphism
  \begin{align*}
    H^i(\Omega^{d,*}_{\bbC,c}(D(0,R)),\bar{\partial})\simeq\begin{cases}
      \mathcal{O}_{\bbC}(D(0,R)),\quad &\text{for } i=1,\\
      0,\ &\text{otherwise},
    \end{cases}
  \end{align*}
  for $d=0,1$.
  Thus, we have an isomorphism
  \begin{align*}
    \bbV(\Obs^{cl}_{\clM})\simeq\bbC[\clJ_\infty(V_0\oplus V_1)]
  \end{align*}
  of vector spaces.
  By \cref{gr_of_q_degenerate}, we have an isomorphism
  \[
  \bbV(\Obs^q_{\clM})\simeq\bbC[\clJ_\infty(V_0\oplus V_1)]\dbr{\hbar}
  \]
  of vector spaces.
  Choose bases $\{v^0_i\}_{i=1,\dots,\dim V_0}$ and $\{v^1_j\}_{j=1,\dots,\dim V_1}$ of $V_0$ and $V_1$, respectively.
  Then the images $a^\alpha\in\bbV(\Obs^q_{\clM})/\hbar\bbV(\Obs^q_{\clM}),\ (\alpha=1,\dots,\dim V_0+\dim V_1)$ of $\{v^0_i\}_i$ and $\{v^1_j\}_j$ under the natural map $\bbC[V_0\oplus V_1]\inj\bbC[\clJ_\infty(V_0\oplus V_1)]\dbr{\hbar}$ form a strongly generator of the $\bbC\dbr{\hbar}$-linear vertex algebra $\bbV(\Obs^q_{\clM})$.
  By assumption, since for $L,L'>0$, the map $\partial_{P(L,L')}$ lowers the degree of a functional on $\clM$, for $n\geq0$, the map $e^{\hbar\partial_{P(L,L')}}$ is given by the following finite sum on $\Sym^n(\clM^*[-1])[\hbar]\subset\Obs^q_{\clM}$.
  \[
  e^{\hbar\partial_{P(L,L')}}=\sum_{k\geq0}^n\frac{1}{k!}(\hbar\partial_{P(L,L')})^k.
  \]
  Thus, for every $\alpha,\beta=1,\dots,\dim V_0+\dim V_1$, the operator product expansion $m_{z,0}(a^\alpha,a^\beta)$ satisfies the operator product expansion \eqref{Ch_assumption_siki}, and hence the $\bbC\dbr{\hbar}$-linear vertex algebra $\bbV(\Obs^q_{\clM})$ satisfies  \cref{assumption_for_Chbar_VA}.
\end{proof}

When considering the quantum action introduced in \cref{ss:2.3}, it is natural to expect that it is induced by the classical action. The following result guarantees this for holomorphic theories.

\begin{lem}[{\cite[Theorem 3.4]{anomaly}}]
\label{holomorphic_counterterm}
Let $(\clM,\Pi,Q_{BV})$ be a holomorphic field theory on $\bbC^n$, $\mathcal{L}$ be a local $L_\infty$-algebra on $\bbC^n$, and $I\in\IA(\mathcal{L},\mathcal{M})$ be a functional.
Then for $L>0$, there exists a limit
\[
\lim_{\epsilon\to0}W(P(\epsilon,L),I).
\]
\end{lem}

The following observation will be useful for the calculations in \cref{ss:3.6}.

\begin{rem}
\label{hol_theory_BV_bracket}
Let $\bigl(\clM=\bigoplus_{l\geq0}\Omega^{l,*}\otimes V_l[-1],\Pi,\bar{\partial}\bigr)$ be a holomorphic free degenerate classical field theory on $\bbC^n$.
Since for $f\in C^\infty(\bbC^n)$, we have
\begin{align}
  \label{Laplacial_on_C}
(e^{-L[\bar{\partial},\bar{\partial}^*]}f)(z)=\int_{w\in\bbC^n}\frac{1}{(2\pi iL)^n}e^{-|z-w|^2/4L}f(w)d^nw,
\end{align}
the heat kernel of the BV Laplacian $[Q_{BV},Q^{GF}]=[\bar{\partial},\bar{\partial}^*]$ is given by
\[
K_L^{H}(z,w)=\frac{1}{(2\pi iL)^{n}}e^{-|z-w|^2/4L}\prod_{i=1}^n(d\bar{z}_i-d\bar{w}_i),
\]
for $z,w\in\bbC^n$.

On the other hand, by Atiyah--Bott lemma, we have a quasi-isomorphism
\[
\Obs^q_{\clM}\qis\Sym\left(\bigoplus_l\Omega^{n-l,*}_{\bbC^n,c}\otimes V_l^\vee[n]\right)\dbr{\hbar}
\]
of precosheaves on $\bbC^n$.
For an open subset $U\subset\bbC^n$, under the above quasi-isomorphism, the cochain complex $\Obs^q_{\clM}(U)$ is generated by the linear observables
\begin{align}
\label{eq:lineobsPhi}
\begin{split}
\Phi^U\colon 
  \Omega^{n-l,*}_{\bbC^n,c}(U)\otimes V_l^\vee
  &\lto \left(\Omega^{l,*}_{\bbC^n}\otimes V_l\right)^!_c(U)\\
  \alpha\otimes\varphi
  &\lmto \left(\beta\otimes v\mapsto\left(\int_U\alpha\wedge\beta\right)\varphi(v)\right),
\end{split}
\end{align}
for $l=0,\dots,n$.
Then for $\alpha,\beta\in\Omega^{0,0}_{\bbC^n,c}(\bbC^n)$, we consider the following integral.
\begin{align*}
I=\lim_{L\to0}&\int_{z,w\in\bbC^n}\alpha(z)\beta(w)\wedge K_L(z,w)\wedge\omega_z\omega_w\\
=\lim_{L\to0}&\int_{z,w\in\bbC^n}\alpha(z)\beta(w)\frac{1}{(2\pi iL)^n}e^{-|z-w|^2/4L}\mathrm{vol}_z\mathrm{vol}_w\\
=\lim_{L\to0}&\int_{u,w\in\bbC^n}\alpha(u+w)\beta(w)\frac{1}{(2\pi iL)^n}e^{-|u|^2/4L}\mathrm{vol}_z\mathrm{vol}_w\\
=\lim_{L\to0}&\int_{u,w\in\bbC^n}\sum_{m_1,\dots,m_n,k_1,\dots,k_n\geq0}\frac{u_1^{m_1}\dots u_n^{m_n}\bar{u}_1^{k_1}\dots\bar{u}_n^{k_n}}{m_1!\dots m_n!k_1!\dots k_n!}\partial_{u_1}^{m_1}\dots\partial_{u_n}^{m_n}\partial_{\bar{u}_1}^{k_1}\dots\partial_{\bar{u}_n}^{k_n}\alpha(w)\beta(w)\\
&\ \times\frac{1}{(2\pi iL)^n}e^{-|u|^2/4L}\mathrm{vol}_z\mathrm{vol}_w
\end{align*}
where $\omega\in\Omega^{n,0}_{\bbC^n}(\bbC^n)$ is a holomorphic volume form on $\bbC^n$ and $\mathrm{vol}\in\Omega^{n,n}_{\bbC^n}(\bbC^n)$ is a volume form on $\bbC^n$.
By the formula of Gaussian integral
\[
\int_{u\in\bbC^n}\prod_{i=1}^nu_i^{p_i}\bar{u}_i^{q_i}e^{-|u|^2/4L}\mathrm{vol}=\delta_{p_1,q_1}\dots\delta_{p_n,q_n}(4\pi L)^n(4L)^{\sum_ip_i}\prod_ip_i!,
\]
for $p_i,q_i\in\bbZ,\ (i=1,\dots,n)$, we have
\begin{align*}
  I=(-2i)^n\int_{w\in\bbC^n}\alpha(w)\beta(w)\mathrm{vol}_w.
\end{align*}
Thus, the scale-$0$ BV bracket $\{-,-\}_0:=\lim_{L\to0}\{-,-\}_L$ is identified with the $P_0$-factorization algebra structure $\{-,-\}$.
\end{rem}

\subsection{\texorpdfstring{$\beta\gamma$}{beta-gamma}-system and universal affine vertex algebra}
\label{ss:3.3}
In this subsection, we briefly review the construction of the $\beta\gamma$-system and universal affine vertex algebra of \cite{CG}.
Here we consider the following classical field theory on $\bbC$:
\begin{align}\label{eq:Mbetagamma}
\begin{split}
    &\left(\clM_{\beta\gamma},\omega_{\beta\gamma},Q_{BV}^{\beta\gamma}\right)\\
    =&\left(\Omega^{0,*}[-1]\oplus\Omega^{1,*}[-1],\omega_{\beta\gamma}(\gamma_0+\beta_0,\gamma_1+\beta_1)=\int_{\bbC}\gamma_0\wedge\beta_1+\beta_0\wedge\gamma_1,Q_{BV}^{\beta\gamma}=\bar{\partial}\right).
\end{split}
\end{align}
Since the triple $(\clM_{\beta\gamma},\omega_{\beta\gamma},Q_{BV}^{\beta\gamma})$ defines a free classical field theory on $\bbC$ in the sense of \cite{CG2}, it is also a free degenerate classical field theory.
Below, we assume that the $S^1$-action on $\clM_{\beta\gamma}$ is induced by the natural action given by the formula \eqref{isometric_action_onCn}.
By \cref{Obs_FA_property}, the factorization algebras $\Obs^q_{\beta\gamma}$ and $\Obs^{cl}_{\beta\gamma}$ of observables of the field theory $(\clM_{\beta\gamma},\omega_{\beta\gamma},Q_{BV}^{\beta\gamma})$ are amenably.

To perform a quantization, we have to choose a gauge fixing operator.
Here, we choose
\[
\bar{\partial}^*=\frac{\partial}{\partial(d\bar{z})}\frac{\partial}{\partial z}.
\]

By definition, the scale $L>0$ heat kernel is a symmetric element $K_L^{\beta\gamma}\in\clM_{\beta\gamma}\otimes\clM_{\beta\gamma}$ that satisfies
\begin{align*}
\omega_{\beta\gamma}(K_L^{\beta\gamma},\varphi)=e^{-L[Q,Q^{GF}]}\varphi  
\end{align*}
for any section $\varphi$ of $\clM_{\beta\gamma}$.
By \eqref{Laplacial_on_C}, the heat kernel $K_L^{\beta\gamma}$ is given by
\begin{align}
  \label{analytic_part_of_heat_kernel_betagamma}
  K^{\beta\gamma}_L(z,w)=\frac{1}{2\pi iL}e^{-|z-w|^2/4L}(d\bar{z}-d\bar{w}).
\end{align}

Also, the propagator is defined by
\begin{align}
\label{propagator_def_betagamma}
\begin{split}
  P^{\beta\gamma}(\epsilon,L)(z,w)
  &=\int_{\epsilon}^L(Q^{GF}\otimes1)K^{\beta\gamma}_t(z,w)dt\\
  &=\int_{\epsilon}^L\frac{1}{2\pi it}\frac{\bar{z}-\bar{w}}{4t}
        e^{-|z-w|^2/4t}dt.
\end{split}
\end{align}
The limit $\epsilon\to0,L\to\infty$ plays an important role, 
which is calculated as 
\begin{align*}
  P^{\beta\gamma}(0,\infty)(z,w)
  &=\lim_{L\to\infty}\lim_{\epsilon\to0}
    \int_{\epsilon}^L\frac{1}{2\pi it}
    \frac{\bar{z}-\bar{w}}{2t}e^{-|z-w|^2/4t}dt\\
  &=\lim_{L\to0}\lim_{\epsilon\to0}\frac{1}{2\pi i}\frac{1}{z-w}
     \int_{|z-w|^2/4L}^{|z-w|^2/4\epsilon}e^{-u}du
   =\frac{1}{2\pi i}\frac{1}{z-w}.
\end{align*}
In the second line, we have made the substitution $u = |z-w|^2/4t$.

For a Stein manifold $X$ and a holomorphic vector bundle $\pi:E\to X$ on $X$, we use the Fr\'{e}chet topology on $\Hol(X,E)$, obtained as a closed subspace of $C^\infty(X,E)$.
The following lemma is useful for identifying the vertex algebras constructed from field theories in the subsequent sections.

\begin{lem}[{\cite[Theorem 5.4.7]{CG}}]
\label{key_lem_Dolbeault}
For a polydisk $D_{r_1,\dots,r_n}\subset\bbC^n$, 
we have a natural isomorphism of differentiable vector spaces
\[
  H^n(\Omega^{0,*}_c(D_{r_1,\dots,r_n})) \simeq 
  \Hol((\bbP^1\setminus D_{r_1})\times\cdots\times(\bbP^1\setminus D_{r_n}),
       \mathcal{O}(-1)^{\boxtimes n}).
\]
Moreover, all the other cohomology groups of 
$\Omega^{0,*}_c(D_{r_1,\dots,r_n})$ are zero as differentiable vector spaces.
\end{lem}

\begin{dfn}[{\cite[Definition 5.4.9]{CG}}]
  \label{def_gerators_of_Obs}
  On any disk $D(x,r)$ centered at the point $x\in\bbC$, let $c_n(x)$ denote the linear classical observables
  \[
  c_n(x):\Omega^{0,0}(D(x,r))\ni\gamma\mapsto\frac{1}{n!}(\partial^n_z\gamma)(x),
  \]
  for $n\geq0$.
  Likewise, for $n>0$, let $b_n(x)$ denote the linear functional
  \[
  b_n(x):\Omega^{1,0}(D(x,r))\ni\beta dz\mapsto\frac{1}{(n-1)!}(\partial^{n-1}_z\beta)(x).
  \]
\end{dfn}

Let introduce some notation to deal with the products of these. 
For multiindices
$K=(k_1,\dots,k_n)\in\bbZ_{>0}^n$ and $L=(l_1,\dots,l_m)\in\bbZ_{\geq0}^m$, let
\begin{align*}
  b_K(x):=b_{k_1}(x)\cdots b_{k_n}(x),\quad
  c_L(x):=c_{l_1}(x)\cdots c_{l_m}(x).
\end{align*}
Note that under the natural $S^1$-action on $H^0(\Obs^{cl}(D(x,r)))$, the elements $b_K(x)$ and $c_L(x)$ are of weight $-|K|=-\sum_{i=1}^nk_i$ and $-|L|=-\sum_{j=1}^ml_j$.

We see that the observables defined in \cref{def_gerators_of_Obs} generate the cohomology.

\begin{lem}[{\cite[Lemma 5.4.10]{CG}}]
\label{Obs_cl_generator}
These observables generate the cohomology of the classical observables in the following sense.
\begin{clist}
\item 
For $i\neq0$ the cohomology group $H^i(\Obs^{cl}_{\beta\gamma}(D(x,r)))$ vanishes as differentiable vector space.

\item 
The monomial $\{b_K(x)c_L(x)\}_{|K|+|L|=n}$ form a basis for the cohomology of the weight $n$ space $H^0(\Obs^{cl}_{\beta\gamma}(D(x,r))_n)$.
Further, they form a basis in the sense of differentiable vector spaces, 
meaning that any smooth section  
$f\colon M \to H^0(\Obs^{cl}_{\beta\gamma}(D(x,r)))_n$
can be expressed uniquely as a sum
\[
  f=\sum_{|K|+|L|=n}f_{KL}b_K(x)c_L(x),
\]
where $f_{KL}\in C^\infty(M)$ and locally on $M$ all but finitely many of the $f_{KL}$ are zero.

\item 
More generally, any smooth section
$f:M\to H^0(\Obs^{cl}_{\beta\gamma}(D(x,r)))$
can be expressed uniquely as a sum
\[
  f=\sum f_{KL}b_K(x)c_L(x)
\]
that is convergent in a natural topology on $C^\infty(M,H^0(\Obs^{cl}_{\beta\gamma}(D(x,r))))$, with the coefficient functions $f_{KL}$ satisfying the following properties:
\begin{itemize}
\item 
Locally on $M$, $f_{k_1,\dots,k_n,l_1,\dots,l_m}=0$ for $n+m\gg0$.
\item 
For fixed $n$ and $m$, the sum
\[
  \sum_{K,L}f_{KL}z_1^{-k_1}\cdots z_n^{k_n}w_1^{-l_1-1}\cdots w_m^{-l_m-1}
\]
is absolutely convergent when $|z_i|\geq r$ and $|w_j|\geq r$, in the natural topology on $C^\infty(M)$.
\end{itemize}
\end{clist}
\end{lem}

%\begin{proof}
%Let $K_{\bbC}$ be the canonical bundle on $\bbC$, and set $x=0$. 
%By definition, $\Obs^{cl}_{\beta\gamma}(D(0,r))$ is a direct sum as differentiable cochain complexes: 
%\[
%  \Obs^{cl}_{\beta\gamma}(D(0,r)) = \bigoplus_{n\ge0}
%  \left(\Omega^{0,*}_c(D(0,r))^n,E^{\boxtimes n}[n]\right)_{S_n},
%\]
%where $E$ is a holomorphic vector bundle $\mathcal{O}_{\bbC}\oplus K_{\bbC}$ on $\bbC$.

%By \cref{key_lem_Dolbeault}, we have an isomorphism
%\[
%  H^i(\Omega^{0,*}_c(D(0,r))^n,E^{\boxtimes n}) \simeq
%  \begin{cases}
%   \Hol\left((\bbP^1\setminus D(0,r))^n,(\mathcal{O}_{\bbC}(-1)\oplus \mathcal{O}_{\bbC}(-1)dz)^{\boxtimes n}\right),\ &\text{for $i=n$},\\
%   0,\ &\text{otherwise},
%  \end{cases}
%\]
%of differentiable vector spaces.
%Therefore, we have an isomorphism
%\[
%  H^0(\Obs^{cl}_{\beta\gamma}(D(0,r))) \simeq \prod_{n\ge0}
%  H^0\left((\bbP^1\setminus D(0,r))^n,(\mathcal{O}_{\bbC}(-1)\oplus \mathcal{O}_{\bbC}(-1)dz)^{\boxtimes n}\right)_{S_n},
%\]
%  of $S^1$-equivariant differentiable vector spaces.
%  Under the above isomorphism, the observable $b_K(0)c_L(0)$ goes to the function
%  \[
% \frac{1}{(2\pi i)^{n+m}}z_1^{-k_1}dz_1\cdots z_n^{-k_n}dz_nw_1^{-l_1-1}\cdots w_m^{-l_m-1}.
 % \]
  %Everything in the statement is now immediate. On $H^0(\Obs^{cl}_{\beta\gamma}(D(0,r)))$ we use the topology that is the colimit of the topologies on holomorphic functions on $(\bbP^1\setminus\ol{D(0,r-\epsilon)})^n$ as $\epsilon\to0$.
%\end{proof}

The proof of the following theorem is essentially the same as Theorem 5.4.2 of \cite{CG}.
The difference from \cite{CG} lies in the definition of the factorization algebra of observables.
In \cite{CG}, observables are defined only for free field theories, and this definition also differs from that of observables in \cite{CG2}.
As mentioned in \cref{rem:factorizationprodisdifferfromvol1}, the observables considered in this note are an extension of the definition in \cite{CG2} to the free case. Therefore, it is not clear that they coincide with observables in the sense of \cite{CG}.
%{red}CGとの違いはfactorization productを今の設定でやることです。これはCG vol2でのfactorization productのeffective interactionが丁度0の時のものと同じになっていて、CG vol1の時点ではこの定義とは異なるはずです。\cref{rem:factorizationprodisdifferfromvol1}これをcrefする
\begin{thm}
  \label{hbar_adic_betagamma_FA}
  There exists an isomorphism
  \[
  \bbV(\Obs^q_{\beta\gamma})\lsto\beta\gamma_\hbar
  \]
  of $\bbC\dbr{\hbar}$-linear vertex algebras.
\end{thm}
\begin{proof}
  We consider the limit as the scale $L\to \infty$.
  By \cref{gr_of_q_degenerate}, we have an isomorphism of prefactorization algebras
  \[
    \gr^\hbar\Obs^q_{\beta\gamma} \simeq 
    \Obs^{cl}_{\beta\gamma}\otimes\bbC\dbr{\hbar}.
  \]
  Thus, by \cref{Obs_cl_generator}, we have an isomorphism
  \[
    \bbV(\Obs^q_{\beta\gamma})\simeq\bbV(\Obs^{cl}_{\beta\gamma})\dbr{\hbar}
  \]
  of vector spaces, and the right hand side of the above isomorphism is spanned by the monomials $b_K(0)c_L(0)$.
  Hence, we have an isomorphism
  \[
  \Phi:\bbV(\Obs^q_{\beta\gamma})\lsto\beta\gamma_\hbar
  \]
  of vector spaces.
  We show that the isomorphism $\Phi$ is an isomorphism of $\bbC\dbr{\hbar}$-linear vertex algebras.
  By Reconstruction \cref{reconstruction_thm_Ch}, it is enough to show that the following OPE satisfies
  \[
  Y(b_1(0),z)Y(c_0(0),w)\sim-\frac{\hbar}{2\pi i}\frac{1}{z-w},\ Y(b_1(0),z)Y(b_1(0),w)\sim0\sim Y(c_0(0),z)Y(c_0(0),w).
  \]
  By definition of the factorization product, we have
  \begin{align*}
    m_{z,0}(b_1(0),c_0(0))&=\exp(-\hbar\partial_{P^{\beta\gamma}(0,\infty)})\left(\exp(\hbar\partial_{P^{\beta\gamma}(0,\infty)})b_1(z)*\exp(\hbar\partial_{P^{\beta\gamma}(0,\infty)})c_0(0)\right)\\
                          &=\exp(-\hbar\partial_{P^{\beta\gamma}(0,\infty)})(b_1(z)*c_0(0))\\
                          &=b_1(z)c_0(0)-\hbar\partial_{P^{\beta\gamma}(0,\infty)}(b_1(z)c_0(0))\\
                          &=b_1(z)c_0(0)-\hbar\frac{1}{2\pi i}\frac{1}{z}.
  \end{align*}
  For $f,g\in\Omega^{0,0}(D(0,1))$, the function
  \[
  (b_1(z)c_0(0))(fdz,g)=f(z)g(0)+g(z)f(0)
  \]
  is smooth at $z=0$. Thus, we have
  \[
  Y(b_1(0),z)c_0(0)\sim-\frac{\hbar}{2\pi i}\frac{1}{z}.
  \]
  Since $\partial_{P(0,\infty)}(b_1(z)b_1(0))=0=\partial_{P(0,\infty)}(c_0(z)c_0(0))$, we have
  \[
  Y(b_1(0),z)b_1(0)\sim0\sim Y(c_0(0),z)c_0(0).
  \]
\end{proof}

In what follows, we describe a factorization algebra realization of the universal affine vertex algebra.

Let $\frg$ be a Lie algebra and $\kappa\in(\Sym^2\frg^*)^\frg$ be a 2-cocycle of $\frg$.
We define a sheaf on $\bbC$ by 
\begin{align*}
  \clG_1\colon \Open(\bbC)^{op} \lto \mathrm{dgLie}, \quad 
  U \lmto \Omega^{0,*}_{\bbC}(U)\otimes\frg.
\end{align*}
For an open subset $U\subset\bbC$, we consider the following cocycle on $(\clG_1)_c(U)=\Omega^{0,*}_c(U)\otimes\frg$
\begin{align*}
  \phi(\kappa)_U(\alpha\otimes x,\beta\otimes y) \ceq \frac{1}{2\pi i}\left(\int_U\partial\alpha\wedge\beta\right)\kappa(x,y),\quad \alpha,\beta\in\Omega^{0,*}_c(U),\ x,y\in\frg.
\end{align*}
Then $\phi(\kappa)=(\phi(\kappa)_U)_{U\in\Open(\bbC)}$ defines a $(-1)$-shifted 2-cocycle of $(\clG_1)_c=\Omega^{0,*}_c\otimes\frg$, and hence we have a twisted factorization enveloping
\begin{align}
\label{def_kac_moody_FA}
  \bbU_{\phi(\kappa)}^\hbar(\clG_1)\colon \Open(\bbC) \lto \mathrm{Ch(DVS)},
\end{align}
and refer to the factorization algebra $\bbU^\hbar_{\phi(\kappa)}(\clG_1)$ on $\bbC$ as a \emph{Kac--Moody factorization algebra}.

\begin{thm}[{\cite[Chapter 5.5]{CG}}]
\label{universal_affine_FA}
There exists an isomorphism
\[
  \bbV\bigl(\bbU^\hbar_{\phi(\kappa)}(\clG_1)\bigr)\lsto V^\kappa(\frg)\dbr{\hbar}
\]
of $\bbC\dbr{\hbar}$-linear vertex algebras.
\end{thm}
\begin{proof}
Fix a positive real number $R>0$.
We define a filtration on $\bbU^\hbar_{\phi(\kappa)}(\clG_1)$ by
\[
  F_p\bbU^\hbar_{\phi(\kappa)}(\clG_1)(D(0,R)) \ceq \Sym^{\leq p}(\Omega^{0,*}_{\bbC,c}(D(0,R))\otimes\frg[1])\dbr{\hbar}.
\]
We then consider the spectral sequence $E^r_{p,q}$ associated with this filtration:
\[
  E^r_{p,q} \Longrightarrow 
  \bbV(\bbU^\hbar_{\phi(\kappa)}(\clG_1)) =
  \bigoplus_{\Delta\in\bbZ}
  H^\bullet(\bbU^\hbar_{\phi(\kappa)}(\clG_1)(D(0,R)_\Delta)).
\]
Since the differential on the $E^0$-page is given by the Dolbeault differential $\bar{\partial}$, we have
\[
  E^1_{p,q}=H^{p+q-p}(\Omega^{0,*}_c(D(0,R))^p,(\mathcal{O}_{\bbC}\otimes\frg)^{\boxtimes p})\dbr{\hbar}.  
\]
By \cref{key_lem_Dolbeault}, we have
\[
  E^1_{p,q}\simeq
  \begin{cases}
    H^0\left((\bbP^1\setminus D(0,R))^p,(\mathcal{O}_{\bbC}(-1)\otimes\frg)^{\boxtimes p}\right)_{S_p}\dbr{\hbar} & p+q=0,\\
    0 &\text{otherwise}.
  \end{cases}
\]
  Therefore, the spectral sequence $E^r$ collapses at the $E^1$-page, and hence we have an isomorphism
  \[
  \bbV(\bbU^\hbar_{\phi(\kappa)}(\clG_1))\simeq\bigoplus_{p\geq0}H^0\left((\bbP^1\setminus D(0,R))^p,(\mathcal{O}_{\bbC}(-1)\otimes\frg)^{\boxtimes p}\right)_{S_p}\dbr{\hbar}
  \]
  of vector spaces.
  For $p\geq0,n_1,\dots,n_p\geq0,x^1,\dots,x^p\in\frg$, the following function forms a basis of the right hand side of the above isomorphism
  \begin{align}
    \label{basis_of_VF_g}
  \frac{1}{(2\pi i)^p}x^1z_1^{-n_1-1}\cdots x^pz_p^{-n_p-1}.
  \end{align}
  Therefore, by sending the basis \eqref{basis_of_VF_g} to the basis $x^1_{(-n_1-1)}\dots x^p_{(-n_p-1)}\ket*{0}$, we have an isomorphism
  \[
  \bbV(\bbU^\hbar_{\phi(\kappa)}(\clG_1))\lsto V^\kappa(\frg)\dbr{\hbar}
  \]
  of vector spaces.
  It can be shown that this is an isomorphism of vertex algebras (see Chapter 5.5 of \cite{CG}).
\end{proof}

\subsection{The BV realization of the Heisenberg vertex algebra}
\label{ss:3.4}
In this subsection, we consider the Heisenberg vertex algebra. 
It is not a classical field theory in the sense of \cite{CG2}, 
but rather a degenerate classical field theory in the sense of \cite{Yoo}.

For $k\in\bbC$, we consider the tuple.
\begin{align}
  \label{Heis_field_space}
\left(\clM_{Heis}^k,\Pi_{Heis}^k,Q_{BV}^{Heis}\right)=\left(\Omega^{1,*}_{\bbC}[-1],\Pi_{Heis}^k=k\partial,Q^{Heis}_{BV}=\bar{\partial}\right).
\end{align}
By \cref{prop: WZW is degenerate classiacal}, the tuple $\bigl(\clM_{Heis}^k,\Pi_{Heis}^k,Q_{BV}^{Heis} \bigr)$ is a free degenerate classical field theory on $\bbC$.
Again, we take the gauge fixing operator $Q^{GF}=\bar{\partial}^*$.
By \cref{heat_kernel_of_degenerate_CLFT} and \eqref{Laplacial_on_C}, the heat kernel is given as follows:
\begin{align}
  \label{heat_kernel_Heis}
  K^{Heis,k}_L(z,w)=k\partial\left(\frac{1}{2\pi iL}e^{-|z-w|^2/4L}(d\bar{z}-d\bar{w})\right).
\end{align}
Also, since $\partial$ commutes with $\bar{\partial}$, the propagator is defined by
\begin{align}
  \label{propagator_Heis}
  \begin{split}
    P^{Heis,k}(\epsilon,L)(z,w)
    &=\int_{\epsilon}^L(Q^{GF}\otimes1)K_t^{Heis,k}(z,w)dt\\
    &=k\partial\int_\epsilon^L\frac{k}{2\pi it}\frac{\bar{z}-\bar{w}}{2t}e^{-|z-w|^2/4t}dt.
  \end{split}
\end{align}
The limit $\epsilon\to0,L\to\infty$ plays an important role, and is given as follows
\begin{align*}
  P^{Heis,k}(0,\infty)(z,w)
  &=\lim_{L\to\infty}\lim_{\epsilon\to0}
  k\partial\int^L_\epsilon\frac{1}{2\pi it}\frac{\bar{z}-\bar{w}}{2t}e^{-|z-w|^2/4t}dt\\
  &=k\partial\frac{1}{2\pi i}\frac{1}{(z-w)}\\
  &=\frac{k}{2\pi i}\frac{1}{(z-w)^2}(dz-dw).
\end{align*}

The factorization algebra of quantum observables of the degenerate classical field theory is defined as follows:
\begin{align*}
 \Obs^q_{Heis:k}
=\Sym((\clM^k_{Heis})^*[-1])
=\Sym((\ol{(\clM^k_{Heis})_c})^![-1]),
\end{align*}
equipped with the differential $\bar{\partial}+\hbar\ev*{-,\Pi^k_{Heis}(-)}$.
By the Atiyah--Bott lemma, we have the following quasi-isomorphism:

\begin{lem}
\label{obs_Heis_quasi-isom}
There exists a quasi-isomorphism
\begin{align*}
\Obs^{cl}_{Heis:k}\qis\left(\Sym(\Omega^{0,*}_c[1])\dbr{\hbar},\bar{\partial}\right).  
\end{align*}
of factorization algebras on $\bbC$.
\end{lem}

Then we have the following assertion.

\begin{thm}
\label{Heis_FA}
There exists an isomorphism
\[
  \bbV\left(\Obs^q_{Heis:k}\right) \lsto \pi^k_\hbar
\]
of $\bbC\dbr{\hbar}$-linear vertex algebras.
\end{thm}
\begin{proof}
We consider the limit as the scale $L\to \infty$.
For $x\in\bbC,r>1,n\geq0$, we define the linear classical observable
\[
  u_n(x)\colon \Omega^{1,0}(D(x,r))\ni\beta dz \lmto 
  \frac{1}{(n-1)!}(\partial^{n-1}_z\beta)(x).
\]
Then similarly to the proof of \cref{Obs_cl_generator}, the vector space
$\bbV(\Obs^{cl}_{Heis,k})$
is spanned by the monomials $u_K(0)$.
By \cref{gr_of_q_degenerate}, we have an isomorphism
\[
  \bbV(\Obs^q_{Heis,k}) \simeq 
  \bbV(\Obs^{cl}_{Heis,k})\otimes\bbC\dbr{\hbar} \simeq 
  \pi_\hbar^k
\]
of vector spaces. We show that the above isomorphism is an isomorphism of $\bbC\dbr{\hbar}$-linear vertex algebras.
By Reconstruction \cref{reconstruction_thm_Ch}, it is enough to show that the following OPE satisfies
\[
  Y(u_1(0),z)Y(u_1(0),w)\sim\frac{\hbar k}{2\pi i}\frac{1}{(z-w)^2}.
\]
Indeed, we have
\begin{align*}
  &m_{z,0}(u_1(0),u_1(0))\\
  &=\exp(\hbar\partial_{P^{Heis,k}(0,\infty)})\left(\exp(-\hbar\partial_{P^{Heis,k}(0,\infty)})u_1(z)*\exp(-\hbar\partial_{P^{Heis,k}(0,\infty)})u_1(0)\right)\\
  &=\exp(\hbar\partial_{P^{Heis,k}(0,\infty)})(u_1(z)*u_1(0))\\
  &=u_1(z)u_1(0)+\hbar\partial_{P^{Heis,k}(0,\infty)}(u_1(z)u_1(0))\\
  &=u_1(z)u_1(0)+\hbar\frac{k}{2\pi i}\frac{1}{z^2}.
\end{align*}
Since, similarly to the proof of \cref{hbar_adic_betagamma_FA}, the first term of the above equality is smooth at $z=0$, we have
\[
  Y(u_1(0),z)u_1(0)\sim\frac{\hbar k}{2\pi i}\frac{1}{z^2}.
\]
\end{proof}

Let $\frh$ be a vector space, and $\kappa\colon\frh\times\frh\to\bbC$ be a symmetric bilinear form.
Take a basis $\{h_i\}_i$ of $\frh$ and a dual basis $\{h_i^*\}_i$ of $\frh^\vee$.
Then we have an isomorphism $(-)^*\colon\frh^\vee\stackrel{\sim}{\to}\frh$ given by $(h_i^*)^*=h_i$.
More generally, we consider the following degenerate free classical field theory on $\bbC$:
\[
(\clM^\kappa_{Heis},\Pi^\kappa_{Heis},Q_{BV}^{Heis}) \ceq \left(\Omega^{1,*}_{\bbC}\otimes\frh[-1],\Pi_{Heis}^\kappa,Q_{BV}^{Heis}=\bar{\partial}\otimes1\right),
\]
where the $(-1)$-shifted Poisson structure $\Pi^\kappa_{Heis}$ is given by
\begin{align*}
  \Pi^\kappa_{Heis}\colon
  \Omega^{0,*}_{\bbC}\otimes\frh^\vee 
  \lto \Omega^{1,*}_{\bbC}\otimes\frh, \quad 
  \omega\otimes\chi \lmto \partial\omega\otimes\chi^*.
\end{align*}

Then we have the following assertion:
\begin{thm}
\label{Ch_kappa_Heis_FA}
There exists an isomorphism
\[
\bbV(\Obs^q_{Heis,\kappa})\lsto\pi^\kappa_\hbar
\]
of $\bbC\dbr{\hbar}$-linear vertex algebras.
\end{thm}
\begin{proof}
The proof of \cref{Heis_FA} works with minor modifications.
\end{proof}

\subsection{The universal affine vertex algebra via Chiral WZW boundary degenerate field theory}
\label{ss:3.5}

In this subsection, let $\frg$ be a semisimple Lie algebra, $(-,-)\colon\frg\otimes\frg\to\bbC$ be a normalized form, and a symmetric bilinear form $\kappa\in(\Sym^2\frg^\vee)^\frg$.

In this subsection, we consider the following tuple:
\begin{align}
  \label{WZW_kappa_def}
  (\clM_{WZW}^\kappa,\Pi_{WZW}^\kappa,Q_{BV}^{WZW})=\left(\Omega^{1,*}_{\bbC}[-1]\otimes\frg^\vee,\Pi_{WZW}^\kappa,Q_{BV}^{WZW}=\bar{\partial}\otimes\mathrm{id}\right),
\end{align}
where the $(-1)$-shifted Poisson structure $\Pi_{WZW}^\kappa$ is given as follows
\begin{align*}
  \Pi_{WZW}^{(0)} \ceq \partial\otimes\kappa(-,-)&\colon 
  \Omega^{0,*}_{\bbC}\otimes\frg\lto
  \Omega^{1,*}_{\bbC}\otimes\frg^\vee,\\
  \Pi_{WZW}^{(1)} \ceq \ad^*(-)(-)&\colon
  \Omega^{0,*}_{\bbC}\otimes\frg\otimes
  \Omega^{1,*}_{\bbC}\otimes\frg^\vee\lto
  \Omega^{1,*}_{\bbC}\otimes\frg^\vee,\\
  \Pi^{(k)}_{WZW} \ceq 0 \qquad & \text{for $k \ne 0,1$},
\end{align*}
where we used the notation \eqref{equ:componentofPoisson}, $\kappa(-,-)\colon\frg\to\frg^*$ is defined by $x\mapsto\kappa(x,-)$, and $\ad^*:\frg\to\mathfrak{gl}(\frg^\vee)$ is the coadjoint representation of $\frg$.
\begin{prop}[{\cite[Proposition 4.4]{Yoo}}]
    \label{prop: WZW is degenerate classiacal}
    The tuple $\bigl(\clM_{WZW}^\kappa,\Pi_{WZW}^\kappa,Q_{BV}^{WZW}\bigr)$ is a degenerate classical field theory on $\bbC$.
\end{prop}
Take a basis $\{J^a\}_{a=1,\dots,\dim\frg}$ of $\frg$ and the dual basis $\{g^a\}_{a=1,\dots,\dim\frg}$ of $\frg^*$, $\{c_{ab}^c\}$ be the corresponding structure constant.
Then for $a,b,c=1,\dots,\dim\frg$, we have
\[
\Pi_{WZW}^{(1)}(\xi\otimes J^a\otimes\omega\otimes g^b)\colon\frg\ni J^c\mapsto\ad^*(J^a)(g^b)(J^c)\xi\wedge\omega=c_{bc}^a\xi\wedge\omega\in\Omega^{1,*}_{\bbC}.
\]
Thus, by \cref{heat_kernel_of_degenerate_CLFT} and \eqref{Laplacial_on_C}, we have
\begin{align*}
  K_L^{(0),WZW}(z,w)&=\partial\left(\frac{1}{2\pi iL}e^{-|z-w|^2/4L}(d\bar{z}-d\bar{w})\right)\otimes\kappa(-,-),\\
  K_L^{(1),WZW}(z,w)\left(\sum_{a}\omega_a\otimes g^a\right)&=\sum_{a,b,c}c_{bc}^a\frac{1}{2\pi iL}e^{-|z-w|^2/4L}\omega_a(z-w)\wedge(d\bar{z}-d\bar{w})\otimes \left(g^b\otimes g^c\right),\\
  K_L^{(k),WZW}&=0,\quad \text{for }k\neq0,1.
\end{align*}
Thus, the propagator is given as follows:
\begin{align*}
  P^{(0),WZW}(0,\infty)(z,w)&=\frac{1}{2\pi i}\frac{1}{(z-w)^2}(dz-dw)\otimes\kappa(-,-),\\
  P^{(1).WZW}(0,\infty)(z,w)\left(\sum_{a}\omega_a\otimes g^a\right)&=\sum_{a,b,c}c_{ab}^c\frac{1}{2\pi i}\frac{1}{z-w}\omega_a(z-w)\otimes \left(g^b\otimes g^c\right),\\
  P^{(k),WZW}(0,\infty)&=0,\quad \text{for }k\neq0,1.
\end{align*}
We denote by $\Obs^q_{WZW,\kappa}$ and $\Obs^{cl}_{WZW,\kappa}$ the factorization algebra of quantum and classical observables, respectively.

\begin{thm}
\label{WZW_VA}
There exists an isomorphism
\[
\bbV(\Obs^q_{WZW,k})\simeq V^{\kappa}_\hbar(\frg)
\]
of $\bbC\dbr{\hbar}$-linear vertex algebras.
In particular, specializing to $\hbar=2\pi i$, we have an isomorphism
\[
\bbV(\Obs^q_{WZW,k})|_{\hbar=2\pi i}\simeq V^k(\frg)
\]
of vertex algebras.
Moreover, there exists a quasi-isomorphism
\[
  \Obs^q_{WZW,\kappa}\qis\bbU^\hbar_{\hbar\phi(\kappa)}(\clG_1^\hbar)
\]
of factorization algebras, where $\clG_1^\hbar=\Omega^{0,*}_{\bbC}\otimes\frg\dbr{\hbar}$ is a sheaf of dgLA on $\bbC$ with differential $\bar{\partial}$ and the Lie bracket given by
\[
  [\alpha\otimes x,\beta\otimes y]=\hbar\alpha\wedge\beta\otimes[x,y],
\]
for sections $\alpha\otimes x,\beta\otimes y$ of $\Omega^{0,*}_{\bbC}\otimes\frg\dbr{\hbar}$.
\end{thm}
\begin{proof}
At first, we show that there exists an isomorphism 
$\bbV(\Obs^q_{WZW,k})\simeq V^{\kappa}_\hbar(\frg)$ 
of $\bbC\dbr{\hbar}$-linear vertex algebras.
We consider the limit as the scale $L\to \infty$.
By the Atiyah--Bott lemma, we have the following quasi-isomorphism
\begin{align}
\label{isom_WZW}
\begin{split}
  \Obs^{cl}_{WZW,\kappa}=\Sym((\clM_{WZW}^k)^*[-1])
  =\Sym(\ol{(\clM^k_{WZW})_c}^![-1])
  \qis\Sym\left(\Omega^{0,*}_c\otimes\frg[1]\right).
\end{split}
\end{align}
For $x\in\bbC$, $r>0$, $n\ge1$, and $a=1,\dotsc,\dim\frg$, 
we define the linear classical observable
\[
  J^a_n(x)\colon 
  \Omega^{1,0}(D(x,r))\otimes\frg^\vee \ni \beta dz\otimes g
  \lmto \frac{1}{(n-1)!}(\partial^{n-1}_z\beta)(x)g(J^a).
\]
Then, similarly to the proof of \cref{Obs_cl_generator}, 
the vector space $\bbV(\Obs^{cl}_{WZW,\kappa})$ is spanned by the monomials
\[
  J^{a_1}_{n_1}(0)\cdots J^{a_m}_{n_m}(0),
\]
for $a_1,\dots,a_m=1,\dots,\dim\frg,\ n_1,\dots,n_m\geq1$.
By \cref{gr_of_q_degenerate}, we have an isomorphism
\[
  \bbV(\Obs^q_{WZW,\kappa}) \simeq 
  \bbV(\Obs^{cl}_{WZW,\kappa}) \otimes 
  \bbC\dbr{\hbar}\simeq V^{\kappa}_{\hbar}(\frg)
\]
of vector spaces.
We show that the above isomorphism is an isomorphism of $\bbC\dbr{\hbar}$-linear vertex algebras.
By \cref{reconstruction_thm_Ch}, it is enough to show that the following OPE satisfies
\[
  Y(J^a_1(0),z)Y(J^b_1(0),w)\sim\frac{\hbar\kappa(J^a,J^b)}{(z-w)^2}+\frac{\hbar\sum_cc_{ab}^cJ^c_1(w)}{z-w},
\]
for $a,b=1,\dots,\dim\frg$.
  
Indeed, we have
\begin{align*}
    &m_{z,0}(J^a_1(0),J^b_1(0))\\
    =&\exp\left(\hbar\partial_{P^{WZW}(0,\infty)}\right)\left(\exp\left(-\hbar\partial_{P^{WZW}(0,\infty)}\right)J^a_1(z)*\exp\left(-\hbar\partial_{P^{WZW}(0,\infty)}\right)J^b_1(0)\right)\\
                              =&\exp\left(\hbar\partial_{P^{WZW}(0,\infty)}\right)(J^a_1(z)*J^b_1(0))\\
                              =&J^a_1(z)J^b_1(0)+\hbar\partial_{P^{WzW}(0,\infty)}(J^a_1(z)J^b_1(0))\\
                              =&J^a_1(z)J^b_1(0)+\hbar\partial_{P^{(0),WZW}(0,\infty)}(J^a_1(z)J^b_1(0))+\hbar\partial_{P^{(1),WZW}(0,\infty)}(J^a_1(z)J^b(0)).
\end{align*}
The first term of the right hand side is a smooth function at the origin.
The second term of the right hand side is a constant observable and is given as follows
\[
  \hbar\frac{\kappa(J^a,J^b)}{2\pi i}\frac{1}{z^2}.
\]
The third term of the right hand side is a linear observable and is given as follows
\begin{align*}
   &\hbar\partial_{P^{(1),WZW}(0,\infty)}(J^a_1(z)J^b_1(0))(fdz\otimes g^c)\\
   &=\hbar\partial_{P^{(1),WZW}(0,\infty)(fdz\otimes g^c)}(J^a_1(z)J^b_1(0))\\
   &=\frac{\hbar}{2\pi i}\frac{c_{ab}^c}{z}f(z)\sim-\frac{\hbar}{2\pi i}\frac{c_{ab}^c}{z}f(0)=-\hbar\frac{\sum_cc_{ab}^cJ^c_1(0)(fdz)}{z},
\end{align*}
for $f\in\Omega^{0,0}(\bbC),\ c=1,\dots,\dim\frg$.
Therefore, we have the following OPE
\[
  Y(J_1^a(0),z)J^b_1(0) \sim 
  \frac{\hbar}{2\pi i}\frac{\kappa(J^a,J^b)}{z^2} + 
  \frac{\hbar}{2\pi i}\frac{\sum_cc_{ab}^cJ^c_1(0)}{z}.
\]

Finally, we show that there exists an isomorphism $\Obs^q_{WZW,\kappa}\simeq\bbU^\hbar_{\phi(\kappa)}(\clG^\hbar_1)$ of factorization algebras.
By the proof of \cref{universal_affine_FA}, we have a quasi-isomorphism
\[
  \Obs^q_{WZW,\kappa}(U)\qis\bbU^\hbar_{\phi(\kappa)}(\clG^\hbar_1)(U),
\]
for each open subset $U\subset\bbC$.
The above proof implies that this quasi-isomorphism is a map of of factorization algebras.
\end{proof}

%%%%%%%%%%%%%%%%%%%%%%%%%%%%%%%%%%%%%%%%%%%%%%%%%%%%%%%%%%%%%%%%%%%%%%%%%%%%%%%%
\subsection{BD center and Vertex algebra center}
\label{ss:3.6}

In this subsection, we introduce the notion of a \emph{BD center} of BD algebra and a BD prefactorization algebra, and consider the relationship between the center of $\bbC\dbr{\hbar}$-linear vertex algebra and the BD center.

\begin{dfn}
\label{BV_center_def}
For a BD algebra $(A,Q,\{-,-\})$, we define
\[
\clZ^{BD}(A)\ceq \{a\in V \mid \{a,b\}=0 \text{ for all } b\in A\},
\]
and call $\clZ^{BD}(A)$ the \emph{BD center} of $A$.

Similarly, for a BD factorization algebra $\clF$ on a manifold $X$, we define a precosheaf
\begin{align*}
  \clZ^{BD}(\clF)\colon \Open(X) \lto \mathrm{pCh(DVS)}, \quad 
  U \lmto \clZ^{BD}(\clF(U)),
\end{align*}
and call $\clZ^{BD}(\clF)$ the \emph{BD center} of $\clF$.
\end{dfn}

For a general BD factorization algebra, it is nontrivial that its BD center is again a factorization algebra. However, for the factorization algebra of observables, the following holds.

\begin{lem}
  \label{hol_center_is_hol}
  Let $(\clM,\Pi,\bar{\partial})$ be a free degenerate holomorphic classical field theory on a manifold $X$.
  Then the BD center $\clZ^{BD}(\Obs^q_{\clM})$ of the factorization algebra $\Obs^q_{\clM}$ of quantum observables is a sub-factorization algebra of $\Obs^q_{\clM}$.
  Moreover, the factorization algebra $\clZ^{BD}(\Obs^q_{\clM})$ is an amenably holomorphic factorization algebra on $\bbC$.
\end{lem}
\begin{proof}
  Let $U,V\subset\bbC$ be disjoint open subsets, $U\sqcup V\subset W\subset\bbC$ be an open subset, and $\{O[L]\}_L\in\clZ^{BD}(\Obs^q_{\clM})(U)$, $\{O'[L]\}_L\in\clZ^{BD}(\Obs^q_{\clM})(V)$, and $\{O''[L]\}_L\in\clZ^{BD}(\Obs^q_{\clM})(W)$ be quantum observables.
  Since the BV brackets $\{(m^{\Obs}_L)^U_W(O[L]),O''[L]\}_L$ and $\{(m^{\Obs}_L(O'[L])),O''[L]\}_L$ vanish, the BV bracket $\{(m^{\Obs}_L)^{U,V}_W(O[L],O'[L]),O''[L]\}_L$ also vanishes.
  Therefore, the factorization product $m^{\Obs}:\Obs^q_{\clM}\times\Obs^q_{\clM}\to\Obs^q_{\clM}$ defines a factorization product in $\clZ^{BD}(\Obs^q_{\clM})$, and therefore the precosheaf $\clZ^{BD}(\Obs^q_{\clM})$ is a sub-factorization algebra of $\Obs^q_{\clM}$.
  Here, the BV bracket is $S^1$-invariant and compatible with the differential $\bar{\partial}+\hbar\Delta_L$, and hence the properties in \cref{FA_to_VA} are preserved.
\end{proof}

In what follows, we compute the vertex algebra associated with the BD center.

\begin{prop}
\label{prop_commutant_betagamma}
There exists an isomorphism
\[
  \clZ^{VA}(\beta\gamma_\hbar)\simeq\bbC\dbr{\hbar}\simeq\bbV(\clZ^{BD}(\Obs^q_{\beta\gamma}))
\]
of $\bbC\dbr{\hbar}$-linear (commutative) vertex algebras.
\end{prop}

\begin{proof}
  We consider the limit as the scale $L\to 0$.
  Recall that for $z\in\bbC, R>0$, we have a quasi-isomorphism
  \[
  \Obs^q_{\beta\gamma}(D(z,R))\qis\Sym\left(\Omega^{0,*}_{\bbC,c}[1]\oplus\Omega^{1,*}_{\bbC,c}[1]\right)(D(z,R))\dbr{\hbar}
  \]
  of cochain complexes.
  By \cref{hol_theory_BV_bracket}, under the above identification, for $\gamma\in\Omega^{0,*}_{\bbC,c}(D(z,R))$ and $\beta\in\Omega^{1,*}_{\bbC,c}(D(z,R))$, we have
  \begin{align*}
    \{\gamma,\beta\}_0=\int_{D(z,R)}\frac{d}{d(d\bar{z})}\gamma\wedge\frac{d}{d(d\bar{z})}\beta\wedge d\bar{z},\quad \{\gamma,\gamma\}_0=0=\{\beta,\beta\}_0.
  \end{align*}
 Dou to the nondegeneracy of the integration pairing, the BD center $\clZ^{BD}(\Obs^{q}_{\beta\gamma})(D(z,R))$ is isomorphic to $\bbC\dbr{\hbar}$.
\end{proof}

\begin{prop}
  \label{prop_commutant_Heis}
  For a vector space $\frh$, and a symmetric bilinear form $\kappa:\frh\times\frh\to\bbC$, there exists an isomorphism
  \[
  \bbV(\clZ^{BD}(\Obs^q_{Heis,\kappa}))\lsto\clZ^{VA}(\pi^\kappa_\hbar)
  \]
  of $\bbC\dbr{\hbar}$-linear (commutative) vertex algebras.
\end{prop}
\begin{proof}
    We consider the limit as the scale $L\to 0$.
    Recall that for $z\in\bbC, R>0$, we have a quasi-isomorphism
    \[
    \Obs^q_{Heis,\kappa}(D(z,R))\qis\Sym\left(\Omega^{0,*}_{\bbC,c}\otimes\frh[1]\right)(D(z,R))\dbr{\hbar}
    \]
    of cochain complexes.
    By \cref{hol_theory_BV_bracket}, under the above identification, for $u\otimes h\in\Omega^{0,*}_{\bbC,c}(D(z,R))$, we have
    \[
    \{u\otimes h,u\otimes h\}_0=\left(\int_{D(z,R)}\frac{d}{d(d\bar{z})}u\wedge\partial \frac{d}{d(d\bar{z})}u\wedge d\bar{z}\right)\kappa(h,h).
    \]
    Again, by the nondegeneracy of the integration pairing and the Leibniz rule of the BV bracket, the BD center $\clZ^{BD}(\Obs^q_{Heis,\kappa})(D(z,R))$ is isomorphic to the subcomplex
    \[
    \Sym\left(\ker{\partial}\otimes\frh[1]+\Omega^{0,*}_{\bbC,c}\otimes\frh^{\perp,\kappa}[1]\right)(D(z,R))\dbr{\hbar},
    \]
    where we denote $\frh^{\perp,\kappa}=\{x\in\frh\mid\kappa(h,x)=0,\text{ for all }h\in\frh\}$.
    Here, we consider the subcomplex
    \[
    \ker\partial\subset\Omega^{0,*}_{\bbC,c}(D(z,R)).
    \]
    The complex $\ker\partial$ consists of Dolbeault forms with compactly supported anti-holomorphic function coefficients.
    For a compactly supported anti-Holomorphic functional, there exists a real number $0<r<R$ such that $\supp{f}\subset \ol{D(z,r)}\subsetneq D(z,R)$.
    Thus, by the identity theorem, the function $f$ vanishes identically, and hence the complex $\ker\partial$ also vanishes.
    Therefore, we have
    \[
    \clZ^{BD}(\Obs^q_{Heis,\kappa})\simeq\Sym(\Omega^{0,*}_{\bbC,c}\otimes\frh^{\perp,\kappa}[1])\dbr{\hbar}.
    \]
    Applying \cref{Ch_kappa_Heis_FA} with $\frh=\frh^{\perp,\kappa}$ and $\kappa=\kappa|_{\frh^{\perp,\kappa}}=0$, we have an isomorphism
    \[
    \bbV(\clZ^{BD}(\Obs^q_{Heis,\kappa}))\simeq\bbC[\clJ_\infty(\frh^{\perp,\kappa})^\vee]\dbr{\hbar}\simeq\clZ^{VA}(\pi^\kappa_\hbar)
    \]
\end{proof}

\begin{prop}
  \label{prop_commutant_WZW}
  For a semisimple Lie algebra $\frg$ and a symmetric invariant bilinear form $\kappa\in(\Sym^2\frg^\vee)^\frg\dbr{\hbar}$, there exists an isomorphism
  \[
  \bbV(\clZ^{BD}(\Obs^q_{WZW,\kappa}))\lsto\clZ^{VA}(V^\kappa_\hbar(\frg))
  \]
  of $\bbC\dbr{\hbar}$-linear (commutative) vertex algebras.
\end{prop}
\begin{proof}
  We consider the limit as the scale $L\to 0$.
  By the Leibniz rule of the BV bracket, for an open subset $U\subset\bbC$, we have
  \[
  \clZ^{BD}(\Obs^q_{WZW,\kappa})(U)\simeq\left\{O\in\Obs^q_{WZW,\kappa}(U) \mid \{\Phi,O\}_0=0 \text{ for all linear observable } \Phi\right\}.
  \]
  Then, by Atiyah--Bott lemma, we have a quasi-isomorphism
  \[
  \clM^*_{WZW,\kappa}(U)\qis\Omega^{0,*}_{\bbC,c}(U)\otimes\frg[2]
  \]
  of cochain complexes.
  Under the quasi-isomorphism, for $\alpha\otimes x\in\Omega^{0,*}_{\bbC,c}(U)\otimes\frg$, we have the associated linear observable
  \begin{align*}
   \Phi^U(\alpha\otimes x)\colon\clM_{WZW,\kappa}(U)
   =\Omega^{1,*}_{\bbC}(U)\otimes\frg^\vee\lto\bbC,\quad 
   \beta\otimes\varphi\lmto\left(\int_U\alpha\wedge\beta\right)\varphi(x).
  \end{align*}

  By \cref{WZW_VA}, the vertex algebra $\bbV(\Obs^q_{WZW,\kappa})$ is spanned by the monomial
  \[
  J^{\vec{a}}_{\vec{n}}(0)\ceq J^{a_1}_{n_1}(0)\cdots J^{a_m}_{n_m}(0)
  \]
  for $a_1,\dots,a_m=1,\dots,\dim\frg,\ n_1,\dots,n_m\geq1$.

  For $n\geq1,\ R>0$, let $f^R\in C^\infty_c(\bbR_{\ge0})$ be a compactly support smooth function on $\bbR_{\ge0}$ such that there exists an open neighborhood of $0\in\bbR_{\ge0}$ such that $f^R$ vanishes on $U$, $\supp{f}^R\subset[0,R)$, and  $\int_{\bbR_{\ge0}}rf^R(r)dr=1$.
  Consider the linear observable $\clJ^{a,R}_n(0):=\Phi^{D(z,R)}(f^R(|z|)z^{-n+1}d\bar{z}\otimes J^a)$.
  We find that it defines a representative of the cohomology class $J^a_{n}(0)\in H^0(\Obs^q_{WZW,\kappa}(D(0,R)))$.
  Indeed, for a holomorphic $(1,0)$-form $\beta dz$ on $D(0,R)$, we have
  \begin{align*}
    \clJ^{a,R}_n(0)(\beta dz\otimes\varphi)&=\left(\int_{D(0,R)}f^R(|z|)z^{-n+1}\beta(z)d^2z\right)\varphi(J^a)\\
   &=\varphi(x)\int_0^R\int^{2\pi}_0r^{-n+1}e^{-i(n-1)\theta}f^R(r)\sum_{k\geq0}\frac{1}{k!}(\partial^k_z\beta)(0)r^ke^{ik\theta}rdrd\theta\\
 &=\varphi(x)\int_0^Rf^R(r)r\frac{1}{(n-1)!}(\partial^{n-1}_z\beta)(0)dr\\
    &=\frac{1}{(n-1)!}(\partial^{n-1}_z\beta)(0)\varphi(x).
  \end{align*}
  Hence, the vertex algebra $\bbV(\clZ^{BD}(\Obs^q_{WZW,\kappa}))$ is spanned by the monomial
  \[
    J^{\vec{a}}_{\vec{n}}(0)=J^{a_1}_{n_1}(0)\cdots J^{a_m}_{n_m}(0),\qquad a_i=1,\dots,\dim\frg,\ n_i\geq1
  \]
  such that the representative $\clJ^{\vec{a},R}_{\vec{n}}(0):=\clJ^{a_1,R}_{n_1}(0)\cdots\clJ^{a_m,R}_{n_m}(0)$ satisfies the following relations
  \[
  \{\Phi^{D(0,R)}(\alpha\otimes x),\clJ^{\vec{a},R}_{\vec{n}}(0)\}^{(0)}_0=0=\{\Phi^{D(0,R)}(\alpha\otimes x),\clJ^{\vec{a},R}_{\vec{n}}(0)\}^{(1)}_0,
  \]
  for every $\alpha\otimes x\in\Omega^{0,*}_{\bbC,c}(D(0,R))\otimes\frg$.

  Then for $\alpha\otimes x\in\operatorname{Im}\bar\partial\subset\Omega^{0,1}_{\bbC,c}(U)$ and $\beta\otimes x\in\ker\bar{\partial}$, by the Stokes lemma, we have
  \[
  \{\Phi^U(\alpha\otimes x),\Phi^U(\beta\otimes y)\}_0=0.
  \]
Furthermore, the subcomplex $\ker\partial\subset\Omega^{0,*}_{\bbC,c}$ vanishes.
Thus, the action
\begin{align*}
  \wt{\{-,-\}}_0\colon 
  \Omega^{0,*}_{\bbC,c}(D(0,R))\otimes\frg\times\bbV(\Obs^q_{WZW,\kappa}) 
  \lto \bbV(\Obs^q_{WZW,\kappa})
\end{align*}
given by
\[
  \wt{\{\alpha\otimes x,J^{\vec{a}}_{\vec{n}}(0)\}}_0 = 
  \{\Phi^{D(0,R)}(\alpha\otimes x),\clJ^{\vec{a},R}_{\vec{n}}(0)\}_0
\]
coincides with the induced action on cohomology:
\begin{align*}
  H^*(\Omega^{0,*}_{\bbC,c}(D(0,R)))\otimes\frg\times\bbV(\Obs^q_{WZW,\kappa}) 
  \lto \bbV(\Obs^q_{WZW,\kappa}).
\end{align*}
Then, by Poincar\'{e}'s lemma, we have
\begin{align*}
  H^1(\Omega^{0,*}_{\bbC,c}(D(0,R)),\bar{\partial}) \simeq
  H^0(\Omega^{0,*}_{\bbC}(D(0,R)),\bar{\partial})d\bar{z}\otimes\frg
  \simeq \mathcal{O}_{\bbC}(D(0,R))d\bar{z}\otimes\frg
  \simeq \frg[z]
\end{align*}
and
\begin{align*}
  H^i(\Omega^{0,*}_{\bbC,c}(D(0,R)),\bar{\partial}) = 0 \quad (i \ne 1).
\end{align*}
Therefore, we have an isomorphism
\[
  \bbV(\clZ^{BD}(\Obs^q_{WZW,\kappa}))\simeq\bbV(\Obs^q_{WZW,\kappa})^{\frg[z]}
\]
of vector spaces.
By \cref{hol_theory_BV_bracket}, note that for $\alpha\otimes x,\beta\otimes y\in\mathcal{O}_{\bbC}(D(0,R))\otimes\frg$, the components of the BV bracket are given as follows
\begin{align*}
  &\{\Phi^{D(0,R)}(\alpha d\bar{z}\otimes x),\Phi^{D(0,R)}(\beta d\bar{z}\otimes y)\}^{(0)}_0\\
  &=\lim_{L\to0}W_0^L\left(\int_{D(0,R)}\alpha(z)\partial_z\beta(z)d^2z\right)
    \kappa(x,y)
   =\left(\int_U\alpha\partial\beta\wedge d\bar{z}\right)\kappa(x,y),\\
  &\{\Phi^{D(0,R)}(\alpha d\bar{z}\otimes x),\Phi^{D(0,R)}(\beta d\bar{z}\otimes y)\}^{(1)}_0(\gamma\otimes\varphi)\\
  &=\left(\int_{D(0,R)}\alpha\beta\wedge\gamma\wedge d\bar{z}\right)\ad^*(x)(\varphi)(y)\\
  &=\left(\int_{D(0,R)}\alpha\beta\gamma\wedge d\bar{z}\right)\varphi([x,y])=\Phi^{D(0,R)}\left(\alpha\beta d\bar{z}\otimes[x,y]\right)(\gamma\otimes\varphi)
\end{align*}
for a field $\gamma\otimes\varphi\in\Omega^{1,*}_{\bbC}(U)\otimes\frg^\vee$.
Thus, by construction of the action of $\frg[z]$ on $\bbV(\Obs^q_{WZW,\kappa})$, the vector space is isomorphic to the center $\clZ^{VA}(V^\kappa_\hbar(\frg))$ of the $\bbC\dbr{\hbar}$-linear vertex algebra $V_\hbar^{\kappa}(\frg)$.
\end{proof}

%%%%%%%%%%%%%%%%%%%%%%%%%%%%%%%%%%%%%%%%%%%%%%%%%%%%%%%%%%%%%%%%%%%%%%%%%%%%%%%%
%%%%%%%%%%%%%%%%%%%%%%%%%%%%%%%%%%%%%%%%%%%%%%%%%%%%%%%%%%%%%%%%%%%%%%%%%%%%%%%%
\section{The BV construction of the Wakimoto realizations}
\label{s.4}

In this section, we construct the Wakimoto realization via BV quantization.
By \cref{quantum_Noether_map}, if we construct a quantum outer-action of the local dgLA $\clG_1$ whose obstruction class solving the inner-equivariant master equation is $\phi(\kappa)$ for some $\kappa\in(\Sym^2\frg^\vee)^\frg$, then we have a map of factorization algebras from the Kac--Moody factorization algebra $\bbU^\hbar_{\phi(\kappa)}(\clG_1)[\hbar^{-1}]$.
In \cite{Cos4}, it is shown that the obstruction class can be computed as the weight of stable graphs (see \cref{equivalent_obstruction_Feynman_diagram_lem}). 
Using this, we compute the obstruction class and construct a factorization algebra analogue of the Wakimoto realization. 

Throughout this section, we assume that $\frg$ is a simple Lie algebra over $\bbC$.
We choose a Cartan subalgebra $\frh\subset\frg$ and a positive part $\Delta_+ \subset \Delta$ of the associated root system.
Then, we have $\Delta=\Delta_+\sqcup(-\Delta_+)$ and the decomposition $\frg=\frn_+\oplus\frh\oplus\frn_-$ of $\frg$.
Let $\{e_\alpha\}_{\alpha\in\Delta_+}$, $\{h_i\}_{i=1,\dots,\rk{\frg}}$, 
and $\{f_\alpha\}_{\alpha\in\Delta_+}$ be a basis of $\frn_+$, $\frh$, and $\frn_-$, respectively.

%%%%%%%%%%%%%%%%%%%%%%%%%%%%%%%%%%%%%%%%%%%%%%%%%%%%%%%%%%%%%%%%%%%%%%%%%%%%%%%%
\subsection{The finite dimensional setting}
\label{ss:4.1}

In this subsection, we review the \emph{vector-field realization} of the Lie algebra $\frg$.

Let $G$ be a connected simply-connected simple Lie group over $\bbC$ such that the associated Lie algebra is $\frg$.
Let $N_+$ be a unipotent subgroup whose associated Lie algebra is $\frn_+$, and $B_-$ be a Borel subgroup whose associated Lie algebra is a Borel subalgebra $\frb_-=\frh\oplus\frn_-$.
Then a flag variety $G/B_-$ is a homogeneous space of the Lie group $G$.

Let $U=N_+\cdot[B_-]\subset G/B_-$ be a big cell.
Since $N_+$ is a unipotent Lie group, the exponential map $\frn_+\to N_+$ is an isomorphism of varieties.
Hence, $U$ is an affine space, and the big cell $U$ is an affine open sub-space of $G/B_-$.

Let $\clD(U)$ be an algebra of differential operators on the affine space $U$, and $\{\clD_{\leq i}(U)\}$ be a natural filtration on $\clD(U)$ given by the order of the differential operator.
Then we have a short exact sequence
\begin{align}
\label{D_leq_1_exact_sequence}
  0 \lto \bbC[U] \lto \clD_{\leq 1}(U) \lto \frX(U) \lto 0,
\end{align}
where $\frX(U)$ is the Lie algebra of vector fields on $U$.
Since it is well-known that the above short exact sequence \eqref{D_leq_1_exact_sequence} has a canonical splitting, we have an injective map $\frX(U)\inj\clD_{\leq1}(U)$.
Since the flag variety $G/B_-$ is homogeneous space of $G$, for every element $a\in\frg=T_eG$, there exists the corresponding vector field $\xi_a$ on $G/B$.
Since the big cell $U$ is open subvariety, the corresponding vector fields define the vector fields on $U$.
Thus, we have the following map
\begin{align*}
  \frg \lto\frX(G/B_-), \quad 
  a \lmto \left(\xi_a:f(y)\mapsto\frac{d}{d\epsilon}f(\exp(-\epsilon a)y)|_{\epsilon=0}\right),
\end{align*}
of Lie algebras.
Therefore, we have an inclusion
\begin{align}
\label{g_to_Vect_U}
  \rho_0\colon 
  \frg \linj \frX(G/B_-) \linj \frX(U)
  \linj \clD_{\leq1}(U)
\end{align}
of Lie algebras.
Taking a coordinate $\{y_\alpha\}_{\alpha\in\Delta_+}$ on $U$, then the explicit formula is well-known and is given as follows: \cite[Chapter 1.5]{Frenkel}
\begin{align}
\label{equ:basismapsto}
\begin{split}
  e_\alpha&\lmto\frac{\partial}{\partial y_\alpha}+\sum_{\beta\in\Delta_+}P^\alpha_\beta(y)\frac{\partial}{\partial y_\beta},\\
  h_i&\lmto-\sum_{\alpha\in\Delta_+}\alpha(h_i)y_\alpha\frac{\partial}{\partial y_\alpha},\\
  f_\alpha&\lmto\sum_{\beta\in\Delta_+}Q^\alpha_\beta(y)\frac{\partial}{\partial y_\beta},
  \end{split}
\end{align}
where $P^\alpha_\beta,Q^\alpha_\beta$ are polynomials on $U\simeq\mathbb{A}^{|\Delta_+|}$.
Since for $x\in\frg$, the differential operator $\rho_0(x)\in\clD_{\leq1}$ is a polynomial of $y_\alpha,\partial/\partial y_\alpha$, we identify the differential operator $\rho_0(x)$ with a polynomial $\rho_0(x)\in\bbC[T^*\frn_+]\simeq\bbC[y_\alpha,\xi_\alpha|\alpha\in\Delta_+]$.

The inclusion \eqref{g_to_Vect_U} is a map of representations of the Lie algebra $\frg$.
Since for $i=1,\dots,\rk{\frg}$ and $\alpha\in\Delta_+$, we have
\begin{align*}
  h_i\cdot y_\alpha
  &=-\sum_{\beta\in\Delta_+}\beta(h_i)
     \left[y_\beta\frac{\partial}{\partial y_\beta},y_\alpha\right]
   =-\alpha(h_i)y_\alpha,\\
  h_i\cdot\frac{\partial}{\partial y_\alpha}
  &=-\sum_{\beta\in\Delta_+}\beta(h_i)
     \left[y_\beta\frac{\partial}{\partial y_\beta},\frac{\partial}{\partial y_\alpha}\right]
   = \alpha(h_i)\frac{\partial}{\partial y_\alpha},
\end{align*}
the weight of $y_\alpha$ and $\frac{\partial}{\partial y_\alpha}$ are given by $-\alpha$ and $\alpha$, respectively.
Furthermore, since the weights of $e_\alpha$ and $f_\alpha$ are given by $\alpha$ and $-\alpha$, respectively, for each $\alpha,\beta\in\Delta_+$, the polynomial $P^\alpha_\beta$ is a homogeneous polynomial of the homogeneous polynomial of  and the polynomial $Q^\alpha_\beta$ is a homogeneous polynomial of weight $\alpha-\beta$.

\begin{lem}
  \label{lem_Q_not_constant}
  For any $\alpha,\beta\in\Delta_+$, the polynomials $P^\alpha_\beta,\ Q^\alpha_\beta$ are not constant.
\end{lem}
\begin{proof}
%最終版では証明は省略してもよいかもしれない
The polynomial $Q^\alpha_\beta$ is a homogeneous polynomial of weight $\alpha-\beta$.
Thus, if the polynomial $Q^\alpha_\beta$ is constant, at least, we must have $\alpha=\beta$.
Suppose that there exists a positive root $\alpha\in\Delta_+$ such that
\[
  \rho_0(f_\alpha)=\frac{\partial}{\partial y_\alpha}+\sum_{\beta\in\Delta_+\setminus\{\alpha\}}Q^\alpha_\beta\frac{\partial}{\partial y_\beta}.
\]
Then since the map $\rho_0$ is a Lie algebra map, we have
  \[
  -\alpha(h_i)\rho_0(f_\alpha)=[\rho_0(h_i),\rho_0(f_\alpha)],
\]
for $i=1,\dots,\rk{\frg}$. The right hand side is
\begin{align*}
  \left[\rho_0(h_i),\rho_0(f_\alpha)\right]
  &=-\sum_{\beta\in\Delta_+}\beta(h_i)\left[y_\beta\frac{\partial}{\partial y_\beta},\frac{\partial}{\partial y_\alpha}+\sum_{\gamma\in\Delta_+\setminus\{\alpha\}}Q^\alpha_\gamma(y)\frac{\partial}{\partial y_\gamma}\right]\\
  &=\alpha(h_i)\frac{\partial}{\partial y_\alpha}+\sum_{\gamma\in\Delta_+\setminus\{\alpha\}}\gamma(h_i)Q^\alpha_\gamma(y)\frac{\partial}{\partial y_\gamma}-\sum_{\beta\in\Delta_+}\sum_{\gamma\in\Delta_+\setminus\{\alpha\}}\left[\frac{\partial}{\partial y_\beta},Q^\alpha_\gamma(y)\right]y_\beta\frac{\partial}{\partial y_\gamma}.
\end{align*}
Since the monomials $y_{\vec{\alpha}}^K\frac{\partial}{\partial y_\beta}$ 
for $n\ge0$, $\vec{\alpha}=(\alpha^{(1)},\dots,\alpha^{(n)})\in\Delta_+^{n}$, 
$K=(k_1,\dots,k_n)\in\bbZ_{\ge0}^n$, and $\beta\in\Delta_+$ are mutually linearly independent, we have
\[
  -\alpha(h_i)\frac{\partial}{\partial y_\alpha}=\alpha(h_i)\frac{\partial}{\partial y_\alpha},
\]
and hence we have $\alpha(h_i)=-\alpha(h_i)$ for all $i=1,\dots,\rk\frg$.
This is a contradiction.
The claim for $P^\alpha_\beta$ follows in the same way.
\end{proof}

On the other hand, we have the short exact sequence \eqref{D_leq_1_exact_sequence}.
We consider an extension of $\frg\to\frX(U)$ by $\frg\to\clD_{\leq1}(U)$.
The extension is classified by first cohomology of $\frg$ with coefficients in $\bbC[U]$.
\begin{lem}
\label{lem_CU_cont_Verma}
There exists an isomorphism
\begin{align*}
  \Hom^{\res}_{U(\frb_-)}(U(\frg),\bbC_{\mathrm{triv}})\simeq\bbC[U]
\end{align*}
of $U(\frg)$-modules, where the $U(\frg)$-module structure on the left hand side is given as follows
\[
(x\cdot f)(P)=f(Px),
\]
for $x\in\frg,\ f\in\Hom^{\res}_{U(\frb_-)}(U(\frg),\bbC_{\mathrm{triv}}),\ P\in U(\frg)$.
\end{lem}
\begin{proof}
  We take a basis $\{e_\alpha,h_i,f_\alpha\}_{\alpha\in\Delta_+,i=1,\dots,\rk\frg}$ of $\frg$.
  By the PBW theorem, any element $P\in U(\frg)$ of $U(\frg)$ is expanded as a linear combination of elements of the form
  \[
  e_{\alpha^{(1)}}^{n_1}\cdots e_{\alpha^{(m)}}^{n_m}h_{i_1}^{k_1}\cdots h_{i_l}^{k_l} f_{\alpha^{(m+1)}}^{p_1}\cdots f_{\alpha^{(m+q)}}^{p_q}
  \]
  for some $n_1,\dots,n_m,k_1,\dots,k_l,p_1,\dots,p_q\geq0,\ \alpha^{(1)},\dots,\alpha^{(m+q)}\in\Delta_+$, and $ i_1,\dots,i_l=1,\dots,\rk\frg$.
  Hence, we have a $U(\frg)$-module isomorphism $\psi:\Hom^{\res}_{U(\frb_-)}(U(\frg),\bbC_{\mathrm{triv}})\stackrel{\simeq}{\to}\bbC[\frn_+]$ such that
  \[
  \psi(\delta_{e_\alpha})=e_\alpha,
  \]
  where $\delta_{e_\alpha}:U(\frg)\to\bbC_{\mathrm{triv}}$ is a $U(\frb_-)$-module map defined by
  \[
  \delta_{e_\alpha}\left(  e_{\alpha^{(1)}}^{n_1}\dots e_{\alpha^{(m)}}^{n_m}h_{i_1}^{k_1}\dots h_{i_l}^{k_l} f_{\alpha^{(m+1)}}^{p_1}\dots f_{\alpha^{(m+q)}}^{p_q}\right)=\begin{cases}
    \delta_{\alpha,\alpha^{(1)}},\quad &\text{for }m=1,n_1=1,\\
    0,\ &\mathrm{otherwise}.
  \end{cases}
  \]
  Since the big cell $U\subset G/B_-$ is identified with the unipotent Lie group $N_+$, the exponential map induces a $U(\frg)$-module isomorphism $(\exp^{-1})^*:\bbC[\frn_+]\stackrel{\sim}{\to}\bbC[N_+]\simeq\bbC[U]$.
  Thus, we have an isomorphism
\[
  \Hom^{\res}_{U(\frb_-)}(U(\frg),\bbC_{\mathrm{triv}}) 
  \xrightarrow{ \ \psi \ } \bbC[\frn_+]
  \xrightarrow{\ (\exp^{-1})^* \ } \bbC[U]
\]
of $U(\frg)$-modules.
\end{proof}

By \cref{lem_CU_cont_Verma} and Shapiro's lemma, we have
\[
H^1(\frg,\bbC[U])\simeq H^1(\frg,\Hom^{\res}_{U(\frb_-)}(U(\frg),\bbC))\simeq H^1(\frb_-,\bbC)\simeq H^1(\frh,\bbC)\simeq\frh^\vee.
\]
We consider the description at the cocycle level.
Since the isomorphism $\frh^\vee\lsto H^1(\frb_-,\bbC)$ is given by the restriction map, we consider the linear map
\begin{align}
\label{equ:htoZ1}
\begin{split}
  \frh^\vee&\lto Z^1(\frb_-,\bbC)\\
  \chi&\lto\begin{dcases}
  h_i \mto \chi(h_i)\\
  f_{\alpha_i} \mto \chi(h_i).
  \end{dcases}
\end{split}
\end{align}
By the construction of the isomorphism of \cref{lem_CU_cont_Verma}, at the cocycle level, the cohomology isomorphisms
\[
\frh^\vee\lsto H^1(\frb_-,\bbC)\lsto H^1(\frg,\Hom^{\res}_{U(\frb_-)}(U(\frg),\bbC))\lsto H^1(\frg,\bbC[U])
\]
can be realized by choosing representatives as follows.
\begin{align*}
  \frh^\vee&  \lto Z^1(\frg,\bbC[U]),\\
  \chi&\lmto
  \begin{cases}
    e_{\alpha_i}\mapsto0\\
    h_i\mapsto\chi(h_i)\\
    f_{\alpha_i}\mapsto\chi(h_i)y_{\alpha_i}.
  \end{cases}
\end{align*}
Thus, for $\chi\in\frh^\vee$, we have a Lie algebra map $\rho_\chi:\frg\to\clD_{\leq1}(U)$ given as follows:
\begin{align}
  \label{extension_g_to_D}
  \begin{split}
  e_{\alpha_i}&\lmto\frac{\partial}{\partial y_{\alpha_i}}+\sum_{\beta\in\Delta_+}P^{\alpha_i}_\beta(y)\frac{\partial}{\partial y_\beta},\\
  h_i&\lmto-\sum_{\alpha\in\Delta_+}\alpha(h_i)y_\alpha\frac{\partial}{\partial y_\alpha}+\chi(h_i),\\
  f_{\alpha_i}&\lmto\sum_{\beta\in\Delta_+}Q^{\alpha_i}_\beta(y)\frac{\partial}{\partial y_\beta}+y_{\alpha_i}\chi(h_i).
    \end{split}
\end{align}
Since for $x\in\frg$, the differential operator $\rho_\chi(x)\in\clD_{\leq1}(U)$ is again a polynomial of $y_\alpha,\partial/\partial y_\alpha$, we identify the differential operator $\rho_\chi(x)$ with a polynomial $\rho_\chi(x)\in\bbC[T^*\frn_+]\simeq\bbC[y_\alpha,\xi_\alpha|\alpha\in\Delta_+]$.
These maps combine into a universal map
\[
\wt{\rho}:\frg\lto\clD_{\leq1}(U)\otimes\frh
\]
such that for each $\chi\in\frh^\vee$ we recover $\rho_\chi$ as the composition of $\wt{\rho}$ and the evaluation at $\chi$ map $\frh\to\bbC$ along the second factor.
We denote by $U(\wt{\rho}):U(\frg)\to\clD(U)\otimes U(\frh)$ by the induced map on the universal enveloping algebra.

Here, we consider the following natural filtration
\[
F_p(\clD(U)\otimes U(\frh))\ceq \bigoplus_{i+j\leq p}\clD_{\leq i}(U)\otimes F^{PBW}_jU(\frh).
\]
on $\clD(U)\otimes U(\frh)$.
Then since the map $U(\wt{\rho})$ preserves the filtrations $\{F^{PBW}_pU(\frg)\}_{p\geq0}$ and $\{F_p(\clD(U)\otimes U(\frh))\}_{p\geq0}$, we have the induced map
\[
  \gr{U}(\wt{\rho})\colon 
  \gr{U}(\frg)\simeq\bbC[\frg^\vee] \lto 
  \gr(\clD(U) \otimes U(\frh)) \simeq \bbC[T^*N_+\times\frh]
\]
of commutative algebras.
Then, we have the induced map
\begin{align}
    \label{eq:Jetinducedmap}
    \clJ_\infty(\gr{U}(\wt{\rho}))\colon\bbC[\clJ_\infty\frg^\vee]\to\bbC[\clJ_\infty(T^*N_+\times\frh)]
\end{align}
on the jet schemes.

The following lemma will be used in the proof of \cref{thm_inj_C}.
\begin{lem}[{\cite[Proposition 5.2]{Frenkel}}]
  \label{lem_jet_tilderho}
  The map \eqref{eq:Jetinducedmap} is injective.
\end{lem}

On the other hand, using the dual basis $\{h_i^*\}_{i=1,\dots,\rk\frg}$ with respect to the basis $\{h_i\}_{i=1,\dots,\rk\frg}$, we have a map $\wt{\rho}_\chi\colon\frg\to\clD_{\leq1}(U)\otimes\frh^\vee$ given by
\begin{align}
\label{tilderho_map}
\begin{split}
  e_{\alpha_i}&\lmto\frac{\partial}{\partial y_i}+\sum_{\beta\in\Delta_+}P^{\alpha_i}_\beta(y)\frac{\partial}{\partial y_\beta},\\
  h_i&\lmto-\sum_{\alpha\in\Delta_+}\alpha(h_i)y_\alpha\frac{\partial}{\partial y_\alpha}+\chi(h_i)h_i^*,\\
  f_{\alpha_i}&\lmto\sum_{\beta\in\Delta_+}Q^{\alpha_i}_{\beta}(y)\frac{\partial}{\partial y_\beta}+y_{\alpha_i}\chi(h_i)h^*_i.
\end{split}
\end{align}

\subsection{Critical level}
\label{ss:4.2}

This subsection is the most important part of this note.
The goal of this subsection is to consider an analogue of the $\rho:\frg\to\mathcal{D}_{\leq1}(U)$ \eqref{equ:basismapsto} introduced in \cref{ss:4.1} for holomorphic field theories on $\bbC$.
By \cite{GMS00,GMS04}, the natural vertex algebra analogue of differential operators is the \emph{chiral differential operator}, which is locally given by the $\beta\gamma$-system.
Since the big cell $U$ is affine, by \cref{hbar_adic_betagamma_FA}, the natural analogue of the algebra of differential operators $\mathcal{D}(U)$ in holomorphic field theory on $\bbC$ is the $\dim U=\dim\frn_+$-fold tensor product of $\clM_{\beta\gamma}$ introduced in \cref{ss:3.3} (see \cite{GGW20}).

The key point here is to consider the lift of the classical action to a quantum action and the corresponding obstruction class.
The proof of the main theorem of this subsection, \cref{Wakimoto_critical}, follows the approach of Chapter 3 of \cite{Frenkel}, while the explicit computation of the obstruction class is based on \cite{GGW20} and \cite{anomaly}.
The method using the bump function in \cref{Noether_map_Wakimoto_critical_level} is based on the approach of \cite{Ni25}.

\subsubsection{The equivariant action}

We consider the free degenerate classical field theory
\[
\left(\clM_{\frg}=\Omega^{0,*}_{\bbC}\otimes\frn_+[-1]\oplus\Omega^{1,*}_{\bbC}\otimes\frn^\vee_+[-1],\omega_{\beta\gamma},Q_{BV}^{\beta\gamma}=\bar{\partial}\right)
\]
on $\bbC$ (see \eqref{eq:Mbetagamma}).
Recall the basis $\{e_\alpha\}_{\alpha\in\Delta_+}$ of $\frn_+$.
Let $\{e_\alpha^*\}_{\alpha\in\Delta_+}$ be the dual basis of $\frn_+^\vee$.
Sections $\gamma\in\Omega^{0,*}_\bbC\otimes\frn_+$ and $\beta\in\Omega^{1,*}_\bbC\otimes\frn_+^\vee$ can be expanded in the following forms.
\begin{align*}
    \gamma=\sum_{\alpha\in\Delta_+}\gamma_\alpha\otimes e_\alpha,\quad \gamma_\alpha\in\Omega
    ^{0,*}_\bbC,\\
    \beta=\sum_{\alpha\in\Delta_+}\beta_\alpha\otimes e_\alpha^*,\quad \beta_\alpha\in\Omega^{1,*}_\bbC.
\end{align*}
We denote a section of $\clM_{\frg}$ by a pair $(\gamma_\alpha,\beta_{\alpha})$.

Similarly, the section $\alpha\in\clG_1$ can be expanded in the following form.
\[
\alpha=\sum_{\alpha\in\Delta_+}\omega_\alpha\otimes e_\alpha+\sum_{i=1}^{\rk{\frg}}\eta_i\otimes h_i+\sum_{\alpha\in\Delta_+}\xi_\alpha\otimes e_\alpha^*,\quad \omega_\alpha,\eta_i,\xi_\alpha\in\Omega^{0,*}_\bbC.
\]
We denote by a section of $\clG_1$ by a tuple $(\omega_\alpha\otimes e_\alpha,\eta_i\otimes h_i,\xi_\alpha\otimes f_\alpha)$.

Below, given $f\in\bbC[T^*\frn_+]=\bbC[y_\alpha,\xi_\alpha\mid\alpha\in\Delta_+]$, we define the Dolbeault form $f(\gamma,\beta)\in\Omega^{*,*}_{\bbC}$ by replacing each $y_\alpha$ and $\xi_\alpha$ with $\gamma_\alpha$ and $\beta_\alpha$, respectively, and the product in the commutative algebra with the wedge product.

Here, we consider a classical action of $\clG_1=\Omega^{0,*}_{\bbC}\otimes\frg$ on the classical field theory $\clM_{\frg}$.
Let $\wt{S}^{Wak}_{\kappa_c}\in\IA(\clG_1,\clM_{\frg})$ be a functional defined as follows:
\begin{align}
  \label{naive_Wakimoto}
  \begin{split}
  \wt{S}^{Wak}_{\kappa_c}&(\omega_\alpha\otimes e_\alpha,\eta_i\otimes h_i,\xi_\alpha\otimes f_\alpha,\gamma_\alpha,\beta_\alpha)\\
 =\sum_{\alpha\in\Delta_+}&\int_{\bbC}\omega_\alpha\wedge\beta_\alpha+\omega_\alpha\wedge\sum_{\beta\in\Delta_+}\beta_\beta\wedge P^\alpha_\beta(\gamma)+\xi_\alpha\wedge\sum_{\beta\in\Delta_+}\beta_\beta\wedge Q^\alpha_\beta(\gamma)-\sum_{i=1}^{\rk\frg}\alpha(h_i)\eta_i\wedge\beta_\alpha\wedge\gamma_\alpha.
  \end{split}
\end{align}
Here we used the polynomials $P^\alpha_\beta,Q^\alpha_\beta$ in \eqref{equ:basismapsto}.
By abuse of notation, we denote by $\wt{S}^{Wak}_{\kappa_c}$ the image of the natural quotient map $\IA(\clG_1,\clM_{\frg})\srj\Act(\clG_1,\clM_{\frg})$.

\begin{lem}
\label{naive_equivariant_classical}
The functional $\wt{S}^{Wak}_{\kappa_c}$ satisfies the equivariant classical master equation.
\end{lem}
\begin{proof}
  Let $\{J^a\}_{a=1,\dots,\dim\frg}$ be a basis of $\frg$.
  The definition \eqref{naive_Wakimoto} of the functional $\wt{S}^{Wak}_{\kappa_c}$ can be rewritten as
  \begin{equation}\label{SWak-action}
  \wt{S}^{Wak}_{\kappa_c}(\eta_a\otimes J^a,\gamma_{\alpha},\beta_{\alpha})=\frac{1}{(2+1)!}\sum_{a=1}^{\dim\frg}\sum_{\alpha\in\Delta_+}\int_\bbC\eta_a\wedge\rho_0(J^a)(\gamma,\beta),
  \end{equation}
  for $\eta_a\in\Omega^{0,*}_{\bbC}$.
  Here we regard $\rho_0(J^a)$ as a polynomial in $\bbC[T^*\frn_+]$.
  The factor $(2+1)!$ comes from the multiplicity in the symmetrization.
  Then, the equivariant classical master equation is equivalent to the following equalities:
  \begin{align}
  \label{eq:equivclassicalmasterfortildeS1}
\bar{\partial}\wt{S}^{Wak}_{\kappa_c}+Q_{BV}^{\beta\gamma}\wt{S}^{Wak}_{\kappa_c}&=0,\\
  \label{eq:equivclassicalmasterfortildeS2}
  (d_{\clG_1}-\bar{\partial})\wt{S}^{Wak}_{\kappa_c}+\frac{1}{2}\{\wt{S}^{Wak}_{\kappa_c},\wt{S}^{Wak}_{\kappa_c}\}&=0.
  \end{align}
  
    The left hand side of \eqref{eq:equivclassicalmasterfortildeS1} has degree 1 as a functional on $\clG_1$.
    For $(\eta_a\otimes J^a)\in\clG_1$, the left hand side of \eqref{eq:equivclassicalmasterfortildeS1} is
    \begin{align*}
        &\frac{1}{(2+1)!}\sum_{a=1}^{\dim\frg}\sum_{\alpha\in\Delta_+}\int_\bbC\bar{\partial}\eta_a\wedge\rho_0(J^a)(\gamma,\beta)+\frac{1}{(2+1)!}\sum_{a=1}^{\dim\frg}\sum_{\alpha\in\Delta_+}\int_\bbC\eta_a\wedge\bar{\partial}(\rho_0(J^a)(\gamma,\beta))\\
      =&\frac{1}{(2+1)!}\sum_{a=1}^{\dim\frg}\sum_{\alpha\in\Delta_+}\int_\bbC\bar{\partial}(\eta_a\wedge\rho_0(J^a)(\gamma,\beta))=0,
    \end{align*}
    where we use the Stokes lemma.
    Thus, we have shown \eqref{eq:equivclassicalmasterfortildeS1}. 

    The left hand side of \eqref{eq:equivclassicalmasterfortildeS2} has degree 2 as a functional on $\clG_1$.
    For $\alpha\otimes x,\alpha'\otimes x'\in\clG_1$, we have
    \begin{align*}
        &\left((d_{\clG_1}-\bar{\partial})\wt{S}^{Wak}_{\kappa_c}+\frac{1}{2}\{\wt{S}^{Wak}_{\kappa_c},\wt{S}^{Wak}_{\kappa_c}\}\right)(\alpha\otimes x,\alpha'\otimes x')\\
      =&\wt{S}^{Wak}_{\kappa_c}(\alpha\wedge\alpha'\otimes[x,x'])-\{\wt{S}^{Wak}_{\kappa_c}(\alpha\otimes x),\wt{S}^{Wak}_{\kappa_c}(\alpha'\otimes x')\}.
    \end{align*}
    By \eqref{analytic_part_of_heat_kernel_betagamma} and \cref{hol_theory_BV_bracket}, the bracket Poisson $\{-,-\}$ is given by 
    \[
    \{I,J\}(\gamma,\beta)=\sum_{\alpha\in\Delta_+}\frac{\partial I(\gamma,\beta)}{\partial\gamma_\alpha}\frac{\partial J(\gamma,\beta)}{\partial\beta_\alpha}+(-1)^{|I||J|}\frac{\partial J(\gamma,\beta)}{\partial\gamma_\alpha}\frac{\partial I(\gamma,\beta)}{\partial\beta_\alpha},\quad I,J\in\operatorname{Obs}^{cl}_{\mathcal{M}_\frg}
    \]
    where the symbols $|I|$ and $|J|$ are cohomological degrees of $I$ and $J$, respectively.
    Thus, for $x,x\in\frg$, we have
    \[
    \rho_0([x,x'])(\gamma,\beta)=\{\rho_0(x),\rho_0(x')\}(\gamma,\beta),
    \]
    and hence the left hand side of \eqref{eq:equivclassicalmasterfortildeS2} vanishes. Thus, the statement is proved.
\end{proof}

\begin{thm}
\label{Wakimoto_critical}
There exists an $S^1$-invariant functional
\[
S^{Wak}_{\kappa_c}\in\IA(\clG_1,\clM_{\frg})[\hbar]
\]
such that the collection of $S^1$-invariant functionals $\{S^{Wak}_{\kappa_c}[L]=\lim_{\epsilon\to0}W(P(\epsilon,L),S^{Wak}_{\kappa_c})\}_{L>0}$ satisfies the equivariant outer-quantum master equation.
\end{thm}

In this subsection, we prove \cref{Wakimoto_critical}.

\subsubsection{The obstruction class}

First, we consider the complex $\IA(\clG_1,\clM_{\frg})$.
By \cref{differential_Inneract_act_classical}, we have a short exact sequence
\begin{align}\label{eq:sesG1}
  0 \lto C^*_{red,loc}(\clG_1) \xrr{i}
  \IA(\clG_1,\clM_{\frg}) \xrr{p}
  \Act(\clG_1,\clM_{\frg}) \lto 0.
\end{align}

\begin{lem}
\label{cohomology_of_red_loc}
The following sequence is exact.
\begin{align}
  \label{long_exact_sequence_inn_act}
  H^1(C^*_{red,loc}(\clG_1))\xrr{H^1(i)} 
  H^1(C^*_{red,loc}(\Omega^{0,*}\otimes\frh) \otimes 
  C^*_{loc}(\Omega^{1,*}_{\bbC}\otimes\frn_+^\vee)) \xrr{H^1(p)}
  H^1(\Act(\clG_1,\clM_{\frg})).
\end{align}
\end{lem}
\begin{proof}

Since the differential of the dgLA $J(\clG_1\ltimes\clM_{\frg})$ is given by 
$\bar{\partial}+p_+$, where 
\[
  p_+:\Omega^{0,*}\otimes\frg \simeq \Omega^{0,*}\otimes(\frn_+\oplus\frh\oplus\frn_-)
  \to
  \Omega^{0,*}\otimes\frn_+[-1]\subset\clM_{\frg}
\]
is the natural projection, we have a quasi-isomorphism
\begin{align*}
  J(\clG_1\ltimes\clM_{\frg})\qis J(\Omega^{0,*}\otimes\frb_-)\oplus J(\Omega^{1,*}_{\bbC}\otimes\frn_+^\vee)[-1]
\end{align*}
of sheaves of dgLAs.
Thus, we have an isomorphism
\begin{align}
  \label{isom_inner_cohomology}
  \begin{split}
  H^p(\IA(\clG_1,\clM_{\frg}))
  &\simeq H^p(C^*_{red,loc}(\Omega^{0,*}\otimes\frb_-))\otimes
  \Sym^{*}(J(\Omega^{1,*}_{\bbC}\otimes\frn_+^\vee[-1])^\vee[-1])\\
  &\simeq 
  H^p(C^*_{red,loc}(\Omega^{0,*}\otimes\frh))\otimes
  C^*_{loc}(\Omega^{1,*}_{\bbC} \otimes \frn_+^\vee[-1]).    
  \end{split}
\end{align}
Then, the short exact sequence \eqref{eq:sesG1}
induces the described \eqref{long_exact_sequence_inn_act}.
\end{proof}

By \cref{holomorphic_counterterm}, there exists a limit
\[
  \wt{S}^{Wak}_{\kappa_c}[L] \ceq \lim_{\epsilon\to0}W(P(\epsilon,L),\wt{S}^{Wak}_{\kappa_c}).
\]
By \cref{lift_classical_to_quantum_obstruction_class}, 
the following functionals are cocycles in the complexes 
$\IA(\clG_1,\clM_{\frg})$ and 
$\Act(\clG_1,\clM_{\frg})$, respectively.
\begin{align*}
  \Theta_{inn}^{Wak,crit} \ceq 
  &\frac{1}{\hbar}\lim_{L\to0}\Bigl(d_{\clG_1}\wt{S}^{Wak}_{\kappa_c}[L]+\bar{\partial}\wt{S}^{Wak}_{\kappa_c}[L] \\
  &\,+\frac{1}{2}\{\wt{S}^{Wak}_{\kappa_c}[L],\wt{S}^{Wak}_{\kappa_c}[L]\}_L+\hbar\Delta_L\wt{S}^{Wak}_{\kappa_c}[L]\Bigr)\ \bmod \hbar,\\
  \Theta_{out}^{Wak,crit} \ceq
  &\frac{1}{\hbar}\lim_{L\to0}\Bigl(d_{\clG_1}\wt{S}^{Wak}_{\kappa_c}[L]+\bar{\partial}\wt{S}^{Wak}_{\kappa_c}[L]\\
  &\,\left.+\frac{1}{2}\{\wt{S}^{Wak}_{\kappa_c}[L],\wt{S}^{Wak}_{\kappa_c}[L]\}_L+\hbar\Delta_L\wt{S}^{Wak}_{\kappa_c}[L]\Bigr)\right|_{\Act(\clG_1,\clM_{\frg})}\ \bmod \hbar,
\end{align*}
where the symbol $(\cdots)|_{\Act(\clG_1,\clM_{\frg})}$ means 
taking the class in the quotient space $\Act(\clG_1,\clM_{\frg})$.
Let $[\Theta^{Wak,crit}_{inn}],[\Theta^{Wak,crit}_{out}]$ be the classes in the first cohomologies $H^1(\IA(\clG_1,\clM_{\frg}))$ and $H^1(\Act(\clG_1,\clM_{\frg}))$, respectively.

For the computation of the obstruction class, the following \cref{equivalent_obstruction_Feynman_diagram_lem} is known.

We review the general method for computing the obstruction classes in \cref{lift_classical_to_quantum_obstruction_class}.
Let $(\clM,\Pi,Q_{BV}=\bar{\partial})$ be a strictly free degenerate classical field theory on $\bbC^d$, $\mathcal{L}$ be a local $L_\infty$-algebra on $\bbC^d$, and $S^\mathcal{L}\in\Act(\mathcal{L},\mathcal{M})$ be an equivariant action functional.
To first order in $\hbar$, the obstruction class solving the equivariant quantum master equation modulo $\hbar^2$ is the following:
\begin{align}
  \Theta^{\clL}[L]:=\hbar^{-1}\left((Q+d_{\clL})S^{\clL}[L]+\hbar\Delta_LS^{\clL}[L]+\frac{1}{2}\{S^{\clL}[L],S^{\clL}[L]\}\right)\ \bmod \hbar.
\end{align}

For a graph $\gamma=(V,E,s,t)$, we denote by $E^{ex}\subset E$ and $E^{in}\subset E$ the set of external edges and internal edges, respectively.
Then $E=E^{ex}\sqcup E^{in}$.
For a graph $\gamma=(V,E,s,t)$, a distinguished internal edge $e\in E^{in}$, $\epsilon,L>0$, a functional $I=\sum_{n}I_n,\ (I_n\in\Sym^n(\mathcal{M}^*))$, and a section $x\in\mathcal{M}$, we denote by $\Gamma=(\gamma,e)$ the pair and define the weight $W_{\Gamma}(P(\epsilon,L),K_\epsilon,I)(x)\in\bbC$ by the following steps (see Chapter 2 of \cite{C11}).
\def\theenumi{\roman{enumi}}
\def\labelenumi{(\theenumi)}
\begin{enumerate}
    \item For a vertex $v\in V$, set $n(v)=\sharp\{e\in E\mid s(e)=v\text{ or }t(e)=v\},\ m(v)=\sharp\{e\in E^{ex}\mid t(e)=v\}$ and put $I_{n(v)}$ on $v$,
    \item For an internal edge $e'\in E^{in}$, put $P(\epsilon,L)$ on $e'$,
    \item For an extrenal edge $e'\in E^{ex}$, put $x$ on $e'$,
    \item For the distinguished internal edge $e\in E^{in}$, put $K_\epsilon$ on $e$,
    \item For each vertex $v\in V$, evaluate the externally placed $x$ to obtain $I_{n(v)}(x)(-)\in\Sym^{n(v)-m(v)}(\mathcal{M}^*)$,
    \item For each internal edge $e'\in E^{\mathrm{in}}$, contract the functionals $I_{n(s(e))}(x)(-)$ and $I_{n(t(e))}(x)(-)$, placed at $s(e)$ and $t(e)$, respectively, using the propagator $P(\epsilon,L)$ or $K_\epsilon$ assigned to $e'$.
\end{enumerate}
The above steps yield a constant, which we denote by $W_{\Gamma}(P(\epsilon,L),K_\epsilon,I)(x)$.
Also, we have a functional $W_\Gamma(P(\epsilon,L),K_\epsilon,I)(-)$ on $\mathcal{M}$, which call the weight of $\Gamma$.

Typically, the sheaf $\clM$ decomposes into several direct summands, the functional $I$ is supported on certain direct summands, and the contractions involving the propagator and the heat kernel are nontrivial only on certain direct summands.
Then, we label each edge by a section of the direct summand of $\clM$.
For a section $y$ of a direct summand $\mathcal{N}\subset\clM$, we call the corresponding edge a $y$-leg.

By Chapter 4.2 of \cite{anomaly}, the obstruction class $\Theta^\clL:=\lim_{L\to0}\Theta^\clL[L]$ can be expanded as the following sum.
\begin{align}
\label{equ:obstructionclassgeneral}
\Theta^\clL=\frac{1}{2\hbar}\sum_{\Gamma=(\gamma,e)}W_{\Gamma}(P(\epsilon,L),S^\clL)\ \bmod\hbar.
\end{align}
Here, the assumption that the functional $S^\mathcal{L}$ satisfies the equivariant classical master equation implies that the coefficient of $\hbar^0$ in the sum on the right-hand side vanishes.

Furthermore, when the Poisson structure $\Pi$ is nondegenerate,\footnote{A $(-1)$-shifted Poisson structure $\Pi\in\PD(\clM^![-1],\clM[1])_{S_2}$ is nondegenerate if $\Pi$ is an isomorphism.} we can characterize more precisely the form of the graphs that contribute to the obstruction class \cref{equivalent_obstruction_Feynman_diagram_lem}.

A wheel graph is a graph $\gamma=(V,E,s,t)$, possibly with external edges such that for each vertex $i\in V$, there uniquely exists two non-external edges $e,e'\in E$ whose end point is $s(e)=i=t(e')$.

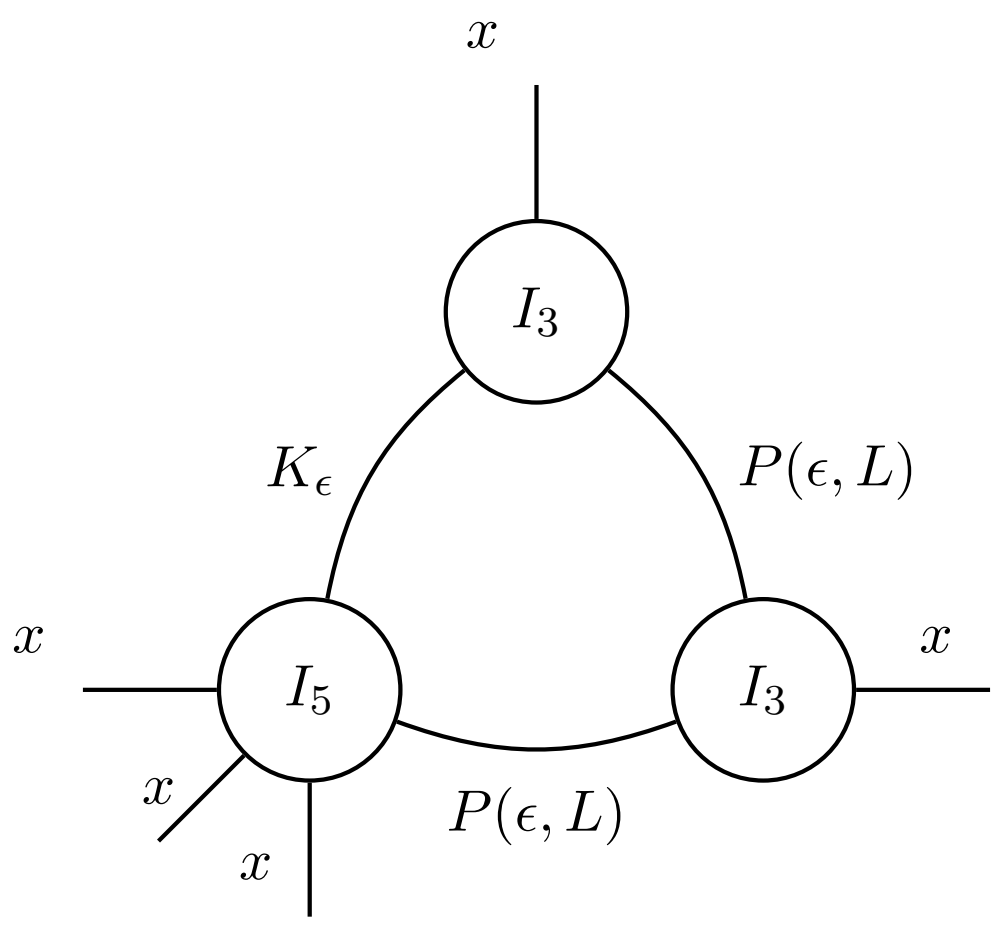
\begin{figure}[htbp]
\centering
\begin{tikzpicture}[node distance=3cm, thick]
 \node [draw, circle, minimum size=12mm](A) at (0,0) {$I_{5}$};
 \node [draw, circle, minimum size=12mm](B) at (3,0) {$I_3$};
 \node [draw, circle, minimum size=12mm](C) at (1.5,2.5) {$I_3$};
 \draw (A) to[bend right=20]  node[midway, below=3pt] {$P(\epsilon,L)$} (B);
 \draw (B) to[bend right=20]  node[midway, right=3pt] {$P(\epsilon,L)$} (C);
 \draw (C) to[bend right=20]  node[midway, left=3pt] {$K_{\epsilon}$} (A);
 \draw (A) -- ++(-1.5,0) node[pos=1, above left=3pt] {$x$};
 \draw (A) -- ++(0,-1.5) node[pos=1, above left=3pt] {$x$};
 \draw (A) -- ++(-1.0,-1.0) node[pos=1, above=3pt] {$x$};
 \draw (B) -- ++(1.5,0) node[pos=1, above left=3pt] {$x$};
 \draw (C) -- ++(0,1.5) node[pos=1, above left=3pt] {$x$};
\end{tikzpicture}
\caption{A example of wheel graph with a distinguished edge.}
\label{fig:triangle}
\end{figure}

It is known that the obstruction class can be computed using the weights defined in this way.

\begin{lem}[{\cite[Corollary 16.0.5]{Cos4}, \cite[Proposition 4.4]{anomaly}}]
\label{equivalent_obstruction_Feynman_diagram_lem}
Let $(\clM,\omega,S)$ be a holomorphic nondegenerate classical field theory on $\bbC^d$, and $\Theta=\lim_{L\to0}\Theta[L]$ be the scale 0 limit of the obstruction class solving the quantum master equation modulo $\hbar^2$.
Then the following holds:
\begin{align}
  \label{equivalent_obstruction_Feynman_diagram_equ}
  \hbar\Theta=\frac{1}{2}\lim_{L\to0}\lim_{\epsilon\to0}\sum_{\Gamma\in Wheel_{d+1,e}}W_\Gamma(P(\epsilon,L),K_\epsilon,I),
\end{align}
where $Wheel_{d+1,e}$ is the set of wheel graphs which has $d$-vertices with one distinguished edge $e$.
\end{lem}

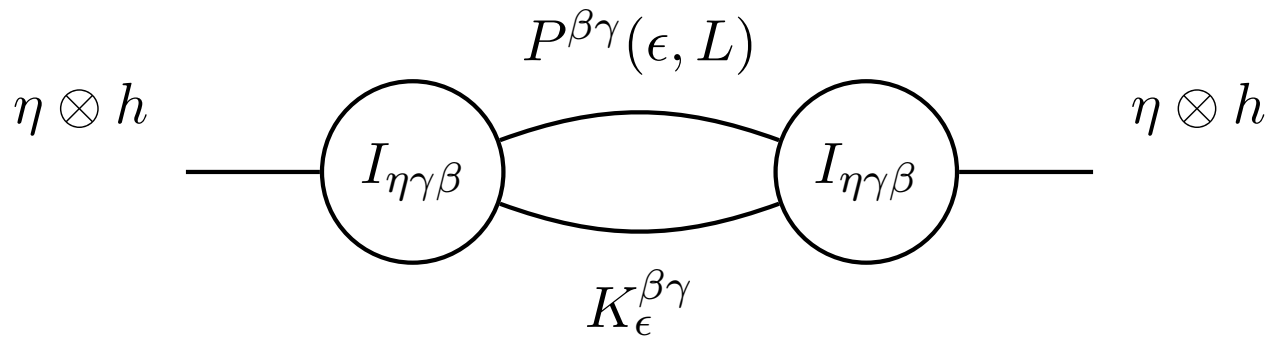
\begin{figure}[htbp]
% --- 1つ目の図 ---
\begin{tikzpicture}[node distance=3cm, thick]

% ノード（丸の中に I）
\node[draw, circle, minimum size=12mm] (A) {$I_{\eta\gamma\beta}$};
\node[draw, circle, minimum size=12mm, right of=A] (B) {$I_{\eta\gamma\beta}$};

% 外側の線とラベル
\draw (A) -- ++(-1.5,0) node[pos=1, above left=3pt] {$\eta\otimes h$};
\draw (B) -- ++(1.5,0)  node[pos=1, above right=3pt] {$\eta\otimes h$};

% 内側の 2 本線
\draw (A) to[bend left=20]  node[midway, above=3pt] {$P^{\beta\gamma}(\epsilon,L)$} (B);
\draw (A) to[bend right=20] node[midway, below=3pt] {$K_\epsilon^{\beta\gamma}$} (B);

\end{tikzpicture}

\caption{The graph which contributes to the obstruction class $[\Theta^{Wak,crit}_{inn}]$.}
\label{obstruction_class_inn}
\end{figure}

Graphically, the homogeneous part of degree $n$ of $\wt{S}^{Wak}_{\kappa_c}$ has $(n-1)$ $\gamma$-legs and one $\beta$-leg.
Since propagators connect only $\gamma$-legs and $\beta$-legs, no $\beta$-leg can appear as an external edge of a two-wheel diagram contributing to $\Theta^{Wak,crit}$.
Thus, the isomorphism \eqref{isom_inner_cohomology} implies that the cohomology class $[\Theta^{Wak,crit}_{inn}]$ belongs to $H^1(C^*_{red,loc}(\Omega^{0,*}_\bbC\otimes\frh))$.
By \cref{cohomology_of_red_loc}, this implies that the cohomology class $[\Theta^{Wak,crit}_{inn}]$ belongs to the image of $H^1(i)$.
Furthermore, a priori, the image of $H^1(i)$ is a functional on $\Omega^{0,*}_\bbC\otimes\frg$.
However, by \eqref{isom_inner_cohomology}, it is determined entirely by its values on $\Omega^{0,*}_\bbC\otimes\frh$.
Thus, by \cref{equivalent_obstruction_Feynman_diagram_lem}, the class $[\Theta^{Wak,crit}_{inn}]$ is a sum of weights of two-wheel graphs whose external edges are all $\eta\otimes h$-legs.
Such a graph is given by \cref{obstruction_class_inn}, 
where the functional $I_{\eta\gamma\beta}$ for $\eta \in\Omega^{0,*}_\bbC$ is defined by 
\[
  I_{\eta\gamma\beta} \ceq
  \sum_{\alpha\in\Delta_+}
  \int_{\bbC}\alpha(h)\beta_\alpha\wedge\eta\wedge\gamma_\alpha.
\]

By the definition of the classical field theory $\clM_\frg$, \eqref{analytic_part_of_heat_kernel_betagamma}, and \eqref{propagator_def_betagamma}, for $\epsilon,L>0$, and $z,w\in\bbC$, the heat kernel $K^{\beta\gamma}_\epsilon$ and the propagator $P^{\beta\gamma}$ of $\clM_\frg$ are
\begin{align}
\label{eq:betagammahetpro}
\begin{split}
      K^{\beta\gamma}_\epsilon(z_1,z_2) &=
  \frac{1}{4\pi \epsilon}e^{-|z_1-z_2|^2/4\epsilon}(d\bar{z}_1-d\bar{z}_2)\otimes\sum_{\alpha\in\Delta_+}e_\alpha\otimes e_\alpha^*,\\
  P^{\beta\gamma}(\epsilon,L)(z_1,z_2) &=
  \int_{\epsilon}^L\frac{1}{4\pi t}\frac{\partial}{\partial z_1}e^{-|z_1-z_2|^2/4t}dt\otimes \sum_{\alpha\in\Delta_+}e_\alpha\otimes e_\alpha^*.
\end{split}
\end{align}
Thus the heat kernel and the propagator decompose as the analytic parts and the algebraic parts:
\[
  K^{\beta\gamma}_\epsilon = 
  K^{\beta\gamma,an}_\epsilon \otimes K^{\beta\gamma,alg}_{\epsilon},\quad 
  P^{\beta\gamma} = P^{an} \otimes P^{alg}.
\]
We consider its value on $a=\eta_i\otimes h_i$, where $\eta_i$ is a Dolbeault form on $\bbC$, and $h_i$ is a basis element of the Cartan subalgebra $\frh$.
Then, according to the decompositions of $K^{\beta\gamma}_\epsilon$ and $P^{\beta\gamma}$, the obstruction class also decomposes as
\[
  \Theta^{Wak,crit}_{inn}(a)= \Theta^{alg}(h_i) \cdot \Theta^{an}(\eta_i),
\]
where $\cdot$ is the multiplication of complex numbers.
Hereafter until \eqref{analytic_part_keisann3}, 
we focus on the analytic part $P^{an}$ and $\Theta^{an}$.
The analytic propagator is defined as 
\[
  P^{an}(\epsilon,L)(z_1,z_2) = 
  \int_\epsilon^L\frac{1}{4\pi t}\frac{\partial}{\partial z_1}e^{-|z_1-z_2|/(4t)}dt.
\]
Let $f,g$ be holomorphic functions on $\bbC$ and $\eta_i=f+g\ d\bar{z}$ be a Dolbeault form. Then,
\begin{align}
\label{analytic_part_keisann1}
  \Theta^{an}(\eta_i)&\\
  =\lim_{L\to0}\lim_{\epsilon\to0}
\notag
  &\int_{\bbC_{z_1}\times\bbC_{z_2}}
    f(z_1)g(z_2)d\bar{z}_2\wedge K_\epsilon^{an}(z_1,z_2)\wedge P^{an}(\epsilon,L)(z_1,z_2) \\
\notag
  %=\lim_{L\to0}\lim_{\epsilon\to0}
  %&\int_{\bbC_{z_1}\times\bbC_{z_2}}\left(f(z_1)g(z_2)
   % \int_{\epsilon}^L\frac{1}{(4\pi)^2\epsilon t}e^{-|z_1-z_2|^2/(4\epsilon)}\frac{\partial}{\partial z_1}e^{-|z_1-z_2|^2/(4t)}dtd^2z_1d^2z_2\right)\\
\notag
  =-\lim_{L\to0}\lim_{\epsilon\to0}
  &\int_{\bbC_{z_1}\times\bbC_{z_2}}\left(f(z_1)g(z_2)
    \int_{\epsilon}^L\frac{\bar{z_1}-\bar{z_2}}{4(4\pi)^2\epsilon t^2}e^{-|z_1-z_2|^2/(4\epsilon)}e^{-|z_1-z_2|^2/(4t)}dtd^2z_1d^2z_2\right)\\
\notag
  =-\lim_{L\to0}\lim_{\epsilon\to0}
  &\int_{\bbC_{w_1}\times\bbC_{w_2}}f\left(\frac{w_1+w_2}{2}\right)g\left(\frac{w_2-w_1}{2}\right)\\
\notag
  &\quad \times\int_{\epsilon}^L\frac{\bar{w_1}}{4(4\pi)^2\epsilon t^2}\exp\left(-\frac{1}{4}(t^{-1}+\epsilon^{-1})|w_1|^2\right)dtd^2z_1d^2z_2,
\end{align}
where in the last equality, we replaced $z_1,z_2$ by $w_1=z_1-z_2,w_2=z_1+z_2$. 

Furthermore, we will use the following trick:
For $t\in\bbR$, define the holomorphic differential operator
\[
D(t)=\left(1-\frac{t}{t+\epsilon}\right)\frac{\partial}{\partial w_1}=\frac{\epsilon}{t+\epsilon}\frac{\partial}{\partial w_1}.
\]
We see that
\begin{align*}
  &D(t)\left(fg\frac{1}{\epsilon t}
   \exp\left(-\frac{1}{4}(t^{-1}+\epsilon^{-1})|w_1|^2\right) \right)\\ 
  &=\left(-\frac{\bar{w}_1}{t}fg+D(t)(fg)\right)
   \frac{1}{4\epsilon t}\exp\left(-\frac{1}{4}(t^{-1}+\epsilon^{-1})|w_1|^2\right).
\end{align*}
Since the left-hand side is a total derivative, 
the $w_1$-integral vanishes by Stokes's theorem, 
and the integral in \eqref{analytic_part_keisann1} is given as follows:
\begin{align}
\label{analytic_part_keisann2}
\begin{split}
  &-\int_{\bbC_{w_1}\times\bbC_{w_2}}\frac{\partial}{\partial w_1}\left(f\left(\frac{w_1+w_2}{2}\right)g\left(\frac{w_2-w_1}{2}\right)\right)\\
  \times&\int_{\epsilon}^L\frac{1}{4(4\pi)^2t(\epsilon+t)}\exp\left(-\frac{1}{4}(t^{-1}+\epsilon^{-1})|w_1|^2\right)dw_1dw_2dt.
\end{split}
\end{align}
Here, for the $w_1$-integration, we take the Taylor expansion of $\frac{\partial}{\partial w_1}\left(f\left(\frac{w_1+w_2}{2}\right)g\left(\frac{w_2-w_1}{2}\right)\right)$ around $w_1=0$ and perform the Gaussian integral.  In the limit $\epsilon\to0$, only the constant of the expansion survives. Therefore, the obstruction class is given as follows:
\begin{align}
\begin{split}
\label{analytic_part_keisann3}
&\Theta^{an}(\eta_i)\\
=&-\lim_{L\to0}\lim_{\epsilon\to0}\frac{1}{(4\pi)^2}\int_{\bbC_{w_2}}\frac{\partial}{\partial w_1}\left(f\left(\frac{w_1+w_2}{2}\right)g\left(\frac{w_2-w_1}{2}\right)\right)(w_1=0)\int_{\epsilon}^L\frac{\epsilon}{(\epsilon+t)^2}dtd^2w_2\\ 
                   =&-\frac{1}{(4\pi)^2}\int_{\bbC_{w_2}}\frac{\partial}{\partial w_1}\left(f\left(\frac{w_1+w_2}{2}\right)g\left(\frac{w_2-w_1}{2}\right)\right)(w_1=0)\lim_{L\to0}\lim_{\epsilon\to0}\left(-\frac{\epsilon}{\epsilon+\epsilon}+\frac{\epsilon}{L+\epsilon}\right)d^2w_2\\
                   =&\frac{1}{2}\frac{1}{(4\pi)^2}\int_{\bbC_{z}}\left(\frac{\partial}{\partial z}f(z)\right)g(z)d^2z.
\end{split}
\end{align}

Let us turn to the algebraic part $\Theta^{alg}(h_i)$.
By \eqref{eq:betagammahetpro}, we find that the contraction operators $\partial_{P^{alg}(\epsilon,L)}$ and $\partial_{K^{alg}_\epsilon}$ by the algebraic parts of the propagator $P^{alg}$ and the heat kernel $K^{alg}_\epsilon$ both become
\[
  \sum_{\alpha\in\Delta_+}
  \frac{\partial}{\partial e_\alpha}\frac{\partial}{\partial e^*_\alpha}.
\]
Thus, the algebraic factor of the weight of the graph in \cref{obstruction_class_inn} is given as follows:
\begin{align*}
  -\sum_{\alpha,\alpha',\alpha''\in\Delta_+}\delta_{\alpha\alpha''}\delta_{\alpha'\alpha''}\alpha(h_i)\alpha'(h_i)=-\sum_{\alpha\in\Delta_+}\alpha(h_i)\alpha(h_i)=-\frac{1}{2}\kappa_{\frg}(h_i,h_i),
\end{align*}
where the $(-1)$-factor comes from the differences in cohomological degrees.

Going back to the total $\Theta^{Wak,crit}_{inn}(a)$, 
the weight of the graph in \cref{obstruction_class_inn} is given by
\[
  \frac{1}{(2\pi i)^2}\frac{1}{2} \times
  \left(-\frac{1}{2}\kappa_{\frg}(h_i,h_i)\right)\int_{\bbC}\eta_i\wedge\partial\eta_i=\phi\left(-\tfrac{1}{2}\kappa_\frg\right)(\eta_i\otimes h_i).
\]
Therefore, the following follows from the fact stated immediately after \cref{equivalent_obstruction_Feynman_diagram_lem}.
\[
  [\Theta^{Wak,crit}_{inn}]=H^1(i)\left(\left[\phi\left(-\tfrac{1}{2}\kappa_{\frg}\right)\right]\right)=H^1(i)([\phi(\kappa_c)]).
\]

\subsubsection{The outer quantum action}

We have determined the cohomology class of $\Theta^{Wak,crit}_{inn}$.
Next, we want to choose a representative of the class
which lifts to a solution of the quantum master equation.

The functional $\phi(\kappa_c)$ is also a cocycle of $C^*_{red,loc}(\clG_1)$.
Since the obstruction class $[\Theta^{Wak,crit}_{inn}]$ coincides with the image $H^1(i)([\phi(\kappa_c)])$, there exists a cohomological degree 0 cochain $\tau\in\IA(\clG_1,\clM_{\frg})$ such that
\begin{align}
  \label{gamma_exact}
i\left(\phi\left(\kappa_c\right)\right)-\Theta^{Wak,crit}_{inn}=d_{\operatorname{Inner}}\tau=d_{\clG_1}\tau+\bar{\partial}\tau+\{\wt{S}^{Wak}_{\kappa_c},\tau\}.  
\end{align}
It is obvious that the functional $\tau$ is not contained in the image of $i$.
Thus, $\tau$ defines a nonzero element of $\Act(\clG_1,\clM_{\frg})$ under the quotient map $\IA(\clG_1,\clM_{\frg})\srj\Act(\clG_1,\clM_{\frg})$.
Since the left hand side of the equality \eqref{gamma_exact} is a degree $2$ functional on $\clG_1$, $\tau$ is a degree $1$ functional on $\clG_1$.
Since in the wheel graph with non-trivial contributions, the internal edges contain one $P^{\beta\gamma}(\epsilon,L)$ and one $K_\epsilon^{\beta\gamma}$, and $\wt{S}^{Wak}_{\kappa_c}$ is linear in $\Omega^{1,*}_\bbC\otimes\frn_+^{\vee}$, the functional $\Theta^{Wak,crit}_{inn}$ is a degree $0$ functional on $\Omega^{1,*}_{\bbC}\otimes\frn^\vee_+[-1]$.
Thus, as a functional on $\Omega^{0,*}_{\bbC}\otimes\frn_+[-1]$, $\tau$ has degree at least $1$, and as a functional on $\Omega^{1,*}_{\bbC}\otimes\frn^\vee_+[-1]$, $\tau$ has degree $0$.
Therefore, there exist polynomials 
\begin{equation}\label{eq:Bpoly}
  B_a\in\bbC[\clJ_\infty\frn_+] \quad 
  (a=1,\dots,\dim\frg)
\end{equation}
such that, for a basis $\{J^a\}_{a=1,\dots,\dim\frg}$ of $\frg$, we have
\[
\tau(\alpha_1\otimes J^1,\dots,\alpha_{\dim\frg}\otimes J^{\dim\frg},\gamma,\beta)=\sum_a\int_{\bbC}\alpha_a\wedge B_a(\gamma)\wedge dz.
\]

We consider the following functional
\begin{equation}\label{eq:SWktau}
  S^{Wak}_{\kappa_c;\tau}=\wt{S}^{Wak}_{\kappa_c}+\hbar\tau\in\Act(\clG_1,\clM_{\frg})[\hbar].
\end{equation}
Here, we again use the same notation for the image under the quotient map $\IA(\clG_1,\clM_{\frg})\srj\Act(\clG_1,\clM_{\frg})$ of the functional $\wt{S}^{Wak}_{\kappa_c}+\hbar\tau\in\operatorname{Inner}(\clG_1,\clM_{\frg})[\hbar]$.
Then by \cref{holomorphic_counterterm}, there exists a limit $S^{Wak}_{\kappa_c;\tau}[L]=\lim_{\epsilon\to0}W(P(\epsilon,L),S^{Wak}_{\kappa_c;\tau})$.
Then by \cref{hol_theory_BV_bracket} and \eqref{gamma_exact}, we have
\begin{align*}
  &\frac{1}{\hbar}\lim_{L\to0}\left((d_{\clG_1}+\bar{\partial})S^{Wak}_{\kappa_c;\tau}[L]+\hbar\Delta_LS^{Wak}_{\kappa_c;\tau}[L]+\frac{1}{2}\{S^{Wak}_{\kappa_c;\tau}[L],S^{Wak}_{\kappa_c;\tau}[L]\}_L\right)\ \bmod \hbar\\
 =&\Theta^{Wak,crit}_{inn}+(d_{\clG_1}+\bar{\partial})\tau+\lim_{L\to0}\{\wt{S}^{Wak}_{\gamma_c},\tau\}_L\ \bmod \hbar\\
 =&\Theta^{Wak,crit}_{inn}+(d_{\clG_1}+\bar{\partial})\tau+\{\wt{S}^{Wak}_{\gamma_c},\tau\}=i\left(\phi\left(\kappa_c\right)\right).
\end{align*}
Therefore, the collection $\{S^{Wak}_{\kappa_c;\tau}[L]\}_{L>0}$ defines an outer quantum action of $\clG_1$ on $\clM_{\frg}$.

\subsubsection{$S^1$-invariance}

Next, choose the functional $S^{Wak}_{\kappa_c;\tau}$ to be $S^1$-invariant.

Recall that the $S^1$-action on $\bbC$ is given by
\[
  \theta\cdot z=e^{-i\theta}z \quad (\theta\in S^1, \, z\in\bbC),
\]
and the $S^1$-action on $\clM_{\frg}$ is given by letting $S^1$ act trivially on $\frn_+$ and $\frn_+^{\vee}$ (see \cref{ss:3.3}).
Since a $(1,0)$-form has $S^1$-weight $1$, $\frg$ must have $S^1$-weight $-1$ for \cref{naive_Wakimoto} to be $S^1$-invariant.

Since $\phi$ and $\wt{S}^{Wak}_{\kappa_c}$ are $S^1$-invariant, the left hand side of \eqref{gamma_exact} is also $S^1$-invariant, and hence we can choose the functional $\tau$ to be $S^1$-invariant.
Note that $\tau$ is $S^1$-invariant if and only if for any $a=1,\dots,\dim\frg$, the polynomial $B_a$ of \eqref{eq:Bpoly} is of the  form
\begin{align}
  \label{B_a_katati}
  B_a(\gamma)(z)=
  \sum_{\alpha\in\Delta_+}R_a^\alpha(\gamma(z))\partial_z\gamma_\alpha(z),
  \quad R_a^\alpha\in\bbC[\frn_+].
\end{align}
%with each $R_a^\alpha\in\bbC[\frn_+]$ a polynomial.

Recall the natural weight module structure on $\clD_{\leq1}(U)$ which is defined by the inclusion \eqref{g_to_Vect_U}.
By the construction of the equivariant action functional $\wt{S}^{Wak}_{\kappa_c}$, the action of $\frg$ on $\clD_{\leq1}(U)$ induces the action of $\frg$ on $\clM_{\frg}$ by replacing the coordinates $y_\alpha,\partial/\partial y_\alpha$ with the sections $\gamma_\alpha,\beta_\alpha$.
Then the weight of the sections $\gamma_\alpha,\beta_\alpha$ are given by
\[
  \operatorname{wt}\gamma_\alpha=-\alpha, \quad 
  \operatorname{wt} \beta_\alpha= \alpha,
\]
for $\alpha\in\Delta_+$.
On the other hand, the adjoint action of $\frg$ on $\frg$ induces the action of $\frg$ on $\clG_1$.
Then, for $a=1,\dots,\dim\frg$ and a section $\alpha\in\Omega^{0,*}_{\bbC}$, 
the weight of the section $\alpha\otimes J^a \in \clG_1$ is given by
\[
  \operatorname{wt}(\alpha\otimes J^a)=-\operatorname{wt}J^a.
\]
By construction, the polynomials $P^\alpha_\beta(\gamma)$ and $Q^\alpha_\beta(\gamma)$ are homogeneous weight polynomials $-\alpha-\beta$ and $\alpha-\beta$ with respect to the above weight, respectively.
Thus, the functional $\wt{S}^{Wak}_{\kappa_c}$ is a homogeneous functional of weight 0.
Since the propagator $P(\epsilon,L)$ and the heat kernel $K_\epsilon$ are also homogeneous functionals of 0 weight, obstruction $\Theta^{Wak,crit}_{inn}$ is also a homogeneous functional of 0 weight.
Indeed, according to \cref{equivalent_obstruction_Feynman_diagram_lem}, the obstruction $\Theta^{Wak,crit}_{inn}$ is the sum of the contraction of the symmetric product $\wt{S}^{Wak}_{\kappa_c}\cdot\wt{S}^{Wak}_{\kappa_c}$ by $P(\epsilon,L)$ and $K_\epsilon$.
Furthermore, it is obvious that the cocycle $\phi(\kappa_c)$ is functional of weight 0 homogeneous.
Therefore, the left hand side of \eqref{gamma_exact} is a homogeneous functional weight 0 and since the differential $d_{\clG_1}+\bar{\partial}+\{\wt{S}^{Wak}_{\kappa_c},-\}$ preserves the weight, the functional $\tau$ can be chosen as a functional weight 0 functional.
Then the polynomial $B_a$ has weight $\operatorname{wt}J^a$ for $a=1,\dots,\dim\frg$.
The equality \eqref{B_a_katati} implies that for $J^a\in\frn_+\oplus\frh$, the polynomial $B_a$ cannot have weight $\operatorname{wt}J^a$, and hence we have $B_a=0$ for $J^a\in\frn_+\oplus\frh$.
On the other hand, for $J^a=f_\alpha$, the polynomial $B_a$ must have the weight $-\alpha$, and hence there exists $c_\alpha\in\bbC$ such that $R^\beta_\alpha=\delta_{\alpha,\beta}c_\alpha$.
Therefore, we have the following assertion:

\begin{prop}
\label{Wakimoto_field_theory_crit}
There exists $c_\alpha\in\bbC\ (\alpha\in\Delta_+)$ such that
\begin{align*}
  \tau(\omega_{\alpha}\otimes e_\alpha,\eta_i\otimes h_i,\xi_\alpha\otimes f_\alpha,\gamma_\alpha,\beta_\alpha)=\sum_{\alpha\in\Delta_+}c_\alpha\int_{\bbC}\xi_\alpha\wedge\partial\gamma_\alpha
\end{align*}
satisfies the equality \eqref{gamma_exact}.
\end{prop}

%Taking $\tau$ as in \cref{Wakimoto_field_theory_crit} and setting $S^{Wak}_{\kappa_c} \ceq S^{Wak}_{\kappa_c;\tau}$ given in \eqref{eq:SWktau}, we obtain the desired functional.

Taking $c_\alpha$ as in \cref{Wakimoto_field_theory_crit} and setting
\[
  S^{Wak}_{\kappa_c} \ceq 
  \wt{S}^{Wak}_{\kappa_c}+\hbar\sum_{\alpha\in\Delta_+}
   c_{\alpha}\int_{\bbC}\xi_\alpha\wedge\partial\gamma_\alpha
\]
given in \eqref{eq:SWktau}, we obtain the desired functional.
The proof of \cref{Wakimoto_critical} is now complete.

\subsubsection{Critical Wakimoto realization}

By \cref{quantum_Noether_map,Wakimoto_critical}, 
we have the following assertion:
\begin{thm}
\label{Noether_map_Wakimoto_critical_level}
The quantum outre-action $\{S^{Wak}_{\kappa_c}[L]\}_{L>0}$ induces a $\bbC\dbr{\hbar}$-linear map of factorization algebras 
\[
  \wt{w}^{fact}_{\kappa_c}\colon
  \bbU^\hbar_{\phi(\kappa_c)}(\clG_1)[\hbar^{-1}]
  \lto \Obs^q_{\clM_{\frg}}[\hbar^{-1}].
\]
The map $\wt{w}^{fact}_{\kappa_c}$ induces a map of $\bbC\dpr{\hbar}$-linear vertex algebras
\[
  \wt{w}^{ch}_{\kappa_c,\hbar}\colon
  V^{\kappa_c}(\frg)\dpr{\hbar} \lto
  (\beta\gamma_\hbar)^{\otimes\dim\frn_+}[\hbar^{-1}]
\]

Specializing to $\hbar=2\pi i$, we obtain the Wakimoto realization at critical level $\kappa_c$:
\[
\wt{w}^{ch}_{\kappa_c}\colon V^{\kappa_c}(\frg)\lto\beta\gamma^{\otimes\dim\frn_+}.
\]
\end{thm}
\begin{proof}
  We consider the limit as the scale $L\to 0$.
  Let $R>0$ and $f\in C^\infty_c(\bbR_{\geq0})$ be a compactly supported smooth function on $\bbR_{\geq0}$ such that $f(0)=0$.
  %Then, there exists an open neighborhood $U\subset\bbR_{\geq0}$ of $0\in\bbR_{\geq0}$ such that $f$ vanishes on $U$. 
  The following elements are representatives of a basis of the cohomology $H^1(\Omega^{0,*}_{\bbC,c}(D(0,R)))$: 
  \[
  z\longmapsto f(|z|)z^{-n}d\bar{z}\ceq\begin{cases}
  f(|z|)z^{-n}d\bar{z} &\text{for } z\neq0\\
  0 &\text{for } z=0    
  \end{cases} \quad (n\geq0).
  \]
  Thus, the cohomology classes of the following elements form a basis of $\bbV(\clF^\kappa_{\frg})$: 
  \[
  D(0,R)^p\ni(z_1,\dots,z_p) \lmto 
  f(|z_1|)z_1^{-n_1}d\bar{z}_1\otimes x_1\cdots f(|z_p|)z_p^{-n_p}d\bar{z}_p\otimes x_p,
  \]
  for $p\geq0,n_1,\dots,n_p\geq0,x_1,\dots,x_p\in\frg$.
  In particular, for $x\in\frg$, under the isomorphism $\bbV(\clF^\kappa_{\frg})\simeq V^\kappa(\frg)$, the base $x_{(-1)}\ket*{0}$ of $V^\kappa(\frg)$ is identified with the cohomology class of the following element
  \[
  D(0,R)\ni z \lmto f(|z|)d\bar{z}\otimes x.
  \]
  Also, similarly to \cref{def_gerators_of_Obs}, we define
  \begin{align*}
    c^\alpha_n(x)\colon\Omega^{0,0}(D(x,R))\otimes\frn_+\ni\gamma\otimes e_\beta&\lmto\frac{1}{n!}(\partial_z\gamma)(x)\delta_{\alpha,\beta},\\
    b^\alpha_m(x)\colon\Omega^{1,0}(D(x,R))\otimes\frn_+^\vee\ni\beta dz\otimes e_\beta^*&\lmto\frac{1}{(m-1)!}(\partial^{m-1}_z\beta)(x)\delta_{\alpha,\beta},
  \end{align*}
  for $\alpha\in\Delta_+$, $n\geq0$, $m>0$, and $x\in\bbC$.
  Then, similarly to \cref{Obs_cl_generator}, the monomials
  \[
  b_K^{\vec{\alpha}}(0)c_L^{\vec{\beta}}(0),
  \]
  form a basis of $\bbV(\Obs^q_{\clM_{\frg}})$.
  Under the isomorphism $\bbV(\Obs^q_{\clM_{\frg}})\simeq(\beta\gamma_\hbar)^{\otimes\dim\frn_+}$, the monomical $  b_K^{\vec{\alpha}}(0)c_L^{\vec{\beta}}(0)$ is identified with the monomial $a_{\vec{\alpha},K}a^*_{\vec{\beta},L}\ket*{0}$ (see \eqref{equ:monomial betagamma}). 
  
By construction of the Noether map \eqref{Noether_map}, for $i=1,\dots,\rk\frg$, we have
\begin{align*}
  \hbar\lim_{L\to0}\wt{w}^{fact}_{\kappa_c}(\omega_{\alpha_i}\otimes e_{\alpha_i})[L](\gamma_\alpha,\beta_\alpha)=\int_{\bbC}\omega_{\alpha_i}\wedge\beta_{\alpha_i}+\omega_{\alpha_i}\wedge\sum_{\beta\in\Delta_+}\beta_\beta\wedge P^{\alpha_i}_{\beta}(\gamma).
\end{align*}
Let $f\in C^\infty_c(\bbR_{\geq0})$ be a compactly supported smooth function on $\bbR_{\geq0}$ with $f(0)=0$ and $\int_{\bbR_{\geq0}}rf(r)dr=\frac{1}{2\pi}$.
%Take an open neighborhood $U$ of $0\in\bbR_{\geq0}$ such that $f$ vanishes on $U$.
Also, let $f_\alpha,g_\alpha$ ($\alpha\in\Delta_+$) be holomorphic functions on $\bbC$.
Using $\int_0^{2\pi}e^{in\theta}d\theta=2\pi\delta_{n,0}$ for $n\in\bbZ$
and the Taylor expansion, we have
\begin{align*}
  &\hbar\lim_{L\to0}\wt{w}^{fact}_{\kappa_c}(f(|z|)d\bar{z}\otimes e_{\alpha_i})[L](f_\alpha,g_\alpha dz)\\
  =&\int_{\bbC}f(|z|)d\bar{z}\wedge g_{\alpha_i}(z)dz+f(|z|)d\bar{z}\wedge\sum_{\beta\in\Delta_+}g_\beta(z)dz\wedge P^{\alpha_i}_{\beta}(f_\alpha(z))\\
=&\left(g_{\alpha_i}(0)+\sum_{\beta\in\Delta_+}g_\beta(0)P^{\alpha_i}_{\beta}(f_\alpha(0))\right)2\pi\int_0^\infty rf(r)dr\\
=&g_{\alpha_i}(0)+\sum_{\beta\in\Delta_+}g_\beta(0)P^{\alpha_i}_{\beta}(f_\alpha(0)).
\end{align*}
Thus, we have
\begin{align}
  \label{e_i_fact_map}
  \hbar\lim_{L\to0}\wt{w}^{fact}_{\kappa_c}(f(|z|)d\bar{z}\otimes e_{\alpha_i})[L]=b_1^{\alpha_i}(0)+\sum_{\beta\in\Delta_+}P^{\alpha_i}_{\beta}(c_0^{\beta}(0))b_1^\beta(0),
\end{align}
and hence we have
\begin{align}
  \label{e_i_VA_map}
  \wt{w}^{ch}_{\kappa_c,\hbar}(e_{\alpha_i,(-1)})=\frac{1}{\hbar}\left(a_{\alpha_i,-1}+\sum_{\beta\in\Delta_+}P^{\alpha_i}_{\beta}(a^*_{\beta,0})a_{\beta,-1}\right).
\end{align}
Similarly we have
\begin{align}
  \label{h_i_f_i_VA_map}
  \begin{split}
    \wt{w}^{ch}_{\kappa_c,\hbar}(h_{i,(-1)})&=-\frac{1}{\hbar}\sum_{\beta\in\Delta_+}\beta(h_i)a^*_{\beta,0}a_{\beta,-1},\\
    \wt{w}^{ch}_{\kappa_c,\hbar}(f_{\alpha_i,(-1)})&=\frac{1}{\hbar}\sum_{\beta\in\Delta_+}Q^{\alpha_i}_{\beta}(a^*_{\beta,0})a_{\beta,-1}+c_{\alpha_i}a^*_{\alpha_i,-1}.
  \end{split}
\end{align}
By \cref{rem:cdprhbarlinear}, the specialization 
\[
\tilde{w}^{ch}_{\kappa_c}:=\tilde{w}^{ch}_{\kappa_c,\hbar}|_{\hbar=2\pi i}:V^{\kappa_c}(\frg)\dpr{\hbar}|_{\hbar=2\pi i}\to(\beta\gamma_\hbar)^{\otimes\dim\frn_+}[\hbar^{-1}]|_{\hbar=2\pi i}
\]
is well-defined.
By \cref{reconstruction_thm_Ch}, the equations \cref{e_i_fact_map,h_i_f_i_VA_map} imply that the map $\wt{w}^{ch}_{\kappa_c}$ is identified with the Wakimoto realization at the critical level $\kappa_c$.
\end{proof}

%%%%%%%%%%%%%%%%%%%%%%%%%%%%%%%%%%%%%%%%%%%%%%%%%%%%%%%%%%%%%%%%%%%%%%%%%%%%%%%%
\subsection{General level}
\label{ss:4.3}
%頂点代数の文脈でこういう話が合って、それを参考にしました。古典的なequivariant functioalは古典論の話から出るべきなので、そこからやって、今のformalismに合わせるために計算をするのがメインパートです。

In this subsection, using \cref{quantum_Noether_map}, we construct a map of factorization algebras from the general-level Kac--Moody factorization algebra $\bbU_{\phi(\kappa)}(\clG_1)$ to a certain factorization algebra of observables.
As in \cref{ss:4.2}, it is natural to construct it from \eqref{tilderho_map} in the classical case.

Recall the Wakimoto realization at general level $\kappa\in\bigl(\Sym^2(\frg^\vee)\bigr)^\frg$ (see \cite[Theorem~5.1]{Frenkel} and \cite[Theorem~6.2.1]{F07}).
It is a map of vertex algebras of the form 
\[
  w^{ch}_\kappa\colon V^\kappa(\frg) \lto 
  \beta\gamma^{\otimes\dim\frn_+}\otimes\pi^{\kappa-\kappa_c}.
\]
We want to lift it to a map of factorization algebras.
By \cref{hbar_adic_betagamma_FA,Ch_kappa_Heis_FA}, it is natural to construct it from the $\clG_1$-action on the strictly free degenerate classical field theory defined by the direct sum of $\clM_{\frg}$ and $\clM_{Heis}^{\kappa-\kappa_c}$.
However, care is needed because $\hbar$ is introduced when carrying out the quantization.

By the explicit formula of $w^{ch}_\kappa$ in \cite{Frenkel,F07}, the terms contributing to the level of the Wakimoto realization arise from the OPEs between pairs of quadratic fields in the $\beta\gamma$-system and pairs of linear fields in the Heisenberg vertex algebra.
In the $\bbC\dbr{\hbar}$-linear setting, we find that the powers of $\hbar$ do not match when considering $\clM_{\frg}$ and $\clM_{Heis}^{\kappa-\kappa_c}$.

The key point here is again to consider the lifting of classical action to a quantum action.
The proof of \cref{Wakimoto_general_level_FA} follows the approach of Theorem 5.1 of \cite{Frenkel}.

In light of the above discussion, we consider a free degenerate classical field theory
\begin{align}
\label{eq:clWdef}
    \left(\clW_\kappa=\clM_{\frg}\oplus\clM^{\hbar(\kappa-\kappa_c)}_{Heis},\Pi=\Pi_{\omega^{\beta\gamma}}\oplus\Pi^{\hbar(\kappa-\kappa_c)}_{Heis},Q_{BV}=\bar{\partial}\right)
\end{align}

on $\bbC$ (see \eqref{eq:Mbetagamma} and \eqref{Heis_field_space}).
Recall the basis $\{h_i\}_{i=1,\dots,\rk\frg}$ of $\frh$.
A section $b\in\clM^{\hbar(\kappa-\kappa_c)}_{Heis}$ can be expanded in the following form.
\[
b=\sum_{i=1}^{\rk{\frg}}b_i\otimes h_i,\quad b_i\in\Omega^{1,*}_\bbC.
\]
We denote a section of $\clW_\kappa$ by a tuple $(\gamma_\alpha,\beta_\alpha,b_i)$.
Here, we consider a classical action of $\clG_1$ on the field theory $\clW_\kappa$.
Let $\wt{S}^{Wak}_\kappa\in\IA(\clG_1,\clW_\kappa)$ be the functional defined by 
\begin{align}
\label{general_level_naive_Wakimoto}
\begin{split}
 &\wt{S}^{Wak}_{\kappa}(\omega_\alpha\otimes e_\alpha,\eta_i\otimes h_i,\xi_\alpha\otimes f_\alpha,\gamma_\alpha,\beta_\alpha,b_i)\\
 \ceq &\widetilde{S}^{Wak}_{\kappa_c}(\omega_\alpha\otimes e_\alpha,\eta_i\otimes h_i,\xi_\alpha\otimes f_\alpha,\gamma_\alpha,\beta_\alpha)+\int_\bbC\sum_{i=1}^{\rk{\frg}}\eta_i\wedge b_i+\sum_{\alpha\in\Delta_+}\xi_{\alpha}\wedge b_{\alpha}\wedge\gamma_{\alpha},
\end{split}
\end{align}
where for a positive root $\alpha=\sum_ic_i\alpha_i\in\Delta_+$, we denote $b_\alpha=\sum_ic_ib_i$.
By abuse of notation, we denote by $\wt{S}^{Wak}_{\kappa}$ the image of the natural quotient map $\IA(\clG_1,\clW_\kappa)\srj\Act(\clG_1,\clW_\kappa)$.

\begin{lem}
\label{lem_equiv_general_one_dim}
The functional $\wt{S}^{Wak}_{\kappa}\in\Act(\clG_1,\clW_\kappa)$ satisfies the equivariant classical master equation. 
\end{lem}
\begin{proof}
The functional $\wt{S}^{Wak}_{\kappa}$ decomposes into the following parts
  \begin{align*}
  \wt{S}^{Wak}_{\kappa}&=\wt{S}^{Wak}_{\kappa_c}+\int_{\bbC}\sum_{i=1}^{\rk\frg}\eta_i\wedge b_i+\sum_{\alpha\in\Delta_+}\xi_{\alpha}\wedge b_{\alpha}\wedge\gamma_{\alpha}\\
  &=\wt{S}^{Wak}_{\kappa_c}+I_{\eta b}+I_{\xi b\gamma}.    
  \end{align*}
  By \cref{naive_equivariant_classical}, it is enough to show that the following equalities
  \begin{align}
    \label{equations_ECMEq}
    \begin{split}
    (d_{\clG_1}+\bar{\partial})(I_{\eta b}+I_{\xi b\gamma})+\{\wt{S}^{Wak}_{\kappa_c},I_{\xi b\gamma}\}_{\clM_{\frg}}&=0,\\
    \{I_{\eta b}+I_{\xi b\gamma},I_{\eta b}+I_{\xi b\gamma}\}_{\clM_{Heis}}&=0,      
    \end{split}
  \end{align}
  where $\{-,-\}_{\clM_{\frg}}$ and $\{-,-\}_{\clM_{Heis}}$ are the $P_0$-factorization algebra structure of the free degenerate classical field theories $\clM_{\frg}$ and $\clM_{Heis}^{\hbar(\kappa-\kappa_c)}$, respectively.
  
  The left hand side of the first equation in \eqref{equations_ECMEq} has degree 1 as a functional on $\clM^{\hbar(\kappa-\kappa_c)}_{Heis}$, and is the sum of its degree 1 and degree 2 parts as a functional on $\clG_1$.
  Then the degree 1 part is given as follows
  \begin{align*}
    &(d_{\clG_1}+\bar{\partial})(I_{\eta b}+I_{\xi b\gamma})(\omega_\alpha\otimes e_\alpha,\eta_i\otimes h_i,\xi_\alpha\otimes f_\alpha)\\
    =&\int_{\bbC}\sum_{i=1}^{\rk\frg}\bar{\partial}(\eta_i\wedge b_\alpha)+\sum_{\alpha\in\Delta_+}\bar{\partial}(\xi_{\alpha}\wedge b_{\alpha}\wedge\gamma_{\alpha})=0,
  \end{align*}
  where we use the Stokes lemma.
  Similarly to the proof of \cref{naive_equivariant_classical}, for $x,x'\in\frg$, we have
  \[
  \tilde{\rho}_{\sum_{\alpha\in\Delta_+}\alpha}([x,x'])(\gamma,\beta,b)=\{\tilde{\rho}_{\sum_{\alpha\in\Delta_+}\alpha}(x),\tilde{\rho}_{\sum_{\alpha\in\Delta_+}\alpha}(x')\}_{\clW_\kappa}(\gamma,\beta,b).
  \]
  Here we regard $\tilde{\rho}_{\sum_{\alpha\in\Delta_+}\alpha}(\gamma,\beta,b)$ as an element of $\Omega^{*,*}_\bbC$ by replacing the coordinate functions $y_\alpha,\xi_\alpha,h_i$ on $\bbC[T^*\frn_+]\otimes\frh^\vee$ with $\gamma_\alpha,\beta_\alpha,b_i$, respectively.
  This implies that the degree 2 part of left hand side of the first equation of \eqref{equations_ECMEq} vanishes.
  
  The left hand side of the second equation in \eqref{equations_ECMEq} has degree 0 as a functional on $\clM^{\hbar(\kappa-\kappa_c)}_{Heis}$, and is given as follows
  \begin{align*}
  &\{I_{\eta b}+I_{\xi b\gamma},I_{\eta b}+I_{\xi b\gamma}\}_{\clM_{Heis}}\\
  =&\hbar\int_{\bbC}\left(\sum_{i=1}^{\rk\frg}(\kappa-\kappa_c)(h_i,h_i)\left(\partial(\eta_i\wedge \eta_i)+\partial(\eta_i\wedge\xi_{\alpha_i}\wedge\gamma_{\alpha_i})\right)\right.\\
  &\left.+\sum_{\alpha\in\Delta_+}\partial(\xi_{\alpha}(\kappa-\kappa_c)(h_\alpha,h_\alpha)\wedge\gamma_{\alpha}\wedge\xi_{\alpha}\wedge\gamma_{\alpha})\right)=0,
  \end{align*}
  where again we use Stokes lemma.
\end{proof}

The following theorem is one of the main theorems in this note.

\begin{thm}
\label{Wakimoto_general_level_FA}
There exists an $S^1$-invariant functional
\[
  S^{Wak}_{\kappa}\in\Act(\clG_1,\clW_\kappa)[\hbar]
\]
such that the collection of $S^1$-invariant functionals 
$\left\{S^{Wak}_\kappa[L] \ceq \lim_{\epsilon\to0}W(P(\epsilon,L),S^{Wak}_\kappa)\right\}_{L>0}$
satisfies the equivariant outer-quantum master equation.
\end{thm}

In the following part of this subsection, we prove \cref{Wakimoto_general_level_FA}.

By \cref{holomorphic_counterterm}, the following limit exists.
\[
  \wt{S}^{Wak}_{\kappa}[L] \ceq 
  \lim_{\epsilon\to0}W(P(\epsilon,L),\wt{S}^{Wak}_{\kappa}).
\]
Similarly to the proof of \cref{Wakimoto_critical}, 
we consider the following obstruction classes
\begin{align*}
  \Theta^{Wak,gen}_{inn}& \ceq 
  \frac{1}{\hbar}\lim_{L\to0}\Bigl(d_{\clG_1}\wt{S}^{Wak}_{\kappa}[L]+\bar{\partial}\wt{S}^{Wak}_{\kappa}[L]&\\
  &+\frac{1}{2}\{\wt{S}^{Wak}_{\kappa}[L],\wt{S}^{Wak}_{\kappa}[L]\}_L+\hbar\Delta_L\wt{S}^{Wak}_{\kappa}[L]\Bigr)\ \bmod \hbar,\\
  \Theta^{Wak,gen}_{out}& \ceq 
  \frac{1}{\hbar}\lim_{L\to0}\Bigl(d_{\clG_1}\wt{S}^{Wak}_{\kappa}[L]+\bar{\partial}\wt{S}^{Wak}_{\kappa}[L]\\
  &+\left.\frac{1}{2}\{\wt{S}^{Wak}_{\kappa}[L],\wt{S}^{Wak}_{\kappa}[L]\}_L+\hbar\Delta_L\wt{S}^{Wak}_{\kappa}[L]\Bigr)\right|_{\Act(\clG_1,\clW_\kappa)}\ \bmod \hbar
\end{align*}
The obstruction classes are 1-cocycles in $\IA(\clG_1,\clW_{\kappa})$ and $\Act(\clG_1,\clW_\kappa)$, respectively.

Now by \eqref{heat_kernel_Heis} and \eqref{propagator_Heis}, the heat kernel of the degenerate classical field theory $\clM^{\hbar(\kappa-\kappa_c)}_{Heis}$ is given as follows:
\begin{align}
\label{equ: Heis heat kernel}
\begin{split}
      K^{Heis,Wak}_L(z,w)&=\frac{\hbar}{2\pi iL}\partial e^{-|z-w|^2/4L}(d\bar{z}-d\bar{w})\otimes c^{\kappa-\kappa_c}\\
                     &=\frac{\hbar}{2\pi iL}\partial_we^{-|z-w|^2/4L}(dz-dw)\wedge(d\bar{z}-d\bar{w})\otimes c^{\kappa-\kappa_c},
\end{split}
\end{align}
where $c^{\kappa-\kappa_c}\in\frh\otimes\frh$ is the image of $\kappa-\kappa_c\in\frh^\vee\otimes\frh^\vee$ under the isomorphism $\frh^\vee\stackrel{\sim}{\to}\frh$ given by $h_i^*\mapsto h_i$.

Since the $(-1)$-shifted Poisson structure of the theory under consideration is degenerate, we cannot apply \cref{equivalent_obstruction_Feynman_diagram_lem} directly.
However, by \cref{equ:obstructionclassgeneral}, the obstruction class $\Theta^{Wak,gen}_{inn}$ can be expanded as the sum of the weights of the graph with a distinguished edge.
The heat kernel and the propagator are decomposed as the sums of the corresponding heat kernels and propagators of $\clM_{\frg}$ and $\clM^{\hbar(\kappa-\kappa_c)}_{Heis}$, respectively.
Consequently, there are no edges connecting a $\beta,\gamma$-leg to a $b_i$-leg.
Since the heat kernel $K_\epsilon^{Heis,Wak}$ is of the order $\hbar$, the only graph that contributes to the obstruction of first-order quantization in $\hbar$ that contains an internal edge labeled by $P^{Heis,Wak}(\epsilon,L)$ or $K^{Heis,Wak}_{\epsilon}$ is the tree graph.
Also, the homogeneous parts of $\wt{S}^{Wak}_{\kappa}$ that has $b_i$-leg, are degree $2$ and $3$ functionals of $\clW_\kappa$ with one $b_i$-leg.
Thus, no graph contributing to the obstruction $\Theta^{Wak,gen}_{inn}$ has an external edge labeled by $b_i$.
This fact implies that the obstruction classes $\Theta^{Wak,gen}_{inn}$ and $\Theta^{Wak,gen}_{out}$ is elements of $\IA(\clG_1,\clM_{\frg})$ and $\Act(\clG_1,\clM_{\frg})$, respectively.

\begin{figure}[htbp]
% --- 1つ目の図 ---
\begin{tikzpicture}[node distance=3cm, thick]
% ノード（丸の中に I）
\node[draw, circle, minimum size=12mm] (A) {$I_{\eta\gamma\beta}$};
\node[draw, circle, minimum size=12mm, right of=A] (B) {$I_{\eta\gamma\beta}$};
% 外側の線とラベル
\draw (A) -- ++(-1.5,0) node[pos=1, left=3pt] {$\eta\otimes h$};
\draw (B) -- ++(1.5,0)  node[pos=1, right=3pt] {$\eta\otimes h$};
% 内側の 2 本線
\draw (A) to[bend left=20]  node[midway, above=3pt] {$P^{\beta\gamma}(\epsilon,L)$} (B);
\draw (A) to[bend right=20] node[midway, below=3pt] {$K_\epsilon^{\beta\gamma}$} (B);
\end{tikzpicture}

% --- 2つ目の図 ---
\begin{tikzpicture}[node distance=3cm, thick]
% ノード（丸の中に I）
\node[draw, circle, minimum size=12mm] (A) {$I_{\eta b}$};
\node[draw, circle, minimum size=12mm, right of=A] (B) {$I_{\eta b}$};
% 外側の線とラベル
\draw (A) -- ++(-1.5,0) node[pos=1, left=3pt] {$\eta\otimes h$};
\draw (B) -- ++(1.5,0)  node[pos=1, right=3pt] {$\eta\otimes h$};
% 内側の 2 本線
\draw (A) to node[midway, above=3pt] {$K^{Heis,Wak}_\epsilon$} (B);
\end{tikzpicture}

\caption{The graphs which contribute to the obstruction class $\Theta^{Wak,gen}_{inn}$.}
\label{obstruction_class_inn_gene}
\end{figure}
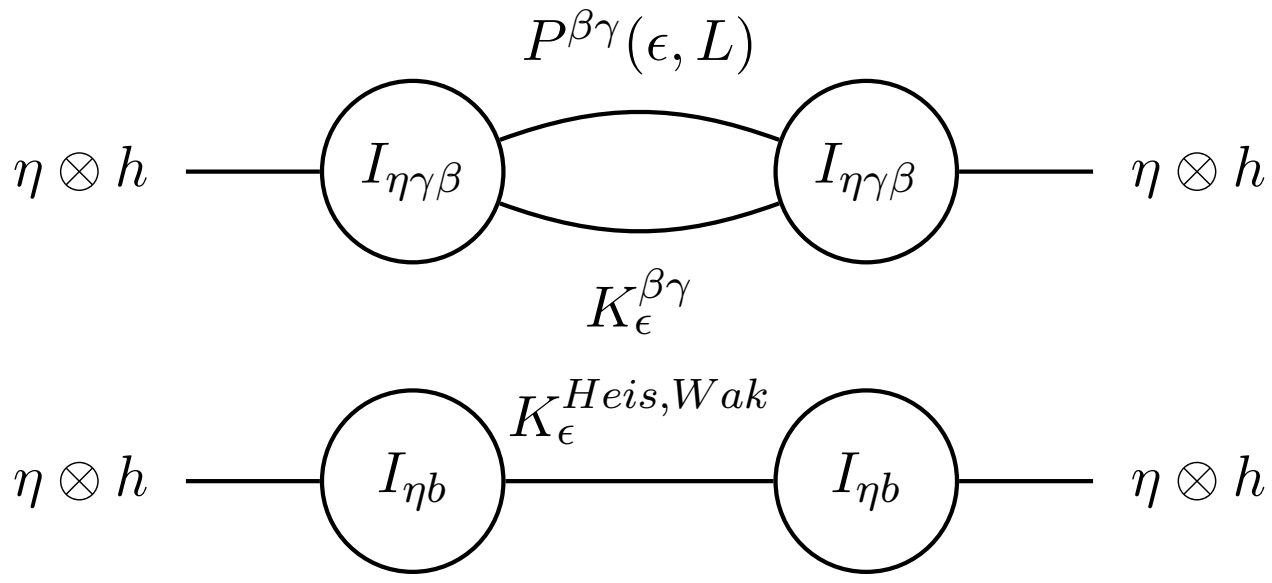

Similarly to the proof of \cref{Wakimoto_critical}, by the isomorphism \eqref{isom_inner_cohomology}, the class $[\Theta^{Wak,gen}_{inn}]\in H^1(\IA(\clG_1,\mathcal{M_{\frg}}))$ is determined by the weight of graphs whose external edges are all $\eta\otimes h$-legs.
Such graphs are given in \cref{obstruction_class_inn_gene}, where the functional $I_{\eta b}$ for $\eta\in\Omega^{0,*}_\bbC$ is defined by
\[
  I_{\eta b} \ceq \int_{\bbC}\eta\wedge b.
\]

The weight of the first graph in \cref{obstruction_class_inn_gene} is already computed in the proof of \cref{Wakimoto_critical}, and is given by $i(\phi(\kappa_c))$.

We consider the weight $\Theta_{2nd}$ of the second graph in \cref{obstruction_class_inn_gene}.
We consider its value on $\eta_i\otimes h_i$, where $\eta_i$ is a Dolbeault form on $\bbC$, and $h_i$ is a basis element of the Cartan subalgebra $\frh$.
As in the typical cases, this decomposes into the analytic part and algebraic part:
\[
\Theta_{2nd}(\eta_i\otimes h_i)=\Theta_{2nd}^{an}(\eta_i)\Theta_{2nd}^{alg}(h_i).
\]
Let $f,g$ be smooth functions on $\bbC$ and $\eta_i=f+gd\bar{z}$ be a differential form.
Then we have
\begin{align*}
  \Theta_{2nd}^{an}(\eta_i)&=\lim_{L\to0}\lim_{\epsilon\to0}\int_{\bbC_{z_1}\times\bbC_{z_2}}f(z_1)f(z_2)\wedge K_\epsilon^{Heis,Wak}(z_1,z_2)\\
                           &=\lim_{L\to0}\lim_{\epsilon\to0}\int_{\bbC_{z_1}\times\bbC_{z_2}}f(z_1)f(z_2)\frac{1}{2\pi i\epsilon}\partial_{z_2}e^{-|z_1-z_2|^2/4\epsilon}d^2z_1d^2z_2\\
                           &=-\lim_{L\to0}\lim_{\epsilon\to0}\int_{\bbC_{z_1}\times\bbC_{z_2}}f(z_1)\partial_{z_2}f(z_2)\frac{1}{2\pi i\epsilon}e^{-|z_1-z_2|^2/4\epsilon}d^2z_1d^2z_2\\
                           &=-\lim_{L\to0}\lim_{\epsilon\to0}\int_{\bbC_{w_1}\times\bbC_{w_2}}f\left(\frac{w_1+w_2}{2}\right)\partial_{w_2}f\left(\frac{w_2-w_1}{2}\right)\frac{1}{2\pi i\epsilon}e^{-|w_1|^2/4\epsilon}d^2w_1d^2w_2\\
                           &=-\lim_{L\to0}\lim_{\epsilon\to0}\int_{\bbC_{w_2}}f(w_2)\left(-\frac{\partial}{\partial w_2}f(w_2)\right)\frac{1}{2\pi i\epsilon}\bigl(\epsilon+O(\epsilon^2)\bigr)d^2w_2\\
                           &=\frac{1}{2\pi i}\int_{\bbC_z}f(z)\frac{\partial}{\partial z}f(z)d^2z,
\end{align*}
where the fifth equality follows from Taylor's theorem and the Gaussian integral.

By \eqref{equ: Heis heat kernel}, the contraction by the algebraic part of the heat kernel $K^{Heis,Wak}_\epsilon$ is nothing but the contraction by $c^{\kappa-\kappa_c}$.
Thus, the algebraic part $\Theta^{alg}_{2nd}$ is given as follows
\[
  \Theta_{2nd}^{alg}(h_i)=(\kappa-\kappa_c)(h_i,h_i).
\]

Therefore, we have
\[
  [\Theta^{Wak,gen}_{inn}] = 
  H^1(i)\left(\left[\phi(\kappa_c)+\phi(\kappa-\kappa_c)\right]\right) =
  H^1(i)([\phi(\kappa)]).
\]
Similarly to the proof of \cref{Wakimoto_critical}, there exists an element $\tau'\in\IA(\clG_1,\clM_{\frg})$ has cohomological degree $0$ such that
\begin{align}
  \label{gen_gamma_exact}
  i(\phi(\kappa))-\Theta^{Wak,gen}_{inn}=d_{\clG_1}\tau'+\bar{\partial}\tau'+\{\wt{S}^{Wak}_{\kappa},\tau'\}.
\end{align}
Then, set
\begin{align*}
    S^{Wak}_{\kappa,\tau'} \ceq \wt{S}^{Wak}_{\kappa}+\hbar\tau'\in\Act(\clG_1,\clW_\kappa)[\hbar].
\end{align*}
By \cref{holomorphic_counterterm}, for $L>0$, there exists a limit
\[
  S^{Wak}_{\kappa,\tau'}[L] \ceq
  \lim_{\epsilon\to0}W(P(\epsilon,L),S^{Wak}_{\kappa,\tau'}),
\]
and by \eqref{gen_gamma_exact}, the obstruction class
\[
  \frac{1}{\hbar}\lim_{L\to0}\left.\Bigl(d_{\clG_1}S^{Wak}_{\kappa,\tau'}[L]+\bar{\partial}S^{Wak}_{\kappa,\tau'}[L]+\frac{1}{2}\{S^{Wak}_{\kappa,\tau'}[L],S^{Wak}_{\kappa,\tau'}[L]\}_L+\hbar\Delta_LS^{Wak}_{\kappa,\tau'}[L]\Bigr)\right|_{\Act(\clG_1,\clW_\kappa)}\ \bmod \hbar
\]
vanishes.
Thus, the collection of the functionals $\{S^{Wak}_{\kappa,\tau'}[L]\}_{L>0}$ defines a quantum outer-action of $\clG_1$ on the quantization of $\clW_\kappa$.

\begin{figure}[htbp]
% --- 1つ目の図 ---
\begin{tikzpicture}[node distance=3cm, thick]

% ノード（丸の中に I）
\node[draw, circle, minimum size=12mm] (A) {$I_{\eta_i\gamma\beta}$};
\node[draw, circle, minimum size=12mm, right of=A] (B) {$I_{\xi_i\partial\gamma}$};

% 外側の線とラベル
\draw (A) -- ++(-1.5,0.5) node[pos=1, left=3pt] {$\eta_i\otimes h_i$};
\draw (A) -- ++(-1.5,-0.5) node[pos=1, left=3pt] {$\gamma_{\alpha_i}$};
\draw (B) -- ++(1.5,0)  node[pos=1, right=3pt] {$\xi_{\alpha_i}\otimes f_{\alpha_i}$};

% 内側の 2 本線
\draw (A) to node[midway, above=3pt] {$K_\epsilon^{\beta\gamma}$} (B);

\end{tikzpicture}

% --- 2つ目の図 ---
\begin{tikzpicture}[node distance=3cm, thick]

% ノード（丸の中に I）
\node[draw, circle, minimum size=12mm] (A) {$I_{\eta_ib_i}$};
\node[draw, circle, minimum size=12mm, right of=A] (B) {$I_{\xi_ib_i\gamma}$};

% 外側の線とラベル
\draw (A) -- ++(-1.5,0) node[pos=1, left=3pt] {$\eta_i\otimes h_i$};
\draw (B) -- ++(1.5,0.5)  node[pos=1, right=3pt] {$\xi_{\alpha_i}\otimes f_{\alpha_i}$};
\draw (B) -- ++(1.5,-0.5)  node[pos=1, right=3pt] {$\gamma_{\alpha_i}$};

% 内側の 2 本線
\draw (A) to node[midway, above=3pt] {$K^{Heis,Wak}_\epsilon$} (B);

\end{tikzpicture}

\caption{The graphs with three external edges labeled by $\eta_i\otimes h_i$, $\xi_{\alpha_i}\otimes f_{\alpha_i}$, and $\gamma_{\alpha_i}$.}

\label{obstruction_class_d_alpha_determine}
\end{figure}
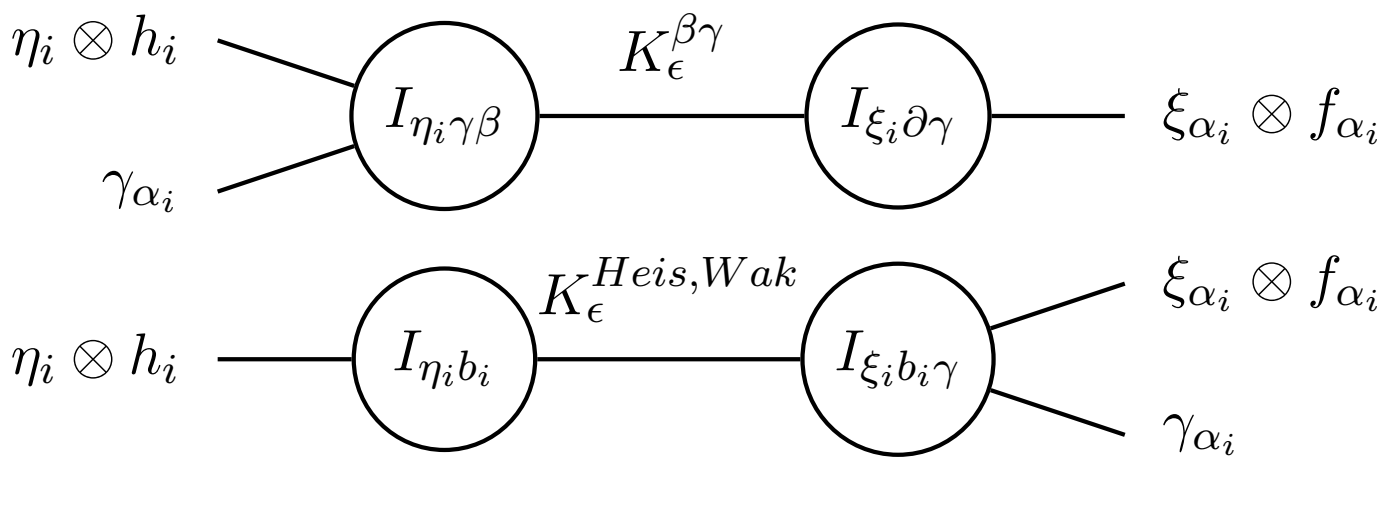

Next, we determine $\tau'$ using an argument similar to that in \cref{Wakimoto_critical}.
We denote the weight of the section $b_i$ by
\[
  \operatorname{wt}b_i=0,\quad i=1,\dots,\rk{\frg}.
\]
Since the functionals $\wt{S}^{Wak}_\kappa$ and $\phi(\kappa)$ are also weight 0 homogeneous functionals, by an argument similar to that of \cref{Wakimoto_critical}, we obtain the following assertion:

\begin{prop}
  \label{gamma_general}
  There exists $c_\alpha'\in\bbC\ (\alpha\in\Delta_+)$ such that
  \[
  \tau'(\omega_{\alpha}\otimes e_\alpha,\eta_i\otimes h_i,\xi_\alpha\otimes f_\alpha,\gamma_\alpha,\beta_\alpha,b_i)=\sum_{\alpha\in\Delta_+}c_\alpha'\int_{\bbC}\xi_\alpha\wedge\partial\gamma_\alpha
  \]
  satisfies the equality \cref{gen_gamma_exact}.
\end{prop}
Finally, we determine the coefficients $c'_\alpha$.
Since the functional $\hbar\tau'$ is of order $\hbar$, the only graph contributing to the obstruction for the first-order quantization in $\hbar$ that contains an vertex labeled by $\hbar\tau'$ is the tree graph with a single vertex labeled by $\hbar\tau$, whose edges are all labeled by either $K^{\beta\gamma}_\epsilon$ or $P^{\beta\gamma}(\epsilon,L)$.
Moreover, writing $c'_\alpha = c_\alpha + d_\alpha$, where $c_\alpha$ is the coefficient appearing in \cref{Wakimoto_field_theory_crit}, the contributions of graphs with neither a $b_i$-leg nor an edge labeled by $P^{Heis,Wak}(\epsilon,L)$ or $K^{Heis,Wak}_\epsilon$ cancel each other.
Therefore, the coefficients $d_{\alpha_i}$ can be determined from the condition that the weights of the graphs in \cref{obstruction_class_d_alpha_determine} cancel, where the functionals $I_{\xi\partial\gamma}$ and $I_{\xi b\gamma}$ for $\xi_{\alpha_i}\in\Omega^{0,*}_\bbC$ are defined by
\begin{align*}
I_{\xi_i\partial\gamma}& \ceq \hbar d_{\alpha_i}\int_{\bbC}\xi_{\alpha_i}\wedge\partial\gamma_{\alpha_i},\\
  I_{\xi_ib_i\gamma}& \ceq \int_{\bbC}\xi_{\alpha_i}\wedge b_i\wedge\gamma_{\alpha_i}.
\end{align*}

It is easy to see that the analytic parts of the weights of the two graphs coincide.
The algebraic part of the weight of the first graph in \cref{obstruction_class_d_alpha_determine} is given as follows
\begin{align*}
  -\sum_{\alpha\in\Delta_+}\alpha(h_i)\delta_{\alpha,\alpha_i}d_{\alpha_i}=-\alpha_i(h_i)d_{\alpha_i}=-2d_{\alpha_i}.
\end{align*}
The algebraic part of the weight of the second graph in \cref{obstruction_class_d_alpha_determine} is given as follows
\begin{align*}
  (\kappa-\kappa_c)(h_i,h_i)=2(\kappa-\kappa_c)(e_{\alpha_i},f_{\alpha_i}).
\end{align*}
Indeed, for a simple Lie algebra $\frg$, the Chevalley generators $\{e_{\alpha_i},h_i,f_{\alpha_i}\}_{i=1,\dots,\rk\frg}$, and the killing form $\kappa_{\frg}$, the following equality holds
\[
\kappa_{\frg}(h_i,h_i)=2\kappa_{\frg}(e_{\alpha_i},f_{\alpha_i}),\quad i=1,\dots,\rk\frg
\]
and for an invariant symmetric bilinear form $\kappa$ on $\frg$, there exists a complex number $k\in\bbC$ such that $\kappa=k\kappa_{\frg}$.
Thus, we have
\begin{align*}
  d_{\alpha_i}=(\kappa-\kappa_c)(e_{\alpha_i},f_{\alpha_i}),
\end{align*}
for $i=1,\dots,\rk\frg$.
For simplicity, we check this only for simple roots.
However, by the definition of the functional \eqref{general_level_naive_Wakimoto}, the same calculation can be carried out for an arbitrary positive root by expanding it as a linear combination of simple roots, and we obtain the following.
\begin{prop}
\label{gamma_Wakimoto_katati}
The functional
  \[
  \tau'(\omega_{\alpha}\otimes e_\alpha,\eta_i\otimes h_i,\xi_\alpha\otimes f_\alpha,\gamma_\alpha,\beta_\alpha,b_i)=\sum_{\alpha\in\Delta_+}\left(c_{\alpha}+(\kappa-\kappa_c)(e_{\alpha},f_{\alpha})\right)\int_{\bbC}\xi_\alpha\wedge\partial\gamma_\alpha(z)
  \]
  satisfies the equality \eqref{gen_gamma_exact}, where $c_{\alpha}$ is the coefficient appearing in \cref{Wakimoto_field_theory_crit}.
\end{prop}

Taking $c_\alpha$ as in \cref{Wakimoto_field_theory_crit} and setting 
\[
S^{Wak}_{\kappa}=\wt{S}^{Wak}_{\kappa}+\hbar\sum\sum_{\alpha\in\Delta_+}\left(c_{\alpha}+(\kappa-\kappa_c)(e_{\alpha},f_{\alpha})\right)\int_{\bbC}\xi_\alpha\wedge\partial\gamma_\alpha,
\]
we obtain the desired functional.
The proof of \cref{Wakimoto_general_level_FA} is now complete.

By \cref{quantum_Noether_map,Wakimoto_general_level_FA}, we have the following assertion:

\begin{cor}
\label{Noether_map_Wakimoto_general_level}
The quantum outer-action $\{S^{Wak}_{\kappa}[L]\}_{L>0}$ induces the $\bbC\dbr{\hbar}$-linear map
\[
w^{fact}_{\kappa}:\bbU^\hbar_{\phi(\kappa)}(\clG_1)[\hbar^{-1}]\to\Obs^q_{\clW_{\kappa}}[\hbar^{-1}]
\]
of factorization algebras.
The map $w^{fact}_{\kappa}$ induces a map
\[
w^{ch}_{\kappa,\hbar}:V^{\kappa}(\frg)\dpr{\hbar}\to\beta\gamma^{\oplus\dim\frn_+}_\hbar\otimes\pi^{\hbar(\kappa-\kappa_c)}_\hbar[\hbar^{-1}]
\]
of $\bbC\dpr{\hbar}$-linear vertex algebras.
Specializing to $\hbar=2\pi i$, we obtain the Wakimoto realization at level $\kappa$:
\[
w^{ch}_{\kappa}\colon V^{\kappa}(\frg)\lto\beta\gamma^{\oplus\dim\frn_+}\otimes\pi^{\kappa-\kappa_c}.
\]
\end{cor}

\begin{proof}
  Similarly to \cref{def_gerators_of_Obs}, we define
  \begin{align*}
    u^i_n(x)\colon\Omega^{1,0}(D(x,R))\otimes\frh^\vee\ni\beta dz\otimes f\lmto\frac{1}{(n-1)!}(\partial^{n-1}_z\beta)(x)f(h_i),
  \end{align*}
  for $i=1,\dots,\rk\frg,n>0,x\in\bbC,R>0$. Then similarly to \cref{Obs_cl_generator}, the monomial
  \[
  u^{\vec{i}}_K(0),
  \]
  for $m\geq0,\ \vec{i}=(i_1,\dots,i_m)\in\{1,\dots,\rk\frg\}^m,\ K=(k_1,\dots,k_m)\in\bbZ_{>0}^{m}$ forms a basis of $\bbV(\Obs^q_{Heis,\hbar(\kappa-\kappa_c)})\simeq\pi^{\hbar(\kappa-\kappa_c)}_\hbar$.

  Let $f\in C^\infty_c(\bbR_{\geq0})$ be a compactly supported smooth function on $\bbR_{\geq0}$ such that $f(0)=0$ and $\int_{\bbR}rf(r)dr=\frac{1}{2\pi}$.
  Then, the same argument as in the proof of  \cref{Noether_map_Wakimoto_critical_level} gives the following.
  \begin{align}
    \label{Wakimoto_gene_siki}
    \begin{split}
          \hbar\lim_{L\to0}w^{fact}_{\kappa}(f(|z|)d\bar{z}\otimes e_{\alpha_i})[L]&=b_1^{\alpha_i}(0)+\sum_{\beta\in\Delta_+}P^{\alpha_i}_\beta(c_0^\beta(0))b_1^\beta(0),\\
          \hbar\lim_{L\to0}w^{fact}_{\kappa}(f(|z|)d\bar{z}\otimes h_i)[L]&=-\sum_{\beta\in\Delta_+}\beta(h_i)c^\beta_0(0)b_1^\beta(0)+u^i_1(0)\\
          \hbar\lim_{L\to0}w^{fact}_{\kappa}(f(|z|)d\bar{z}\otimes f_{\alpha_i})[L]&\\
          =\sum_{\beta\in\Delta_+}Q^{\alpha_i}_\beta(c_0^\beta(0))b_1^\beta(0)+&\hbar(c_{\alpha_i}+(\kappa-\kappa_c)(e_{\alpha_i},f_{\alpha_i}))c^{\alpha_i}_1(0)+u^i_1(0)c^{\alpha_i}_0(0),
    \end{split}
  \end{align}
  for $i=1,\dots,\rk\frg$.
  By \cref{rem:cdprhbarlinear}, the specialization $w^{ch}_{\kappa}:=w^{ch}_{\kappa,\hbar}|_{\hbar=2\pi i}$ is well-defined.
  By \cref{reconstruction_thm_Ch} and the equation \eqref{Wakimoto_gene_siki}, we obtain the desired result.
\end{proof}

\subsection{The injectivity of the map \texorpdfstring{$H^*(w^{fact}_\kappa)$}{H*(w-fact-kappa)}}
\label{ss:4.4}

In this subsection, for any $\kappa\in(\Sym^2\frg^\vee)^\frg$, we show that the map $H^*(w^{fact}_\kappa)$ of factorization algebras given in \cref{Noether_map_Wakimoto_general_level} is injective.
This is the factorization algebra analogue of the well-known injectivity of the Wakimoto realization for vertex algebras \cite[Proposition 5.2]{Frenkel}.
Using a filtration analogous to the one employed in the proof of Proposition 5.2 of \cite{Frenkel}, we show that the induced map on the associated graded spaces is injective.
%, i.e. we show the following Theorem.

\begin{thm}
\label{thm_inj_C}
For any $\kappa\in(\Sym^2\frg^\vee)^\frg$, the induced map
\[
  H^*(w^{fact}_\kappa)\colon 
  H^*(\bbU^\hbar_{\phi(\kappa)}(\clG_{1}))[\hbar^{-1}] \lto 
  H^*(\Obs^q_{\clW_\kappa})[\hbar^{-1}]
\]
is an injective map of factorization algebras.
\end{thm}

Below, we show the above theorem.

We denote
\[
  \frB \ceq \{D(z_1,r_1)\sqcup\cdots\sqcup D(r_l,z_l) \mid 
  l\in\bbZ_{\geq1},\ r_i\in\bbR_{>0},\ z_i\in\bbC,\ D(z_i,r_i)\cap D(z_j,r_j)=\emptyset,\ (i\neq j)\}.
\]
Then $\frB$ is an open basis of $\bbC$, and for each open subset $U\subset\bbC$, the family
\[
  \frB(U) \ceq 
  \{D(z_1,r_1)\sqcup\dots\sqcup D(r_l,z_l)\in\frB \mid 
  D(z_i,r_i)\subset U \ (i=1,\dots,l)\}
\]
is a Weiss cover of $U$.
Thus, to show that the map
\[
  H^*(w^{fact}_\kappa)\colon 
  H^*(\bbU^\hbar_{\phi(\kappa)}(\clG_{1}))[\hbar^{-1}] \lto 
  H^*(\Obs^q_{\clW_\kappa})[\hbar^{-1}]
\]
is an injective map of factorization algebras, it is enough to show that for every $z\in\bbC,\ R>0$, the map $H^*(w^{fact}_\kappa)(D(z,R))$ of differentiable vector spaces is injective.

For $p\geq0$, set
\begin{align*}
   F_p\bbU^\hbar_{\phi(\kappa)}(\clG_1)(D(z,R)) \ceq 
   \Sym^{\leq p}\Bigl(\Omega^{0,*}_{\bbC,c}(D(z,R))\otimes\frg[1]\Bigr)\dbr{\hbar}
\end{align*}
and
\begin{align*}
 &F_p\Obs^q_{\clW_\kappa}(D(z,R)) \ceq 
  \Sym\left((\Omega^{0,*}_{\bbC}\otimes\frn_+[-1])^*(D(z,R))[-1]\right)\\
 &\otimes\prod_{k+l\leq p}
  \Sym^k\left((\Omega^{1,*}_{\bbC}\otimes\frn_+^\vee[-1])^*(D(z,R))[-1]\right)
  \otimes
  \Sym^l\left((\Omega^{1,*}_{\bbC}\otimes\frh^\vee[-1])^*(D(z,R))[-1]\right)\dbr{\hbar}.
\end{align*}
Then $\gr\bbU^\hbar_{\phi(\kappa)}(\clG_1)(D(z,R))$ and $\gr\Obs^q_{\clW_1}(D(z,R))$ denote the graded spaces associated to the filtrations $\{F^p\bbU^\hbar_{\phi(\kappa)}(\clG_1)(D(z,R))\}_{p\geq0}$ and $\{F^p\Obs^q_{\clW_1}(D(z,R))\}_{p\geq0}$, respectively.
Since the map $w^{fact}_\kappa$ and differentials preserve these filtrations, we have the associated map
\[
  \gr{w}^{fact}_\kappa(D(z,R))\colon 
  \gr\bbU^\hbar_{\phi(\kappa)}(\clG_1)(D(z,R))[\hbar^{-1}] \lto
  \gr\Obs^q_{\clW_\kappa}(D(z,R))[\hbar^{-1}]
\]
of cochain complexes.
Thus, it is enough to show that the induced map
\[
  H^*(\gr{w}^{fact}_\kappa(D(z,R)))\colon 
  H^*(\gr\bbU^\hbar_{\phi(\kappa)}(\clG_1)(D(z,R)))[\hbar^{-1}] \lto 
  H^*(\gr\Obs^q_{\clW_\kappa}(D(z,R)))[\hbar^{-1}]
\]
is injective.
Take a compactly supported smooth function $f\in C^\infty_c(\bbR_{\ge0})$ such that $f(0)=0$ and $\int_{\bbR_{\geq0}}rf(r)dr=1$.
Then we have isomorphisms
\begin{align*}
  &H^i(\gr\bbU^\hbar_{\phi(\kappa)}(\clG_1)(D(z,R)))\\
  \simeq&\begin{dcases}
    \operatorname{Span}_{\bbC}\left\{f(|z_1|)z_1^{-n_1}d\bar{z}_1\otimes x^1\cdots f(|z_p|)z_p^{-n_p}d\bar{z}_p\otimes x^p\mid \substack{p\ge0,\ x^i\in\frg,\\ n_1,\dots,n_p\ge0}\right\}\dbr{\hbar}\\
    \simeq\bbC[\clJ_\infty\frg^*]\dbr{\hbar},  &\text{for }i=0,\\
    0, &\text{otherwise},
  \end{dcases}\\
  &H^i(\gr\Obs^q_{\clW_\kappa}(\clG_1)(D(z,R)))\\
  \simeq&\begin{dcases}
    \operatorname{Span}_{\bbC}\left\{\left.b^{\vec{\alpha}}_K(0)c^{\vec{\beta}}_L(0)u^{\vec{i}}_{M}(0)\right|\substack{n,m,p\geq0,K\in\bbZ^n_{>0},L\in\bbZ_{\ge0}^m,M\in\bbZ_{>0}^p,\\ \vec{\alpha}\in\Delta_+^n,\vec{\beta}\in\Delta_+^m,\vec{i}\in\{1,\dots,\rk{\frg}\}^p}\right\}\dbr{\hbar}\\
    \simeq\bbC[\clJ_\infty(\frn_+\oplus\frn_+^\vee\oplus\frh)]\dbr{\hbar}\simeq\bbC[\clJ_\infty(T^*U\times\frh)]\dbr{\hbar},  &\text{for }i=0,\\
    0, &\text{otherwise},
  \end{dcases}
\end{align*}
Similarly to the proof of \cref{Noether_map_Wakimoto_critical_level}, for $i=1,\dots,\rk{\frg}$, we obtain
\begin{align}
\label{Wakimoto_gene_another_gr_siki}
\begin{split}
  H^0(\gr{w}^{fact}_\kappa)(f(|z|)d\bar{z}\otimes e_{\alpha_i})
  &=\frac{1}{\hbar}\left(b_1^{\alpha_i}(0)+\sum_{\beta\in\Delta_+}P^{\alpha_i}_\beta(c_0^\beta(0))b_1^\beta(0)\right),\\
  H^0(\gr{w}^{fact}_\kappa)(f(|z|)d\bar{z}\otimes h_i)
  &=\frac{1}{\hbar}\left(-\sum_{\beta\in\Delta_+}\beta(h_i)c^\beta_0(0)b_1^\beta(0)+u^i_1(0)\right)\\
  H^0(\gr{w}^{fact}_\kappa)(f(|z|)d\bar{z}\otimes f_{\alpha_i})
  &=\frac{1}{\hbar}\left(\sum_{\beta\in\Delta_+}Q^{\alpha_i}_\beta(c_0^\beta(0))b_1^\beta(0)+u^i_1(0)c^{\alpha_i}_0(0)\right).
\end{split}
\end{align}

Thus, the map $H^0(\gr{w}^{fact}_\kappa)$ is identified with the tensor product of the map $\clJ_\infty(\gr{U}(\wt{\rho})):\bbC[\clJ_\infty\frg^*]\to\bbC[\clJ_\infty(T^*N_+\times\frh)]$ with $\bbC\dpr{\hbar}$, and hence by \cref{lem_jet_tilderho}, it is injective.
Thus, the map $H^0(\gr{w}^{fact}_\kappa)(D(z,R))$ is injective, and hence the map $H^0(w^{fact}_\kappa)(D(z,R))$ is injective.

The proof of \cref{thm_inj_C} is now complete.

The following assertion follows immediately from \cref{Noether_map_Wakimoto_general_level,thm_inj_C}.

\begin{cor}
\label{cor_Wakimoto_inj}
For any $\kappa\in(\Sym^2\frg^\vee)^\frg$, the maps $w^{ch}_{\kappa,\hbar}$ and $w^{ch}_\kappa$ are injective.
\end{cor}

This corollary recovers the well-known injectivity of the Wakimoto realization (see Proposition 5.2 of \cite{Frenkel}).

%%%%%%%%%%%%%%%%%%%%%%%%%%%%%%%%%%%%%%%%%%%%%%%%%%%%%%%%%%%%%%%%%%%%%%%%%%%%%%%%
\subsection{The BD center and critical level}
\label{ss:4.5}

It is known that the Wakimoto realization induces a free field realization of the Feigin–Frenkel center (see Chapter 9.1 of \cite{Frenkel}).
In this subsection, we consider a factorization algebra analogue of this result. In \cref{ss:3.6}, we introduced the notion of the BD center and computed the BD centers of the factorization algebras of observables associated with several field theories. By \cref{prop_commutant_WZW}, the vertex algebra corresponding to the BD center of the factorization algebra of observables of $\clM_{WZW,\hbar\kappa_c}$ is the Feigin–Frenkel center.

In this subsection, we construct a map of BD factorization algebras between the factorization algebras of observables from the map of factorization algebras constructed in \cref{Noether_map_Wakimoto_general_level}.
We also show that this map induces a map between their BD centers \cref{BV_center_def}, which gives a free field realization of the Feigin--Frenkel center.

The proof of \cref{BD_center_intertwine} follows the approach of the proof of Lemma 9.1 of \cite{Frenkel}.

Let $\frg$ be a simple Lie algebra and $\kappa\in(\Sym^2\frg^\vee)^\frg$ be a 2-cocycle.

\begin{prop}
\label{BD_FA_morphism_Wakimoto}
There exists a map
\begin{align*}
  \scW^{fact}_\kappa\colon 
  \Obs^q_{WZW,\hbar\kappa} \lto \Obs^q_{\clW_\kappa}
\end{align*}
of BD factorization algebras on $\bbC$ such that the induced map
\[
  H^*(\scW^{fact}_\kappa)\colon 
  H^*(\Obs^q_{WZW,\kappa}) \lto H^*(\Obs^q_{\clW_\kappa})
\]
is injective.
\end{prop}

\begin{proof}
  By \cref{WZW_VA}, we have a quasi-isomorphism
  \[
  \Obs^q_{WZW,\hbar\kappa}\qis\bbU^\hbar_{\hbar^2\phi(\kappa)}(\clG_1^\hbar)
  \]
  of factorization algebras.
  %By the proof of \cref{prop_commutant_WZW}, the map $H^*(\Delta_{WZW,\infty}):H^*(\Obs^q_{WZW,\hbar\kappa})\to H^*(\Obs^q_{WZW,\hbar\kappa})$ on the cohomology which is induced by the BV-Laplacian $\Delta_{WZW,\infty}$ on the field theory $\clM_{WZW,\hbar\kappa}$ is identified with the Lie bracket on $H^*(\Obs^q_{WZW,\hbar\kappa})\simeq\Sym(\mathcal{O}_{\bbC}\otimes\frg)\dbr{\hbar}$ given by
  %\[
  %[f\otimes x,g\otimes y]=\hbar fg\otimes[x,y]+\hbar^2\left(\int f\partial_zgd^2z\right)\kappa(x,y),
  %\]
  %for sections $f\otimes x,g\otimes y$ of $\mathcal{O}_{\bbC}\otimes\frg\dbr{\hbar}$.
  %Thus, we have a quasi-isomorphism of cosheaves of cochain complexes
  %\[
  %\Obs^q_{WZW,\hbar\kappa}\simeq\bbU^\hbar_{\hbar^2\phi(\kappa)}(\clG_1^\hbar),
  %\]
  %where $\clG_1^\hbar=\Omega^{0,*}_{\bbC}\otimes\frg\dbr{\hbar}$ is a sheaf of dgLA on $\bbC$ with differential $\bar{\partial}$ and Lie bracket given by
  %\[
  %[\alpha\otimes x,\beta\otimes y]=\hbar\alpha\wedge\beta\otimes[x,y],
  %\]
  %for sections $\alpha\otimes x,\beta\otimes y$ of $\Omega^{0,*}_{\bbC}\otimes\frg\dbr{\hbar}$.
  
  Here, for an open subset $U\subset\bbC$, we consider the following linear map
  \begin{align*}
    F\colon \clG_1^\hbar \lto\clG_1\dbr{\hbar}, \quad 
    \alpha\otimes x \lmto\hbar\alpha\otimes x.
  \end{align*}
  Then the linear map $F$ induces a linear map $\Sym(F):\bbU^\hbar_{\hbar^2\phi(\kappa)}(\clG_1^\hbar)\to\bbU^\hbar_{\phi(\kappa)}(\clG_1)$ of precosheaves on $\bbC$.
  Then it is easy to show that the map $\Sym(F)$ is a dgLA map, and hence the composition
  \[
  \Obs^q_{WZW,\hbar\kappa} \lsto 
  \bbU^\hbar_{\hbar^2\phi(\kappa)}(\clG_1^\hbar) \xrr{\Sym(F)} 
  \bbU^\hbar_{\phi(\kappa)}(\clG_1)
  \]
  of factorization algebras.
  We define a map $\scW^{fact}_{\kappa}:\Obs^q_{WZW,\hbar\kappa}\to\Obs^q_{\clW_\kappa}[\hbar^{-1}]$ by the following composition
  \begin{align*}
    \Obs^q_{WZW,\hbar\kappa} \sto \bbU^\hbar_{\hbar^2\phi(\kappa)}
    (\clG_1^\hbar)\xrightarrow{\Sym(F)}
    \bbU^\hbar_{\phi(\kappa)}(\clG_1) \inj 
    \bbU^\hbar_{\phi(\kappa)}(\clG_1)[\hbar^{-1}] 
    \xrightarrow{w^{fact}_\kappa} \Obs^q_{\clW_\kappa}[\hbar^{-1}].
  \end{align*}
  Explicitly, for an open subset $U\subset\bbC$ and a monomial $m=\alpha_1\otimes x_1\dots\alpha_n\otimes x_n\in \Sym^n(\Omega^{0,*}_{\bbC,c}(U)\otimes\frg[1])$, it is given as follows.
  \begin{align*}
    (\scW^{fact}_\kappa)_U(m) = \hbar^nw^{fact}_\kappa(m)
    =\left.\hbar^n\ev*{m,\exp\bigl(S^{Wak}_\kappa[L]/\hbar\bigr)}\right|_{0\times\clW_\kappa}.
  \end{align*}
  where $\ev*{-,-}$ is the duality pairing between the symmetric products $\Sym(\Omega^{0,*}_{\bbC,c}(U)\otimes\frg[1])$ and $\Sym((\Omega^{0,*}_{\bbC,c}\otimes\frg[-1])^*(U)[-1])$.
  Since the functional $S^{Wak}_\kappa[L]$ has degree at least one as a functional on $\clG_1$, the observable
  \[
  \left.\ev*{ m,\exp(S^{Wak}_\kappa[L]/\hbar) }\right|_{0\times\clW_\kappa}
  \in\Obs^q_{\clW_\kappa}(U)
  \]
  has degree at least $-n$ as a rational function in $\hbar$.
  Thus, the observable $(\scW^{fact}_\kappa)_U(m)$ is a polynomial in $\hbar$, and hence the image of $\scW^{fact}_\kappa$ belongs to the subfactorization algebra $\Obs^q_{\clW_\kappa}\subset\Obs^q_{\clW_\kappa}[\hbar^{-1}]$.

  Finally, the map $\scW^{fact}_\kappa$ is a map of BD factorization algebras.
  Since the map $\scW^{fact}_\kappa$ is a composition of maps of factorization algebras, it is also a map of factorization algebras.
  Take an open subset $U\subset\bbC$ of $\bbC$.
  Then the map $(\scW^{fact}_\kappa)_U$ is a map of cochain complexes, and hence we have
  \[
  (\bar{\partial}+\hbar\Delta_L^{Wak})\circ(\scW^{fact}_\kappa)_U=(\scW^{fact}_\kappa)_U\circ(\bar{\partial}+\hbar\Delta_L^{WZW}),
  \]
  where $\Delta^{Wak}_L$ and $\Delta^{WZW}_L$ are the BV Laplacians of the field theories $\clW_\kappa$ and $\clM_{WZW,\hbar\kappa}$, respectively.
  It is obvious that the map $(\scW^{fact}_\kappa)_U$ commutes with the Dolbeault differential $\bar{\partial}$, and hence the map $(\scW^{fact}_\kappa)_U$ commutes with the BV Laplacians.
  Thus, we have
  \begin{align*}
    \left\{(\scW^{fact}_\kappa)_U(-),(\scW^{fact}_\kappa)_U(-)\right\}_L^{Wak}=(\scW^{fact}_\kappa)_U\left(\left\{-,-\right\}^{WZW}_L\right),
  \end{align*}
  where $\{-,-\}^{Wak}_L$ is the BV bracket on $\clW_\kappa$ and $\{-,-\}^{WZW}_L$ is the BV bracket on $\clM_{WZW,\hbar\kappa}$.
  Therefore, the map $\scW^{fact}_\kappa$ is a map of BD factorization algebras.
  The injectivity of the map $H^*(\scW^{fact}_\kappa)$ follows from the injectivity of $H^*(w^{fact}_\kappa)$ and $\Sym(F)$.
\end{proof}

\begin{rem}
\label{inj_hbar_filtration}
Considering the filtration induced by the powers of $\hbar$, 
the injectivity of $\scW^{fact}_\kappa$ follows in the same way as  \cref{thm_inj_C}.
This method seems to provide a more natural proof.
\end{rem}

\begin{cor}
\label{cor_map_hbar_VA_Wakimoto}
The map $\scW^{fact}_{\kappa}$ defines an injective map
\[
  V^{\hbar\kappa}_\hbar(\frg) \linj
  \beta\gamma_\hbar^{\oplus\dim\frn_+}\otimes\pi^{\hbar(\kappa-\kappa_c)}_\hbar
\]
of $\bbC\dbr{\hbar}$-linear vertex algebras.
Specializing to $\hbar=2\pi i$, we obtain the Wakimoto realization at the level $\kappa$:
\[
  V^{\kappa}(\frg) \linj
  \beta\gamma^{\oplus\dim\frn_+}\otimes\pi^{\kappa-\kappa_c}
\]
\end{cor}

Recall the notion of BD center of a BD factorization algebra (\cref{BV_center_def}).
By \cref{prop_commutant_betagamma,prop_commutant_Heis}, we have
\begin{align*}
  \clZ^{BD}(\Obs^q_{\clW_\kappa})\simeq\begin{cases}
    \bbC\dbr{\hbar},\quad &\text{for }\kappa\neq\kappa_c,\\
    \Obs^q_{Heis,0},\quad &\text{for } \kappa=\kappa_c.
  \end{cases}
\end{align*}
Then we have the following assertion.
\begin{lem}
\label{BD_center_intertwine}
For any $\kappa\in(\Sym^2\frg^\vee)^\frg$, the image of the BD center $\clZ^{BD}(\Obs^q_{WZW,\hbar\kappa})$ under the map $\scW^{fact}_\kappa$ is contained in the BD center $\clZ^{BD}(\Obs^q_{\clW_\kappa})$.
In particular, for the non-critical level $\kappa\neq\kappa_c$, the BD center $\clZ^{BD}(\Obs^q_{WZW,\hbar\kappa})$ is isomorphic to $\bbC\dbr{\hbar}$.
\end{lem}
\begin{proof}
  We consider the limit as the scale $L\to 0$.
  Similarly to the proof of \cref{thm_inj_C}, it is enough to show that for an open disk $D(z,R)$, the image of $\clZ^{BD}(\Obs^q_{WZW,\hbar\kappa})(D(z,R))$ under the map $(\scW^{fact}_\kappa)_{D(z,R)}$ is contained in the BD center $\clZ^{BD}(\Obs^q_{\clW_\kappa})(D(z,R))$.
  Then we define the maps
  \begin{align*}
    \Phi_{\beta}:\Omega^{0,*}_{\bbC,c}(D(z,R))\otimes\frn_+&\to\left(\Omega^{1,*}_{\bbC}\otimes\frn_+^\vee\right)^!_c(D(z,R))\\
    %\alpha\otimes x&\mapsto\left(\beta\otimes\varphi\mapsto\left(\int_{D(z,R)}\alpha\wedge\beta\right)\varphi(x)\right),\\
    \Phi_{\gamma}:\Omega^{1,*}_{\bbC,c}(D(z,R))\otimes\frn_+^\vee&\to\left(\Omega^{0,*}_{\bbC}\otimes\frn_+\right)^!_c(D(z,R))\\
    %\alpha\otimes x&\mapsto\left(\beta\otimes\varphi\mapsto\left(\int_{D(z,R)}\alpha\wedge\beta\right)\varphi(x)\right),\\
    \Phi_{Heis}:\Omega^{0,*}_{\bbC,c}(D(z,R))\otimes\frh&\to\left(\Omega^{1,*}_{\bbC}\otimes\frh^\vee\right)^!_c(D(z,R))\\
    %\alpha\otimes x&\mapsto\left(\beta\otimes\varphi\mapsto\left(\int_{D(z,R)}\alpha\wedge\beta\right)\varphi(x)\right),\\
   \Phi_{WZW}:\Omega^{0,*}_{\bbC,c}(D(z,R))\otimes\frg&\to\left(\Omega^{1,*}_{\bbC}\otimes\frg^\vee\right)^!_c(D(z,R))\\
    %\alpha\otimes x&\mapsto\left(\beta\otimes\varphi\mapsto\left(\int_{D(z,R)}\alpha\wedge\beta\right)\varphi(x)\right)
  \end{align*}
  by integration pairing
  \begin{align*}
      \Omega^{*,*}_{\bbC,c}(D(z,R))&\lto\left(\Omega^{*,*}_\bbC\right)^!_c(D(z,R)),\\
      \alpha&\lmto\left(\beta\mapsto\int_{(D(z,R))}\alpha
      \wedge\beta\right)
  \end{align*}
  and the natural evaluation pairing (see \eqref{eq:lineobsPhi}).
  Note that for $\beta\otimes x\in\Omega^{0,*}_{\bbC,c}(D(z,R))\otimes\frg$, we have
  \begin{align*}
    \lim_{L\to0}(\scW^{fact}_\kappa)_{D(z,R)}(\Phi_{WZW}(\beta\otimes x))[L]=\hbar\lim_{L\to0}(w^{fact}_\kappa)_{D(z,R)}(\beta\otimes x)[L],
  \end{align*}
  Since the image of $\beta\otimes e_\alpha\in\Omega^{0,*}_{\bbC,c}(D(z,R))\otimes\frn_+$ contains a linear observable $\Phi_{\beta}(\alpha\otimes e_\alpha)$ for some $\alpha\in\Omega^{1,*}_{\bbC,c}(D(z,R))$, the following monomials generate the cochain complex $\Obs^q_{\clW_\kappa}(D(z,R))$.
  \begin{align}
    \label{monomial_generator_Obs_W}
    \prod_{j=1}^J\Phi_{Heis}(b_{j}\otimes h_{i_j})\prod_{k=1}^K(\scW^{fact}_\kappa)_{D(z,R)}(\Phi_{WZW}(\omega_{k}\otimes e_{\alpha^{(c_k)}}))\prod_{l=1}^L\Phi_{\gamma}(\gamma_{l}\otimes e_{\alpha^{(d_l)}}),
  \end{align}
  for $b_{j},\in\Omega^{0,*}_{\bbC,c}(D(z,R)),\gamma_l,\omega_k\in\Omega^{1,*}_{\bbC,c}(D(z,R)),\ h_{i_j}\in\frh,\alpha^{(c_k)},\alpha^{(d_l)}\in\Delta_+$.
  Since we have
  \begin{align*}
      (\scW^{fact}_\kappa)_{D(z,R)}&(\clZ^{BD}(\Obs^q_{WZW,\hbar\kappa}))\\
      \subset&\bigcap_{\alpha\otimes x\in\Omega^{0,*}_{\bbC,c}(D(z,R))\otimes\frg}\ker\left(  \{(\scW^{fact}_\kappa)_{D(z,R)}(\Phi_{WZW}(\alpha\otimes x)),-\}^{Wak}_0\right),
  \end{align*}
  it is enough to show that if the monomial \eqref{monomial_generator_Obs_W} belongs to the intersection of the kernels of the map
  \[
  \{(\scW^{fact}_\kappa)_{D(z,R)}(\Phi_{WZW}(\alpha\otimes x)),-\}^{Wak}_0
  \]
  for $\alpha\otimes x\in\Omega^{0,*}_{\bbC,c}(D(z,R))\otimes\frg$, then it belongs to the BD center $\clZ^{BD}(\Obs^q_{\clW_\kappa})$.
  First, for a section $\omega\otimes e_\alpha\in\Omega^{0,*}_{\bbC,c}(D(z,R))\otimes\frg,\ (\alpha\in\Delta_+)$, we consider the map
  \[
  \{(\scW^{fact}_\kappa)_{D(z,R)}(\Phi_{WZW}(\omega\otimes e_\alpha)),-\}^{Wak}_0.
  \]
  Then we have
  \begin{align*}
  &\{(\scW^{fact}_\kappa)_{D(z,R)}(\Phi_{WZW}(\omega\otimes e_\alpha)),(\scW^{fact}_\kappa)_{D(z,R)}(\Phi_{WZW}(\omega'\otimes e_{\alpha'}))\}^{Wak}_0\\
  &=(\scW^{fact}_\kappa)_{D(z,R)}\left(\{\Phi_{WZW}(\omega\otimes e_\alpha),\Phi_{WZW}(\omega'\otimes e_{\alpha'})\}^{WZW}_0\right)\\
  &=(\scW^{fact}_\kappa)_{D(z,R)}\left(\Phi_{WZW}(\omega\wedge\omega'\otimes[e_\alpha,e_{\alpha'}])+\kappa(e_\alpha,e_{\alpha'})\int_{D(z,R)}\frac{d\omega}{d(d\bar{z})}\partial\frac{d\omega'}{d(d\bar{z})}\wedge d\bar{z}\right)=0,
  \end{align*}
  for $\omega'\otimes e_{\alpha'}\in\Omega^{0,*}_{\bbC,c}(D(z,R))\otimes\frg,\ (\alpha'\in\Delta_+)$.
  Since the functional $S^{Wak}_\kappa[L](\omega\otimes e_\alpha)$ is a functional of degree $0$ on $\clM^{\hbar(\kappa-\kappa_c)}_{Heis}$, for $b\otimes h\in\Omega^{0,*}_{\bbC,c}(D(z,R))\otimes\frh$, we have
  \[
  \{(\scW^{fact}_\kappa)_{D(z,R)}(\Phi_{WZW}(\omega\otimes e_\alpha)),\Phi_{Heis}(b\otimes h)\}^{Wak}_0=0.
  \]
  By \cref{lem_Q_not_constant}, the functional $S^{Wak}_\kappa[L](\omega\otimes e_\alpha)$ is the sum of a linear functional on $\Omega^{1,*}_{\bbC}\otimes\frn_+^\vee$ and a functional that is of degree $1$ as functional on $\Omega^{1,*}_{\bbC}\otimes\frn_+^\vee$ and of degree at least $1$ as functional on $\Omega^{0,*}_{\bbC}\otimes\frn_+$.
  Therefore, considering the degree of functional on $\Omega^{0,*}_{\bbC}(D(z,R))\otimes\frn_+$, we have the fact that the monomial 
  \[
  \prod_{l=1}^L\Phi_{\gamma}(\gamma_{l}\otimes e_{\alpha^{(d_l)}}) 
  \]
  belongs to the interaction of the kernels of the maps $\{(\scW^{fact}_\kappa)_{D(z,R)}(\Phi_{WZW}(\omega\otimes e_\alpha)),-\}^{Wak}_0$ if and only if the monomial is constant.

  Next, for a section $\eta\otimes h\in\Omega^{0,*}_{\bbC,c}(D(z,R))\otimes\frg$, we consider the map
  \[
  \{(\scW^{fact}_\kappa)_{D(z,R)}(\Phi_{WZW}(\eta\otimes h)),-\}^{Wak}_0.
  \]
  Then we have
  \begin{align*}
  &\{(\scW^{fact}_\kappa)_{D(z,R)}(\Phi_{WZW}(\eta\otimes h)),(\scW^{fact}_\kappa)_{D(z,R)}(\Phi_{WZW}(\omega\otimes e_{\alpha}))\}^{Wak}_0\\
  &=(\scW^{fact}_\kappa)_{D(z,R)}\left(\{\Phi_{WZW}(\eta\otimes h),\Phi_{WZW}(\omega\otimes e_{\alpha})\}^{WZW}_0\right)\\
  &=(\scW^{fact}_\kappa)_{D(z,R)}\left(\Phi_{WZW}(\eta\wedge\omega\otimes[h,e_\alpha])+\kappa(h,e_\alpha)\int_{D(z,R)}\frac{d\eta}{d(d\bar{z})}\partial\frac{d\omega}{d(d\bar{z})}\wedge d\bar{z}\right)\\
  &=\alpha(h)(\scW^{fact}_\kappa)_{D(z,R)}\left(\Phi_{WZW}(\eta\wedge\omega\otimes e_\alpha)\right),
  \end{align*}
  for a section $\omega\otimes e_\alpha\in\Omega^{0,*}_{\bbC,c}(D(z,R))\otimes\frg,\ (\alpha\in\Delta_+)$, and hence we have the fact that the monomial 
  \[
  \prod_{k=1}^K(\scW^{fact}_\kappa)_{D(z,R)}(\Phi_{WZW}(\omega_{k}\otimes e_{\alpha^{(c_k)}}))
  \]
  belongs to the interaction of the kernels of the maps $\{(\scW^{fact}_\kappa)_{D(z,R)}(\Phi_{WZW}(\eta\otimes h)),-\}^{Wak}_0$ if and only if the monomial is constant.
  Therefore, the image of the BD center $\clZ^{BD}(\Obs^q_{WZW,\hbar\kappa})$ under the map $\scW^{fact}_\kappa$ is contained in the sub-factorization algebra $\Obs^q_{Heis,\hbar(\kappa-\kappa_c)}$.
  Furthermore, the functional $S^{Wak}_\kappa[L](\eta\otimes h)$ contains a linear functional on $\clM^{\hbar(\kappa-\kappa_c)}_{Heis}$.
  Thus, considering the degree of the functional on $\clM^{\hbar(\kappa-\kappa_c)}_{Heis}$, we have the fact that the monomial 
  \[
  \prod_{j=1}^J\Phi_{Heis}(b_{j}\otimes h_{i_j})
  \]
  belongs to the interaction of the kernels of the maps $\{(\scW^{fact}_\kappa)_{D(z,R)}(\Phi_{WZW}(\eta\otimes h)),-\}^{Wak}_0$ if and only if the monomial belongs to the BD center $\clZ^{BD}(\Obs^q_{Heis,\hbar(\kappa-\kappa_c)})$.
\end{proof}

\begin{cor}
\label{cor_center_critical_level}
The map $\scW^{fact}_{\kappa_c}$ define an injective map
\[
  \clZ^{BD}(\scW^{fact}_{\kappa_c}): 
  \clZ^{BD}(\Obs^q_{WZW,\hbar\kappa_c})\linj\Obs^q_{Heis,0}
\]
of factorization algebras. 
The above map $\clZ^{BD}(\scW^{fact}_{\kappa_c})$ induces an injective map
\[
  \frz(\widehat{\frg})\dbr{\hbar} \linj \pi_\hbar^0
\]
of $\bbC\dbr{\hbar}$-linear vertex algebras.
Specializing to $\hbar=2\pi i$, we obtain the free field realization of Feigin--Frenkel center:
\[
  \frz(\widehat{\frg}) \linj \pi^0,
\]
where $\frz(\hat{\frg})$ is the Feigin--Frenkel center (see \cref{example_VA} (iv)).
\end{cor}

\begin{proof}
By \cref{prop_commutant_WZW}, we have an isomorphism
\[
  \bbV(\clZ^{BD}(\Obs^q_{WZW,\hbar\kappa_c})) \simeq 
  \clZ^{VA}(V^{\hbar\kappa_c}_\hbar(\frg)) \simeq
  \frz(\widehat{\frg})\dbr{\hbar}
\]
of $\bbC\dbr{\hbar}$-linear vertex algebras.
This implies that the described assertion.
\end{proof}

The following \cref{prop_BD_FA_WZW_to_Wak} at the critical level
can be shown in the same way as \cref{BD_FA_morphism_Wakimoto}.

\begin{prop}
\label{prop_BD_FA_WZW_to_Wak}
The following map is a map of BD factorization algebras.
\begin{align*}
  \wt{\scW}^{fact}_\kappa=\wt{w}^{fact}_{\kappa_c}\circ\Sym(F)\colon 
  \Obs^q_{WZW,\hbar\kappa_c} \lto \Obs^q_{\clM_{\frg}}
\end{align*} 
\end{prop}

\begin{prop}
\label{not_inj_critical_level}
The induced map
\[
  H^*(\wt{\scW}^{fact}_{\kappa_c})\colon
  H^*(\Obs^q_{WZW,\hbar\kappa_c}) \lto H^*(\Obs^q_{\clM_{\frg}})
\]
  is not injective.
\end{prop}

\begin{proof}
By the proof of \cref{BD_center_intertwine}, the image of $\clZ^{BD}(\Obs^q_{WZW,\hbar\kappa_c})$ under the map $\wt{\scW}^{fact}_{\kappa_c}$ is contained in the BD center $\clZ^{BD}(\Obs^q_{\clM_{\frg}})$.
On the other hand, by \cref{prop_commutant_betagamma}, we have $\clZ^{BD}(\Obs^q_{\clM_{\frg}})\simeq\bbC\dbr{\hbar}$.
By \cref{cor_center_critical_level}, the factorization algebra $\clZ^{BD}(\Obs^q_{WZW,\hbar\kappa_c})$ has nontrivial generators and the cohomology $H^*(\clZ^{BD}(\Obs^q_{WZW,\hbar\kappa_c}))$ also has nontrivial generators.
Thus, the kernel of the map $H^*(\wt{\scW}^{fact}_{\kappa_c})$ is nontrivial and hence the map $H^*(\wt{\scW}^{fact}_{\kappa_c})$ is not injective.
\end{proof}

%%%%%%%%%%%%%%%%%%%%%%%%%%%%%%%%%%%%%%%%%%%%%%%%%%%%%%%%%%%%%%%%%%%%%%%%%%%%%%%%
%%%%%%%%%%%%%%%%%%%%%%%%%%%%%%%%%%%%%%%%%%%%%%%%%%%%%%%%%%%%%%%%%%%%%%%%%%%%%%%%
\section{The higher dimensional analogue}
%The higher dimensional Wakimoto realizations 
%for the higher Kac--Moody factorization algebras
\label{s.5}

In this section, we consider the higher dimensional analogues of the results in previous section.

By \cref{equivalent_obstruction_Feynman_diagram_lem}, the obstruction class for the holomorphic field theory on $\bbC^d$ is given by the sum of the weights of wheel graphs with $(d+1)$ vertices.
This is a key difference from holomorphic field theories on $\bbC$. Therefore, to treat the general-level case, we need to construct an equivariant action functional from a “higher extension''.
In \cref{ss:5.3}, we return to the finite-dimensional setting, introduce the higher extension, and show that it carries an $L_\infty$-algebra structure. In \cref{ss:5.4}, as an analogue of \cref{ss:4.3}, we show that level deformations can be constructed using the equivariant action functional determined by the $L_\infty$-algebra structure arising from the higher extension.

\subsection{The higher Kac--Moody factorization algebras}
\label{ss:5.1}
In this subsection, we review the construction of the higher Kac--Moody factorization algebra.
This is introduced in \cite{GW}.

Fix a simple Lie algebra $\frg$, and we call the following element $\theta_c\in(\Sym^{d+1}\frg^\vee)^\frg$ defined by 
\[
\theta_c(x_0,\dots,x_{d})=-\sum_{\alpha\in\Delta_+}\alpha(x_0)\dots\alpha(x_d)
\]
the \emph{critical $(d+1)$-form} of $\frg$.
When $d=1$, we have $\theta_1=\kappa_c=-\frac{1}{2}\kappa_{\frg}$, which coincides with the usual critical level.

By an abuse of notation, for $\theta\in(\Sym^{d+1}\frg^\vee)^\frg$ and an open subset $U\subset\bbC^d$, we define the following cochain on $(\clG_d)_c(U)=\Omega^{0,*}_{\bbC^d,c}(U)\otimes\frg$.
\begin{align*}
  (\phi(\theta))_U(\alpha^{(0)}\otimes x^0,\dots,\alpha^{(d)}\otimes x^{d})=\frac{1}{(2\pi i)^{d+1}}\left(\int_U\alpha^{(0)}\wedge\partial\alpha^{(1)}\wedge\dots\wedge\partial \alpha^{(d)}\right)\theta(x^{0},\dots,x^{d}),
\end{align*}
for $\alpha^{(i)}\otimes x^i\in\Omega^{0,*}_{\bbC^d,c}(U)\otimes\frg,\ (i=0,\dots,d)$.
By Lemma 1.21 of \cite{GW}, the map of sheaves $\phi(\theta)=((\phi(\theta))_U)_{U\in\Open(\bbC^{d})}$ defines a $(-1)$-shifted 2-cocycle of $(\clG_d)_c=\Omega^{0,*}_{\bbC^d,c}\otimes\frg$.

\begin{dfn}[{\cite[Definition 1.24]{GW}}]
  \label{higher_Kac_Moody_FD_def}
  Let $\frg$ be a Lie algebra and let $\theta\in(\Sym^{d+1}\frg^\vee)^\frg$.
  A \emph{higher Kac--Moody factorization algebra} $\bbU_{\phi(\theta)}(\clG_d)$ on $\bbC^d$ is a twisted factorization enveloping
  \[
  \bbU_{\phi(\theta)}(\clG_d)\colon \Open(\bbC^d)\to\mathrm{pCh(DVS)}.
  \]
\end{dfn}

 It is obvious that for a 2-cocycle $\kappa\in(\Sym^2\frg^\vee)^\frg$, the higher Kac--Moody factorization algebra $\bbU_{\phi(\kappa)}(\clG_1)$ on $\bbC$ is nothing but the Kac--Moody factorization algebra (see \cref{def_kac_moody_FA}).

\subsection{Critical level}
\label{ss:5.2}
Throughout the following, we assume that $\frg$ is a simple Lie algebra.

Here, as a two-dimensional analogue, we consider the following classical field theory on $\bbC^d$:
\begin{align}
\label{eq:clMddef}
\begin{split}
    &(\clM^d_{\frg},\omega_{\beta\gamma},Q_{BV}^{\beta\gamma})\\
=&\left(\Omega^{0,*}_{\bbC^d}\otimes\frn_+[-1]\oplus\Omega^{d,*}_{\bbC^d}\otimes\frn_+^\vee[d-2],\omega_{\beta\gamma}=I_{\bbC^d}\circ\ev*{-,-},Q_{BV}^{\beta\gamma}=\bar{\partial}\right)
\end{split}
\end{align}
where $I_{\bbC^d}$ is an integral pairing and $\ev*{-,-}:\frn_+\otimes\frn^\vee_+\to\bbC$ is a natural evaluation pairing.
Recall that  $(S^1)^d$ acts on $\bbC^d$ by rotations in each component, i,e. for $(\theta_1,\dots,\theta_d)\in(S^1)^d$ and $(z_1,\dots,z_d)\in\bbC^d$ the action is given as follows.
\[
(\theta_1,\dots,\theta_d)\cdot(z_1,\dots,z_d)=(e^{-i\theta_1}z_1,\dots,e^{-i\theta_d}z_d).
\]

As in \cref{ss:4.2}, we denote a section of $\mathcal{M}_\frg^d$ by $(\gamma_\alpha,\beta_\alpha)$.
For a polynomial $f\in\bbC[T^*\frn_+]$, we define an element $f(\gamma,\beta)\in\Omega^{*,*}_{\bbC^d}$ by replacing coordinate functions $y_\alpha,\xi_\alpha$ by $\gamma_\alpha,\beta_\alpha$, respectively.

Let $\wt{S}^{hiWak}_{d,\theta_c}\in\IA(\clG_d,\clM^d_{\frg})$ be a functional given as follows:
\begin{align}
  \label{naive_Wakimoto_d_dim}
  \begin{split}
  \wt{S}^{hiWak}_{d,\theta_c}&(\omega_\alpha\otimes e_\alpha,\eta_i\otimes h_i,\xi_\alpha\otimes f_\alpha,\gamma_\alpha,\beta_\alpha)\\
 =\sum_{\alpha\in\Delta_+}&\int_{\bbC^d}\omega_\alpha\wedge\beta_\alpha+\omega_\alpha\wedge\sum_{\beta\in\Delta_+}\beta_\beta\wedge P^\alpha_\beta(\gamma)+\xi_\alpha\wedge\sum_{\beta\in\Delta_+}\beta_\beta\wedge Q^\alpha_\beta(\gamma)-\sum_{i=1}^{\rk\frg}\alpha(h_i)\eta_i\wedge\beta_\alpha\wedge\gamma_\alpha.
  \end{split}
\end{align}
Again, by abuse of notation, we denote by $\wt{S}^{hiWak}_{d,\theta_c}$ the image of the natural quotient map $\IA(\clG_d,\clM_\frg^d)\srj\Act(\clG_d,\clM^d\frg)$.

\begin{lem}
\label{naive_Wakimoto_higher_dim}
The functional $\wt{S}^{hiWak}_{d,\theta_c}$ satisfies the equivariant classical master equation.
\end{lem}
\begin{proof}
  The proof of \cref{naive_equivariant_classical} works with minor modifications.
\end{proof}

\begin{thm}
\label{higher_dim_Wakimoto_critical}
There exists an $(S^1)^d$-invariant functional  $S^{hiWak}_{d,\theta_c}\in\Act(\clG_d,\clM_{\frg}^d)[\hbar]$ such that the collection of $(S^1)^d$-invariant functionals
\[
\left\{S^{hiWak}_{d,\theta_c}[L] \ceq \lim_{\epsilon\to0}W(P(\epsilon,L),S^{hiWak}_{d,\theta_c})\right\}_{L>0}
\]
satisfies the equivariant outer-quantum master equation.
\end{thm}
\begin{proof}
  The proof of \cref{Wakimoto_critical} works with minor modifications.
  
  Similarly to \cref{cohomology_of_red_loc}, we have the following exact sequence
  \begin{align}
    \label{long_exact_sequence_inn_act_higherdim}
  H^1(C^*_{red,loc}(\clG_d))\stackrel{H^1(i)}{\to} H^1(C^*_{red,loc}(\Omega^{0,*}_{\bbC^d}\otimes\frh))\otimes\mathscr{O}_{loc}(\Omega^{d,*}_{\bbC^d}\otimes\frn_+^\vee[-1])\stackrel{H^1(p)}{\to} H^1(\Act(\clG_d,\clM_{\frg}^d)).
  \end{align}
  By \cref{holomorphic_counterterm}, there exists a limit
  \[
  \wt{S}^{hiWak}_{d,\theta_c}[L] \ceq \lim_{\epsilon\to0}W(P(\epsilon,L),\wt{S}^{hiWak}_{d,\theta_c}).
  \]
  Similarly to the proof of \cref{Wakimoto_critical}, we consider the following obstruction classes
  \begin{align*}
    \Theta^{hiWak,crit}_{inn}& \ceq 
    \frac{1}{\hbar}\lim_{L\to0}(d_{\clG_d}\wt{S}^{hiWak}_{d,\theta_c}[L]+\bar{\partial}\wt{S}^{hiWak}_{d,\theta_c}[L]\\
    &+\frac{1}{2}\{\wt{S}^{hiWak}_{d,\theta_c}[L],\wt{S}^{hiWak}_{d,\theta_c}[L]\}_L+\hbar\Delta_L\wt{S}^{hiWak}_{d,\theta_c}[L])\ \bmod \hbar,\\
  \Theta^{hiWak,crit}_{out}&\ceq
  \frac{1}{\hbar}\lim_{L\to0}(d_{\clG_d}\wt{S}^{hiWak}_{d,\theta_c}[L]+\bar{\partial}\wt{S}^{hiWak}_{d,\theta_c}[L]\\
  &+\frac{1}{2}\{\wt{S}^{hiWak}_{d,\theta_c}[L],\wt{S}^{hiWak}_{d,\theta_c}[L]\}_L+\hbar\Delta_L\wt{S}^{hiWak}_{d,\theta_c}[L])|_{\Act(\clG_d,\clM^d_{\frg})}\ \bmod \hbar
  \end{align*}
  The obstruction classes are 1-cocycles in $\IA(\clG_d,\clM^d_{\frg})$ and $\Act(\clG_d,\clM^d_{\frg})$, respectively.\\
  Let $[\Theta^{hiWak,crit}_{inn}],\ [\Theta^{hiWak,crit}_{out}]$ be the classes in the first cohomologies $H^1(\IA(\clG_d,\clM^d_{\frg}))$ and $H^1(\Act(\clG_d,\clM^d_{\frg}))$, respectively.
  Then, similarly to the proof of \cref{Wakimoto_critical}, the cohomology class $[\Theta^{hiWak,crit}_{inn}]$ is determined by the weight of the $(d+1)$-Wheel graph whose external edges are all $\eta\otimes h$-leg. Such a graph is given in \cref{obstruction_class_inn_higher_Wak}, where the functional $I_{\eta\gamma\beta}^{hi}$ for $\eta\in\Omega^{0,*}_{\bbC^d}$ is defined by
  \[
  I^{hi}_{\eta\gamma\beta}:=\sum_{\alpha\in\Delta_+}\int_{\bbC^d}\alpha(h)\beta_\alpha\wedge\eta\wedge\gamma_\alpha.
  \]
  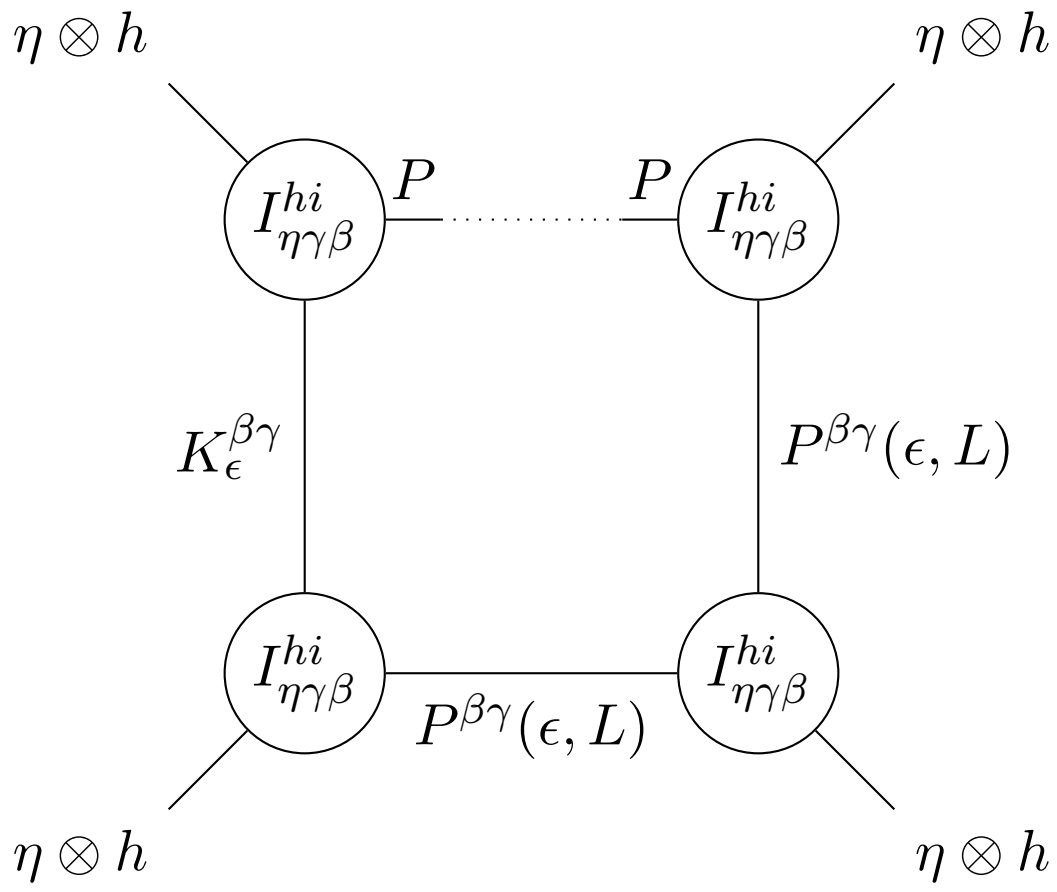
\begin{figure}[htbp]
\begin{tikzpicture}[
  vertex/.style={
    circle,
    draw,
    fill=white,
    inner sep=2pt,
    minimum size=6mm
  }
]

% vertices
\node[vertex] (v1) at (-1.5,1.5) {$I^{hi}_{\eta\gamma\beta}$};
\node[vertex] (v2) at (1.5,1.5) {$I^{hi}_{\eta\gamma\beta}$};
\node[vertex] (v3) at (1.5,-1.5) {$I^{hi}_{\eta\gamma\beta}$};
\node[vertex] (v4) at (-1.5,-1.5) {$I^{hi}_{\eta\gamma\beta}$};

% top edge: solid -- dotted -- solid
\draw (v1) -- node[above] {$P$} (-0.6,1.5);
\draw[dotted] (-0.6,1.5) -- (0.6,1.5);
\draw (0.6,1.5) -- node[above] {$P$} (v2);

% right edge
\draw (v2) -- node[right] {$P^{\beta\gamma}(\epsilon,L)$} (v3);

% bottom edge
\draw (v3) -- node[below] {$P^{\beta\gamma}(\epsilon,L)$} (v4);

% left edge
\draw (v4) -- node[left] {$K^{\beta\gamma}_\epsilon$} (v1);

% external edges
\draw (v1) -- ++(-0.9,0.9) node[above left] {$\eta\otimes h$};
\draw (v2) -- ++(0.9,0.9) node[above right] {$\eta\otimes h$};
\draw (v3) -- ++(0.9,-0.9) node[below right] {$\eta\otimes h$};
\draw (v4) -- ++(-0.9,-0.9) node[below left] {$\eta\otimes h$};

\end{tikzpicture}
\caption{The $(d+1)$-wheel graph which contributes to the obstruction class $[\Theta^{hiWak,crit}_{inn}]$.}
\label{obstruction_class_inn_higher_Wak}
\end{figure}

By a similar computation as in \cref{Wakimoto_critical}, we find that the weight of the graph in \cref{obstruction_class_inn_higher_Wak} is given as follows: (see Appendix B of \cite{GW})
\[
\left(-\sum_{\alpha\in\Delta_+}\alpha(h_i)^{d+1}\right)\frac{1}{(2\pi i)^d}\frac{1}{(d+1)!}\int_{\bbC^d}\eta_i\wedge\partial\eta_i\wedge\dots\wedge\partial\eta_i=\phi(\theta_c^{d+1})(\eta_i\otimes h_i,\dots,\eta_i\otimes h_i).
\]
Hence, we have
\[
[\Theta^{hiWak,crit}_{inn}]=H^1(i)\left[\phi(\theta_c^{d+1})\right].
\]
This implies that there exists a cohomological degree 0 element $\tau_d\in\IA(\clG_d,\clM_{\frg}^d)$ such that
\begin{align}
  \label{gamma_d_exact}
  i(\phi(\theta_c^{d+1}))-\Theta^{hiWak,crit}_{inn}=d_{\clG_d}\tau_d+\bar{\partial}\tau_d+\{\wt{S}^{hiWak}_{d,\theta_c},\tau_d\}.
\end{align}
It is obvious that the functional $\tau_d$ is not contained in the image of $i$.
Thus, under the quotient map, $\tau_d$ defines a nonzero element of $\Act(\clG_d,\clM_{\frg}^d)$.
Since the left hand side of the equality \eqref{gamma_d_exact} is a degree $(d+1)$ functional on $\clG_d$, $\tau_d$ is a degree $d$ functional on $\clG_d$.

Since $\phi$ and $\wt{S}^{hiWak}_{d,\theta_c}$ are $(S^1)^d$-invariant, the left hand side of \eqref{gamma_d_exact} is also $(S^1)^{d}$-invariant, and hence we can choose $\tau_d$ to be $(S^1)^d$-invariant.
 Similarly to the proof of \cref{Wakimoto_critical}, the $(S^1)^d$-action on $\frg$ is given by weight $(-1,\dots,-1)$.

Then, similarly to \cref{Wakimoto_field_theory_crit}, we have the following assertion:
\begin{prop}
  \label{higher_Wakimoto_crit_gamma}
  There exists $c(d)_{\vec{\alpha}}\in\bbC\ (\vec{\alpha}\in\Delta_+^{d})$ such that
  \[
  \tau_d(\omega_\alpha\otimes e_\alpha,\eta_i\otimes h_i,\xi_\alpha\otimes f_\alpha,\gamma_\alpha,\beta_\alpha)=\sum_{\vec{\alpha}\in\Delta_+^d}c(d)_{\vec{\alpha}}\int_{\bbC^d}\xi_{\alpha_1}\wedge\dots\wedge\xi_{\alpha_d}\wedge\partial\gamma_{\alpha_1}\wedge\dots\wedge\partial\gamma_{\alpha_d}
  \]
  satisfies the equality \eqref{gamma_d_exact}.
\end{prop}

Taking $c(d)_{\vec{\alpha}}$ as in \cref{higher_Wakimoto_crit_gamma} and setting 
\[
S^{hiWak}_{d,\theta_c} \ceq \wt{S}^{hiWak}_{d,\theta_c}+\hbar\sum_{\vec{\alpha}\in\Delta_+^d}c(d)_{\vec{\alpha}}\int_{\bbC^d}\xi_{\alpha_1}\wedge\dots\wedge\xi_{\alpha_d}\wedge\partial\gamma_{\alpha_1}\wedge\dots\wedge\partial\gamma_{\alpha_d},
\]
we obtain the described functional.
The proof of \cref{higher_dim_Wakimoto_critical} is now complete.
\end{proof}

Again, by \cref{quantum_Noether_map}, we have the following assertion:
\begin{cor}
\label{Noether_map_higher_Wakimoto_critical_level}
The quantum outre-action $\{S^{hiWak}_{d,\theta_c}[L]\}_{L>0}$ induces the $\bbC\dbr{\hbar}$-linear map
\[
  \wt{\frw}^{fact}_{d,\theta_c}\colon\bbU^\hbar_{\phi(\theta_c)}(\clG_d)[\hbar^{-1}]\to\Obs^q_{\clM_{\frg}^d}[\hbar^{-1}]
\]
of factorization algebras.
\end{cor}

\begin{rem}
\label{d_1_case_Wakimoto}
In the case of $d=1$, the above map $\wt{\frm}^{fact}_{1,\theta_c}$ coincides with the map $\wt{w}^{fact}_{\kappa_c}$ constructed in \cref{Noether_map_Wakimoto_critical_level}.
\end{rem}

%%%%%%%%%%%%%%%%%%%%%%%%%%%%%%%%%%%%%%%%%%%%%%%%%%%%%%%%%%%%%%%%%%%%%%%%%%%%%%%%
\subsection{The finite dimensional setting :revisit}
\label{ss:5.3}

In considering higher dimensional analogues, it is natural to study some sense of  ``higher extensions'' of \eqref{g_to_Vect_U}.
Recall that the Lie algebra map \eqref{extension_g_to_D} is classified by the first cohomology $H^1(\frg,\bbC[U])$.
The natural higher analogue should be a map that is classified by the higher cohomology $H^*(\frg,\bbC[N_+])$.
Here, we consider the extensions as $L_\infty$-algebras.
\begin{prop}
\label{ext_n_Lie}
Let $n\geq0$ be a non-negative integer, and $\mu\in Z^{n+1}(\frg,\bbC[U])$ be an $n+1$-cocycle.
Then the quadruple $(\frg\oplus\frX(U)\oplus\bbC[U][n-1],l_k^\mu)$ is an $L_\infty$-algebra, where
\begin{itemize}
  \item $l_2^\mu$ is defined by the Lie bracket on $\frg$, the $\frg$ action on $\frX(U)$ defined by $\rho_0:\frg\inj\frX(U)$, and the natural action of $\frg\oplus\frX(U)$ on $\bbC[U]$,
  \item $l_{n+1}^\mu$ is defined by $\mu:\frg^{\otimes(n+1)}\to\bbC[U]$,
  \item For $k\neq 2,n+1$, the $l_k^\mu$-bracket is zero.
\end{itemize}
\end{prop}
\begin{proof}
  By Theorem 55 of \cite{JBS}, the subspace $\frg\oplus\bbC[U][n-1]$ is an $L_\infty$-algebra.
  By construction of the action of $\frg$ on $\bbC[U]$, the action factors through the natural action of $\frX(U)$ on $\bbC[U]$, i.e. the diagram commutes
  \[
  \begin{tikzcd}
      \frg\times\bbC[U] \arrow[r] \arrow[d, "\rho_0\times\mathrm{id}"] 
        & \bbC[U] \arrow[d] \\
      \frX(U)\times\bbC[U] \arrow[r]
        & \bbC[U].
  \end{tikzcd}
  \]
  Therefore, the vector space $\frg\oplus\frX(U)\oplus\bbC[U][n-1]$ is also an $L_\infty$-algebra.
\end{proof}
By \cref{lem_CU_cont_Verma} and Shapiro's lemma, we have
\[
  H^{n+1}(\frg,\bbC[U]) \simeq 
  H^{n+1}(\frg,\Hom^{\res}_{U(\frb_-)}(U(\frg),\bbC))\simeq 
  H^{n+1}(\frb_-,\bbC)\simeq 
  H^{n+1}(\frh,\bbC)\simeq\bigwedge{}^{n+1}\frh^\vee.
\]
As in \cref{ss:4.1}, let us consider the inverse map 
\begin{align}\label{eq:5.3:inv}
  \bigwedge{}^{n+1}\frh^\vee \lsto 
  H^{n+1}(\frb_-,\bbC) \lsto 
  H^{n+1}(\frg,\Hom^{\res}_{U(\frb_-)}(U(\frg),\bbC))
  \lsto H^{n+1}(\frg,\bbC[U]).
\end{align}
By the construction of the isomorphisms of \cref{lem_CU_cont_Verma} and \eqref{equ:htoZ1}, the isomorphism \eqref{eq:5.3:inv} can be realized by choosing representatives as
\begin{align}
\nonumber
  &S\colon \bigwedge{}^{n+1}\frh^\vee \lto Z^{n+1}(\frg,\bbC[U]),\\
\label{eq:5.3:S}
  &\begin{cases}
    S(\theta)(h_{i_1}\wedge\dots\wedge h_{i_{n+1}}) \ceq
    \theta(h_{i_1},\dots,h_{i_{n+1}}),\\
    S(\theta)(f_{\alpha_{i_1}}\wedge\dots\wedge f_{\alpha_{i_{n+1}}}) \ceq 
    \theta(h_{i_1},\dots,h_{i_{n+1}})y_{\alpha_{i_1}}\dots y_{\alpha_{i_{n+1}}},\\
    S(\theta)(m) \ceq 0 \quad (\text{$m$: the other basis elements}).
  \end{cases}
\end{align}

Using this observation, in \cref{ss:5.4} following, we consider a higher-dimensional level deformation.

\subsection{General level}
\label{ss:5.4}
In this subsection, we consider the higher Wakimoto realization of the general level.

As a natural higher dimensional analogue of the free degenerate classical field theory $\clM_{Heis}$, we consider the following free degenerate classical field theory on $\bbC^{2n+1}$.
\begin{align}
  \label{higher_Heis}
  \begin{split}
      (\clM^{\hbar}_{n,Heis},\Pi^{\hbar}_{n,Heis},Q_{BV}^{Heis})
      \ceq\left(\Omega^{n+1,*}_{\bbC^{2n+1}}\otimes\frh[n-1],\Pi_{n,Heis},Q_{BV}^{Heis}=\bar{\partial}\otimes\mathrm{id}\right),
  \end{split}
\end{align}
where the $(-1)$-shifted Poisson structure $\Pi^{\hbar}_{n,Heis}$ is given by
\begin{align*}
  &\Pi^{(0)}_{n,Heis}\colon
   \Omega^{n,*}_{\bbC^{2n+1}}\otimes\frh^\vee
   \lto \Omega^{n+1,*}_{\bbC^{2n+1}}\otimes\frh, \quad
   \omega\otimes f \lmto \hbar\partial\omega\otimes f^*,\\
  &\Pi^{(k)}_{n,Heis} \ceq 0 \quad (k \ne 0),
\end{align*}
where $f^*\in\frh$ is the image of $f\in\frh^\vee$ under the isomorphism $\frh^\vee\stackrel{\sim}{\to}\frh$.

Then we have a free degenerate classical field theory $\clW_{2n+1}=\clM_{\frg}^d\oplus\clM^{\hbar}_{n,Heis}$, where $\clM^d_\frg$ is the degenerate classical field theory defined in \eqref{eq:clMddef}.
As in \cref{ss:4.3}, we denote a section of $\clW_{2n+1}$ by $(\gamma_\alpha,\beta_\alpha,b_i)$.

%We generalize the functional \eqref{general_level_C1_Wakimoto} to higher dimensions and consider the following functional:

Recalling the cocycle $S(\theta)$ in \eqref{eq:5.3:S}, 
it is natural to consider the following functional
\begin{align}
  \label{general_level_naive_hiWakimoto}
  \begin{split}
&\wt{S}^{hiWak}_{n,k}(\omega_\alpha\otimes e_\alpha,\eta_i\otimes h_i,\xi_\alpha\otimes f_\alpha,\gamma_\alpha,\beta_\alpha,b_i)\\
 \ceq&\wt{S}^{hiWak}_{2n+1,\theta_c}(\omega_\alpha\otimes e_\alpha,\eta_i\otimes h_i,\xi_\alpha\otimes f_\alpha,\gamma_\alpha,\beta_\alpha)\\
 &+\sqrt{\frac{k+h^\vee}{h^\vee}}\sum_{\alpha\in\Delta_+}\int_{\bbC^{2n+1}}\left(\sum_{i=1}^{\rk\frg}\alpha(h_i)^{n+1}\eta_i\wedge(\partial\eta_{i})^{\wedge n}\wedge b_\alpha\right.\\
 &\left.+\sum_{\beta\in\Delta_+}\alpha([e_\beta,f_\beta])^{n+1}\xi_\beta\wedge(\partial\xi_\beta)^{\wedge n}\wedge b_{\alpha}\wedge\gamma_{\beta}^{\wedge(n+1)}\right)\in\IA(\clG_{2n+1},\clW_{2n+1}),
  \end{split}
\end{align}
where for positive root $\alpha=\sum_ic_i\alpha_i\in\Delta_+$, we denote $b_\alpha=\sum_ic_ib_i$.
Again, by abuse of notation, we denote by $\tilde{S}^{hiWak}_{n,k}$ by the image of the natural quotient map $\IA(\clG_{2n+1},\clW_{2n+1})\srj\Act(\clG_{2n+1},\clW_{2n+1})$.

\begin{lem}
  \label{naive_hiWaki_equivariant}
  The functional $\wt{S}^{hiWak}_{n,k}\in\Act(\clG_{2n+1},\clW_{2n+1})$ satisfies the equivariant classical master equation.
\end{lem}
\begin{proof}
      The proof of \cref{lem_equiv_general_one_dim} works with minor modifications.
\end{proof}

\begin{thm}
\label{general_level_higher_Wakimoto}
There exists a functional $S^{hiWak}_{n,k}\in\Act(\clG_{2n+1},\clW_{2n+1})[\hbar]$
such that the collection of $(S^1)^{2n+1}$-invariant functionals 
\[
\left\{S^{hiWak}_{n,k}[L]=\lim_{\epsilon\to0}W(P(\epsilon,L),S^{hiWak}_{n,k})\right\}_{L>0}
\]
satisfies the equivariant outer-quantum master equation.
\end{thm}
\begin{proof}
    The proof of \cref{Wakimoto_general_level_FA} works with minor modifications.

    Similarly to the proof of \cref{Wakimoto_general_level_FA}, the following obstruction classes
    \begin{align*}
      \Theta^{hiWak,gen}_{inn}&:=\frac{1}{\hbar}\lim_{L\to0}(d_{\clG_{2n+1}}\wt{S}^{hiWak}_{n,k}[L]+\bar{\partial}\wt{S}^{hiWak}_{n,k}[L]\\
      &+\frac{1}{2}\{\wt{S}^{hiWak}_{n,k}[L],\wt{S}^{hiWak}_{n,k}[L]\}_L+\hbar\Delta_L\wt{S}^{hiWak}_{n,k}[L])\ \bmod \hbar,\\
  \Theta^{hiWak,gen}_{out}:=&\frac{1}{\hbar}\lim_{L\to0}(d_{\clG_{2n+1}}\wt{S}^{hiWak}_{n,k}[L]+\bar{\partial}\wt{S}^{hiWak}_{n,k}[L]\\
  &\left.+\frac{1}{2}\{\wt{S}^{hiWak}_{n,k}[L],\wt{S}^{hiWak}_{n,k}[L]\}_L+\hbar\Delta_L\wt{S}^{hiWak}_{n,k}[L]\right)|_{\Act(\clG_{2n+1},\clW_{2n+1})}\ \bmod \hbar.
    \end{align*}
    Similarly to the proof of \cref{Wakimoto_general_level_FA}, the cohomology class $[\Theta^{hiWak,gen}_{inn}]$ belongs to the subspace $ H^1(\IA(\clG_{2n+1},\clM_{\frg}^{2n+1}))\subset H^1(\IA(\mathcal{G}_{2n+1},\clW_{2n+1}))$ and is determined by the weight of the graph whose external edges are all $\eta\otimes h$-legs.
 Again, such graphs are given in \cref{obstruction_class_inn_gene_higher_dim}, where the functional $I_{\eta b}$ for $\eta\in\Omega^{0,*}_{\bbC^{2n+1}}$ is defined by
\[
I^{hi}_{\eta b}:=\sqrt{\frac{k+h^\vee}{h^\vee}}\sum_{\alpha\in\Delta_+}\int_{\bbC^{2n+1}}\alpha(h)^{n+1}\eta\wedge(\partial\eta)^{\wedge n}\wedge b_\alpha.
\]
    \begin{figure}[htbp]
% --- 1つ目の図 ---

\begin{tikzpicture}[
  vertex/.style={
    circle,
    draw,
    fill=white,
    inner sep=2pt,
    minimum size=6mm
  }
]

% vertices
\node[vertex] (v1) at (-1.5,1.5) {$I^{hi}_{\eta\gamma\beta}$};
\node[vertex] (v2) at (1.5,1.5) {$I^{hi}_{\eta\gamma\beta}$};
\node[vertex] (v3) at (1.5,-1.5) {$I^{hi}_{\eta\gamma\beta}$};
\node[vertex] (v4) at (-1.5,-1.5) {$I^{hi}_{\eta\gamma\beta}$};

% top edge: solid -- dotted -- solid
\draw (v1) -- node[above] {$P$} (-0.6,1.5);
\draw[dotted] (-0.6,1.5) -- (0.6,1.5);
\draw (0.6,1.5) -- node[above] {$P$} (v2);

% right edge
\draw (v2) -- node[right] {$P^{\beta\gamma}(\epsilon,L)$} (v3);

% bottom edge
\draw (v3) -- node[below] {$P^{\beta\gamma}(\epsilon,L)$} (v4);

% left edge
\draw (v4) -- node[left] {$K^{\beta\gamma}_\epsilon$} (v1);

% external edges
\draw (v1) -- ++(-0.9,0.9) node[above left] {$\eta\otimes h$};
\draw (v2) -- ++(0.9,0.9) node[above right] {$\eta\otimes h$};
\draw (v3) -- ++(0.9,-0.9) node[below right] {$\eta\otimes h$};
\draw (v4) -- ++(-0.9,-0.9) node[below left] {$\eta\otimes h$};

\end{tikzpicture}

% --- 2つ目の図 ---
\begin{tikzpicture}[node distance=3cm, thick]

% ノード（丸の中に I）
\node[draw, circle, minimum size=12mm] (A) {$I^{hi}_{\eta b}$};
\node[draw, circle, minimum size=12mm, right of=A] (B) {$I^{hi}_{\eta b}$};

% 外側の線とラベル
\draw (A) -- ++(-1.5,0.5)  node[pos=1, left=3pt] {$\eta\otimes h$};
\draw (A) -- ++(-1.5,-0.5)  node[pos=1, left=3pt] {$\eta\otimes h$};
\draw (A) -- ++(-1.5,-0.5)  node[pos=1, above right=5pt] {$\vdots$};
\draw (B) -- ++(1.5,0.5)  node[pos=1, right=3pt] {$\eta\otimes h$};
\draw (B) -- ++(1.5,-0.5)  node[pos=1, right=3pt] {$\eta\otimes h$};
\draw (B) -- ++(1.5,-0.5)  node[pos=1, above left=5pt] {$\vdots$};

% 内側の 2 本線
\draw (A) to node[midway, above=3pt] {$K^{n,Heis}_\epsilon$} (B);

\end{tikzpicture}

\caption{The graph which contributes to the obstruction class $\Theta^{hiWak,gen}_{inn}$.}
\label{obstruction_class_inn_gene_higher_dim}
\end{figure}

The weight of the first graph in \cref{obstruction_class_inn_gene} is already computed in the proof of \cref{higher_dim_Wakimoto_critical}, and is given by $i(\phi(\theta_c))$.

Next, we consider the weight $\Theta_{2nd}$ of the second graph in \cref{obstruction_class_inn_gene}.
Similarly to the above cases, this decomposes into an analytic part and an algebraic part:
\[
\Theta_{2nd}(\eta_i\otimes h_i)=\Theta_{2nd}^{an}(\eta_i)\Theta^{alg}_{2nd}(h_i).
\]
Then we have
\begin{align*}
  \Theta^{an}_{2nd}(f)
=&\lim_{L\to0}\lim_{\epsilon\to0}\int_{\bbC^{2n+1}_{z_1}\times\bbC^{2n+1}_{z_2}}f(z_1)f(z_2)(\partial f(z_1))^{\wedge n}\wedge(\partial f(z_2))^{\wedge n}\wedge K^{n,Heis}_\epsilon(z_1,z_2)\\
                      %=\lim_{L\to0}\lim_{\epsilon\to0}&\int_{\bbC^{2n+1}_{z_1}\times\bbC^{2n+1}_{z_2}}f(z_1)f(z_2)(\partial f(z_1))^{\wedge n}\wedge(\partial f(z_2))^{\wedge n}\wedge\frac{1}{(2\pi i\epsilon)^{2n+1}}\partial e^{-|z_1-z_2|^2/4\epsilon}\\
                      %&\wedge\bigwedge{}_{i=1}^{2n+1}(d\bar{z}_1^{(i)}-d\bar{z}_2^{(i)})\\
                      %=\lim_{L\to0}\lim_{\epsilon\to0}&\int_{\bbC^{2n+1}_{z_1}\times\bbC^{2n+1}_{z_2}}\partial(f(z_1)f(z_2))\wedge(\partial f(z_1))^{\wedge n}\wedge(\partial f(z_2))^{\wedge n}\wedge\frac{1}{(2\pi i\epsilon)^{2n+1}} e^{-|z_1-z_2|^2/4\epsilon}\\
                      %&\wedge\bigwedge{}_{i=1}^{2n+1}(d\bar{z}_1^{(i)}-d\bar{z}_2^{(i)})\\
                      =&\frac{1}{(2\pi i)^{2n+1}}\frac{1}{(2n+1)!}\int_{\bbC_w^{2n+1}}f(w)(\partial f(w))^{\wedge 2n+1}d^{2n+1}\bar{w}.
\end{align*}

The algebraic part is given as follows
\begin{align*}
  \Theta^{alg}_{2nd}(h_i)=\frac{k+h^\vee}{h^\vee}\sum_{\alpha\in\Delta_+}\alpha(h_i)^{2n+2}=-\left(\frac{k}{h^\vee}+1\right)\theta_c(h_i,\dots,h_i).
\end{align*}
Therefore, we have
\begin{align*}
  [\Theta^{hiWak,gen}_{inn}]=H^1(i)\left([\phi(\theta_c)]-\left(\frac{k}{h^\vee}+1\right)[\phi(\theta_c)]\right)=H^1(i)\left(-\frac{k}{h^\vee}[\phi(\theta_c)]\right).
\end{align*}
This implies that there exists a cohomological degree 0 element $\tau_{2n+1}'\in\IA(\clG_d,\clM_{\frg}^{2n+1})$ such that
\begin{align}
  \label{higher_Wakimoto_gamma_d_exact}
  -\frac{k}{h^\vee}i(\phi(\theta_c))-\Theta^{hiWak,gen}_{inn}=d_{\clG_{2n+1}}\tau_{2n+1}'+\bar{\partial}\tau_{2n+1}'+\{\wt{S}^{hiWak}_{n,k},\tau_{2n+1}'\}.
\end{align}

Noting that the functional $\wt{S}^{hiWak}_{n,k}$ contains the degree $n+1$ part of $\clG_{2n+1}$, similarly to \cref{higher_Wakimoto_crit_gamma}, there exists $c'(2n+1)_{\vec{\alpha}}\in\bbC,\ (\vec{\alpha}\in\Delta_+^{2n+1})$ and $d(2n+1)_{\vec{\beta}},\ (\beta\in\Delta_+^{n+1})$ such that
\begin{align*}
&\tau_{2n+1}'(\omega_\alpha\otimes e_\alpha,\eta_i\otimes h_i,\xi_\alpha\otimes f_\alpha,\gamma_\alpha,\beta_\alpha)\\
=&\sum_{\vec{\alpha}\in\Delta_+^{2n+1}}c'(2n+1)_{\vec{\alpha}}\int_{\bbC^{2n+1}}\xi_{\alpha_1}\wedge\dots\wedge\xi_{\alpha_{2n+1}}\wedge\partial\gamma_{\alpha_1}\wedge\dots\wedge\partial\gamma_{\alpha_{2n+1}}\\
  +&\sum_{\vec{\beta}\in\Delta_+^{n+1}}d(2n+1)_{\vec{\beta}}\int_{\bbC^{2n+1}}\xi_{\beta_1}\wedge\partial\xi_{\beta_2}\wedge\dots\wedge\partial\xi_{\beta_{n+1}}\wedge\partial\gamma_{\beta_1}\wedge\dots\wedge\partial\gamma_{\beta_{n+1}}.
\end{align*}

Similarly to the proof of \cref{Wakimoto_general_level_FA}, the coefficient $d(2n+1)_{\vec{\beta}}$ can be determined from the condition that the weights of the graphs in \cref{obstruction_class_d_alpha_determine_higher_Wakimoto} cancel, where the functionals $I_{\xi\partial\gamma}$ and $I_{\xi b\gamma}$ for $\xi_{\beta_i}\in\Omega^{0,*}_{\bbC^{2n+1}}$ are defined by
\begin{align*}
  I_{\xi\partial\gamma}&=\hbar d(2n+1)_{\vec{\beta}}\int_{\bbC^{2n+1}}\xi_{\beta_1}\wedge\partial\xi_{\beta_2}\wedge\dots\wedge\partial\xi_{\beta_{n+1}}\wedge\partial\gamma_{\beta_1}\wedge\dots\wedge\partial\gamma_{\beta_{n+1}},\\
  I_{\xi b\gamma}&=\sqrt{\frac{k+h^\vee}{h^\vee}}\sum_{\alpha\in\Delta_+}\alpha([e_{\alpha_i},f_{\alpha_i}])^{n+1}\int_{\bbC^{2n+1}}\xi_{\alpha_i}\wedge (\partial\xi_{\alpha_i})^{\wedge n}\wedge b_\alpha\wedge\gamma_{\alpha_i}^{\wedge(n+1)}.
\end{align*}
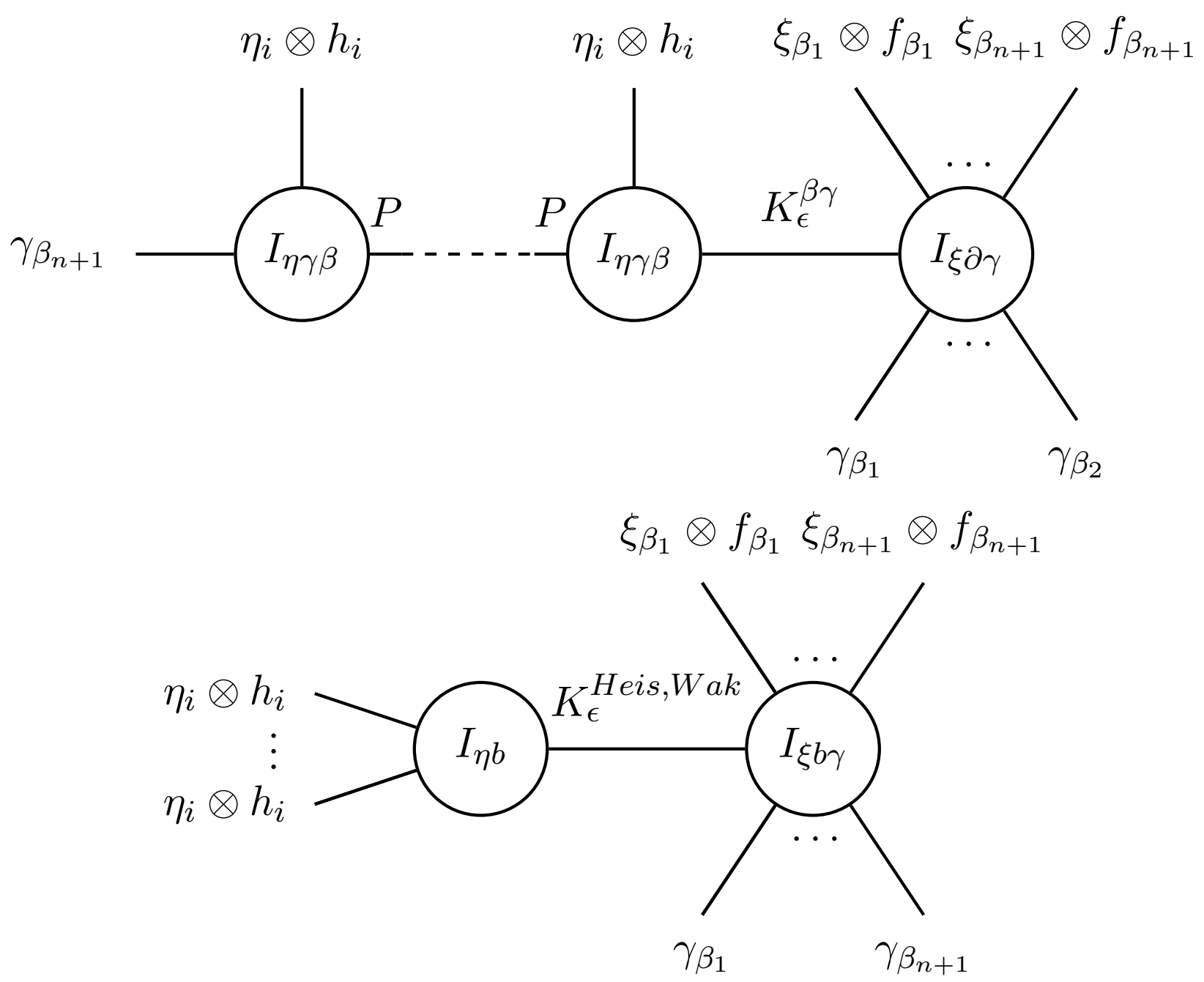
\begin{figure}[htbp]
% --- 1つ目の図 ---
\begin{tikzpicture}[node distance=3cm, thick]

% ノード（丸の中に I）
\node[draw, circle, minimum size=12mm] (A) {$I_{\eta\gamma\beta}$};
\node[draw, circle, minimum size=12mm, right of=A] (B) {$I_{\xi\partial\gamma}$};
\node[draw, circle, minimum size=12mm, left of=A] (C) {$I_{\eta\gamma\beta}$};

% 外側の線とラベル
\draw (A) -- ++(0,1.5) node[pos=1, above=3pt] {$\eta_{i}\otimes h_{i}$};
\draw (B) -- ++(-1.0,1.5)  node[pos=1, above=3pt] {$\xi_{\beta_1}\otimes f_{\beta_1}$};
\draw (B) -- ++(-1.0,1.5)  node[pos=1, below right=15pt] {$\ \dots$};
\draw (B) -- ++(1.0,1.5)  node[pos=1, above=3pt] {$\xi_{\beta_{n+1}}\otimes f_{\beta_{n+1}}$};
\draw (B) -- ++(-1.0,-1.5)  node[pos=1, below=3pt] {$\gamma_{\beta_1}$};
\draw (B) -- ++(-1.0,-1.5)  node[pos=1, above right=15pt] {$\ \dots$};
\draw (B) -- ++(1.0,-1.5)  node[pos=1, below=3pt] {$\gamma_{\beta_2}$};
\draw (C) -- ++(0,1.5) node[pos=1, above=3pt] {$\eta_i\otimes h_i$};
\draw (C) -- ++(-1.5,0.0) node[pos=1, left=3pt] {$\gamma_{\beta_{n+1}}$};

% 内側の 2 本線
\draw (C) -- node[midway, above=3pt] {$P$} ($(C)!0.3!(A)$);
\draw[dashed] ($(C)!0.3!(A)$) -- ($(C)!0.7!(A)$);
\draw ($(C)!0.7!(A)$) -- node[midway, above=3pt] {$P$} (A);
\draw (A) to node[midway, above=3pt] {$K_\epsilon^{\beta\gamma}$} (B);

\end{tikzpicture}

% --- 2つ目の図 ---
\begin{tikzpicture}[node distance=3cm, thick]

% ノード（丸の中に I）
\node[draw, circle, minimum size=12mm] (A) {$I_{\eta b}$};
\node[draw, circle, minimum size=12mm, right of=A] (B) {$I_{\xi b\gamma}$};

% 外側の線とラベル
\draw (A) -- ++(-1.5,0.5) node[pos=1, left=3pt] {$\eta_i\otimes h_i$};
\draw (A) -- ++(-1.5,-0.5) node[pos=1, left=3pt] {$\eta_i\otimes h_i$};
\draw (A) -- ++(-1.5,-0.5) node[pos=1, above left=5pt] {$\vdots$};
\draw (B) -- ++(-1.0,1.5)  node[pos=1, above=3pt] {$\xi_{\beta_1}\otimes f_{\beta_1}$};
\draw (B) -- ++(-1.0,1.5)  node[pos=1, below right=15pt] {$\ \dots$};
\draw (B) -- ++(1.0,1.5)  node[pos=1, above=3pt] {$\xi_{\beta_{n+1}}\otimes f_{\beta_{n+1}}$};
\draw (B) -- ++(-1.0,-1.5)  node[pos=1, below=3pt] {$\gamma_{\beta_1}$};
\draw (B) -- ++(-1.0,-1.5)  node[pos=1, above right=15pt] {$\ \dots$};
\draw (B) -- ++(1.0,-1.5)  node[pos=1, below=3pt] {$\gamma_{\beta_{n+1}}$};

% 内側の 2 本線
\draw (A) to node[midway, above=3pt] {$K^{Heis,Wak}_\epsilon$} (B);

\end{tikzpicture}

\caption{The graphs with $3n+2$ external edges labeled by $(n+1)$-$\eta\otimes h$, $(n+1)$-$\xi\otimes f$, and $n$-$\gamma$.}

\label{obstruction_class_d_alpha_determine_higher_Wakimoto}
\end{figure}

By performing integration by parts, we find that the analytic parts of the weights of the two graphs in \cref{obstruction_class_d_alpha_determine_higher_Wakimoto} coincide.
The algebraic part of the weight of the first graph in \cref{obstruction_class_d_alpha_determine_higher_Wakimoto} is given as follows
\begin{align*}
  -d(2n+1)_{\vec{\beta}}\sum_{j=1}^{n+1}\beta_j(h_i)^{n+1}.
\end{align*}
The algebraic part of the weight of the second graph in \cref{obstruction_class_d_alpha_determine_higher_Wakimoto} is given as follows
\begin{align*}
  \frac{k+h^\vee}{h^\vee}\sum_{\alpha\in\Delta_+}\alpha(h_i)^{n+1}\alpha([e_{\beta_1},f_{\beta_1}])\dots\alpha([e_{\beta_{n+1}},f_{\beta_{n+1}}]).
\end{align*}
Therefore, we have
\begin{align}
    \label{eq:d2n+1explicit}
    d(2n+1)_{\vec{\beta}}=\left(\sum_{j=1}^{n+1}\beta_j(h_i)^{n+1}\right)^{-1}\frac{k+h^\vee}{h^\vee}\sum_{\alpha\in\Delta_+}\alpha(h_i)^{n+1}\alpha([e_{\beta_1},f_{\beta_1}])\dots\alpha([e_{\beta_{n+1}},f_{\beta_{n+1}}]).
\end{align}

Therefore, we obtain the following.
\begin{prop}
\label{gamma_higher_Wakimoto_katati}
There exist $c(2n+1)_{\vec{\alpha}}'\in\bbC$ such that the functional
\begin{align*}
  &\tau_{2n+1}'(\omega_\alpha\otimes e_\alpha,\eta_i\otimes h_i,\xi_\alpha\otimes f_\alpha,\gamma_\alpha,\beta_\alpha)\\
=&\sum_{\vec{\alpha}\in\Delta_+^{2n+1}}c(2n+1)'_{\vec{\alpha}}\int_{\bbC^{2n+1}}\xi_{\alpha_1}\wedge\dots\wedge\xi_{\alpha_{2n+1}}\wedge\partial\gamma_{\alpha_1}\wedge\dots\wedge\partial\gamma_{\alpha_{2n+1}}\\  +&\sum_{\vec{\beta}\in\Delta_+^{n+1}}d(2n+1)_{\vec{\beta}}\int_{\bbC^{2n+1}}\xi_{\beta_1}\wedge\partial\xi_{\beta_2}\wedge\dots\wedge\partial\xi_{\beta_{n+1}}\wedge\partial\gamma_{\beta_1}\wedge\dots\wedge\partial\gamma_{\beta_{n+1}}
\end{align*}
  satisfies the equality \eqref{gen_gamma_exact}.
\end{prop}
Taking $c(2n+1)_{\vec{\alpha}}'$ as in \cref{gamma_higher_Wakimoto_katati} and setting
\begin{align*}
S^{hiWak}_{n,k} \ceq &\wt{S}^{hiWak}_{n,k}+\hbar\left(\sum_{\vec{\alpha}\in\Delta_+^{2n+1}}c(2n+1)'_{\vec{\alpha}}\int_{\bbC^{2n+1}}\xi_{\alpha_1}\wedge\dots\wedge\xi_{\alpha_{2n+1}}\wedge\partial\gamma_{\alpha_1}\wedge\dots\wedge\partial\gamma_{\alpha_{2n+1}}\right.\\  +&\left.\sum_{\vec{\beta}\in\Delta_+^{n+1}}d(2n+1)_{\vec{\beta}}\int_{\bbC^{2n+1}}\xi_{\beta_1}\wedge\partial\xi_{\beta_2}\wedge\dots\wedge\partial\xi_{\beta_{n+1}}\wedge\partial\gamma_{\beta_1}\wedge\dots\wedge\partial\gamma_{\beta_{n+1}}\right),    
\end{align*}
where $d(2n+1)_{\vec{\beta}}\in\bbC$ is the constant defined in \eqref{eq:d2n+1explicit}, 
we obtain the described functional.
The proof of \cref{general_level_higher_Wakimoto} is now complete.
\end{proof}

Again by \cref{quantum_Noether_map,general_level_higher_Wakimoto}, we have the following assertion:
\begin{cor}
  \label{Noether_map_higher_Wakimoto_general_level}
  The quantum outer-action $\{S^{hiWak}_{n,k}[L]\}_{L>0}$ induces the $\bbC\dbr{\hbar}$-linear map
  \[
  \frw^{fact}_{n,k}\colon\bbU^\hbar_{\phi(-k\theta_c/h^\vee)}(\clG_{2n+1})[\hbar^{-1}]\lto\Obs^q_{\clW_{2n+1}}[\hbar^{-1}]
  \]
  of factorization algebras.
\end{cor}

\begin{rem}
\label{d_1_case_general_Wakimoto}
We denote by $(-,-)\colon\frg\otimes\frg\to\bbC$ the normalized form on $\frg$.
For $k\neq -h^\vee$, under the isomorphism
\begin{align*}
    \clM^{\hbar(k+h^\vee)(-,-)}_{Heis}&\lto\clM^\hbar_{0,Heis}\\
    b_i&\mapsto\sqrt{\frac{k+h^\vee}{h^\vee}}\frac{1}{h^\vee}\sum_{\alpha\in\Delta_+}\alpha(h_i)b_\alpha,
\end{align*}
the above map $\frm^{fact}_{0,k}$ is identified with the map $w^{fact}_{k(-,-)}$ constructed in \cref{Noether_map_Wakimoto_general_level}.
\end{rem}

\begin{rem}
    \label{rem_even_dim}
    Since the natural finite-dimensional $L_\infty$-algebra setting is exhausted by \cref{ss:5.3}, it is not possible to construct a Wakimoto-type realization in even complex dimensions as a natural analogue of the finite-dimensional setting.
    
\end{rem}

\subsection{Kernel of the higher Wakimoto realization}
\label{ss:5.5}
In this subsection, we consider the injectivity of the higher Wakimoto realization given by \cref{Noether_map_higher_Wakimoto_critical_level,Noether_map_higher_Wakimoto_general_level}.
In this subsection, let $d\geq2$, and $\frg$ be a simple Lie algebra.
By \cref{not_inj_critical_level}, at the critical level, the induced map of factorization algebras is not injective. Moreover, its kernel contains elements that define elements of the Feigin–Frenkel center.
A notable example is the Sugawara vector (see Section 18.4.3 of \cite{FBZ}).
It is natural to expect that the higher-dimensional analogue of the Sugawara vector lies in the kernel of the induced map $H^*(\wt{\frw}^{fact}_{d,\theta_c})$.
In fact, this expectation is correct (see \cref{kernel_of_higher_Wakimoto_critical}).
Moreover, for $n\geq2$, we can show that the higher-dimensional analogue lies in the kernel of the induced map $H^*(\wt{\frw}^{fact}_{n,k})$ (see \cref{not_inj_ngeq2}).

To consider the injectivity of the maps $H^*(\wt{\frm}^{fact}_{d,\theta_c})$ and $H^*(\frw^{fact}_{n,k})$, it is useful to consider the BD centers of the BD factorization algebras $\Obs^q_{\clM_{\frg}^d}$ and $\Obs^q_{\clW_{2n+1}}$.

\begin{prop}
  \label{BD_center_of_higher_beta_gamma}
  There exists an isomorphism
  \[
  \clZ^{BD}(\Obs^q_{\clM_{\frg}})\simeq\bbC\dbr{\hbar}
  \]
  of factorization algebras.
\end{prop}
\begin{proof}
  The proof of \cref{prop_commutant_betagamma} works with minor modifications.
\end{proof}
%\begin{prop}
%  For $n\geq0$, there exists an isomorphism
%  \[
%  \clZ^{BD}\left(\Obs^q_{\clM^\hbar_{n,Heis}}\right)\simeq\bbC\dbr{\hbar}
%  \]
%  of factorization algebras on $\bbC^{2n+1}$.
%\end{prop}
%\begin{proof}
%    The proof of \cref{prop_commutant_Heis} works with minor modifications.
%\end{proof}

\begin{prop}
  \label{kernel_of_higher_Wakimoto_critical}
  The induced map
  \[
  H^*(\wt{\frw}^{fact}_{d,\theta_c})\colon H^*\left(\bbU^\hbar_{\phi(\theta_c)}(\clG_d)\right)[\hbar^{-1}]\lto H^*\left(\Obs^q_{\clM_{\frg}^d}\right)[\hbar^{-1}]
  \]
  is not injective.
\end{prop}
\begin{proof}

Here we denote
\[
\frB\ceq\{D(z_1,r_1)\sqcup\dots\sqcup D(z_l,r_l)|l\in\bbZ_{\geq1},\ z_i\in\bbC^d,\ D(z_i,r_i)\cap D(z_j,r_j)=\emptyset,\ (i\neq j)\},
\]
then $\frB$ is an open basis of $\bbC^d$ such that for each open subset $U\subset\bbC^d$, the family
\[
\frB(U)\ceq\{D(z_1,r_1)\sqcup\dots\sqcup D(z_l,r_l)\in\frB|l\in\bbZ_{\geq1},\ D(z_i,r_i)\subset U,\ (i=1,\dots,l)\}
\]
is a Weiss cover of $U$. Thus, to show that the induced map $H^*(\wt{\frw}^{fact}_{d,\theta_c})$ is not injective, it is enough to show that for every $z\in\bbC^d,\ R>0$, the induced map $H^*(\wt{\frw}^{fact}_{d,\theta_c})(D(z,R))$ is not injective.

Take a compactly supported smooth function $f\in C^\infty_c(\bbR_{\geq0})$ on $\bbR_{\geq0}$ such that $f(0)=0$ and $\int_{\bbR_{\geq0}}r^{2d-1}f(r)fr=(d-1)!/2\pi^d$.
Then for $n_1,\dots,n_d\in\bbZ_{\geq0}$ and $x\in\frg$, we define a compactly supported Dolbeault form $x_{n_1,\dots,n_d}(z)$ whose cohomological degree is $1-d$ by
\[
(x_{n_1,\dots,n_d}(z))(w) \ceq f(|w-z|)(w_1-z_1)^{-n_1}\dots (w_d-z_d)^{-n_d}d(\bar{w}_1-\bar{z}_1)\wedge\dots\wedge d(\bar{w}_d-\bar{z}_d)\otimes x.
\]
Then the compactly supported Dolbeault form $\{x_{n_1,\dots,n_d}(z)\}_{n_1,\dots,n_d\geq0,\ x\in\frg}$ generates the cohomology $H^*(\bbU^\hbar_{\phi(\theta_c)}(\clG_d)(D(z,R)))$.
Also, for $n_1,\dots,n_d\geq0,\ m_1,\dots,m_d\geq1,\ \alpha\in\Delta_+$, we define observables
\begin{align*}
  c_{(n_1,\dots,n_d)}^\alpha(z):\Omega^{0,0}_{\bbC^d}(D(z,R))&\otimes\frn_+\ni\gamma\otimes x\mapsto\frac{1}{(n_1)!\dots(n_d)!}(\partial^{n_1}_{z_1}\dots\partial^{n_d}_{z_d}\gamma)(z)e_\alpha^*(x),\\
  b_{(m_1,\dots,m_d)}^\alpha(z):\Omega^{d,0}_{\bbC^d}(D(z,R))&\otimes\frn_+^\vee\ni\beta dz_1\wedge\dots\wedge dz_d\otimes\varphi\\
  &\ \mapsto\frac{1}{(m_1-1)!\dots(m_d-1)!}(\partial^{m_1-1}_{z_1}\dots\partial^{m_d-1}_{z_d}\beta)(z)\varphi(e_\alpha).
\end{align*}
The cohomological degrees of the observables $c^\alpha_{n_1,\dots,n_d}(z)$ and $b^\alpha_{m_1,\dots,m_d}(z)$ are $0$ and $1-d$, respectively.

For $\alpha\in\Delta_+$, we have
\begin{align*}
  &\hbar\lim_{L\to0}\wt{\frw}^{fact}_{\theta_c}((e_\alpha)_{(0,\dots,0)}(z))[L](f_\alpha,g_\alpha dz_1\wedge\dots\wedge dz_d)\\
 =&\int_{w\in\bbC^d}\left(f(|w-z|)g_\alpha(w)+f(|w-z|)\sum_{\beta\in\Delta_+}g_\beta(w)P^\alpha_\beta(f_\alpha(w))\right)d^dw\\
 =&\left(g_{\alpha}(z)+\sum_{\beta\in\Delta_+}g_\beta(z)P^\alpha_\beta(f_\alpha(z))\right)\int_{S^{2d-1}}\mathrm{vol}_{S^{2d-1}}\int_0^\infty r^{2d-1}f(r)dr\\
 =&g_{\alpha}(z)+\sum_{\beta\in\Delta_+}g_\beta(z)P^\alpha_\beta(f_\alpha(z)),
\end{align*}
for holomorphic functions $f_\beta,g_\beta,\ (\beta\in\Delta_+)$ on $\bbC^d$.
Thus, we have
\begin{align*}
  \hbar\lim_{L\to0}\wt{\frw}^{fact}_{d,\theta_c}((e_\alpha)_{(0,\dots,0)}(z))[L]=b_{(1,\dots,1)}^{\alpha}(z)+\sum_{\beta\in\Delta_+}P^\alpha_\beta(c_{(0,\dots,0)}^\alpha(z))b_{(1,\dots,1)}^\beta(z).
\end{align*}
Similarly, we have
\begin{align*}
  \hbar\lim_{L\to0}\wt{\frw}^{fact}_{d,\theta_c}((h_i)_{(0,\dots,0)}(z))[L]=-\sum_{\beta\in\Delta_+}\beta(h_i)c_{(0,\dots,0)}^\beta(z)b_{(1,\dots,1)}^\beta(z),\\
  \hbar\lim_{L\to0}\wt{\frw}^{fact}_{d,\theta_c}((f_\alpha)_{(0,\dots,0)}(z))[L]=\sum_{\beta\in\Delta_+}Q^\alpha_\beta(c_{(0,\dots,0)}^\alpha(z))b_{(1,\dots,1)}^\beta(z).
\end{align*}
Then we define the following observable
\begin{align*}
  \omega_{(0,\dots,0)}^{hiSug}(z) \ceq &\sum_{a=1}^{\dim\frg}(J_a)_{(0,\dots,0)}(z)(J^a)_{(0,\dots,0)}(z)\\
 =&\sum_{\alpha\in\Delta_+}\left((e_\alpha)_{(0,\dots,0)}(z)(f_\alpha)_{(0,\dots,0)}(z)+(f_{\alpha})_{(0,\dots,0)}(z)(e_\alpha)_{(0,\dots,0)}(z)\right)\\
 &+\sum_{i=1}^{\rk\frg}(h_i)_{(0,\dots,0)}(z)(h_i)_{(0,\dots,0)}(z),
\end{align*}
where $\{J^a\}_{a=1,\dots,\dim\frg}$ is a basis of $\frg$, and $\{J_a\}_{a=1,\dots,\dim\frg}$ be a dual basis of $\frg$ with respect to the normalized form.
Since the observable $\omega^{hiSug}_{(0,\dots,0)}(z)$ contains no observable $(f_\alpha)_{(0,\dots,0)}(z)$ associated with terms of degree $2$ or higher in the Chevalley basis $f_\alpha$, we have
\begin{align*}
&\lim_{L\to0}\wt{\frw}^{fact}_{d,\theta_c}(\omega^{hiSug}_{(0,\dots,0)}(z))[L]\\
=&\sum_{a=1}^{\dim\frg}\left(\lim_{L\to0}\wt{\frw}^{fact}_{d,\theta_c}((J_a)_{(0,\dots,0)}(z))[L]\right)\left(\lim_{L\to0}\wt{\frw}^{fact}_{d,\theta_c}((J_a)_{(0,\dots,0)}(z))[L]\right).
\end{align*}
Then the observable $\wt{\frw}^{fact}_{d,\theta_c}(\omega^{hiSug}_{(0,\dots,0)}(z))$ contained in the BD center $\clZ^{BD}(\Obs^q_{\clM_{\frg}}(D(z,R)))\simeq\bbC\dbr{\hbar}$.
Therefore, the element $\omega^{hiSug}_{(0,\dots,0)}(z)$ belongs to the kernel of the map $H^*(\wt{\frw}^{fact}_{d,\theta_c})(D(z,R))$, and hence the map $H^*(\wt{\frw}^{fact}_{d,\theta_c})(D(z,R))$ is not injective.

\end{proof}

Furthermore, $\omega^{hiSug}_{(0,\dots,0)}$ also defines an element of the kernel at a general level.

\begin{prop}
  \label{not_inj_ngeq2}
  For $n\geq2$ and $k\in\bbC$, the induced map
  \[
  H^*(\frw^{fact}_{n,k})\colon H^*(\bbU^\hbar_{\phi(-k\theta_c/h^\vee)}(\clG_{2n+1}))[\hbar^{-1}]\lto H^*(\Obs^q_{\clW_{2n+1}})[\hbar^{-1}]
  \]
  is not injective.
\end{prop}

\begin{proof}
  Similarly to the proof of \cref{kernel_of_higher_Wakimoto_critical}, it is enough to show that for $z\in\bbC^d,R>0$, the induced map $H^*(\frw^{fact}_{n,k})(D(z,R))$ is not injective.
  By construction, the functional $S^{hiWak}_{n,k}-S^{hiWak}_{2n+1,\theta_c}$ is a functional of degree at least $3$ in $\Omega^{0,*}_{\bbC^{2n+1},c}\otimes\frh\subset(\clG_{2n+1})_c$, degree $0$ in $\Omega^{0,*}_{\bbC^{2n+1},c}\otimes\frn_+\subset(\clG_{2n+1})_c$, and degree at least $3$ in $\Omega^{0,*}_{\bbC^{2n+1},c}\otimes\frn_-$.
  Thus, we have
  \[
  \lim_{L\to0}\wt{\frw}^{fact}_{n,k}(\omega^{hiSug}_{(0,\dots,0)}(z))[L]=\lim_{L\to0}\wt{\frw}^{fact}_{2n+1,\theta_c}(\omega^{hiSug}_{(0,\dots,0)}(z))[L].
  \]
  Hence, the image of $\omega^{hiSug}_{(0,\dots,0)}(z)$ is contained in the BD center
  \[
  \clZ^{BD}(\Obs^q_{\clM_{\frg}^{2n+1}}(D(z,R)))\subset\Obs^q_{\clW_{2n+1}}(D(z,R)).
  \]
  Therefore,the element $\omega^{hiSug}_{(0,\dots,0)}(z)$ is contained the kernel of the map $H^*(\frw^{fact}_{n,k})(D(z,R))$, and hence the map $H^*(\frw^{fact}_{n,k})(D(z,R))$ is not injective.
\end{proof}

On the other hand, when $n=1$, $\frw^{fact}_{1,k}$ contains a quadratic functional on $\clG_{3}$, so the above proof does not go through.
\begin{conjecture}
\label{conj_n_1_inj}
For $n=1$ and $k\neq-h^\vee$, the induced map
\[
H^*(\frw^{fact}_{1,k})\colon H^*(\bbU^\hbar_{\phi(-k\theta_c/h^\vee)}(\clG_3))[\hbar^{-1}]\lto H^*(\Obs^q_{\clW_3})[\hbar^{-1}]
\]
of factorization algebras is injective.
\end{conjecture}

%Finally we consider the case of $n=1$.
%Explicitly, the functional $S^{hiWak}_{1,k}$ is given as follows:
%\begin{align*}
%  &\wt{S}^{hiWak}_{1,k}(\omega_\alpha\otimes e_\alpha,\eta_i\otimes h_i,\xi_\alpha\otimes f_\alpha,\gamma_\alpha,\beta_\alpha,b_i)\\
%=&\int_{\bbC^{3}}\sum_{\alpha\in\Delta_+}\left(\omega_\alpha\wedge\beta_\alpha+\omega_\alpha\wedge\sum_{\beta\in\Delta_+}\beta_\beta\wedge P^\alpha_\beta(\gamma)+\xi_\alpha\wedge\sum_{\beta\in\Delta_+}\beta_\beta\wedge Q^\alpha_\beta(\gamma)-\sum_{i=1}^{\rk\frg}\alpha(h_i)\eta_i\wedge\beta_\alpha\wedge\gamma_\alpha\right.\\
% &\left.+\sqrt{\frac{k+h^\vee}{h^\vee}}\sum_{i=1}^{\rk\frg}\alpha(h_i)^{2}\eta_i\wedge\partial\eta_{i}\wedge b_\alpha+\sqrt{\frac{k+h^\vee}{h^\vee}}\sum_{\beta\in\Delta_+}\alpha([e_\beta,f_\beta])^{2}\xi_\beta\wedge\partial\xi_\beta\wedge b_{\alpha}\wedge\gamma_{\beta}\wedge\gamma_\beta\right)+\hbar\gamma'_{3}.
%\end{align*}
%たぶん1次元のときと同様のfiltrationを考えることによりgrで$\gamma_3'$の項を落として考えるとよいはずである。そのあとで微分の次数についてのfiltrationのgrの各degree partを見る。このときすべてのgraded partで消えるものがkernelになるが、たぶんすべてで消えるものは取れないと思われる。
%まず$\wt{\frw}^{fact}_{3,\theta_c}$のkernelは$\bbC^1$の時のアナログを考えるとFF-centerの$\bbC^3$versionによって生成されるものになるはずである。
%これは$U(\frg)$のcenterの元の``$3$-fold jet''みたいな概念を考えることにより定式化可能なはずである。
%このとき微分作用素の次数が変わらないパートが消えるための必要十分条件はinputがこの``3d FF-center''により生成されるmodの元であることである。
%したがって3d FF-centerの元に対して$\frw^{fact}_{1,k}$の高次の微分のパートが消えないことを示せばよい。これもたぶん微分の次数を細かく見るとたぶんいけるはずである。つまり、levelをつける前のcenterの記述さえできればよい。しかしそれが全然わからない。

\section{Future directions}

Let us close this note with some remarks and comments.
\def\theenumi{\roman{enumi}}
\def\labelenumi{(\theenumi)}

\begin{enumerate}
    \item The construction in this note should admit a natural SUSY analogue.
    T. Iwane \cite{I26} introduced the SUSY analogue of the Wakimoto realization of $N=1$ SUSY affine Vertex algebras.
    Furthermore, S. Yanagida \cite{Yan26} introduced the notion of \emph{$N=1$ SUSY (pre)factorization algebra}, the SUSY analogue of the functor $\bbV$, and constructed an example of a SUSY factorization algebra using free field theory on a $N=1$ SUSY curve.
    Their construction should be understood using factorization algebras by considering a natural SUSY analogue of the construction in this note, which should also admit a natural higher-dimensional generalization.
    
    \item Recall subsections \cref{ss:3.6} and \cref{ss:4.5}.
    We introduced the notion of BD center of BD factorization algebra and showed that the isomorphism
    \[
    \bbV(\clZ^{BD}(\Obs^q_{WZW,\hbar\kappa_c}))\simeq\frz(\widehat{\frg})\dbr{\hbar}
    \]
    of $\bbC\dbr{\hbar}$-linear commutative vertex algebras.
    This should be related to studies of the Geometric Langlands conjecture using the Feigin–Frenkel center. We expect that two-dimensional holomorphic field theory will provide a framework for applying physical methods to this problem (see \cite{ABCFGLRR24,Frenkel,F07}).
    Also, this should also be related to studies of the Geometric Langlands conjecture using four-dimensional $\mathcal{N}=4$ supersymmetric Yang–Mills theory (see \cite{EY18,EY19,EY25}).
    
    By relating four-dimensional and two-dimensional field theories via dimensional reduction (see \cite{ESW22,Yoo25}), we expect to describe the relationship between the two approaches mentioned above in the language of factorization algebras.
    
    \item By \cref{quasi_isom}, we have a $P_0$-factorization algebra structure on the factorization algebra $\Obs^{cl}$ of classical observables.
    The relationship between integrable systems and vertex algebras is described in \cite{FF99}.
    We aim to formulate this relationship in the language of factorization algebras and to investigate its higher-dimensional generalization.
\end{enumerate}

\subsection*{Acknowledgements}
The author would like to thank his advisor, S. Yanagida, for his continued guidance and encouragement.
The author used ChatGPT (GPT-6, OpenAI) for
language editing and for drafting a small amount of explanatory text.


\begin{thebibliography}{GWW25}

\bibitem[ABC+24]{ABCFGLRR24}
D.~Arinkin, D.~Beraldo, J.~Campbell, L.~Chen, J.~Faergeman, D.~Gaitsgory, 
K.~Lin, S.~Raskin, and N.~Rozenblyum, 
\emph{Proof of the geometric Langlands conjecture II: %
 Kac--Moody localization and the FLE}
preprint (2024), 448~pp.; arXiv:2405.03648.

\bibitem[BC04]{JBS}
J.~C.~Baez and A.~S.~Crans,
\emph{Higher-dimensional algebra {VI}: {Lie} 2-algebras},
Theory Appl.\ Categ.\ \textbf{12} (2004), 492--538.

\bibitem[BLL+15]{BLL+15}
C.~Beem, M.~Lemos, P.~Liendo, W.~Peelaers, L.~Rastelli, and B.~C.~van Rees,
\emph{Infinite Chiral Symmetry in Four Dimensions}, 
Commun.\ Math.\ Phys., \textbf{336} no. 3, (2015) 1359--1433,

\bibitem[BMR19]{BMR19}
C.~Beem, C.~Meneghelli, and L.~Rastelli.
 \emph{Free field realizations from the Higgs branch}, 
 J.\ High Energy Phys., (9):058, 2019.

\bibitem[BY16]{Yoo}
D.~Butson and P.~Yoo,
\emph{Degenerate classical field theories and boundary theories},
preprint (2016), 38~pp.; arXiv:1611.00311. %[math-ph].

\bibitem[Cos11a]{Cos4}
K.~Costello,
\emph{A geometric construction of the {Witten} genus, {II}},
preprint (2011), 92~pp.; arXiv:1112.0816. %[math.QA].

\bibitem[Cos11b]{C11}
K.~Costello, 
\emph{Renormalization and Effective Field Theory}, 
Math.\ Surv.\ Monog. \textbf{170}, Amer.\ Math.\ Soc., 2011. 

\bibitem[CGa19]{CG19}
K.~Costello and D.~Gaiotto.
\emph{Vertex Operator Algebras and 3d $N = 4$ gauge theories}, 
J.\ High Energy Phys., (5):018, 2019. 

\bibitem[CG17]{CG}
K.~Costello, O.~Gwilliam, 
\emph{Factorization Algebras in Quantum Field Theory volume 1}, 
new mathematical monographs \textbf{31}, Cambridge Univ.\ Press, 2017. 

\bibitem[CG21]{CG2}
K.~Costello, O.~Gwilliam, 
\emph{Factorization Algebras in Quantum Field Theory volume 2}, 
new mathematical monographs \textbf{41}, Cambridge Univ.\ Press, 2021. 

\bibitem[ESW22]{ESW22}
C.~Elliott, P.~Safronov and B.~R.~Williams,
\emph{A taxonomy of twists of supersymmetric Yang--Mills theory},
Selecta Math., \textbf{28} (2022), no.~4, 73.
%doi:10.1007/s00029-022-00786-y

\bibitem[EY18]{EY18}
C.~Elliott, P.~Yoo,
\emph{Geometric Langlands twists of $N=4$ gauge theory from derived algebraic geometry},
preprint (2018) 55~pp.; arXiv:1507.03048.

\bibitem[EY19]{EY19}
C.~Elliott and P.~Yoo,
\emph{A Physical Origin for Singular Support Conditions in Geometric Langlands Theory},
Commun.\ Math.\ Phys., \textbf{368} (2019), no.~3, 985--1050.
%doi:10.1007/s00220-019-03438-z


\bibitem[EY25]{EY25}
C.~Elliott and P.~Yoo,
\emph{Quantum geometric Langlands categories from $\mathcal{N}=4$ super Yang-Mills theory},
Commun. Num. Theor. Phys. \textbf{19}, no.1, 75-134 (2025)
%doi:10.4310/cntp.250215003117


\bibitem[FF90]{FF90}
B.~Feigin and E.~Frenkel,
\emph{Representations of affine Kac--Moody algebras and bosonization},
in \emph{Physics and Mathematics of Strings},
World Sci., 1990,  pp.~271--316.
 
\bibitem[FF92]{FF92}
B.~Feigin and E.~Frenkel,
\emph{ Affine Kac--Moody algebras at the critical level and
  Gelfand-Dikii algebras},
Int.\ J.\ Mod.\ Phys.\ A \textbf{7}, Suppl.\ 1A (1992), 197--215.

\bibitem[FF99]{FF99}
B.~Feigin and E.~Frenkel,
\emph{Integrable hierarchies and Wakimoto modules},
preprint (1999), 36~pp.; arXiv:math/9902103.

\bibitem[FBZ04]{FBZ}
E.~Frenkel, D.~Ben-Zvi,
\emph{Vertex Algebras and Algebraic Curves}, 2nd ed.,
Math.\ Surv.\ Monog., \textbf{88}, Amer.\ Math.\ Soc., 2004. 

\bibitem[Fre02]{Frenkel}
E.~Frenkel,
\emph{Lectures on Wakimoto modules, opers and the center at the critical level},
preprint (2002), 85~pp., 
arXiv:math/0210029. %[math.QA]

\bibitem[Fre07]{F07}
E. Frenkel. 
\emph{Langlands Correspondence for Loop Groups}, 
volume 103 of Cambridge Studies in Advanced Mathematics. 
Cambridge University Press, Cambridge. 2007.

\bibitem[GGW20]{GGW20}
V.~Gorbounov, O.~Gwilliam, B.~Williams,
\emph{Chiral Differential Operators via Quantization of the Holomorphic $\sigma$
Model}, Ast\'{e}risque, 419 (2020), 210~pp

\bibitem[GLL17]{RGLi}
R.~E.~Grady, Q.~Li, and S.~Li,
\emph{{Batalin--Vilkovisky} quantization and the algebraic index},
Adv.\ Math., \textbf{317} (2017), 575--639.

\bibitem[GMS00]{GMS00}
V.~Gorbounov, F.~Malikov, V.~Schechtman,
\emph{Gerbes of chiral differential operators}, 
Math.\ Res.\ Lett., \textbf{7} (2000), 55–66.

\bibitem[GMS04]{GMS04}
 V. Gorbounov, F. Malikov, V. Schechtman,
 \emph{Gerbes of chiral differential operators II. Vertex algebroids}, Invent.\ 
Math., 155 (2004), 605–680.

\bibitem[GWW25]{GWW25}
Z.~Gui, M.~Wang, B.~R.~Williams, 
\emph{Higher-dimensional Chiral Algebras in the Jouanolou Model and free-field realization}, preprint (2025) 81~pp.; arXiv:2510.26608. %[math.QA].

\bibitem[GW25]{GW}
O.~Gwilliam and B.~R.~Williams,
\emph{Higher {Kac--Moody} algebras and symmetries of holomorphic field theories},
Adv.\ Theor.\ Math.\ Phys., \textbf{25} (2021), no.~1, 129--239.
%preprint (2018); arXiv:1810.06534 [math.QA].

\bibitem[Iwa26]{I26}
T.~Iwane,
\emph{Free field realizations of $ N= 1$ SUSY affine vertex algebras}, 
preprint (2026) 24~pp.; arXiv:2608.22340.


\bibitem[Li02]{Li02}
H.~Li,
\emph{Vertex algebras and vertex Poisson algebras},
Comm.\ Contemp.\ Math., \textbf{6} (2004), no.~1, 61--110.
%arXiv:math/0209310 [math.QA], 2002.

\bibitem[Ni25]{Ni25}
Y.~Nishinaka,
\emph{Factorization envelopes and enveloping vertex algebras}, 
preprint (2025), 49~pp; arXiv:2512.07635.

\bibitem[Yan26]{Yan26}
S.~Yanagida,
\emph{BV construction of SUSY vertex algebras from SUSY factorization algebras}, preprint (2026), 26~pp.; arXiv:2606.05706.


\bibitem[Yoo25]{Yoo25}
P.~Yoo,
\emph{Twists of Supersymmetric Yang--mills Theory in Topological String Theory},
J.\ Korean Math.\ Soc., \textbf{62} (2025), no.~6, 1539--1572.
%doi:10.4134/JKMS.j250056


\bibitem[Wak86]{Wak86}
M.~Wakimoto,
\emph{Fock representations of the affine Lie algebra $A^{(1)}_1$},
Comm.\ Math.\ Phys., \textbf{104} (1986), no.~4, 605--609.

\bibitem[Wil20]{anomaly}
B.~R.~Williams,
\emph{Renormalization for holomorphic field theories},
Comm.\ Math.\ Phys., \textbf{374} (2020), no.~3, 1693--1742.

\end{thebibliography}
\end{document}